
\ifx\shlhetal\undefinedcontrolsequence\let\shlhetal\relax\fi


\input amstex
\expandafter\ifx\csname mathdefs.tex\endcsname\relax
  \expandafter\gdef\csname mathdefs.tex\endcsname{}
\else \message{Hey!  Apparently you were trying to
  \string\input{mathdefs.tex} twice.   This does not make sense.} 
\errmessage{Please edit your file (probably \jobname.tex) and remove
any duplicate ``\string\input'' lines}\endinput\fi




\catcode`\X=12\catcode`\@=11

\def\n@wcount{\alloc@0\count\countdef\insc@unt}
\def\n@wwrite{\alloc@7\write\chardef\sixt@@n}
\def\n@wread{\alloc@6\read\chardef\sixt@@n}
\def\r@s@t{\relax}\def\v@idline{\par}\def\@mputate#1/{#1}
\def\l@c@l#1X{\firstpart.#1}\def\gl@b@l#1X{#1}\def\t@d@l#1X{{}}

\def\RED#1{\pfeilsw
            #1%
				              \pfeilso}

\def\crossrefs#1{\ifx\all#1\let\tr@ce=\all\else\def\tr@ce{#1,}\fi
   \n@wwrite\cit@tionsout\openout\cit@tionsout=\jobname.cit 
   \write\cit@tionsout{\tr@ce}\expandafter\setfl@gs\tr@ce,}
\def\setfl@gs#1,{\def\@{#1}\ifx\@\empty\let\next=\relax
   \else\let\next=\setfl@gs\expandafter\xdef
   \csname#1tr@cetrue\endcsname{}\fi\next}
\def\m@ketag#1#2{\expandafter\n@wcount\csname#2tagno\endcsname
     \csname#2tagno\endcsname=0\let\tail=\all\xdef\all{\tail#2,}
   \ifx#1\l@c@l\let\tail=\r@s@t\xdef\r@s@t{\csname#2tagno\endcsname=0\tail}\fi
   \expandafter\gdef\csname#2cite\endcsname##1{\expandafter
     \ifx\csname#2tag##1\endcsname\relax\RED{##1}\else\csname#2tag##1\endcsname\fi
     \expandafter\ifx\csname#2tr@cetrue\endcsname\relax\else
     \write\cit@tionsout{#2tag ##1 cited on page \folio.}\fi}
   \expandafter\gdef\csname#2page\endcsname##1{\expandafter
     \ifx\csname#2page##1\endcsname\relax\RED{##1}\else\csname#2page##1\endcsname\fi
     \expandafter\ifx\csname#2tr@cetrue\endcsname\relax\else
     \write\cit@tionsout{#2tag ##1 cited on page \folio.}\fi}
   \expandafter\gdef\csname#2tag\endcsname##1{\expandafter
      \ifx\csname#2check##1\endcsname\relax
      \expandafter\xdef\csname#2check##1\endcsname{}%
      \else\immediate\write16{Warning: #2tag ##1 used more than once.}\fi
      \multit@g{#1}{#2}##1/X%
      \write\t@gsout{#2tag ##1 assigned number \csname#2tag##1\endcsname\space
      on page \number\count0.}%
   \csname#2tag##1\endcsname}}

\def\multit@g#1#2#3/#4X{\def\t@mp{#4}\ifx\t@mp\empty%
      \global\advance\csname#2tagno\endcsname by 1 
      \expandafter\xdef\csname#2tag#3\endcsname
      {#1\number\csname#2tagno\endcsnameX}%
   \else\expandafter\ifx\csname#2last#3\endcsname\relax
      \expandafter\n@wcount\csname#2last#3\endcsname
      \global\advance\csname#2tagno\endcsname by 1 
      \expandafter\xdef\csname#2tag#3\endcsname
      {#1\number\csname#2tagno\endcsnameX}
      \write\t@gsout{#2tag #3 assigned number \csname#2tag#3\endcsname\space
      on page \number\count0.}\fi
   \global\advance\csname#2last#3\endcsname by 1
   \def\t@mp{\expandafter\xdef\csname#2tag#3/}%
   \expandafter\t@mp\@mputate#4\endcsname
   {\csname#2tag#3\endcsname\lastpart{\csname#2last#3\endcsname}}\fi}
\def\t@gs#1{\def\all{}\m@ketag#1e\m@ketag#1s\m@ketag\t@d@l p
\let\realscite\scite
\let\realstag\stag
   \m@ketag\gl@b@l r \n@wread\t@gsin
   \openin\t@gsin=\jobname.tgs \re@der \closein\t@gsin
   \n@wwrite\t@gsout\openout\t@gsout=\jobname.tgs }
\outer\def\localtags{\t@gs\l@c@l}
\outer\def\globaltags{\t@gs\gl@b@l}
\outer\def\newlocaltag#1{\m@ketag\l@c@l{#1}}
\outer\def\newglobaltag#1{\m@ketag\gl@b@l{#1}}

\newif\ifpr@ 
\def\m@kecs #1tag #2 assigned number #3 on page #4.%
   {\expandafter\gdef\csname#1tag#2\endcsname{#3}
   \expandafter\gdef\csname#1page#2\endcsname{#4}
   \ifpr@\expandafter\xdef\csname#1check#2\endcsname{}\fi}
\def\re@der{\ifeof\t@gsin\let\next=\relax\else
   \read\t@gsin to\t@gline\ifx\t@gline\v@idline\else
   \expandafter\m@kecs \t@gline\fi\let \next=\re@der\fi\next}
\def\pretags#1{\pr@true\pret@gs#1,,}
\def\pret@gs#1,{\def\@{#1}\ifx\@\empty\let\n@xtfile=\relax
   \else\let\n@xtfile=\pret@gs \openin\t@gsin=#1.tgs \message{#1} \re@der 
   \closein\t@gsin\fi \n@xtfile}

\newcount\sectno\sectno=0\newcount\subsectno\subsectno=0
\newif\ifultr@local \def\ultralocal{\ultr@localtrue}
\def\firstpart{\number\sectno}
\def\lastpart#1{\ifcase#1 \or a\or b\or c\or d\or e\or f\or g\or h\or 
   i\or k\or l\or m\or n\or o\or p\or q\or r\or s\or t\or u\or v\or w\or 
   x\or y\or z \fi}

\def\resetall{\global\advance\sectno by 1\subsectno=0
   \gdef\firstpart{\number\sectno}\r@s@t}
\def\resetsub{\global\advance\subsectno by 1
   \gdef\firstpart{\number\sectno.\number\subsectno}\r@s@t}
\def\newsection#1\par{\resetall\vskip0pt plus.3\vsize\penalty-250
   \vskip0pt plus-.3\vsize\bigskip\bigskip
   \message{#1}\leftline{\bf#1}\nobreak\bigskip}
\def\subsection#1\par{\ifultr@local\resetsub\fi
   \vskip0pt plus.2\vsize\penalty-250\vskip0pt plus-.2\vsize
   \bigskip\smallskip\message{#1}\leftline{\bf#1}\nobreak\medskip}


\newdimen\marginshift

\newdimen\margindelta
\newdimen\marginmax
\newdimen\marginmin

\def\margininit{       
\marginmax=3 true cm                  
				      
\margindelta=0.1 true cm              
\marginmin=0.1true cm                 
\marginshift=\marginmin
}    

\def\t@gsjj#1,{\def\@{#1}\ifx\@\empty\let\next=\relax\else\let\next=\t@gsjj
   \def\@@{p}\ifx\@\@@\else
   \expandafter\gdef\csname#1cite\endcsname##1{\citejj{##1}}
   \expandafter\gdef\csname#1page\endcsname##1{\RED{##1}}
   \expandafter\gdef\csname#1tag\endcsname##1{\tagjj{##1}}\fi\fi\next}
\newif\ifshowstuffinmargin
\showstuffinmarginfalse
\def\jjtags{\ifx\shlhetal\relax 
  \else
\ifx\shlhetal\undefinedcontrolseq
\else
\showstuffinmargintrue
\ifx\all\relax\else\expandafter\t@gsjj\all,\fi\fi \fi
}

\def\tagjj#1{\realstag{#1}\oldmginpar{\zeigen{#1}}}
\def\citejj#1{\rechnen{#1}\mginpar{\zeigen{#1}}}     

\def\rechnen#1{\expandafter\ifx\csname stag#1\endcsname\relax\RED{#1}\else
                           \csname stag#1\endcsname\fi}

\newdimen\theight

\def\marginfont{\sevenrm}

\def\trymarginbox#1{\setbox0=\hbox{\marginfont\hskip\marginshift #1}%
		\global\marginshift\wd0 
		\global\advance\marginshift\margindelta}

\def \oldmginpar#1{%
\ifvmode\setbox0\hbox to \hsize{\hfill\rlap{\marginfont\quad#1}}%
\ht0 0cm
\dp0 0cm
\box0\vskip-\baselineskip
\else 
             \vadjust{\trymarginbox{#1}%
		\ifdim\marginshift>\marginmax \global\marginshift\marginmin
			\trymarginbox{#1}%
                \fi
             \theight=\ht0
             \advance\theight by \dp0    \advance\theight by \lineskip
             \kern -\theight \vbox to \theight{\rightline{\rlap{\box0}}%
\vss}}\fi}

\newdimen\upordown
\global\upordown=8pt
\font\tinyfont=cmtt8 
\def\mginpar#1{\smash{\hbox to 0cm{\kern-10pt\raise7pt\hbox{\tinyfont #1}\hss}}}
\def\mginpar#1{{\hbox to 0cm{\kern-10pt\raise\upordown\hbox{\tinyfont #1}\hss}}\global\upordown-\upordown}


\def\t@gsoff#1,{\def\@{#1}\ifx\@\empty\let\next=\relax\else\let\next=\t@gsoff
   \def\@@{p}\ifx\@\@@\else
   \expandafter\gdef\csname#1cite\endcsname##1{\zeigen{##1}}
   \expandafter\gdef\csname#1page\endcsname##1{?}
   \expandafter\gdef\csname#1tag\endcsname##1{\zeigen{##1}}\fi\fi\next}
\def\verbatimtags{\showstuffinmarginfalse
\ifx\all\relax\else\expandafter\t@gsoff\all,\fi}
\def\zeigen#1{\hbox{$\scriptstyle\langle$}#1\hbox{$\scriptstyle\rangle$}}


\def\margintag#1{\ifshowstuffinmargin\oldmginpar{\zeigen{#1}}\fi}

\def\marginplain#1{\ifshowstuffinmargin\mginpar{{#1}}\fi}
\def\marginbf#1{\marginplain{{\bf \ \ #1}}}

\def\(#1){\edef\dot@g{\ifmmode\ifinner(\hbox{\noexpand\etag{#1}})
   \else\noexpand\eqno(\hbox{\noexpand\etag{#1}})\fi
   \else(\noexpand\ecite{#1})\fi}\dot@g}

\newif\ifbr@ck
\def\eat#1{}
\def\[#1]{\br@cktrue[\br@cket#1'X]}
\def\br@cket#1'#2X{\def\temp{#2}\ifx\temp\empty\let\next\eat
   \else\let\next\br@cket\fi
   \ifbr@ck\br@ckfalse\br@ck@t#1,X\else\br@cktrue#1\fi\next#2X}
\def\br@ck@t#1,#2X{\def\temp{#2}\ifx\temp\empty\let\neext\eat
   \else\let\neext\br@ck@t\def\temp{,}\fi
   \def\teemp{#1}\ifx\teemp\empty\else\rcite{#1}\fi\temp\neext#2X}
\def\resetbr@cket{\gdef\[##1]{[\rtag{##1}]}}
\def\references{\resetbr@cket\newsection References\par}

\newtoks\symb@ls\newtoks\s@mb@ls\newtoks\p@gelist\n@wcount\ftn@mber
    \ftn@mber=1\newif\ifftn@mbers\ftn@mbersfalse\newif\ifbyp@ge\byp@gefalse
\def\defm@rk{\ifftn@mbers\n@mberm@rk\else\symb@lm@rk\fi}
\def\n@mberm@rk{\xdef\m@rk{{\the\ftn@mber}}%
    \global\advance\ftn@mber by 1 }
\def\rot@te#1{\let\temp=#1\global#1=\expandafter\r@t@te\the\temp,X}
\def\r@t@te#1,#2X{{#2#1}\xdef\m@rk{{#1}}}
\def\b@@st#1{{$^{#1}$}}\def\str@p#1{#1}
\def\symb@lm@rk{\ifbyp@ge\rot@te\p@gelist\ifnum\expandafter\str@p\m@rk=1 
    \s@mb@ls=\symb@ls\fi\write\f@nsout{\number\count0}\fi \rot@te\s@mb@ls}
\def\byp@ge{\byp@getrue\n@wwrite\f@nsin\openin\f@nsin=\jobname.fns 
    \n@wcount\currentp@ge\currentp@ge=0\p@gelist={0}
    \re@dfns\closein\f@nsin\rot@te\p@gelist
    \n@wread\f@nsout\openout\f@nsout=\jobname.fns }
\def\m@kelist#1X#2{{#1,#2}}
\def\re@dfns{\ifeof\f@nsin\let\next=\relax\else\read\f@nsin to \f@nline
    \ifx\f@nline\v@idline\else\let\t@mplist=\p@gelist
    \ifnum\currentp@ge=\f@nline
    \global\p@gelist=\expandafter\m@kelist\the\t@mplistX0
    \else\currentp@ge=\f@nline
    \global\p@gelist=\expandafter\m@kelist\the\t@mplistX1\fi\fi
    \let\next=\re@dfns\fi\next}
\def\symbols#1{\symb@ls={#1}\s@mb@ls=\symb@ls} 
\def\bigsymbol{\textstyle}
\symbols{\bigsymbol\ast,\dagger,\ddagger,\sharp,\flat,\natural,\star}
\def\ftnumbers{\ftn@mberstrue} \def\ftsymbols{\ftn@mbersfalse}
\def\paginal{\byp@ge} \def\resetftnumbers{\ftn@mber=1}
\def\ftnote#1{\defm@rk\expandafter\expandafter\expandafter\footnote
    \expandafter\b@@st\m@rk{#1}}

\long\def\jump#1\endjump{}
\def\ssum{\mathop{\lower .1em\hbox{$\textstyle\Sigma$}}\nolimits}

\def\qed{\nobreak\kern 1em \vrule height .5em width .5em depth 0em}
\def\newneq{\hbox{\rlap{\hbox to 1\wd9{\hss$=$\hss}}\raise .1em 
   \hbox to 1\wd9{\hss$\scriptscriptstyle/$\hss}}}
\def\subsetne{\setbox9 = \hbox{$\subset$}\mathrel{\hbox{\rlap
   {\lower .4em \newneq}\raise .13em \hbox{$\subset$}}}}
\def\supsetne{\setbox9 = \hbox{$\subset$}\mathrel{\hbox{\rlap
   {\lower .4em \newneq}\raise .13em \hbox{$\supset$}}}}

\def\vbar{\mathchoice{\vrule height6.3ptdepth-.5ptwidth.8pt\kern-.8pt}
   {\vrule height6.3ptdepth-.5ptwidth.8pt\kern-.8pt}
   {\vrule height4.1ptdepth-.35ptwidth.6pt\kern-.6pt}
   {\vrule height3.1ptdepth-.25ptwidth.5pt\kern-.5pt}}
\def\f@dge{\mathchoice{}{}{\mkern.5mu}{\mkern.8mu}}
\def\b@c#1#2{{\rm \mkern#2mu\vbar\mkern-#2mu#1}}
\def\b@b#1{{\rm I\mkern-3.5mu #1}}
\def\b@a#1#2{{\rm #1\mkern-#2mu\f@dge #1}}
\def\bb#1{{\count4=`#1 \advance\count4by-64 \ifcase\count4\or\b@a A{11.5}\or
   \b@b B\or\b@c C{5}\or\b@b D\or\b@b E\or\b@b F \or\b@c G{5}\or\b@b H\or
   \b@b I\or\b@c J{3}\or\b@b K\or\b@b L \or\b@b M\or\b@b N\or\b@c O{5} \or
   \b@b P\or\b@c Q{5}\or\b@b R\or\b@a S{8}\or\b@a T{10.5}\or\b@c U{5}\or
   \b@a V{12}\or\b@a W{16.5}\or\b@a X{11}\or\b@a Y{11.7}\or\b@a Z{7.5}\fi}}

\catcode`\X=11 \catcode`\@=12




\let\thischap\jobname

\def\partof#1{\csname returnthe#1part\endcsname}
\def\chapof#1{\csname returnthe#1chap\endcsname}

\def\setchapter#1,#2,#3;{%
  \expandafter\def\csname returnthe#1part\endcsname{#2}%
  \expandafter\def\csname returnthe#1chap\endcsname{#3}%
}

\def\setprevious#1 #2 {%
  \expandafter\def\csname set#1page\endcsname{\input page-#2}
}

\setchapter  E53,B,N;       \setprevious E53 null
\setchapter  300z,B,A;       \setprevious 300z E53
\setchapter  88r,B,I;       \setprevious 88r 300z
\setchapter 300a,A,II.A;      \setprevious 300a 88r
\setchapter 300b,A,II.B;       \setprevious 300b 300a
\setchapter 300c,A,II.C;       \setprevious 300c 300b
\setchapter 300d,A,II.D;       \setprevious 300d 300c
\setchapter 300e,A,II.E;       \setprevious 300e 300d
\setchapter 300f,A,II.F;       \setprevious 300f 300e
\setchapter 300g,A,II.G;       \setprevious 300g 300f
\setchapter  600,B,III;       \setprevious  600 300g
\setchapter  705,B,IV;       \setprevious   705 600
\setchapter  734,B,V;        \setprevious   734 705
\setchapter  E46,B,VI;      \setprevious    E46 734
\setchapter  838,B,VII;      \setprevious   838 E46
\setchapter  300x,B,;      \setprevious   300x 838

\newwrite\pageout
\def\rememberpagenumber{\let\setpage\relax
\openout\pageout page-\jobname  \relax \write\pageout{\setpage\the\pageno.}}

\def\recallpagenumber{\csname set\jobname page\endcsname
\def\headmark##1{\rightheadtext{\chapof{\jobname}.##1}}\WRITETOC}
\def\setupchapter#1{\leftheadtext{\chapof{\jobname}. #1}}

\def\setpage#1.{\pageno=#1\relax\advance\pageno1\relax}

\def\cprefix#1{
\edef\theotherpart{\partof{#1}}\edef\theotherchap{\chapof{#1}}%
\ifx\theotherpart\thispart
   \ifx\theotherchap\thischap 
    \else 
     \theotherchap%
    \fi
   \else 
     \theotherchap\fi}

\def\sectioncite[#1]#2{%
     \cprefix{#2}#1}

\def\chaptercite#1{Chapter \cprefix{#1}}

\edef\thispart{\partof{\thischap}}
\edef\thischap{\chapof{\thischap}}

\def\lastpage of '#1' is #2.{\expandafter\def\csname lastpage#1\endcsname{#2}}

\def\xCITE#1{\chaptercite{#1}}
\def\yCITE[#1]#2{\cprefix{#2}.\scite{#2-#1}}

\newwrite\writetoc
\immediate\openout\writetoc \jobname.toc
\def\addcontents#1{\def\WRITETOC{\immediate\write\writetoc{\noexpand\tocentry{\chapof{\jobname}}{#1}{\number\pageno}}}}


\def\spuriousreset{}


\expandafter\ifx\csname citeadd.tex\endcsname\relax
\expandafter\gdef\csname citeadd.tex\endcsname{}
\else \message{Hey!  Apparently you were trying to
\string\input{citeadd.tex} twice.   This does not make sense.} 
\errmessage{Please edit your file (probably \jobname.tex) and remove
any duplicate ``\string\input'' lines}\endinput\fi

\sectno=-1   
\localtags
\jjtags
\define\rest{\restriction}

\newbox\noforkbox \newdimen\forklinewidth
\forklinewidth=0.3pt   
\setbox0\hbox{$\textstyle\bigcup$}
\setbox1\hbox to \wd0{\hfil\vrule width \forklinewidth depth \dp0
                        height \ht0 \hfil}
\wd1=0 cm
\setbox\noforkbox\hbox{\box1\box0\relax}
\def\unionstick{\mathop{\copy\noforkbox}\limits}
\def\nonfork#1#2_#3{#1\unionstick_{\textstyle #3}#2}
\def\nonforkin#1#2_#3^#4{#1\unionstick_{\textstyle #3}^{\textstyle #4}#2}     
%
\setbox0\hbox{$\textstyle\bigcup$}
\setbox1\hbox to \wd0{\hfil{\sl /\/}\hfil}
\setbox2\hbox to \wd0{\hfil\vrule height \ht0 depth \dp0 width
                                \forklinewidth\hfil}
\wd1=0cm
\wd2=0cm
\newbox\doesforkbox
\setbox\doesforkbox\hbox{\box1\box0\relax}
\def\nunionstick{\mathop{\copy\doesforkbox}\limits}

\def\fork#1#2_#3{#1\nunionstick_{\textstyle #3}#2}
\def\forkin#1#2_#3^#4{#1\nunionstick_{\textstyle #3}^{\textstyle #4}#2}     
\NoBlackBoxes

\define\ortp{\text{\bf tp}}

\define\sftp{\text{\rm tp}}
\define\qftp{\hbox{$\sftp_{\text{\rm qf}}$}}

\define\mr{\medskip\roster}
\define\sn{\smallskip\noindent}
\define\mn{\medskip\noindent}
\define\bn{\bigskip\noindent}
\define\ub{\underbar}
\define\wilog{\text{without loss of generality}}
\define\ermn{\endroster\medskip\noindent}

\define\dbcu{\dsize\bigcup}
\define \nl{\newline}

\define\Sfor{\text{\rm Sfr}}
\define\Sav{\text{\rm Sav}}

\magnification=\magstep 1
\documentstyle{amsppt}

{    
\catcode`@11

\ifx\alicetwothousandloaded@\relax
  \endinput\else\global\let\alicetwothousandloaded@\relax\fi

\gdef\subjclass{\let\savedef@\subjclass
 \def\subjclass##1\endsubjclass{\let\subjclass\savedef@
   \toks@{\def\usualspace{{\rm\enspace}}\eightpoint}%
   \toks@@{##1\unskip.}%
   \edef\thesubjclass@{\the\toks@
     \frills@{{\noexpand\rm2000 {\noexpand\it Mathematics Subject
       Classification}.\noexpand\enspace}}%
     \the\toks@@}}%
  \nofrillscheck\subjclass}
} 


\expandafter\ifx\csname alice2jlem.tex\endcsname\relax
  \expandafter\xdef\csname alice2jlem.tex\endcsname{\the\catcode`@}
\else \message{Hey!  Apparently you were trying to
\string\input{alice2jlem.tex}  twice.   This does not make sense.}
\errmessage{Please edit your file (probably \jobname.tex) and remove
any duplicate ``\string\input'' lines}\endinput\fi

\expandafter\ifx\csname bib4plain.tex\endcsname\relax
  \expandafter\gdef\csname bib4plain.tex\endcsname{}
\else \message{Hey!  Apparently you were trying to \string\input
  bib4plain.tex twice.   This does not make sense.}
\errmessage{Please edit your file (probably \jobname.tex) and remove
any duplicate ``\string\input'' lines}\endinput\fi

\def\renewcommand{\newcommand}	       
\edef\cite{\the\catcode`@}%
\catcode`@ = 11
\let\@oldatcatcode = \cite
\chardef\@letter = 11
\chardef\@other = 12
%
%
%
%
\def\@innerdef#1#2{\edef#1{\expandafter\noexpand\csname #2\endcsname}}%
%
%
\@innerdef\@innernewcount{newcount}%
\@innerdef\@innernewdimen{newdimen}%
\@innerdef\@innernewif{newif}%
\@innerdef\@innernewwrite{newwrite}%
%
%
%
\def\@gobble#1{}%
%
%
%
\ifx\inputlineno\@undefined
   \let\@linenumber = \empty 
\else
   \def\@linenumber{\the\inputlineno:\space}%
\fi
%
%
%
\def\@futurenonspacelet#1{\def\cs{#1}%
   \afterassignment\@stepone\let\@nexttoken=
}%
\begingroup 
\def\\{\global\let\@stoken= }%
\\ 
\endgroup
\def\@stepone{\expandafter\futurelet\cs\@steptwo}%
\def\@steptwo{\expandafter\ifx\cs\@stoken\let\@@next=\@stepthree
   \else\let\@@next=\@nexttoken\fi \@@next}%
\def\@stepthree{\afterassignment\@stepone\let\@@next= }%
%
%
%
\def\@getoptionalarg#1{%
   \let\@optionaltemp = #1%
   \let\@optionalnext = \relax
   \@futurenonspacelet\@optionalnext\@bracketcheck
}%
%
%
\def\@bracketcheck{%
   \ifx [\@optionalnext
      \expandafter\@@getoptionalarg
   \else
      \let\@optionalarg = \empty
      \expandafter\@optionaltemp
   \fi
}%
\def\@@getoptionalarg[#1]{%
   \def\@optionalarg{#1}%
   \@optionaltemp
}%
%
%
%
\def\@nnil{\@nil}%
\def\@fornoop#1\@@#2#3{}%
\def\@for#1:=#2\do#3{%
   \edef\@fortmp{#2}%
   \ifx\@fortmp\empty \else
      \expandafter\@forloop#2,\@nil,\@nil\@@#1{#3}%
   \fi
}%
\def\@forloop#1,#2,#3\@@#4#5{\def#4{#1}\ifx #4\@nnil \else
       #5\def#4{#2}\ifx #4\@nnil \else#5\@iforloop #3\@@#4{#5}\fi\fi
}%
\def\@iforloop#1,#2\@@#3#4{\def#3{#1}\ifx #3\@nnil
       \let\@nextwhile=\@fornoop \else
      #4\relax\let\@nextwhile=\@iforloop\fi\@nextwhile#2\@@#3{#4}%
}%
%
%
%
\@innernewif\if@fileexists
\def\@testfileexistence{\@getoptionalarg\@finishtestfileexistence}%
\def\@finishtestfileexistence#1{%
   \begingroup
      \def\extension{#1}%
      \immediate\openin0 =
         \ifx\@optionalarg\empty\jobname\else\@optionalarg\fi
         \ifx\extension\empty \else .#1\fi
         \space
      \ifeof 0
         \global\@fileexistsfalse
      \else
         \global\@fileexiststrue
      \fi
      \immediate\closein0
   \endgroup
}%
%
%
%
%
\def\bibliographystyle#1{%
   \@readauxfile
   \@writeaux{\string\bibstyle{#1}}%
}%
\let\bibstyle = \@gobble
%
%
\let\bblfilebasename = \jobname
\def\bibliography#1{%
   \@readauxfile
   \@writeaux{\string\bibdata{#1}}%
   \@testfileexistence[\bblfilebasename]{bbl}%
   \if@fileexists
      \nobreak
      \@readbblfile
   \fi
}%
\let\bibdata = \@gobble
%
%
\def\nocite#1{%
   \@readauxfile
   \@writeaux{\string\citation{#1}}%
}%
\@innernewif\if@notfirstcitation
%
%
\def\cite{\@getoptionalarg\@cite}%
%
%
\def\@cite#1{%
   \let\@citenotetext = \@optionalarg
   \printcitestart
   \nocite{#1}%
   \@notfirstcitationfalse
   \@for \@citation :=#1\do
   {%
      \expandafter\@onecitation\@citation\@@
   }%
   \ifx\empty\@citenotetext\else
      \printcitenote{\@citenotetext}%
   \fi
   \printcitefinish
}%
\newif\ifweareinprivate
\weareinprivatetrue
\ifx\shlhetal\undefinedcontrolseq\weareinprivatefalse\fi
\ifx\shlhetal\relax\weareinprivatefalse\fi
\def\@onecitation#1\@@{%
   \if@notfirstcitation
      \printbetweencitations
   \fi
   \expandafter \ifx \csname\@citelabel{#1}\endcsname \relax
      \if@citewarning
         \message{\@linenumber Undefined citation `#1'.}%
      \fi
     \iftrue 
      \expandafter\gdef\csname\@citelabel{#1}\endcsname{%
\strut 
\vadjust{\vskip-\dp\strutbox
\vbox to 0pt{\vss\parindent0cm \leftskip=\hsize 
\advance\leftskip3mm
\advance\hsize 4cm\strut\openup-4pt 
\rightskip 0cm plus 1cm minus 0.5cm ?  #1 ?\strut}}
         {\tt
            \escapechar = -1
            \nobreak\hskip0pt\pfeilsw
            \expandafter\string\csname#1\endcsname
             \pfeilso
            \nobreak\hskip0pt
         }%
      }%
     \else  
      \expandafter\gdef\csname\@citelabel{#1}\endcsname{%
            {\tt\expandafter\string\csname#1\endcsname}
      }%
     \fi  
   \fi
   \csname\@citelabel{#1}\endcsname
   \@notfirstcitationtrue
}%
%
%
\def\@citelabel#1{b@#1}%
%
%
\def\@citedef#1#2{\expandafter\gdef\csname\@citelabel{#1}\endcsname{#2}}%
%
%
%
\def\@readbblfile{%
   \ifx\@itemnum\@undefined
      \@innernewcount\@itemnum
   \fi
   \begingroup
      \def\begin##1##2{%
         \setbox0 = \hbox{\biblabelcontents{##2}}%
         \biblabelwidth = \wd0
      }%
      \def\end##1{}
      %
      %
      \@itemnum = 0
      \def\bibitem{\@getoptionalarg\@bibitem}%
      \def\@bibitem{%
         \ifx\@optionalarg\empty
            \expandafter\@numberedbibitem
         \else
            \expandafter\@alphabibitem
         \fi
      }%
      \def\@alphabibitem##1{%
         \expandafter \xdef\csname\@citelabel{##1}\endcsname {\@optionalarg}%
         \ifx\biblabelprecontents\@undefined
            \let\biblabelprecontents = \relax
         \fi
         \ifx\biblabelpostcontents\@undefined
            \let\biblabelpostcontents = \hss
         \fi
         \@finishbibitem{##1}%
      }%
      \def\@numberedbibitem##1{%
         \advance\@itemnum by 1
         \expandafter \xdef\csname\@citelabel{##1}\endcsname{\number\@itemnum}%
         \ifx\biblabelprecontents\@undefined
            \let\biblabelprecontents = \hss
         \fi
         \ifx\biblabelpostcontents\@undefined
            \let\biblabelpostcontents = \relax
         \fi
         \@finishbibitem{##1}%
      }%
      \def\@finishbibitem##1{%
         \biblabelprint{\csname\@citelabel{##1}\endcsname}%
         \@writeaux{\string\@citedef{##1}{\csname\@citelabel{##1}\endcsname}}%
         \ignorespaces
      }%
      %
      %
      \let\em = \bblem
      \let\newblock = \bblnewblock
      \let\sc = \bblsc
      \frenchspacing
      \clubpenalty = 4000 \widowpenalty = 4000
      \tolerance = 10000 \hfuzz = .5pt
      \everypar = {\hangindent = \biblabelwidth
                      \advance\hangindent by \biblabelextraspace}%
      \bblrm
      \parskip = 1.5ex plus .5ex minus .5ex
      \biblabelextraspace = .5em
      \bblhook
      \input \bblfilebasename.bbl
   \endgroup
}%
%
%
\@innernewdimen\biblabelwidth
\@innernewdimen\biblabelextraspace
%
%
%
\def\biblabelprint#1{%
   \noindent
   \hbox to \biblabelwidth{%
      \biblabelprecontents
      \biblabelcontents{#1}%
      \biblabelpostcontents
   }%
   \kern\biblabelextraspace
}%
%
%
%
\def\biblabelcontents#1{{\bblrm [#1]}}%
%
%
\def\bblrm{\rm}%
%
%
\def\bblem{\it}%
%
%
\def\bblsc{\ifx\@scfont\@undefined
              \font\@scfont = cmcsc10
           \fi
           \@scfont
}%
%
%
\def\bblnewblock{\hskip .11em plus .33em minus .07em }%
%
%
\let\bblhook = \empty
%
%
%
\def\printcitestart{[}
\def\printcitefinish{]}
\def\printbetweencitations{, }
\def\printcitenote#1{, #1}
%
%
%
\let\citation = \@gobble
%
%
%
\@innernewcount\@numparams
%
%
\def\newcommand#1{%
   \def\@commandname{#1}%
   \@getoptionalarg\@continuenewcommand
}%
%
%
\def\@continuenewcommand{%
   \@numparams = \ifx\@optionalarg\empty 0\else\@optionalarg \fi \relax
   \@newcommand
}%
%
%
\def\@newcommand#1{%
   \def\@startdef{\expandafter\edef\@commandname}%
   \ifnum\@numparams=0
      \let\@paramdef = \empty
   \else
      \ifnum\@numparams>9
         \errmessage{\the\@numparams\space is too many parameters}%
      \else
         \ifnum\@numparams<0
            \errmessage{\the\@numparams\space is too few parameters}%
         \else
            \edef\@paramdef{%
               \ifcase\@numparams
                  \empty  No arguments.
               \or ####1%
               \or ####1####2%
               \or ####1####2####3%
               \or ####1####2####3####4%
               \or ####1####2####3####4####5%
               \or ####1####2####3####4####5####6%
               \or ####1####2####3####4####5####6####7%
               \or ####1####2####3####4####5####6####7####8%
               \or ####1####2####3####4####5####6####7####8####9%
               \fi
            }%
         \fi
      \fi
   \fi
   \expandafter\@startdef\@paramdef{#1}%
}%
%
%
%
%
\def\@readauxfile{%
   \if@auxfiledone \else 
      \global\@auxfiledonetrue
      \@testfileexistence{aux}%
         \begingroup
            \endlinechar = -1
            \catcode`@ = 11
\citation{Sh:c}
\citation{Sh:1}
\citation{Eh57}
\citation{BoNe94}
\citation{Mo65}
\citation{Sh:a}
\citation{Sh:a}
\citation{Sh:c}
\citation{HHL00}
\citation{Sh:10}
\citation{Sh:a}
\citation{Sh:93}
\citation{KiPi98}
\citation{GIL02}
\citation{Sh:225}
\citation{Sh:225a}
\citation{BLSh:464}
\citation{LwSh:489}
\citation{HyTu91}
\citation{HShT:428}
\citation{HySh:474}
\citation{HySh:529}
\citation{HySh:602}
\citation{CoSh:919}
\citation{BlSh:156}
\citation{BlSh:156}
\citation{BlSh:156}
\citation{Bl85}
\citation{Sh:197}
\citation{Sh:205}
\citation{Sh:284c}
\citation{LwSh:560}
\citation{LwSh:687}
\citation{LwSh:871}
\citation{Sh:840}
\citation{Sh:e}
\citation{GbTl06}
\citation{Sh:c}
\citation{Sh:800}
\citation{Sh:e}
\citation{Sh:800}
\citation{Ke71}
\citation{BaFe85}
\citation{Bal0x}
\citation{Mw85}
\citation{KM67}
\citation{Ke71}
\citation{Sh:3}
\citation{Sh:54}
\citation{Hy98}
\citation{HySh:629}
\citation{HySh:632}
\citation{HySh:629}
\citation{GrLe02}
\citation{GrLe0x}
\citation{GrLe00a}
\citation{Le0x}
\citation{Le0y}
\citation{ChKe62}
\citation{ChKe66}
\citation{He74}
\citation{Str76}
\citation{HeIo02}
\citation{Sh:54}
\citation{BY0y}
\citation{BeUs0x}
\citation{Pi0x}
\citation{ShUs:837}
\citation{Ke71}
\citation{Fr75}
\citation{Sh:48}
\citation{Sh:48}
\citation{Sh:87a}
\citation{Sh:87b}
\citation{GrHa89}
\citation{MkSh:366}
\citation{HaSh:323}
\citation{Zi0xa}
\citation{Zi0xb}
\citation{Sh:300}
\citation{ShVi:635}
\citation{Va02}
\citation{GrSh:222}
\citation{GrSh:238}
\citation{GrSh:259}
\citation{Sh:394}
\citation{Gr91}
\citation{BlSh:330}
\citation{BlSh:360}
\citation{BlSh:393}
\citation{GrVa0xa}
\citation{GrVa0xb}
\citation{BKV0x}
\citation{MaSh:285}
\citation{KlSh:362}
\citation{Sh:472}
\citation{Bal0x}
\citation{Sh:88}
\citation{Jn56}
\citation{Jn60}
\citation{Sh:88}
\citation{Jn56}
\citation{Jn60}
\citation{Jn56}
\citation{Jn60}
\citation{Jn56}
\citation{Jn60}
\citation{MoVa62}
\citation{MoVa62}
\citation{Sh:88}
\citation{Sh:87b}
\citation{Sh:87a}
\citation{Sh:87a}
\citation{Sh:87b}
\citation{Sh:300}
\citation{Sh:155}
\citation{Sh:300}
\citation{Sh:576}
\citation{Sh:88}
\citation{Sh:e}
\citation{Sh:576}
\citation{Sh:576}
\citation{Sh:576}
\citation{Sh:c}
\citation{Sh:200}
\citation{Bal88}
\citation{Sh:3}
\citation{Sh:394}
\citation{Sh:394}
\citation{Sh:56}
\citation{Sh:1}
\citation{Sh:93}
\citation{Sh:576}
\citation{Sh:576}
\citation{Sh:868}
\citation{Sh:576}
\citation{Sh:576}
\citation{Sh:87a}
\citation{Sh:87b}
\citation{Sh:87b}
\citation{Sh:87a}
\citation{Sh:87b}
\citation{Sh:87a}
\citation{Sh:87b}
\citation{Sh:F735}
\citation{Sh:576}
\citation{Sh:603}
\citation{BlSh:862}
\citation{Sh:603}
\citation{Sh:576}
\citation{Sh:576}
\citation{Sh:576}
\citation{Sh:322}
\citation{Sh:87a}
\citation{Sh:87b}
\citation{MaSh:285}
\citation{KlSh:362}
\citation{Sh:472}
\citation{Sh:576}
\citation{Sh:576}
\citation{HaSh:323}
\citation{ShVi:648}
\citation{MaSh:285}
\citation{KlSh:362}
\citation{Sh:472}
\citation{Sh:394}
\citation{Bal0x}
\citation{Sh:842}
\citation{Sh:842}
\citation{Sh:842}
\citation{Sh:c}
\citation{Sh:c}
\citation{Sh:a}
\citation{Sh:c}
\citation{Sh:c}
\citation{Sh:429}
\citation{He92}
\citation{Sh:839}
\citation{Sh:576}
\citation{JrSh:875}
\citation{BlSh:862}
\citation{Zi0xa}
\citation{Zi0xb}
\citation{Sh:F709}
\citation{Sh:E54}
\citation{Sh:F735}
\citation{Sh:F782}
\citation{Sh:F888}
\citation{Sh:E56}
\citation{Sh:F841}
\citation{Ke70}
\citation{Sh:c}
\citation{Sh:842}
\citation{Sh:3}
\citation{HySh:676}
\citation{GrLe0x}
\citation{Sh:839}
\citation{Sh:300}
\citation{Sh:87a}
\citation{Sh:87b}
\citation{Sh:842}
\citation{Sh:300}
\citation{Sh:c}
\citation{Sh:576}
\citation{Sh:F888}
\citation{Sh:842}
\citation{Sh:576}
\citation{Sh:842}
\citation{ignore-this-bibtex-warning}
\citation{Sh:48}
\citation{Sh:c}
\citation{Sh:839}
\citation{Sh:576}
\citation{Sh:576}
\citation{Sh:87a}
\citation{Sh:87b}
\citation{Sh:48}
\citation{Sh:576}
\citation{Sh:E36}
\citation{Sh:56}
\citation{Sh:603}
\citation{Sh:603}
\citation{ignore-this-bibtex-warning}
\citation{Sh:48}
\citation{Sh:87a}
\citation{Sh:87b}
\citation{Sh:87a}
\citation{Sh:87b}
\citation{Sh:48}
\citation{Sh:c}
\citation{Sh:87a}
\citation{Sh:87b}
\citation{Sh:48}
\citation{Sh:87a}
\citation{Sh:87b}
\citation{Sh:868}
\citation{Sh:868}
\citation{Sh:87a}
\citation{Sh:87b}
\citation{Sh:87a}
\citation{Sh:87b}
\citation{Sh:87a}
\citation{Sh:87b}
\citation{Sh:576}
\citation{Sh:46}
\citation{Sh:576}
\citation{Sh:576}
\citation{Sh:88}
\citation{Sh:87a}
\citation{Sh:87b}
\citation{Sh:88}
\citation{Sh:48}
\citation{Sh:88}
\citation{Mw85a}
\citation{Sh:88}
\citation{Sh:F709}
\citation{Sh:88}
\citation{Ke71}
\citation{KM67}
\citation{Sh:87a}
\citation{Sh:87b}
\citation{Sh:3}
\citation{Sh:c}
\citation{Sh:48}
\citation{Ke70}
\citation{Mo70}
\citation{Sh:a}
\citation{Ke70}
\citation{Sh:43}
\citation{Sc76}
\citation{DvSh:65}
\citation{Sh:f}
\citation{Sh:E45}
\citation{Sh:420}
\citation{Sh:E12}
\citation{Sh:c}
\citation{Sh:c}
\citation{Sh:c}
\citation{Sh:394}
\citation{Jo56}
\citation{Jo60}
\citation{KM67}
\citation{Sh:c}
\citation{Sh:c}
\citation{Sh:c}
\citation{Sh:868}
\citation{Sh:c}
\citation{Sh:f}
\citation{Ha61}
\citation{Sh:48}
\citation{Sh:c}
\citation{GrSh:259}
\citation{Sh:48}
\citation{Mo70}
\citation{Mo70}
\citation{Bg}
\citation{Sh:202}
\citation{Sh:43}
\citation{BKM78}
\citation{Sh:48}
\citation{Sh:87a}
\citation{Sh:87b}
\citation{Ke71}
\citation{Sh:48}
\citation{Sh:87a}
\citation{Sh:87b}
\citation{Sh:c}
\citation{Sh:48}
\citation{Sh:87a}
\citation{Sh:87b}
\citation{GrSh:174}
\citation{Sh:87a}
\citation{Sh:87b}
\citation{Sh:87a}
\citation{Sh:87b}
\citation{Sh:87a}
\citation{Sh:87b}
\citation{ignore-this-bibtex-warning}
\citation{Sh:155}
\citation{Sh:c}
\citation{Sh:c}
\citation{Sh:c}
\citation{Sh:c}
\citation{Sh:E54}
\citation{Sh:300}
\citation{Sh:a}
\citation{Sh:16}
\citation{GrSh:222}
\citation{GrSh:238}
\citation{GrSh:259}
\citation{McSh:55}
\citation{GrSh:174}
\citation{Sh:a}
\citation{Sh:c}
\citation{Sh:715}
\citation{Sh:c}
\citation{Sh:c}
\citation{Sh:c}
\citation{Sh:3}
\citation{Sh:e}
\citation{Sh:c}
\citation{MaSh:285}
\citation{Sh:300}
\citation{Sh:e}
\citation{Sh:220}
\citation{Sh:e}
\citation{Sh:16}
\citation{GrSh:222}
\citation{GrSh:259}
\citation{Sh:c}
\citation{Sh:300}
\citation{Sh:e}
\citation{Sh:e}
\citation{Sh:220}
\citation{Sh:e}
\citation{Sh:220}
\citation{Sh:e}
\citation{Sh:e}
\citation{Sh:c}
\citation{GrSh:222}
\citation{GrSh:259}
\citation{Sh:c}
\citation{Sh:c}
\citation{Sh:c}
\citation{Sh:c}
\citation{Sh:c}
\citation{Sh:715}
\citation{Mo65}
\citation{Sh:c}
\citation{ignore-this-bibtex-warning}
\citation{BlSh:330}
\citation{BlSh:360}
\citation{BlSh:393}
\citation{Sh:839}
\citation{Sh:c}
\citation{Sh:c}
\citation{Sh:e}
\citation{Sh:300}
\citation{Sh:3}
\citation{Sh:3}
\citation{Sh:54}
\citation{Sh:3}
\citation{Sh:54}
\citation{Sh:3}
\citation{Sh:3}
\citation{Sh:3}
\citation{Sh:3}
\citation{Sh:c}
\citation{Sh:3}
\citation{Sh:3}
\citation{Sh:3}
\citation{Sh:54}
\citation{Sh:3}
\citation{Sh:3}
\citation{Sh:c}
\citation{HySh:629}
\citation{HySh:676}
\citation{Sh:3}
\citation{ignore-this-bibtex-warning}
\citation{Sh:c}
\citation{Sh:132}
\citation{Sh:c}
\citation{Sh:48}
\citation{Sh:e}
\citation{Sh:420}
\citation{Sh:108}
\citation{Sh:88a}
\citation{Sh:420}
\citation{Sh:351}
\citation{Sh:420}
\citation{Sh:420}
\citation{Sh:351}
\citation{Sh:E54}
\citation{Sh:e}
\citation{DvSh:65}
\citation{Sh:b}
\citation{Sh:87b}
\citation{Sh:c}
\citation{Sh:300}
\citation{Sh:e}
\citation{Sh:300}
\citation{Sh:e}
\citation{Sh:e}
\citation{Sh:300}
\citation{Sh:e}
\citation{ignore-this-bibtex-warning}
\citation{Sh:c}
\citation{Sh:c}
\citation{Sh:a}
\citation{Sh:c}
\citation{Sh:c}
\citation{Sh:a}
\citation{Sh:c}
\citation{Sh:c}
\citation{Sh:c}
\citation{Sh:c}
\citation{Sh:c}
\citation{Sh:c}
\citation{Sh:c}
\citation{Sh:c}
\citation{Sh:c}
\citation{Sh:a}
\citation{Sh:c}
\citation{Sh:c}
\citation{Sh:E54}
\citation{ignore-this-bibtex-warning}
\citation{Sh:c}
\citation{Sh:87a}
\citation{Sh:87b}
\citation{Sh:c}
\citation{Sh:E54}
\citation{Sh:300}
\citation{Sh:e}
\citation{Sh:e}
\citation{Sh:E54}
\citation{Sh:E54}
\citation{Sh:c}
\citation{Sh:c}
\citation{Sh:c}
\citation{Sh:a}
\citation{Sh:c}
\citation{Sh:c}
\citation{Sh:c}
\citation{Sh:c}
\citation{ignore-this-bibtex-warning}
\citation{RuSh:117}
\citation{Sh:f}
\citation{Sh:E54}
\citation{Sh:e}
\citation{Sh:300}
\citation{Sh:e}
\citation{Sh:e}
\citation{Sh:300}
\citation{Sh:E54}
\citation{Sh:E54}
\citation{Sh:E54}
\citation{ignore-this-bibtex-warning}
\citation{Sh:E54}
\citation{Sh:E54}
\citation{KoSh:796}
\citation{ignore-this-bibtex-warning}
\citation{Sh:576}
\citation{Sh:576}
\citation{Sh:576}
\citation{Sh:48}
\citation{Sh:576}
\citation{JrSh:875}
\citation{Sh:F888}
\citation{Sh:576}
\citation{Sh:48}
\citation{Sh:576}
\citation{Sh:87a}
\citation{Sh:87b}
\citation{Sh:87a}
\citation{Sh:87b}
\citation{Sh:48}
\citation{Sh:88}
\citation{Sh:576}
\citation{Sh:576}
\citation{Sh:c}
\citation{Sh:87a}
\citation{Sh:87b}
\citation{Sh:576}
\citation{Sh:87a}
\citation{Sh:576}
\citation{Sh:576}
\citation{Sh:576}
\citation{Sh:576}
\citation{Sh:576}
\citation{Sh:576}
\citation{Sh:576}
\citation{Sh:48}
\citation{Mw85a}
\citation{Sh:88}
\citation{Sh:300}
\citation{Sh:576}
\citation{Sh:842}
\citation{Sh:576}
\citation{Sh:576}
\citation{Sh:87a}
\citation{Sh:87b}
\citation{Sh:87a}
\citation{Sh:87b}
\citation{Sh:87b}
\citation{JrSh:875}
\citation{Sh:842}
\citation{Sh:88}
\citation{Sh:576}
\citation{Sh:48}
\citation{Sh:87a}
\citation{Sh:87b}
\citation{Sh:c}
\citation{Sh:88}
\citation{Sh:87a}
\citation{Sh:87b}
\citation{Sh:48}
\citation{Sh:576}
\citation{Sh:c}
\citation{Sh:c}
\citation{HuSh:342}
\citation{Sh:c}
\citation{Sh:88}
\citation{Sh:88}
\citation{Sh:48}
\citation{Sh:48}
\citation{Sh:48}
\citation{Sh:48}
\citation{Sh:48}
\citation{Ke71}
\citation{Sh:48}
\citation{Sh:48}
\citation{Sh:48}
\citation{Sh:48}
\citation{Sh:48}
\citation{Sh:48}
\citation{Sh:48}
\citation{Sh:48}
\citation{Sh:48}
\citation{Sh:48}
\citation{Sh:48}
\citation{Sh:48}
\citation{Sh:48}
\citation{Sh:48}
\citation{Sh:576}
\citation{Sh:576}
\citation{Sh:576}
\citation{Sh:576}
\citation{Sh:576}
\citation{Sh:576}
\citation{Sh:576}
\citation{Sh:576}
\citation{Sh:576}
\citation{Sh:576}
\citation{Sh:576}
\citation{Sh:576}
\citation{Sh:576}
\citation{Sh:576}
\citation{Sh:576}
\citation{Sh:576}
\citation{Sh:576}
\citation{Sh:576}
\citation{Sh:576}
\citation{Sh:576}
\citation{Sh:576}
\citation{Sh:576}
\citation{Sh:576}
\citation{Sh:482}
\citation{Sh:576}
\citation{Sh:576}
\citation{Sh:832}
\citation{Sh:576}
\citation{Sh:576}
\citation{Sh:576}
\citation{Sh:576}
\citation{Sh:576}
\citation{Sh:576}
\citation{Sh:576}
\citation{Sh:576}
\citation{Sh:576}
\citation{Sh:842}
\citation{Sh:576}
\citation{Sh:48}
\citation{ignore-this-bibtex-warning}
\citation{Sh:F735}
\citation{Sh:c}
\citation{Sh:F735}
\citation{Sh:576}
\citation{Sh:576}
\citation{Sh:48}
\citation{Sh:c}
\citation{Sh:c}
\citation{Sh:225a}
\citation{Sh:429}
\citation{Sh:c}
\citation{ShHM:158}
\citation{Las88}
\citation{Sh:F735}
\citation{Sh:c}
\citation{Sh:576}
\citation{Sh:576}
\citation{Sh:31}
\citation{Sh:576}
\citation{Sh:c}
\citation{Sh:576}
\citation{Sh:F735}
\citation{Sh:c}
\citation{Sh:F735}
\citation{Sh:F735}
\citation{Sh:842}
\citation{Sh:842}
\citation{Sh:F735}
\citation{Sh:c}
\citation{Sh:839}
\citation{Sh:839}
\citation{Sh:300}
\citation{Sh:e}
\citation{Sh:576}
\citation{Sh:87b}
\citation{Sh:c}
\citation{Sh:c}
\citation{Sh:87b}
\citation{Sh:87b}
\citation{Sh:c}
\citation{Sh:87a}
\citation{Sh:87b}
\citation{Sh:87a}
\citation{Sh:87b}
\citation{Sh:842}
\citation{Sh:394}
\citation{Sh:842}
\citation{MaSh:285}
\citation{KlSh:362}
\citation{Sh:472}
\citation{Sh:394}
\citation{Sh:842}
\citation{ignore-this-bibtex-warning}
\citation{Sh:F782}
\citation{Sh:842}
\citation{Sh:842}
\citation{Sh:F820}
\citation{Sh:F782}
\citation{Sh:394}
\citation{Di}
\citation{Sh:g}
\citation{Sh:394}
\citation{Di}
\citation{Sh:300}
\citation{Sh:e}
\citation{Sh:88}
\citation{Sh:F782}
\citation{Kin66}
\citation{Sh:F820}
\citation{Sh:F782}
\citation{Lv71}
\citation{Sh:e}
\citation{Sh:e}
\citation{MiRa65}
\citation{Sh:620}
\citation{Sh:F782}
\citation{Sh:394}
\citation{Sh:394}
\citation{Sh:394}
\citation{Sh:394}
\citation{Sh:394}
\citation{Sh:394}
\citation{Sh:589}
\citation{Sh:g}
\citation{Sh:g}
\citation{Sh:394}
\citation{Sh:394}
\citation{Sh:394}
\citation{Sh:394}
\citation{Sh:394}
\citation{Sh:F782}
\citation{KlSh:362}
\citation{Sh:472}
\citation{Sh:394}
\citation{Sh:394}
\citation{Sh:394}
\citation{ignore-this-bibtex-warning}
\citation{Sh:576}
\citation{Sh:603}
\citation{Sh:576}
\citation{Sh:603}
\citation{Sh:705}
\citation{Sh:c}
\citation{MaSh:285}
\citation{KlSh:362}
\citation{Sh:472}
\citation{Sh:576}
\citation{Sh:603}
\citation{Sh:576}
\citation{Sh:576}
\citation{Sh:603}
\citation{Sh:576}
\citation{Sh:576}
\citation{Sh:576}
\citation{Sh:603}
\citation{Sh:576}
\citation{Sh:F888}
\citation{Sh:300}
\citation{BlSh:862}
\citation{Sh:351}
\citation{Sh:F888}
\citation{Sh:E45}
\citation{Sh:g}
\citation{Sh:g}
\citation{KjSh:409}
\citation{Sh:603}
\citation{Sh:430}
\citation{Sh:F888}
\citation{Sh:922}
\citation{Sh:576}
\citation{ignore-this-bibtex-warning}
\citation{Sh:87b}
\citation{Sh:88r}
\citation{Sh:88}
\citation{Sh:576}
\citation{Sh:603}
\citation{Sh:88r}
\citation{Sh:705}
\citation{Sh:576}
\citation{Sh:87b}
\citation{Sh:88}
\citation{Sh:576}
\citation{Sh:603}
\citation{Sh:87b}
\citation{Sh:842}
\citation{Sh:576}
\citation{Sh:576}
\citation{Sh:576}
\citation{Sh:576}
\citation{Sh:576}
\citation{Sh:576}
\citation{Sh:603}
\citation{Sh:576}
\citation{Sh:f}
\citation{Sh:b}
\citation{Sh:f}
\citation{Sh:f}
\citation{DvSh:65}
\citation{Sh:b}
\citation{Sh:f}
\citation{Sh:E45}
\citation{Sh:576}
\citation{Sh:576}
\citation{Sh:576}
\citation{Sh:603}
\citation{Sh:603}
\citation{Sh:576}
\citation{Sh:576}
\citation{Sh:576}
\citation{Sh:E45}
\citation{Sh:576}
\citation{Sh:576}
\citation{Sh:87b}
\citation{Sh:48}
\citation{Sh:F888}
\citation{Sh:F888}
\citation{Sh:576}
\citation{Sh:F841}
\citation{Sh:F841}
\citation{GiSh:577}
\citation{Sh:576}
\citation{Sh:842}
\citation{Sh:E45}
\citation{Sh:f}
\citation{Sh:460}
\citation{Sh:829}
\citation{Sh:576}
\citation{Sh:603}
\citation{Sh:g}
\citation{Sh:g}
\citation{Sh:430}
\citation{Sh:g}
\citation{Sh:g}
\citation{Sh:f}
\citation{Sh:g}
\citation{Sh:g}
\citation{Sh:g}
\citation{Sh:g}
\citation{Sh:g}
\citation{Sh:g}
\citation{Sh:g}
\citation{Sh:g}
\citation{Sh:g}
\citation{J}
\citation{Sh:460}
\citation{Sh:829}
\citation{Ha61}
\citation{Sh:g}
\citation{Sh:420}
\citation{Sh:g}
\citation{Sh:g}
\citation{Ha61}
\bibstyle{lit-plain}
\bibdata{lista,listb,listx,listf,liste}
\@citedef{Bal88}{Bal88}
\@citedef{Bal0x}{Bal0x}
\@citedef{BlSh:862}{BlSh 862}
\@citedef{Bl85}{Bl85}
\@citedef{BKV0x}{BKV0x}
\@citedef{BLSh:464}{BLSh 464}
\@citedef{BlSh:156}{BlSh 156}
\@citedef{BlSh:330}{BlSh 330}
\@citedef{BlSh:360}{BlSh 360}
\@citedef{BlSh:393}{BlSh 393}
\@citedef{BaFe85}{BaFe85}
\@citedef{BKM78}{BKM78}
\@citedef{BY0y}{BY0y}
\@citedef{BeUs0x}{BeUs0x}
\@citedef{BoNe94}{BoNe94}
\@citedef{Bg}{Bg}
\@citedef{ChKe66}{ChKe66}
\@citedef{ChKe62}{ChKe62}
\@citedef{CoSh:919}{CoSh:919}
\@citedef{DvSh:65}{DvSh 65}
\@citedef{Di}{Di}
\@citedef{Eh57}{Eh57}
\@citedef{Fr75}{Fr75}
\@citedef{GiSh:577}{GiSh 577}
\@citedef{GbTl06}{GbTl06}
\@citedef{Gr91}{Gr91}
\@citedef{GrHa89}{GrHa89}
\@citedef{GIL02}{GIL02}
\@citedef{GrLe0x}{GrLe0x}
\@citedef{GrLe00a}{GrLe00a}
\@citedef{GrLe02}{GrLe02}
\@citedef{GrSh:259}{GrSh 259}
\@citedef{GrSh:174}{GrSh 174}
\@citedef{GrSh:238}{GrSh 238}
\@citedef{GrSh:222}{GrSh 222}
\@citedef{GrVa0xa}{GrVa0xa}
\@citedef{GrVa0xb}{GrVa0xb}
\@citedef{Ha61}{Ha61}
\@citedef{HHL00}{HHL00}
\@citedef{HaSh:323}{HaSh 323}
\@citedef{He74}{He74}
\@citedef{HeIo02}{HeIo02}
\@citedef{He92}{He92}
\@citedef{HuSh:342}{HuSh 342}
\@citedef{Hy98}{Hy98}
\@citedef{HySh:474}{HySh 474}
\@citedef{HySh:529}{HySh 529}
\@citedef{HySh:632}{HySh 632}
\@citedef{HySh:602}{HySh 602}
\@citedef{HySh:629}{HySh 629}
\@citedef{HySh:676}{HySh 676}
\@citedef{HShT:428}{HShT 428}
\@citedef{HyTu91}{HyTu91}
\@citedef{JrSh:875}{JrSh 875}
\@citedef{J}{J}
\@citedef{Jn56}{Jn56}
\@citedef{Jo56}{Jo56}
\@citedef{Jn60}{Jn60}
\@citedef{Jo60}{Jo60}
\@citedef{KM67}{KM67}
\@citedef{Ke70}{Ke70}
\@citedef{Ke71}{Ke71}
\@citedef{KiPi98}{KiPi98}
\@citedef{Kin66}{Kin66}
\@citedef{KjSh:409}{KjSh 409}
\@citedef{KlSh:362}{KlSh 362}
\@citedef{KoSh:796}{KoSh 796}
\@citedef{Las88}{Las88}
\@citedef{LwSh:871}{LwSh 871}
\@citedef{LwSh:489}{LwSh 489}
\@citedef{LwSh:560}{LwSh 560}
\@citedef{LwSh:687}{LwSh 687}
\@citedef{Lv71}{Lv71}
\@citedef{Le0x}{Le0x}
\@citedef{Le0y}{Le0y}
\@citedef{McSh:55}{McSh 55}
\@citedef{MaSh:285}{MaSh 285}
\@citedef{Mw85a}{Mw85a}
\@citedef{Mw85}{Mw85}
\@citedef{MkSh:366}{MkSh 366}
\@citedef{MiRa65}{MiRa65}
\@citedef{MoVa62}{MoVa62}
\@citedef{Mo65}{Mo65}
\@citedef{Mo70}{Mo70}
\@citedef{Pi0x}{Pi0x}
\@citedef{RuSh:117}{RuSh 117}
\@citedef{Sc76}{Sc76}
\@citedef{Sh:88r}{Sh 88r}
\@citedef{Sh:E12}{Sh:E12}
\citation{Sh:g}
\@citedef{Sh:F888}{Sh:F888}
\@citedef{Sh:322}{Sh 322}
\@citedef{Sh:482}{Sh 482}
\@citedef{Sh:922}{Sh:922}
\@citedef{Sh:F820}{Sh:F820}
\@citedef{Sh:832}{Sh 832}
\@citedef{Sh:840}{Sh 840}
\@citedef{Sh:e}{Sh:e}
\@citedef{Sh:F782}{Sh:F782}
\@citedef{Sh:800}{Sh 800}
\@citedef{Sh:F841}{Sh:F841}
\@citedef{Sh:F735}{Sh:F735}
\@citedef{Sh:842}{Sh 842}
\@citedef{Sh:839}{Sh 839}
\@citedef{Sh:705}{Sh 705}
\@citedef{Sh:868}{Sh 868}
\@citedef{Sh:1}{Sh 1}
\@citedef{Sh:3}{Sh 3}
\@citedef{Sh:10}{Sh 10}
\@citedef{Sh:16}{Sh 16}
\@citedef{Sh:31}{Sh 31}
\@citedef{Sh:48}{Sh 48}
\@citedef{Sh:46}{Sh 46}
\@citedef{Sh:43}{Sh 43}
\@citedef{Sh:54}{Sh 54}
\@citedef{Sh:56}{Sh 56}
\@citedef{Sh:a}{Sh:a}
\@citedef{Sh:108}{Sh 108}
\@citedef{Sh:93}{Sh:93}
\@citedef{Sh:b}{Sh:b}
\@citedef{Sh:132}{Sh 132}
\@citedef{Sh:87a}{Sh 87a}
\@citedef{Sh:87b}{Sh 87b}
\@citedef{Sh:202}{Sh 202}
\@citedef{Sh:200}{Sh 200}
\@citedef{Sh:205}{Sh 205}
\@citedef{Sh:197}{Sh 197}
\@citedef{Sh:155}{Sh 155}
\@citedef{Sh:88a}{Sh 88a}
\@citedef{Sh:88}{Sh 88}
\@citedef{Sh:220}{Sh 220}
\@citedef{Sh:225}{Sh 225}
\@citedef{Sh:300}{Sh 300}
\@citedef{Sh:225a}{Sh 225a}
\@citedef{Sh:c}{Sh:c}
\@citedef{Sh:284c}{Sh 284c}
\@citedef{Sh:429}{Sh 429}
\@citedef{Sh:351}{Sh 351}
\@citedef{Sh:420}{Sh 420}
\@citedef{Sh:g}{Sh:g}
\@citedef{Sh:430}{Sh 430}
\@citedef{Sh:f}{Sh:f}
\@citedef{Sh:394}{Sh 394}
\@citedef{Sh:620}{Sh 620}
\@citedef{Sh:589}{Sh 589}
\@citedef{Sh:460}{Sh 460}
\@citedef{Sh:576}{Sh 576}
\@citedef{Sh:472}{Sh 472}
\@citedef{Sh:603}{Sh 603}
\@citedef{Sh:715}{Sh 715}
\@citedef{Sh:829}{Sh 829}
\@citedef{ShHM:158}{ShHM 158}
\@citedef{ShUs:837}{ShUs 837}
\@citedef{ShVi:648}{ShVi 648}
\@citedef{ShVi:635}{ShVi 635}
\@citedef{Sh:E45}{Sh:E45}
\@citedef{Sh:E54}{Sh:E54}
\@citedef{Sh:E56}{Sh:E56}
\@citedef{Sh:F709}{Sh:F709}
\@citedef{Sh:E36}{Sh:E36}
\@citedef{Str76}{Str76}
\@citedef{Va02}{Va02}
\@citedef{Zi0xa}{Zi0xa}
\@citedef{Zi0xb}{Zi0xb}

          \endgroup
      \immediate\openout\@auxfile = \jobname.aux
   \fi
}%
%
%
\newif\if@auxfiledone
\ifx\noauxfile\@undefined \else \@auxfiledonetrue\fi
%
%
%
%
\@innernewwrite\@auxfile
\def\@writeaux#1{\ifx\noauxfile\@undefined \write\@auxfile{#1}\fi}%
%
%
%
\ifx\@undefinedmessage\@undefined
   \def\@undefinedmessage{No .aux file; I won't give you warnings about
                          undefined citations.}%
\fi
%
%
\@innernewif\if@citewarning
\ifx\noauxfile\@undefined \@citewarningtrue\fi
%
%
%
\catcode`@ = \@oldatcatcode

\def\pfeilso{\leavevmode
            \vrule width 1pt height9pt depth 0pt\relax
           \vrule width 1pt height8.7pt depth 0pt\relax
           \vrule width 1pt height8.3pt depth 0pt\relax
           \vrule width 1pt height8.0pt depth 0pt\relax
           \vrule width 1pt height7.7pt depth 0pt\relax
            \vrule width 1pt height7.3pt depth 0pt\relax
            \vrule width 1pt height7.0pt depth 0pt\relax
            \vrule width 1pt height6.7pt depth 0pt\relax
            \vrule width 1pt height6.3pt depth 0pt\relax
            \vrule width 1pt height6.0pt depth 0pt\relax
            \vrule width 1pt height5.7pt depth 0pt\relax
            \vrule width 1pt height5.3pt depth 0pt\relax
            \vrule width 1pt height5.0pt depth 0pt\relax
            \vrule width 1pt height4.7pt depth 0pt\relax
            \vrule width 1pt height4.3pt depth 0pt\relax
            \vrule width 1pt height4.0pt depth 0pt\relax
            \vrule width 1pt height3.7pt depth 0pt\relax
            \vrule width 1pt height3.3pt depth 0pt\relax
            \vrule width 1pt height3.0pt depth 0pt\relax
            \vrule width 1pt height2.7pt depth 0pt\relax
            \vrule width 1pt height2.3pt depth 0pt\relax
            \vrule width 1pt height2.0pt depth 0pt\relax
            \vrule width 1pt height1.7pt depth 0pt\relax
            \vrule width 1pt height1.3pt depth 0pt\relax
            \vrule width 1pt height1.0pt depth 0pt\relax
            \vrule width 1pt height0.7pt depth 0pt\relax
            \vrule width 1pt height0.3pt depth 0pt\relax}

\def\pfeilsw{ \leavevmode 
            \vrule width 1pt height0.3pt depth 0pt\relax
            \vrule width 1pt height0.7pt depth 0pt\relax
            \vrule width 1pt height1.0pt depth 0pt\relax
            \vrule width 1pt height1.3pt depth 0pt\relax
            \vrule width 1pt height1.7pt depth 0pt\relax
            \vrule width 1pt height2.0pt depth 0pt\relax
            \vrule width 1pt height2.3pt depth 0pt\relax
            \vrule width 1pt height2.7pt depth 0pt\relax
            \vrule width 1pt height3.0pt depth 0pt\relax
            \vrule width 1pt height3.3pt depth 0pt\relax
            \vrule width 1pt height3.7pt depth 0pt\relax
            \vrule width 1pt height4.0pt depth 0pt\relax
            \vrule width 1pt height4.3pt depth 0pt\relax
            \vrule width 1pt height4.7pt depth 0pt\relax
            \vrule width 1pt height5.0pt depth 0pt\relax
            \vrule width 1pt height5.3pt depth 0pt\relax
            \vrule width 1pt height5.7pt depth 0pt\relax
            \vrule width 1pt height6.0pt depth 0pt\relax
            \vrule width 1pt height6.3pt depth 0pt\relax
            \vrule width 1pt height6.7pt depth 0pt\relax
            \vrule width 1pt height7.0pt depth 0pt\relax
            \vrule width 1pt height7.3pt depth 0pt\relax
            \vrule width 1pt height7.7pt depth 0pt\relax
            \vrule width 1pt height8.0pt depth 0pt\relax
            \vrule width 1pt height8.3pt depth 0pt\relax
            \vrule width 1pt height8.7pt depth 0pt\relax
            \vrule width 1pt height9pt depth 0pt\relax
      }


\def\widestnumber#1#2{}

\def\citewarning#1{\ifx\shlhetal\relax 
    \else
    \par{#1}\par
    \fi
}

\def\rm{\fam0 \tenrm}

\def\fakesubhead#1\endsubhead{\bigskip\noindent{\bf#1}\par}



%
%
%

%

\font\textrsfs=rsfs10
\font\scriptrsfs=rsfs7
\font\scriptscriptrsfs=rsfs5

\newfam\rsfsfam
\textfont\rsfsfam=\textrsfs
\scriptfont\rsfsfam=\scriptrsfs
\scriptscriptfont\rsfsfam=\scriptscriptrsfs

\edef\oldcatcodeofat{\the\catcode`\@}
\catcode`\@11

\def\Cal@@#1{\noaccents@ \fam \rsfsfam #1}

\catcode`\@\oldcatcodeofat


\expandafter\ifx \csname margininit\endcsname \relax\else\margininit\fi

\long\def\red#1\endred{}
\long\def\green#1\endgreen{}
\long\def\blue#1\endblue{}
\long\def\private#1\endprivate{}

\def\endred{ \unmatched endred! }
\def\endgreen{ \unmatched endgreen! }
\def\endblue{ \unmatched endblue! }
\def\endprivate{ \unmatched endprivate! }

\ifx\latexcolors\undefinedcs\def\latexcolors{}\fi

\def\emptycs{}
\def\evaluatelatexcolors{%
        \ifx\latexcolors\emptycs\else
        \expandafter\xxevaluate\latexcolors\xxfertig\evaluatelatexcolors\fi}
\def\xxevaluate#1,#2\xxfertig{\setupthiscolor{#1}%
        \def\latexcolors{#2}}


\font\smallfont=cmsl7
\def\rutgerscolor{\ifmmode\else\endgraf\fi\smallfont
\advance\leftskip0.5cm\relax}
\def\setupthiscolor#1{\edef\tmptmpcs{\noexpand\bgroup\noexpand\rutgerscolor
\noexpand\def\noexpand\currentcolor{#1}%
\noexpand}%
\expandafter\let\csname#1\endcsname\tmptmpcs
\def\tmptmpcs{\checkColorUnmatched{#1}\popthecolor}
\expandafter\let\csname end#1\endcsname\tmptmpcs}

\def\checkColorUnmatched#1{\def\expectcolor{#1}%
    \ifx\expectcolor\currentcolor   
    \else \edef\failhere{\noexpand\tryingToClose '\currentcolor' with end\expectcolor}\failhere\fi}

\def\currentcolor{???}

\def\popthecolor{\ifmmode\else\endgraf\fi\egroup}

\expandafter\def\csname#1\endcsname{}

\evaluatelatexcolors

 \let\outerhead\head
 \def\head{\innerhead}
 \let\innerhead\outerhead

 \let\outersubhead\subhead
 \def\subhead{\innersubhead}
 \let\innersubhead\outersubhead

 \let\outersubsubhead\subsubhead
 \def\subsubhead{\innersubsubhead}
 \let\innersubsubhead\outersubsubhead

 \let\outerproclaim\proclaim
 \def\proclaim{\innerproclaim}
 \let\innerproclaim\outerproclaim

 %
 %
 %
 %

\def\demo#1{\medskip\noindent{\it #1.\/}}
\def\enddemo{\smallskip}

\def\remark#1{\medskip\noindent{\it #1.\/}}
\def\endremark{\smallskip}

%
%

\pageheight{8.5truein}
 \topmatter
 \title{Categoricity and solvability of a.e.c., quite highly} \endtitle

\author {Saharon Shelah \thanks {\null\newline 
The author would like to thank the Israel Science Foundation for
partial support of this research (Grant No. 242/03). Publication 734.
 \null\newline
I would like to thank Alice Leonhardt for the beautiful typing. 
} \endthanks} \endauthor 

 \affil{The Hebrew University of Jerusalem \\
  Einstein Institute of Mathematics\\
  Edmond J. Safra Campus, Givat Ram \\ 
 Jerusalem 91904, Israel
  \medskip
  Department of Mathematics \\
  Hill Center-Busch Campus \\
  Rutgers, The State University of New Jersey \\
  110 Frelinghuysen Road \\
  Piscataway, NJ 08854-8019 USA} \endaffil
 
 \endtopmatter
\document
 

\head {\S0 Introduction} \endhead  \resetall \sectno=0
 \spuriousreset
\bigskip

The hope which motivates this work is
\sn
\margintag{734-0n.0}\ub{\stag{734-0n.0} Conjecture}:  If ${\frak K}$ is an a.e.c. \ub{then} 
either for every
large enough cardinal $\mu,{\frak K}$ is categorical in $\mu$ or for
every large enough cardinal $\mu,{\frak K}$ is not categorical in $\mu$.
\nl
Why do we consider this a good dream?  See \chaptercite{E53}.

Our main result is \scite{734-f.20}, it says that if ${\frak K}$ is categorical
in $\mu$ (ignoring few exceptional $\mu$'s) and 
$\lambda \in [\text{LS}({\frak K}),\mu)$ has countable
cofinality and is a fix point of the sequence of the $\beth_\alpha$'s,
(moreover a limit of such cardinals) \ub{then} there is a superlimit $M \in
K_\lambda$ for which ${\frak K}_{[M]} = {\frak K}_\lambda \restriction
\{M':M' \cong M\}$ has the amalgamation property
(and a good $\lambda$-frame ${\frak s}$ with 
${\frak K}_{\frak s} = {\frak K}_{[M]}$).
Note that \chaptercite{705} seems to give a strong indication that finding good
$\lambda$-frames is a significant advance.  This may be considered an
unsatisfactory evidence of an advance, being too much phrased in the work's
own terms.  So we prove in \S5 - \S7 that for a 
restrictive context we make a clear cut advance: assuming
amalgamation and enough instances of $2^\lambda < 2^{\lambda^+}$
occurs, much more than the conjecture holds, see \chaptercite{E53} on
background.

Note that as we try to get results on $\lambda = \beth_\lambda >
\text{ LS}({\frak K})$, clearly it does not particularly matter if for
$\kappa \in (\text{LS}({\frak K}),\lambda)$ we use, e.g. $\kappa_1 =
\kappa^+$ or $\kappa_1 = \beth_{(2^\kappa)^+}(= \beth_{1,1}(\kappa))$
or even $\beth_{1,7}(\kappa)$.

After \scite{734-f.20} the next natural step is to show that ${\frak
s}_\lambda$ has the better properties dealt with in \chaptercite{600},
\chaptercite{705}, see \cite{Sh:F782}.
Note that
if we strengthen the assumption on $\mu$ in \S4 (to $\mu = \mu^{<
\lambda}$), then it relies on \S1 only.  Without this we need \S2 
(hence \scite{734-am3.2.3}(1),(4)).
\nl
Originally we have used here categoricity assumptions but lately it
seems desirable to use a weaker one: (variants of) solvability.
About being solvable, see \sectioncite[\S4]{E53}(B), \cite{Sh:842}.  This
seems better as it is a candidate for being an ``outside"
generalization of being superstable (rather than of being categorical).  
\nl
Here we use solvable when it does not require much change; for more on
it see \cite{Sh:842}, \cite{Sh:F820} and 
on material delayed from here see \cite{Sh:F782}.
\nl
Note we can systematically use $K^{\text{sc}(\theta)\text{-lin}}$,
say with $\theta = \aleph_0$ or $\theta = \text{\rm LS}({\frak K})$ instead of
$K^{\text{lin}}$; see Definition \scite{734-0n.2}(8).  
In several respects this is better, but not enough
to make us use it. Also working more it seemed we can get rid of
``wide", ``wide over", see Definition \scite{734-0n.2}(1),(2),(3).  
If instead proving the existence of a 
good $\lambda$-frame it suffices for us to prove the existence of
almost good $\lambda$-frame, \ub{then} the assumption on
$\lambda$ can be somewhat weaker (fixed point instead limit of fix
points of the sequence of the $\beth_\alpha$'s).  In \S7 we sometimes give
alternative quotations in \cite{Sh:394} but do not rely on it.
\nl
We thank Mor Doron, Esther Gruenhut,  
Aviv Tatarski and Alex Usvyatsov for their help in
proofreading.
\nl
Basic knowledge on infinitary logics is assumed, see e.g. \cite{Di};
though the reader may just read the definition here in
\sectioncite[\S5]{E53} and believe some quoted results.
\bigskip

\demo{\stag{734-0n.1} Notation}  Let $\beth_{0,\alpha}(\lambda) =
 \beth_\alpha(\lambda) := \lambda + \Sigma\{\beth_\beta(\lambda):\beta
 < \alpha\}$.  Let 
$\beth_{1,\alpha}(\lambda)$ be defined by induction on
$\alpha:\beth_{1,0}(\lambda) = \lambda$, for limit $\beta$ we let
$\beth_{1,\beta} = \dsize \sum_{\gamma < \beta}
\beth_{1,\gamma}$ and $\beth_{1,\beta + 1}(\lambda) = \beth_\mu$ where
$\mu = (2^{\beth_{1,\beta}(\lambda)})^+$.
\enddemo
\bigskip

\remark{\stag{734-0n.1A} Remark}  1) For our purpose, usually 
$\beth_{1,\beta +1}(\lambda) = \beth_{\delta(\mu)}$ where $\mu = 
\beth_{1,\beta}(\lambda)$ suffice, see e.g. \sectioncite[\S1]{300a} in
particular on $\delta(-)$.  Generally 
$\mu = (\beth_{1,\beta}(\lambda))^+$ is a more natural definition,
but:
\mr
\item "{$(a)$}"  the difference is not significant, e.g. for $\alpha$
limit we get the same value
\sn
\item "{$(b)$}"   our use of omitting types makes our 
choice more natural.
\ermn
2)   We do not use
but it is natural to define $\beth_{\gamma +1,0}(\lambda) = \lambda,
\beth_{\gamma +1,\beta +1}(\lambda) = 
\beth_{\gamma,\mu}(\lambda)$ with $\mu = (2^{\beth_{\gamma
+1,\beta}(\lambda)})^+,
\beth_{\gamma +1,\delta}(\lambda) = \dsize \sum_{\beta < \delta}
\beth_{\gamma +1,\beta}(\lambda)$ and $\beth_{\delta,0}(\lambda) =
\sup\{\beth_{\gamma,0}(\lambda):\gamma < \delta\} = \lambda,
\beth_{\delta,\beta +1}(\lambda) =
\beth_{\delta,\beta}(\beth_{\delta,\beta}(\lambda)),\beth_{\delta,\delta_1}
= \sup\{\beth_{\delta,\alpha}(\lambda):\alpha < \delta_1\}$; 
this is used, e.g. in \cite[ChV]{Sh:g}.
\endremark
\bigskip

\definition{\stag{734-0n.4} Definition}  Assume $M$ is a model, $\tau =
\tau_M$ is its vocabulary and $\Delta$ is a language (or just a 
set of formulas) in some logic, in the vocabulary $\tau$.

For any set $A \subseteq M$ and set $\Delta$ of formulas in the
vocabulary $\tau_M$, let $\Sfor^\alpha_\Delta(A,M)$ which we call the set of
formal $(\Delta,\alpha)$-types over $A$ in $M$, be the set of 
$p$ such that 
\mr
\item "{$(a)$}"  $p$ a set of formulas of the form $\varphi(\bar x,\bar
a)$ where $\varphi(\bar x,\bar y) \in \Delta,\bar x = \langle x_i:i <
\alpha \rangle$ and $\bar a \in {}^{\ell g(\bar y)}A$
\sn
\item "{$(b)$}"  if $\Delta$ is closed under negation (which is the case we
use here) then for any $\varphi(\bar x,\bar y) \in \Delta$ with $\bar x$ as
above and $\bar a \in {}^{\ell g(\bar y)}A$ we have $\varphi(\bar x,\bar a)
\in p$ or $\neg \varphi(\bar x,\bar a) \in p$.
\endroster
\enddefinition
\bn
Recall
\definition{\stag{734-0n.5} Definition}  1) For ${\frak K}$ an a.e.c. we
say $M \in {\frak K}_\theta$ is a superlimit (model in ${\frak K}$ or
in ${\frak K}_\theta$) when:
\mr
\item "{$(a)$}"  $M$ is universal
\sn
\item "{$(b)$}"    if $\delta$ is a limit ordinal $< \theta^+$ and
$\langle M_\alpha:\alpha \le \delta \rangle$ is 
$\le_{{\frak K}_\theta}$-increasing continuous and $\alpha < \delta \Rightarrow
M_\alpha \cong M$ then $M_\delta \cong M$ (equivalently, ${\frak
K}^{[M]}_\theta = {\frak K} \restriction \{N:N \cong M\}$ is a
$\theta$-a.e.c.)
\sn
\item "{$(c)$}"  there is $N$ such that $M <_{\frak K} N \in {\frak
K}_\theta$ and $N$ is isomorphic to $M$.
\ermn
2) We say $M \in {\frak K}_\theta$ is locally superlimit 
when we weaken clause (a) to
\mr
\item "{$(a)^-$}"  if $N \in {\frak K}_\theta$ 
is a $\le_{\frak K}$-extension of $M$ then
$N$ can be $\le_{\frak K}$-embedded into $M$.
\ermn
3) We say that $M$ is pseudo superlimit when in part (1) clauses
(b),(c) hold (but we omit clause (a)); see \scite{734-0n.6}(7) below.
\nl
3A) For $M \in K_\lambda$ let ${\frak K}_{[M]} = {\frak K}^{[M]}_\lambda$ be
${\frak K} \restriction \{N:N \cong M\}$.
\nl
4)  In (1) we may say globally superlimit.
\enddefinition
\bigskip

\demo{\stag{734-0n.6} Observation}  Assume (${\frak K}$ is an a.e.c. and)
${\frak K}_\lambda \ne \emptyset$.
\nl
1) If ${\frak K}$ is categorical in $\lambda$ and there are $M
<_{{\frak K}_\lambda} N$ \ub{then} every 
$M \in {\frak K}_\lambda$ is superlimit.
\nl
2) If every/some $M \in {\frak K}_\lambda$ is superlimit \ub{then}
every/some $M \in K_\lambda$ is locally superlimit.
\nl
3) If every/some $M \in {\frak K}_\lambda$ is locally superlimit
\ub{then} every/some $M \in {\frak K}_\lambda$ is pseudo superlimit.
\nl
4) If some $M \in {\frak K}_\lambda$ is superlimit \ub{then} every
locally superlimit $M' \in {\frak K}_\lambda$ is isomorphic to $M$.
\nl
5) If $M$ is superlimit in ${\frak K}$ \ub{then} $M$ is locally
superlimit in ${\frak K}$.  If $M$ is locally superlimit in ${\frak
K}$, \ub{then} $M$ is pseudo superlimit in ${\frak K}$.  If $M$ is
locally superlimit in ${\frak K}_\theta$ \ub{then} ${\frak K}_\theta$ has
the joint embedding property \ub{iff} $M$ is superlimit.
\nl
6) In Definition \scite{734-0n.5}(1), clause (c)
follows from
\mr
\item "{$(c)^-$}"  LS$({\frak K}) \le \theta$ and $K_{\ge \theta^+}
\ne \emptyset$.
\ermn
7)  $M \in K_\lambda$ is pseudo-superlimit \ub{iff} ${\frak K}_{[M]}$ is a
$\lambda$-a.e.c. and $\le_{{\frak K}_{[M]}}$ is not the equality.  Also
Definition \scite{734-0n.5}(3A) is compatible with \marginbf{!!}{\cprefix{600}.\scite{600-0.33}}.
\enddemo
\bigskip

\definition{\stag{734-0n.7} Definition}  For an a.e.c. ${\frak K}$, let
${\frak K}^{\text{sl}}_\mu,{\frak K}^{\text{ls}}_\mu,
{\frak K}^{\text{pl}}_\mu$ be the
class of $M \in {\frak K}_\mu$ which are superlimit, locally superlimit, pseudo
superlimit respectively with the partial order 
$\le_{{\frak K}^{\text{sl}}_\mu},\le_{{\frak K}^{\text{ls}}_\mu},
\le_{{\frak K}^{\text{pl}}_\mu}$ 
being $\le_{\frak K} \restriction K^{\text{sl}}_\mu,\le_{\frak K}
\restriction K^{\text{pl}}_\mu$ respectively.
\enddefinition
\bigskip

\definition{\stag{734-11.1.3A} Definition}  1) $\Phi$ is proper for linear
orders \ub{when}:
\mr
\item "{$(a)$}"  for some vocabulary $\tau = \tau_\Phi =
\tau(\Phi),\Phi$ is an $\omega$-sequence, the $n$-th element a complete
quantifier free $n$-type in the vocabulary $\tau$
\sn
\item "{$(b)$}"  for every linear order $I$ there is a $\tau$-model
$M$ denoted by EM$(I,\Phi)$, generated by $\{a_t:t \in I\}$ such that
$s \ne t \Rightarrow a_s \ne a_t$ for $s,t \in I$ and
 $\langle a_{t_0},\dotsc,a_{t_{n-1}}\rangle$ realizes the quantifier free
$n$-type from clause (a) whenever $n < \omega$ and
$t_0 <_I \ldots <_I t_{n-1}$; so
really $M$ is determined only up to isomorphism but we may ignore this and
use $I_1 \subseteq J_1 \Rightarrow \text{ EM}(I_1,\Phi) \subseteq
\text{ EM}(I_2,\Phi)$.  We call $\langle a_t:t \in I\rangle$ ``the"
skeleton of $M$; of course again ``the" is an abuse of notation as it is not
necessarily unique.
\ermn
1A) If $\tau \subseteq \tau(\Phi)$ then we let EM$_\tau(I,\Phi)$ be the
$\tau$-reduct of EM$(I,\Phi)$.
\nl
2) $\Upsilon^{\text{or}}_\kappa[{\frak K}]$ is the class of $\Phi$
proper for linear orders satisfying clauses 
$(a)(\alpha),(b),(c)$ of Claim \scite{734-0n.8}(1) below and
$|\tau(\Phi)| \le \kappa$.  The default value of $\kappa$ is
LS$({\frak K})$ and then we may write $\Upsilon^{\text{or}}_{\frak K}$
or $\Upsilon^{\text{or}}[{\frak K}]$ and for simplicity always $\kappa \ge 
\text{ LS}({\frak K})$ (and so $\kappa \ge |\tau_{\frak K}|$).
\nl
3) We define ``$\Phi$ proper for $K$" similarly when in clause (b) of
part (1) we demand $I \in K$, so $K$ is a class of $\tau_K$-models,
i.e.
\mr
\item "{$(a)$}"  $\Phi$ is a function, giving for a quantifier free
$n$-type in $\tau_K$, a quantifier free $n$-type in $\tau_\Phi$
\sn
\item "{$(b)'$}"  in clause (b) of part (1), the quantifier free type
which $\langle a_{t_0},\dotsc,a_{t_{n-1}}\rangle$ realizes in $M$ is
$\Phi(\text{tp}_{\text{qf}}(\langle
t_0,\dotsc,t_{n-1}\rangle,\emptyset,M))$ for $n <
\omega,t_0,\dotsc,t_{n-1} \in I$.
\endroster
\enddefinition
\bigskip

\proclaim{\stag{734-0n.8} Claim}  1) Let ${\frak K}$ be an a.e.c. and
$M \in K$ be of cardinality $\ge
\beth_{1,1}(\text{\rm LS}({\frak K}))$ recalling we naturally assume
$|\tau_{\frak K}| \le \text{\rm LS}({\frak K})$ as usual.

\ub{Then} there is a $\Phi$ such that $\Phi$ is proper for linear orders and:
\mr
\item "{$(a)$}"  $(\alpha) \quad \tau_{\frak K} 
\subseteq \tau_\Phi$,
\sn
\item "{${{}}$}"  $(\beta) \quad |\tau_\Phi| = 
{ \text{\rm LS\/}}({\frak K}) + |\tau_{\frak K}|$ 
\sn
\item "{$(b)$}"  for any linear order $I$ the model 
{\rm EM}$(I,\Phi)$ has cardinality $|\tau(\Phi)|+ |I|$ and 
we have {\rm EM}$_{\tau({\frak K})}(I,\Phi) \in K$
\sn
\item "{$(c)$}"  for any linear orders $I \subseteq J$ we have 
{\rm EM}$_{\tau({\frak K})}(I,\Phi) \le_{\frak K}$ 
{\rm EM}$_{\tau({\frak K})}(J,\Phi)$
\sn
\item "{$(d)$}"  for every finite linear order $I$, the model
{\rm EM}$_{\tau({\frak K})}(I,\Phi)$ can be $\le_{\frak K}$-embedded
into $M$.
\ermn
2) If we allow {\rm LS}$({\frak K}) < |\tau_{\frak K}|$ 
and there is $M \in {\frak K}$ of
cardinality $\ge \beth_{1,1}(\text{\rm LS}({\frak K}) + |\tau_{\frak
K}|)$, \ub{then} there is $\Phi \in \Upsilon^{\text{or}}_{\text{LS}({\frak
K})+|\tau(\Phi)|} [{\frak K}]$ such that {\rm EM}$(I,\Phi)$ has cardinality
$\le \text{\rm LS}({\frak K})$ for $I$ finite.  Hence ${\Cal E}$ has $\le
2^{\text{LS}({\frak K})}$ equivalence classes where ${\Cal E} =
\{(P_1,P_2):P_1,P_2 \in \tau_\Phi$ and
$P^{\text{EM}(I,\Phi)}_1 = P^{\text{EM}(I,\Phi)}_2$ for every linear
order $I\}$. 
\nl
3) Actually having a model of cardinality $\ge \beth_\alpha$ for every
$\alpha < (2^{\text{LS}({\frak K})+|\tau({\frak K})|})^+$ suffice
(in part (2)). 
\endproclaim
\bigskip

\demo{Proof}  Follows from the existence of a representation of
${\frak K}$ as a PC$_{\mu,2^\mu}$-class when $\mu = 
\text{ LS}({\frak K}) + |\tau({\frak K})|$ in
\marginbf{!!}{\cprefix{88r}.\scite{88r-1.4}}(3),(4),(5) and \marginbf{!!}{\cprefix{88r}.\scite{88r-1.8}} 
(or see \cite[0.6]{Sh:394}).  \hfill$\square_{\scite{734-0n.8}}$
\enddemo
\bigskip

\remark{\stag{734-0n.8.13} Remark}  Note that some of the 
definitions and claims below will be used
only in remarks: $K^{\text{sc}(\kappa)}_\theta$ from \scite{734-0n.2}(8), in
\scite{734-11.1.6T}; and some only in \S6,\S7 (and part of \S5 needed for
it): $\Upsilon^{\text{lin}}_\kappa[2]$ from \scite{734-X1.2}(5) (and even less
$\Upsilon^{\text{lin}}_\kappa[\alpha(*)]$ from Definition
\scite{734-0n.2}(9)).  Also the use of
$\le^\otimes_\kappa,\le^{\text{ie}}_\kappa,\le^\oplus_\kappa$ is
marginal.  
\endremark
\bigskip

\definition{\stag{734-X1.2} Definition}   We define partial orders
$\le^\oplus_\kappa,\le^{\text{ie}}_\kappa$ and $\le^\otimes_\kappa$ on 
$\varUpsilon^{\text{or}}_\kappa[{\frak K}]$ 
(for $\kappa \ge \text{ LS}({\frak K}))$ as follows: 
\newline
1)  $\Psi_1 \le^\oplus_\kappa \Psi_2$ \underbar{if} 
$\tau(\Psi_1) \subseteq \tau (\Psi_2)$ and EM$_{\tau({\frak K})}(I,\Psi_1) 
\le_{\frak K} \text{ EM}_{\tau({\frak K})}(I,\Psi_2)$ 
and EM$(I,\Psi_1) = \text{ EM}_{\tau(\Psi_1)}(I,\Psi_1) \subseteq 
\text{ EM}_{\tau(\Psi_1)}(I,\Psi_2)$ for any linear order $I$. \newline
Again for $\kappa = \text{ LS}({\frak K})$ we may drop the $\kappa$. \newline
2) For $\Phi_1,\Phi_2 \in \Upsilon^{\text{or}}_\kappa[{\frak K}]$, we say 
$\Phi_2$ is an inessential extension of $\Phi_1$ and write 
$\Phi_1 \le^{\text{ie}}_\kappa \Phi_2$ \ub{if} 
$\Phi_1 \le^\oplus_\kappa \Phi_2$ and for every linear order $I$,
we have (note: there may be more function symbols in $\tau(\Phi_2)!)$

$$
\text{ EM}_{\tau({\frak K})}(I,\Phi_1) = 
\text{ EM}_{\tau({\frak K})}(I,\Phi_2).
$$
\mn
3) Let $\Upsilon^{\text{lin}}_\kappa$ be the class of 
$\Psi$ proper for linear order and (producing a linear order extending
the original one, i.e.) such that:
\mr
\item "{$(a)$}"  $\tau(\Psi)$  has cardinality $\le \kappa$ and the
two-place predicate $<$ belongs to $\tau(\Psi)$ 
\sn
\item "{$(b)$}"  EM$_{\{<\}}(I,\Psi)$ is a linear order which is an 
extension of $I$ in the sense that EM$(I,\Phi) \models ``a_s < a_t"$
iff $I \models ``s < t"$; in fact we usually stipulate
$[t \in  I \Rightarrow a_t = t]$.
\ermn
4) $\Phi_1 \le^\otimes_\kappa \Phi_2$ \ub{iff} there is $\Psi$ such
that
\mr
\item "{$(a)$}"  $\Psi \in \Upsilon^{\text{lin}}_\kappa$
\sn
\item "{$(b)$}"  $\Phi_\ell \in \Upsilon^{\text{or}}_\kappa[{\frak K}]$ 
for $\ell=1,2$
\sn
\item "{$(c)$}"  $\Phi'_2 \le^{\text{ie}}_\kappa \Phi_2$ where
$\Phi'_2 = \Psi \circ \Phi_1$, i.e. for every linear order $I$ we have
\endroster

$$
\text{ EM}(I,\Phi'_2) = \text{ EM}(\text{EM}_{\{<\}}(I,\Psi),\Phi_1).
$$
\mn
5) $\Upsilon^{\text{lin}}_\kappa[2]$ is the class of $\Psi$ proper
for $K^{\text{lin}}_{\tau^*_2}$ and producing structures from
$K^{\text{lin}}_{\tau^*_2}$ extending the originals, i.e.
\mr
\item "{$(a)$}"  $\tau^*_2 = \{<,P_0,P_1\}$ where $P_0,P_1$ are unary
predicates, $<$ a binary predicate
\sn 
\item "{$(b)$}"   $K^{\text{lin}}_{\tau^*_2} = \{M:M$ a
$\tau^*_2$-model, $<^M$ a linear order, $\langle P^M_0,P^M_1\rangle$ a
partition of $M\}$
\sn
\item "{$(c)$}"  the two-place predicate $<$ and the one place
predicates $P_0,P_1$ belong to $\tau(\Psi)$
\sn
\item "{$(d)$}"  if $I \in K^{\text{lin}}_{\tau^*_2}$ then $M 
= \text{ EM}_{\tau^*_2}(I,\Phi)$ belongs to $K^{\text{lin}}_{\tau^*_2}$ and
$<^M$ is a linear
order and $I\models s <t \Rightarrow M \models a_s < a_t$ and $t \in
P^I_\ell \Rightarrow a_\ell \in P^M_\ell$.
\ermn
6) Similarly $\Upsilon^{\text{lin}}_\kappa[\alpha(*)]$ using
$K^{\text{lin}}_{\tau^*_{\alpha(*)}}$ (see below in \scite{734-0n.2}(9)). 
\enddefinition
\bigskip

\proclaim{\stag{734-0n.em.7} Claim}  Assume $\Phi \in
\Upsilon^{\text{or}}_{\frak K}$.
\nl
1) If $\pi$ is an isomorphism from the
linear order $I_1$ onto the linear order $I_2$ \ub{then} it induces a
unique isomorphism $\hat \pi$ from $M_1 = \text{\rm EM}(I_1,\Phi)$
onto $M_2 = \text{\rm EM}(I_2,\Phi)$ such that:
\mr
\item "{$(a)$}"  $\hat\pi(a_t) = a_{\pi(t)}$ for $t \in I$
\sn
\item "{$(b)$}"  $\hat\pi(\sigma^{M_1}(a_{t_0},
\dotsc,a_{t_{n-1}})) =
\sigma^{M_2}(a_{\pi(t_0)},\dotsc,a_{\pi(t_{n-1})})$ where
$\sigma(x_0,\dotsc,x_{n-1})$ is a $\tau_\Phi$-term and
$t_0,\dotsc,t_{n-1} \in I_1$.
\ermn
2) If $\pi$ is an automorphism of the linear order $I$ \ub{then} it
induces a unique automorphism $\hat \pi$ of {\rm EM}$(I,\Phi)$ (as
above with $I_1 = I = I_2$).
\endproclaim
\bigskip

\remark{\stag{734-0n.em.14} Remark}  1) So in \scite{734-X1.2}(2) 
we allow further expansion 
by functions definable from earlier ones (composition or even 
definition by cases), as long as the number is $\le \kappa$. 
\nl
2) Of course, in \scite{734-0n.em.7} is true for trivial ${\frak K}$.
\endremark
\bn
So we may be interested in some classes of linear orders; below
\scite{734-0n.2}(1) is used much more than the others and also
\scite{734-0n.2}(5),(6) are used not so few times, in particular parts
(8),(9) are not used till \S5.
\definition{\stag{734-0n.2} Definition}  1) A linear order $I$ is
$\kappa$-wide when for every $\theta < \kappa$ there is a monotonic
sequence of lenth $\theta^+$ in $I$.
\nl
2) A linear order $I$ is $\kappa$-wider if $|I| \ge
\beth_{1,1}(\kappa)$.
\nl
3) $I_2$ is $\kappa$-wide over $I_1$ if $I_1 \subseteq I_2$ and for every
$\theta < \kappa$ there is a convex subset of $I_2$ disjoint to $I_1$ which
is $\theta^+$-wide.  We say ``$I_2$ is wide over $I_1$" if ``$I_2$ is
$|I_1|$-wide over $I_2$".
\nl
4) $K^{\text{lin}}[K^{\text{lin}}_\lambda]$ is the class of 
linear orders [of cardinality $\lambda$].
\nl
5)  Let $K^{\text{flin}}$ be the class of infinite linear order $I$
such that every interval has cardinality $|I|$ and is with neither first
nor last elements.
\nl
6) Let the two-place relation $\le_{K^{\text{flin}}}$ on
$K^{\text{flin}}$ be defined by: $I \le_{K^{\text{flin}}} J$ iff
$I,J \in K^{\text{flin}}$ and $I \subseteq J$ and either $I=J$ or
$J \backslash I$ is a dense subset of $J$ and for every $t \in J
\backslash I,I$ can be embedded into $J \restriction \{s \in J
\backslash I:(\forall r \in I)(s <_J r \equiv t <_J r)\}$.
\nl
6A) Let the two-place relation $\le^*_{K^{\text{flin}}}$ on
$K^{\text{lin}}$ be defined similarly omitting ``$I \in
K^{\text{flin}}$" (but not $J \in K^{\text{flin}}$).
\nl
7) $K^{\text{flin}}_\theta = \{I \in K^{\text{flin}}:|I| = \theta\}$
 and $\le_{K^{\text{flin}}_\theta} = \le_{K^{\text{flin}}}
 \restriction K^{\text{flin}}_\theta$.
\nl
8) $K^{\text{sc}(\kappa)-\text{lin}}_\theta$
 is the class of linear orders of cardinality
$\theta$ which are the union of $\le \kappa$ scattered linear orders
(recalling $I$ is scattered when there is no $J \subseteq I$
isomorphic to the rationals).
 If $\kappa = \aleph_0$ we may omit it (i.e. write $K^{\text{sc-lin}}_\theta$).
\nl
9) Let $\tau^*_{\alpha(*)} = \{<\}
\cup \{P_i:i < \alpha(*)\},P_i$ a monadic predicate,
$K^{\text{lin}}_{\tau^*_{\alpha(*)}} = \{I:I$ a $\tau^*_{\alpha(*)}$-model,
$<^I$ a linear order and $\langle P^I_i:i < \alpha(*) \rangle$ a
partition of $I\}$.
If $\alpha(*)=1$ we may omit $P^I_0$, so $I$ is a linear order, so any
ordinal can be treated as a member of $K^{\text{lin}}_{\tau^*_1}$. 
\enddefinition
\bigskip

\demo{\stag{734-0n.3} Observation}  1) If $|I| > 2^\theta$ \ub{then} $I$ is 
$\theta^+$-wide.
\nl
2) If $|I| \ge \lambda$ and $\lambda$ is a 
strong limit cardinal \ub{then} $I$ is $\lambda$-wide. 
\nl
3) $(K^{\text{flin}}_\theta,\le_{K^{\text{flin}}_\theta})$ almost is a
 $\theta$-a.e.c., only smoothness may fail.
\nl
4) If $I_1 \in K^{\text{lin}}$ \ub{then} for some $I_2 \in
K^{\text{flin}}$ we have: $|I_2| = |I_1| + \aleph_0$ and 
$I_1 \le^*_{K^{\text{flin}}} I_2$; and
$(\forall I_0)[I_0 \subseteq I_1 \wedge 
I_0 \in K^{\text{flin}} \Rightarrow  I_0 \le_{K^{\text{flin}}} I_2]$.
\nl
5) If $I_1$ is $\kappa$-wide and $I_1 <_{K^{\text{flin}}} I_2$ \ub{then}
$I_2$ is $\kappa$-wide over $I_2$.
\enddemo
\bigskip

\remark{Remark}   If in the definition of $\le_{K^{\text{flin}}}$
in \scite{734-0n.2}(6) we can add ``$(\forall t \in I)
(\exists t' \in J)[t' <_J t \wedge (\forall s \in I)(s <_I t
\rightarrow s <_J t')]$" (and its dual, i.e. inverting the order).
So we can strengthen \scite{734-0n.2}(6) by the demand above. 
\endremark
\bigskip

\demo{Proof}  1) By Erd\"os-Rado Theorem, i.e., by $(2^\theta)^+
\rightarrow (\theta^+)^2_2$.
\nl
2) Follows by part (1). 
\nl
3),4),5) Easy.   \hfill$\square_{\scite{734-0n.3}}$
\enddemo
\bigskip

\proclaim{\stag{734-X1.3} Claim}  1)  
$(\varUpsilon^{\text{or}}_{\kappa[{\frak K}]},\le^\otimes_\kappa),
(\Upsilon^{\text{or}}_\kappa[{\frak K}],<^{\text{ie}}_\kappa)$
and $(\varUpsilon^{\text{or}}_{\kappa[{\frak K}]},\le^\oplus)$ are 
partial orders (and $\le^\otimes_\kappa,
\le^{\text{ie}}_\kappa \subseteq \le^\oplus_\kappa$). 
\newline
2) If $\Phi_i \in \Upsilon^{\text{or}}_\kappa[{\frak K}]$ and the sequence
$\langle \Phi_i:i < \delta \rangle$ is a 
$\le^\otimes_\kappa$-increasing sequence, $\delta < \kappa^+$, \ub{then} it 
has a $<^\otimes_\kappa$-{\rm l.u.b.} 
$\Phi \in \Upsilon^{\text{or}}_\kappa[{\frak K}]$, and 
{\rm EM}$(I,\Phi) = \dsize 
\bigcup_{i<\delta} { \text{\rm EM\/}}(I,\Phi_i)$ for every linear
order $I$, i.e. $\tau(\Phi) = \cup\{\tau(\Phi_i):i < \delta\}$ and for
every $j < \delta$ we have {\rm EM}$_{\tau(\Phi_j)}(I,\Phi) =
\cup\{\text{\rm EM}_{\tau(\Phi_i)}(I,\Phi):i \in [j,\delta)\}$.  
\nl
3) Similarly for $<^\oplus_\kappa$ and $\le^{\text{ie}}_\kappa$.
\nl
4) If $\Phi \in \Upsilon^{\text{lin}}_\kappa$ and $I \in
K^{\text{lin}}$ \ub{then} $I \subseteq \,
\text{\rm EM}_{\{<\}}(I,\Phi)$ as linear orders stipulating (as in
\scite{734-X1.2}(3)) that $a_t = t$.
\endproclaim
\bigskip

\demo{Proof}  Easy.    \hfill$\square_{\scite{734-X1.3}}$
\enddemo
\bn
Recall various well known facts on $\Bbb L_{\infty,\theta}$.
\proclaim{\stag{734-0n.13} Claim}   1) If $M,N$ are $\tau$-models of
cardinality $\lambda$, {\rm cf}$(\lambda) = \aleph_0$ and $M \equiv_{\Bbb
L_{\infty,\lambda}} N$ \ub{then} $M \cong N$. 
\nl
2) If $M,N$ are $\tau$-models then $M \equiv_{\Bbb L_{\infty,\theta}} N$
   \ub{iff} there is ${\Cal F}$ such that
\mr
\item "{$\circledast$}"  $(a)(\alpha) \quad$ each $f \in {\Cal F}$ is a
   partial isomorphism from $M$ to $N$
\sn
\item "{${{}}$}"  $\quad (\beta) \quad {\Cal F} \ne \emptyset$
\sn
\item "{${{}}$}"  $\quad (\gamma) \quad$ if $f \in {\Cal F}$ and $A
\subseteq \text{\rm Dom}(f)$ then $f \restriction A \in {\Cal F}$
\sn
\item "{${{}}$}"  $(b) \quad$ if $f \in {\Cal F},A \in [M]^{< \theta}$
and $B \in [N]^{< \theta}$ then for some $g \in {\Cal F}$ we have
\nl

\hskip25pt  $f \subseteq g,A \subseteq \text{\rm Dom}(g),
B \subseteq \text{\rm Rang}(g)$.
\ermn
2A) If $M \subseteq N$ are $\tau$-models, \ub{then} $M \prec_{\Bbb
L_{\infty,\theta}} N$ iff for some ${\Cal F}$ clauses
$\circledast(a),(b)$ hold together with
\mr
\item "{${{}}$}"  $(c) \quad$ if $A \in [M]^{< \theta}$ then for some $f
\in {\Cal F}$ we have {\rm id}$_A \subseteq f$.
\ermn
2B) In part (2) (and part (2A)), we can omit subclause $(\gamma)$ of
clause (a), and if ${\Cal F}$ satisfies $(a)(\alpha),(\beta)+(b)$ (and
$(c)$), \ub{then} also
${\Cal F}' = \{f \restriction A:f \in {\Cal F}$ and $A \subseteq
\text{\rm Dom}(f)\}$ satisfies the demands.
\nl
2C) Let $M,N$ be $\tau$-models and define ${\Cal F} = \{f$: for some
$\bar a \in {}^{\theta >}M,f$ is a function from {\rm Rang}$(\bar a)$ to $N$
such that $(M,\bar a) \equiv_{\Bbb L_{\infty,\theta}} 
(N,f(\bar a))\}$ \ub{then}
$M \equiv_{\Bbb L_{\infty,\theta}} N$ \ub{iff} ${\Cal F} \ne
\emptyset$ iff ${\Cal F}$ satisfies clauses (a),(b) of $\circledast$. 
\nl
3) If $M$ is a $\tau$-model, $\theta = \text{\rm cf}(\theta)$ and $\mu =
\|M\|^{< \theta}$ \ub{then} for some $\gamma < \mu^+$ and $\Delta
 \subseteq \Bbb L_{\mu^+,\theta}(\tau)$ of cardinality $\le \mu$ such that
each $\varphi(\bar x) \in \Delta$ is of quantifier depth $< \gamma$, we have
\mr
\item "{$(a)$}"  for $\bar a,\bar b \in {}^{\theta >}M$ we have $(M,\bar
 a) \equiv_{\Bbb L_{\infty,\theta}}(M,\bar b)$ iff $\sftp_\Delta(\bar
 a,\emptyset,M) = \sftp_\Delta(\bar a,\emptyset,M)$
\sn
\item "{$(b)$}"  for any $\tau$-model $N$ we have $N \equiv_{\Bbb
 L_{\infty,\theta}} M$ \ub{iff} $\{\sftp_\Delta(\bar
 a,\emptyset,N):\bar a \in {}^{\theta >}N\} =\{\sftp_\Delta(\bar
 a,\emptyset,M):\bar a \in {}^{\theta >}M\}$.
\ermn
4) Assume $\chi > \mu = \mu^{< \kappa}$ and $x \in {\Cal H}(\chi)$.  
There is ${\frak B}$ such
that (in fact clauses (d)-(g) follow from clauses (a),(b),(c))
\mr
\item "{$(a)$}"  ${\frak B} \prec ({\Cal H}(\chi),\in)$ has
cardinality $\mu$,
\sn
\item "{$(b)$}"  $\mu +1 \subseteq {\frak B}$ and $[{\frak B}]^{<
\kappa} \subseteq {\frak B}$ and $x \in {\frak B}$
\sn
\item "{$(c)$}"  ${\frak B} \prec_{\Bbb L_{\kappa,\kappa}} 
({\Cal H}(\chi),\in)$
\sn
\item "{$(d)$}"  if ${\frak K}$ is an a.e.c. with {\rm LS}$({\frak K}) +
|\tau({\frak K})| \le \mu$ and ${\frak K} \in {\frak B}$ (which means
$\{(M,N):M \le_{\frak K} N$ has universes $\subseteq 
\text{\rm LS}({\frak K})\} \in {\frak B}$) \ub{then}
{\roster 
\itemitem{ $(\alpha)$ }  $M \in {\frak K} \cap {\frak B} \Rightarrow M
\restriction {\frak B} := M \restriction ({\frak B} \cap M) \le_{\frak
K} M$
\sn
\itemitem{ $(\beta)$ }  if $M \le_{\frak K} N$ belongs to ${\frak B}$
then $M \restriction {\frak B} \le_{\frak K} N \restriction {\frak B}$
\endroster}  
\item "{$(e)$}"  if ${\frak K}$ is as in (d), $\Phi \in
\Upsilon^{\text{or}}_{\le \mu} [{\frak K}] \cap {\frak B}$ and $I \in
{\frak B}$ is a linear order and so $M = \,\text{\rm EM}(I,\Phi) \in {\frak
B}$ \ub{then} $I' = I \restriction {\frak B} \subseteq I$ and $M
\restriction {\frak B} = \text{\rm EM}(I',\Phi)$ so $(M \restriction
\tau({\frak K})) \restriction {\frak B} = \text{\rm EM}_{\tau({\frak
K})}(I',\Phi) \le_{\frak K} M \restriction \tau({\frak K})$
\sn 
\item "{$(f)$}"  if $|\tau| \le \mu,\tau \in {\frak B}$ and $M,N \in
{\frak B}$ are $\tau$-models, \ub{then}
{\roster
\itemitem{ $(\alpha)$ }  $M \restriction {\frak B} 
\prec_{\Bbb L_{\kappa,\kappa}[\tau]} M$
\sn
\itemitem{ $(\beta)$ }   $M \not\equiv_{\Bbb L_{\infty,\kappa}[\tau]} N$
\ub{iff} $(M \restriction {\frak B}) 
\not\equiv_{\Bbb L_{\infty,\kappa}[\tau]}(N \restriction {\frak B})$
\sn
\itemitem{ $(\gamma)$ }  if $M \subseteq N$ \ub{then} 
$(M \prec_{\Bbb L_{\infty,\kappa}(\tau)} N)$ iff $(M \restriction
{\frak B}) \prec_{\Bbb L_{\infty,\kappa}(\tau)} (N \restriction {\frak
B})$; this applies also to $(M,\bar a),(N,\bar a)$ for $\bar a \in
{}^{\kappa >} M$
\endroster}
\item "{$(g)$}"  if $I \in K^{\text{flin}}$ then $I_1 \cap {\frak B}
\in K^{\text{flin}}$ and if $I_1 <^*_{K^{\text{flin}}} I_2$ then $(I_1
\cap {\frak B}) <^*_{K^{\text{flin}}} (I_2 \cap {\frak B})$. 
\endroster
\endproclaim
\bigskip

\demo{Proof}  1)-3) and 4)(a),(b),(c)  Well known, e.g. see \cite{Di}.
\nl
4) Clauses (d),(e),(f): as in \scite{734-0n.8}(1), i.e. by absoluteness.
Also clause (g) should be clear.   \hfill$\square_{\scite{734-0n.13}}$
\enddemo
\bigskip

\remark{\stag{734-0n.16} Remark}  1) We will be able to add, in \scite{734-0n.13}(4):
\mr
\item "{$(g)$}"  if ${\frak K}$ is as in clause (d) and $\tau =
\tau_{\frak K}$ \ub{then} in clause (f) we can replace $\Bbb
L_{\infty,\kappa}(\tau)$ by $\Bbb L_{\infty,\kappa}[{\frak K}]$ and
$\Bbb L_{\kappa,\kappa}(\tau)$ by $\Bbb L_{\kappa,\kappa}[{\frak K}]$,
see Definition \scite{734-11.1A} and Fact \scite{734-11.1B}(5).
\ermn
2) We use part (4) in \scite{734-11.10}(3).
\endremark
\bigskip

\definition{\stag{734-0n.17} Definition}  For a model $M$ and for a set
$\Delta$ of formulas in the vocabulary of $M,
\bar x = \langle x_i:i < \alpha\rangle,A \subseteq M$ and 
$\bar a \in {}^\alpha M$ \ub{let} the $\Delta$-type of 
$\bar a$ over $A$ in $M$ be tp$_\Delta(\bar a,A,M) 
= \{\varphi(\bar x,\bar b):M \models \varphi[\bar a,\bar b]$
where $\varphi = \varphi(\bar x,\bar y) \in \Delta$ 
and $\bar b \in {}^{\ell g(\bar y)}A\}$.
\enddefinition
\goodbreak

\head {\S1 Amalgamation in $K^*_\lambda$} \endhead  \resetall \sectno=1
 \spuriousreset
\bigskip

Our aim is to investigate what is implied by \scite{734-11.0.3} below but
instead of assuming it we shall shortly assume only some of its
consequences.  For our purpose here, for $\theta \in [\text{LS}({\frak
K}),\lambda),\lambda = \beth_\lambda$ it does not really matter if we use
$\kappa = \beth_{1,1}(\theta)$ or $\kappa =
\beth_{1,1}(\beth_n(\theta))$ or $\beth_{1,n}(\theta)$, as we are
trying to analyze models in $K_\lambda$.
\bigskip

\remark{\stag{734-11n.0.7} Remark}  1) We can in our claims use only $\Phi \in 
\Upsilon^{\text{or}}_{\frak K} =
\Upsilon^{\text{or}}_{\text{LS}({\frak K})}[{\frak K}]$ 
because for every $\theta \ge 
\text{ LS}({\frak K})$ we can replace ${\frak K}$ by ${\frak K}_{\ge \theta}$
as LS$({\frak K}_{\ge \theta}) = \theta$ when ${\frak K}_{\ge \theta}
\ne \emptyset$, of course.
\nl
2) As usual we assume $|\tau_{\frak K}| \le \text{ LS}({\frak K})$ just 
for convenience, otherwise we should just replace LS$({\frak K})$ by
LS$({\frak K}) + |\tau_{\frak K}|$.  
\endremark
\bigskip

\demo{\stag{734-11.0} Hypothesis}
\mr
\item "{$(a)$}"  ${\frak K} = (K,\le_{\frak K})$ is an a.e.c. with vocabulary 
$\tau = \tau({\frak K})$ (and we can assume $|\tau| \le \text{
LS}({\frak K})$ for notational simplicity)
\sn
\item "{$(b)$}"   ${\frak K}$ has arbitrarily large models
(equivalently has a model of cardinality 
$\ge \beth_{1,1}(\text{LS}({\frak K})))$, not used, e.g. in 
\scite{734-11.1B}, \scite{734-11.1.21} but from \scite{734-11.2} on it is
used extensively.
\endroster
\enddemo
\bigskip

\definition{\stag{734-11.0.3} Definition}  We say $(\mu,\lambda)$ or really
$(\mu,\lambda,\Phi)$ is a weak/strong/pseudo ${\frak K}$-candidate
when (weak is the default value):
\mr
\item "{$(a)$}"  $\mu > \lambda = \beth_\lambda > \text{ LS}({\frak
K})$ (e.g. the first beth fix point 
$> \text{ \rm LS}({\frak K})$, see \scite{734-e.2}; in the main case 
$\lambda$ has cofinality $\aleph_0$)
\sn
\item "{$(b)$}"  ${\frak K}$ categorical in $\mu$ and $\Phi \in
\Upsilon^{\text{or}}_{\frak K}$ 
\nl
or just
\sn
\item "{$(b)^-$}"  ${\frak K}$ is weakly/strongly/pseudo
solvable in $\mu$ and $\Phi \in \Upsilon^{\text{or}}_{\frak K}$ 
witnesses it; see below.
\endroster
\enddefinition
\bigskip

\definition{\stag{734-11.1} Definition}  1) We say ${\frak K}$ is weakly
$(\mu,\kappa)$-solvable \ub{when} $\mu \ge \kappa \ge 
\text{ LS}({\frak K})$ and there is $\Phi \in
\Upsilon^{\text{or}}_\kappa[{\frak K}]$ witnessing it, which means that
$\Phi \in \Upsilon^{\text{or}}_\kappa[{\frak K}]$
and {\rm EM}$_{\tau({\frak K})}(I,\Phi)$ is a locally superlimit
member of ${\frak K}_\mu$ for every linear order $I$ of
cardinality $\mu$.  We may say $({\frak K},\Phi)$ is weakly
$(\mu,\kappa)$-solvable and we may say 
$\Phi$ witness that ${\frak K}$ is weakly $(\mu,\kappa)$-solvable.

If $\kappa = \text{ LS}({\frak K})$ we may omit it,
saying ${\frak K}$ or $({\frak K},\Phi)$ is weakly $\mu$-solvable in $\mu$.  
\nl
2) ${\frak K}$ is strongly $(\mu,\kappa)$-solvable \ub{when} 
$\mu \ge \kappa \ge \text{ LS}({\frak K})$ and some $\Phi \in
\Upsilon^{\text{or}}_\kappa[{\frak K}]$ witness it which means that
if $I \in K^{\text{lin}}_\lambda$ then EM$_{\tau[{\frak K}]}(I,\Phi)$
is superlimit (for ${\frak K}_\lambda$).  We use the conventions from
part (1).
\nl
3) We say ${\frak K}$ is pseudo $(\mu,\kappa)$-solvable 
when $\mu \ge \kappa \ge \text{ LS}({\frak K})$ and
there is $\Phi \in \Upsilon^{\text{or}}_\kappa[{\frak K}]$ witnessing it 
which means that for some $\mu$-a.e.c. ${\frak K}'$ with no
$\le_{{\frak K}'}$-maximal member, we have $M \in {\frak K}'$
\ub{iff} $M \cong \text{ EM}_{\tau({\frak K})}(I,\Phi)$ for some $I
\in K^{\text{lin}}_\mu$ iff $M \cong \,\text{\rm EM}_{\tau({\frak
K})}(I,\Phi)$ for every $I \in K^{\text{lin}}_\mu$.  
We use the conventions from part (1).
\nl
4) Let $(\mu,\kappa)$-solvable mean weakly $(\mu,\kappa)$-solvable,
 etc., (including \scite{734-11.0.3})
\enddefinition
\bigskip

\proclaim{\stag{734-11.1.2} Claim}  1) In Definition \scite{734-11.0.3}, clause
(b) implies clause (b)$^-$.  Also in Definition \scite{734-11.1} ``${\frak
K}$ is strongly $(\mu,\kappa)$-solvable" implies ``${\frak K}$ is weakly
$(\mu,\kappa)$-solvable" which implies ``${\frak K}$ is pseudo
$(\mu,\kappa)$-solvable".  Similarly for $({\frak K},\Phi)$.
\nl
2) Assume $ \Phi \in \Upsilon^{\text{or}}_\kappa[{\frak K}]$;
if clause $(b)^-$ of \scite{734-11.0.3} or just 
$\dot I(\mu,{\frak K}) < 2^\mu$, or just
$2^\mu > \dot I(\mu,\{\text{\rm EM}_{\tau({\frak K})}(I,\Phi):
I \in K^{\text{lin}}_\mu\})$ for some $\mu$ satisfying 
{\rm LS}$({\frak K}) < \kappa^+ < \mu$ \ub{then} we can deduce that
\mr
\item "{$(*)$}"  $\Phi$, really $({\frak K},\Phi)$ has the $\kappa$-non-order 
property, where the $\kappa$-non-order property means that:
\nl
\ub{if} $I$ is a linear order of cardinality $\kappa,\bar t^1,\bar t^2
\in {}^\kappa I$ form a $\Delta$-system pair (see below) and $\langle
\sigma_i(\bar x):i <\kappa \rangle$ lists the $\tau(\Phi)$-terms (with
the sequence $\bar x$ of variables being 
$\langle x_i:i < \kappa\rangle)$ and $\langle a_t:t \in I\rangle$ is ``the"
indiscernible sequence generating {\rm EM}$(I,\Phi)$ 
(i.e. as usual ``$\langle a_t:t
\in I\rangle$" is ``the" skeleton of {\rm EM}$(I,\Phi)$, so generating
it, see Definition \scite{734-11.1.3A}) \ub{then} for some $J \supseteq I$
there is an automorphism of {\rm EM}$_{\tau({\frak K})}(J,\Phi)$ which
exchanges $\langle \sigma_i(\langle a_{t^1_i}:i < \kappa \rangle):i <
\kappa\rangle$ and 
$\langle \sigma_i(\langle a_{t^2_i}:i < \kappa \rangle):i < \kappa\rangle$.  
\nl
where
{\roster
\itemitem{ $\boxtimes$ }  $\bar t^1,\bar t^2 \in {}^\alpha I$ is a
$\Delta$-system pair when for some $J \supseteq I$ there are $\bar
t^\zeta \in {}^\alpha J$ for $\zeta \in \kappa \backslash \{1,2\}$
such that $\langle \bar t^\alpha:\alpha < \kappa \rangle$ is an
indiscernible sequence for quantifier free formulas in the linear order $J$.
\endroster}
\endroster
\endproclaim
\bigskip

\demo{Proof}  1) The first sentence holds by 
Claim \scite{734-0n.8}(1) and Definition \scite{734-11.1.3A} (and Claim
\scite{734-0n.6}).  The second and third sentences follows by \scite{734-0n.6}.
\nl
2) Otherwise we get a contradiction by \cite[Ch.III]{Sh:300} or better 
\cite[III]{Sh:e}.   \hfill$\square_{\scite{734-11.1}}$
\enddemo
\bigskip

\definition{\stag{734-11.1.3} Definition}  1) If ${\Cal M}'$ is a class of
linear orders and $\Phi \in \Upsilon^{\text{or}}_\kappa[{\frak K}]$
then we let $K[{\Cal M}',\Phi] 
= \{\text{\rm EM}_{\tau({\frak K})}(I,\Phi):I \in {\Cal M}'\}$.
\nl
2) Let $K^{u(\kappa)\text{-lin}}_\theta$ be the class of linear orders
$I$ of cardinality $\theta$ such that for some scattered
\footnote{i.e. one into which the rational order cannot be embedded} linear
order $J$ and $\Phi$ proper for $K^{\text{lin}}$ such that $<$ belongs to
$\tau_\Phi,|\tau_\Phi| \le \kappa$ we have $I$ is embeddable into
EM$_{\{<\}}(J,\Phi)$.  If we omit $\kappa$ we mean LS$({\frak K})$.
If $\kappa = \aleph_0$ we may omit it.
\enddefinition
\bigskip

\remark{\stag{734-11.1.6T} Remark}  1) Note that in Definition \scite{734-11.1}(1) we
\ub{can} restrict ourselves to $I \in 
K^{\text{sc}(\theta)\text{-lin}}_\lambda$, see
\scite{734-0n.2}(8) and even $I \in K^{u(\theta)\text{-lin}}$ see
\scite{734-11.1.3}(2), 
i.e., assume $2^\mu > \dot I(\mu,K[{\Cal M}',\Phi])$, 
for ${\Cal M}' = K^{\text{sc}(\theta)\text{-lin}}_\lambda$ or ${\Cal M}' =
K^{u(\theta)\text{-lin}}_\lambda$ and restrict the
conclusion $(*)$ to $I \in K^{\text{sc}(\theta)\text{-lin}}$.
A gain is that, if $\lambda > \theta$, every $I \in
K^{\text{sc}(\theta)\text{-lin}}_\lambda$ is $\lambda$-wide so later
$K^* = K^{**}$, and being solvable is a weaker demand.  But it is less
natural.  Anyhow we presently do not deal with this.
\nl
1A)  Note that $K^{\text{sc}(\theta)-\text{lin}}_\lambda \supseteq
K^{u(\theta)-\text{lin}}_\lambda$. 
\nl
2) An aim of \scite{734-11.1.7} below is to show that: by 
changing $\Phi$ instead of assuming $I_1 \subset I_2 \wedge (I_2$
is $\kappa$-wide over $I_1$) it suffices to assume 
$I_1 \subset I_2 \wedge (I_2$ is $\kappa$-wide).
\endremark
\bigskip

\proclaim{\stag{734-11.1.7} Claim}  For every $\Phi_1 \in
\Upsilon^{\text{or}}_\kappa[{\frak K}]$ there is $\Phi_2$ such that
\mr
\item "{$(a)$}"  $\Phi_2 \in \Upsilon^{\text{or}}_\kappa[{\frak K}]$
and if $\Phi_1$ witnesses ${\frak K}$ is weakly/strongly/pseudo
$(\lambda,\kappa)$-solvable then so does $\Phi_2$
\sn
\item "{$(b)$}"  $\tau_{\Phi_1} \subseteq \tau_{\Phi_2}$ and
$|\tau_{\Phi_2}| = |\tau_{\Phi_1}| + \aleph_0$
\sn
\item "{$(c)$}"  for any $I_2 \in K^{\text{lin}}$ there are $I_1$ and
$h$ such that:
{\roster
\itemitem{ $(\alpha)$ }  $I_1 \in K^{\text{lin}}$ and even $I_1 \in
K^{\text{flin}}$, see \scite{734-0n.2}(5) 
\sn
\itemitem{ $(\beta)$ }  $h$ is an embedding of $I_2$ into $I_1$
\sn
\itemitem{ $(\gamma)$ }   there is an isomorphism $f$ from {\rm
EM}$_{\tau(\Phi_1)}(I_2,\Phi_2)$ onto {\rm EM}$(I_1,\Phi_1)$ such
that $f(a_t) = a_{h(t)}$ for $t \in I_2$
\sn
\itemitem{ $(\delta)$ }  if $J_1 = I_1 \restriction \text{\rm Rang}(h)$
 and we let ${\Cal E} = \{(t_1,t_2):t_1,t_2 \in I_1 \backslash J_1$ and
$(\forall s \in J_1)(s < t_1 \equiv s < t_2)\}$ \ub{then}: ${\Cal E}$ is an
equivalence relation and each equivalence class has $\ge |I_2|$
members and $J_1 \le_{K^{\text{flin}}} I_1$, see \scite{734-0n.2}(6)
\sn
\itemitem{ $(\varepsilon)$ }  [not used] if $\emptyset \ne J_2
\subseteq I_2$, $J_1 = \{t \in I_1$: for some $\tau(\Phi_2)$-term
$\sigma(x_0,\,\dotsc,\,x_{n-1})$ and some $t_0, \,\dots\,,\, t_{n-1} \in J_2$ we have
$f^{-1}(a_t) = 
\sigma^{\text{EM}(I_2,\Phi_2)}(a_{t_0},\dotsc,a_{t_{n-1}})\}$
and $J'_1
\subseteq \text{\rm Rang}(h) \backslash J_1$ and $t \in J'_1$ then $\{s
\in t/{\Cal E}:f^{-1}(a_s)$ belongs to the Skolem hull of 
$\{f^{-1}(a_r):r \in J'_1\}$ in {\rm EM}$(I_2,\Phi)\}$ has
cardinality $\ge |J'_1|$ and $J'_1$ and its inverse can be embedded
into it; in fact, $I_1$ and its inverse are embeddable into any 
interval of $I_2$.
\endroster}
\endroster
\endproclaim
\bigskip

\remark{Remark}  1) We can express it by $\le^\otimes_\kappa$, see
\scite{734-X1.2}(4).   So for some $\Psi$ proper for linear orders such
that $\tau_\Psi$ is countable, the two-place predicate $<$ belongs to
$\tau_\Psi$ and above EM$_{\{<\}}(I_2,\Psi)$ is $I_1$.
\nl
2) In fact, $J_2 \subset I_2 \Rightarrow \text{ EM}_{\{<\}}(J_2,\Psi)
<_{K^{\text{flin}}} \text{ EM}_{\{<\}}(I_2,\Psi)$ and $I_2
<^*_{K^{\text{flin}}} \text{ EM}_{\{<\}}(I_2,\Phi)$ when we identify $t
\in I_2$ with $a_t$.  
\endremark
\bigskip

\demo{Proof}  For $I_2 \in K^{\text{lin}}$ let the set of elements of
$I_1$ be $\{\eta:\eta$ is
a finite sequence of elements from $(\Bbb Z \backslash \{0\}) \times
I_2\}$.  For $\eta
\in I_1$ let $(\ell_{\eta,k},t_{\eta,k})$ be $\eta(k)$ for $k < \ell
g(\eta)$.

Lastly, $I_1$ is ordered by: 
$\eta_1 < \eta_2$ iff for some $n$ one of the following occurs
\mr
\item "{$\circledast$}"  $(a) \quad \eta_1 \restriction n = \eta_2
\restriction n,\ell g(\eta_1) > n,\ell g(\eta_2) > n$ and
$\ell_{\eta_1,n} < \ell_{\eta_2,n}$
\sn
\item "{${{}}$}"  $(b) \quad \eta_1 \restriction n = \eta_2 \restriction
n,\ell g(\eta_1) > n,\ell g(\eta_2) > n,\ell_{\eta_1,n} =
\ell_{\eta_2,n} > 0$
and
 \nl

\hskip25pt $t_{\eta_1,n} <_{I_2} t_{\eta_2,n}$
\sn
\item "{${{}}$}"  $(c) \quad \eta_1 \restriction n = \eta_2
\restriction n,\ell g(\eta_1) > n,\ell g(\eta_2) > n,\ell_{\eta_1,n} =
 \ell_{\eta_2,n} < 0$ and
\nl

\hskip25pt  $t_{\eta_2,n} <_{I_2} t_{\eta_1,n}$
\sn
\item "{${{}}$}"  $(d) \quad \eta_1 \restriction n = \eta_2
\restriction n,\ell g(\eta_1) = n,\ell g(\eta_2) > n$ and
$\ell_{\eta_2,n} > 0$
\sn
\item "{${{}}$}"  $(e) \quad \eta_1 \restriction n = \eta_2
\restriction n,\ell g(\eta_1) > n,\ell g(\eta_2) = n$ and
$\ell_{\eta_1,n} < 0$.  
\ermn
We identify $t \in I_1$ with the pair $(1,t)$.  Now check.  
\hfill$\square_{\scite{734-11.1.7}}$
\enddemo
\bigskip

\definition{\stag{734-11.1A} Definition}  1) Let the language
$\Bbb L_{\theta,\partial}
[{\frak K}]$ or $\Bbb L_{\theta,\partial,{\frak K}}$ where 
$\theta \ge \partial \ge \aleph_0$ and $\theta$ is possibly 
$\infty$, be defined like the infinitary logic $\Bbb L_{\theta,\partial}
(\tau_{\frak K})$, except that we deal only with models from $K$ and
we add for $i^* < \partial$ the 
atomic formula ``$\{x_i:i < i^*\}$ is the universe
of a $\le_{\frak K}$-submodel", with obvious syntax and semantics.
Of course, it is interesting normally only for $\partial > 
\text{ LS}({\frak K})$ and
recall that any formula has $< \partial$ free variables.
\nl
2) For $M$ a $\tau_{\frak K}$-model and 
$N \in K$ let $M \prec_{\Bbb L_{\theta,\partial}[{\frak K}]} N$
means that $M \subseteq N$ and if $\varphi(\bar x,\bar y)$ is a
formula from $\Bbb L_{\theta,\partial}[{\frak K}]$ and 
$N \models (\exists \bar
x)\varphi(\bar x,\bar b)$ where $\bar b \in {}^{\ell g(\bar y)}M$,
\ub{then} for some $\bar a \in {}^{\ell g(\bar x)}M$ we have $N \models
\varphi[\bar a,\bar b]$.
\enddefinition
\bn
\margintag{734-11.1B}\ub{\stag{734-11.1B} Fact}:  1) If $\theta \ge \partial > \text{\rm LS}({\frak K})$ 
and $M,N$ are $\tau_{\frak K}$-models and 
$N \in K$ and $M \prec_{{\Bbb L}_{\theta,\partial} [{\frak K}]} N$,
\ub{then} $M \le_{\frak K} N$ and $M \in K$.
\nl
2) The relation $\prec_{\Bbb L_{\theta,\partial}[{\frak K}]}$ can also
be defined as
usual: $M \prec_{\Bbb L_{\theta,\partial}[{\frak K}]} N$ \ub{iff} $M,N
\in K,M \subseteq N$ and for every $\varphi(\bar x) \in 
\Bbb L_{\theta,\partial}[{\frak K}]$ and $\bar a \in {}^{\ell g(\bar x)}M$
we have $M \models \varphi[\bar a]$ iff $N \models \varphi[\bar a]$. 
\nl
3) If $N \in {\frak K}$ and $M$ is a $\tau_K$-model satisfying
 $M \prec_{\Bbb L_{\infty,\kappa}} N$ and $\kappa > \text{ LS}({\frak K})$ 
\ub{then} $M \in K,M \le_{\frak K} N$ and 
$M \prec_{\Bbb L_{\infty,\kappa}[{\frak K}]} N$.
\nl
4) If $N \in K,M$ a $\tau_K$-model and $M 
\equiv_{\Bbb L_{\infty,\kappa}} N$ where $\kappa > \text{LS}({\frak K})$
\ub{then} $M \in K$ and $M \equiv_{\Bbb L_{\infty,\kappa}[{\frak K}]}
N$.
\nl
5) The parallel of \scite{734-0n.13}(2) holds for 
$\Bbb L_{\infty,\kappa}[{\frak K}]$, i.e. there is ${\Cal F}$
satisfying clauses (a),(b) there and
\mr
\item "{$(d)$}"  if $f \in {\Cal F}$ then
{\roster
\itemitem{ $(\alpha)$ }   $M \restriction \text{ Dom}(f) 
\le_{\frak K} M$
\sn
\itemitem{ $(\beta)$ }  $N \restriction \text{ Rang}(f) \le_{\frak K} M$.
\endroster}
\ermn
6) Also the parallel of \scite{734-0n.13}(2A) holds for $\Bbb
L_{\infty,\kappa}[{\frak K}]$.
\nl
7) The parallel of \scite{734-0n.13}(4) holds for $\Bbb
L_{\infty,\kappa}[{\frak K}]$.
\bigskip

\demo{Proof}  Part (1) is straight (knowing \sectioncite[\S1]{88r} or
\cite[\S1]{Sh:88}).  Part (2) is proved as in the Tarski-Vaught
criterion and parts (5),(6),(7) are proved as in \scite{734-0n.13}.

Toward proving parts (3),(4) we first assume just
\mr
\item "{$\boxtimes_1$}"  $M,N$ are $\tau_K$-models, $N \in K$ and 
$M \equiv_{\Bbb L_{\infty,\kappa}} N$ and
$\kappa > \text{ LS}({\frak K})$ and 
$\lambda \in [\text{LS}({\frak K}),\kappa)$
\ermn
and we define:
\mr
\item "{$\boxdot$}"  $(a) \quad I = I_\lambda = \{(f,M',N'):M' \subseteq M 
\text{ and } N' \subseteq N \text{ and } f \text{ is an isomorphism}$
\nl

\hskip25pt from $M' \text{ onto } N' \text{ and } \|M'\| \le \lambda 
\text{ and letting } \bar a \text{ list } M' \text{ we have}$
\nl

\hskip25pt $(M,\bar a) \equiv_{\Bbb L_{\infty,\kappa}} (N,f(\bar a))\}$
\sn
\item "{${{}}$}"  $(b) \quad$ for $t \in I$ let $t = (f_t,M_t,N_t)$
\sn
\item "{${{}}$}"  $(c) \quad$ for $\ell=0,1,2$ we define the two-place
relation $\le^\ell_I$ on $I$: 
\nl

\hskip25pt let $s \le^\ell_I t$ hold \ub{iff}
{\roster
\itemitem{ ${{}}$ }  $(\alpha) \quad \ell=0$ and $M_s \subseteq M_t 
\wedge N_s \subseteq N_t$
\sn
\itemitem{ ${{}}$ }  $(\beta) \quad \ell=1$ and 
$M_s \le_{\frak K} M_t \wedge N_s \le_{\frak K} N_t$
\sn
\itemitem{ ${{}}$ }  $(\gamma) \quad \ell=2$ and $f_s \subseteq f_t$
\endroster}
\item "{${{}}$}"   $(d) \quad I_1 = I^1_\lambda := \{t \in I_0:N_t
\le_{\frak K} N\}$ and let $\le^\ell_{I_1} = \le^\ell_I \restriction
I_1$ for $\ell=0,1,2$.
\ermn
Now easily
\mr
\item "{$(*)_0$}"  $(\alpha) \quad I \ne \emptyset$ is 
partially ordered by $\le^\ell_I$ for $\ell=0,1,2$
\sn
\item "{${{}}$}"  $(\beta) \quad s \le^1_I t \Rightarrow s \le^0_I t$
\sn
\item "{${{}}$}"  $(\gamma) \quad s \le^2_I t \Rightarrow s \le^0_I t$.
\ermn
[Why?  Straight, e.g. $I \ne \emptyset$ by \scite{734-0n.13}(1).]
\mr
\item "{$(*)_1$}"  if $t \in I_1$ \ub{then} $M_t \in K_{\le \lambda}$ and
$N_t \in K_{\le \lambda}$.
\ermn
[Why?  As $t \in I_1$ by the definition of $I$ we have
$N_t \in K_{\le \lambda}$ (because $N_t \le_{\frak K} N$) and 
$M_t \in K_{\le \lambda}$ as $f_t$ is an isomorphism from $M_t$ onto $N_t$.]
\mr
\item "{$(*)_2$}"  if $s \in I,A \in [M]^{\le \lambda}$ and $B \in 
[N]^{< \lambda}$ \ub{then} for some $t$ we have $s \le^2_I t$ and $A
\subseteq M_t$ and $B \subseteq N_t$.
\ermn
[Why?  By the properties of $\equiv_{\Bbb L_{\infty,\kappa}}$, see
\scite{734-0n.13}(2C) as $\kappa > \lambda,M \equiv_{\Bbb
L_{\infty,\kappa}}$ and the definition of $I$.]
\mr
\item "{$(*)_3$}"  if $s \le^2_{I_1} t$ \ub{then} $s \le^1_I t$, i.e.
$M_s \le_{\frak K} M_t$ and $N_s \le_{\frak K} N_t$.
\ermn
[Why?  As $s,t \in I_1$ we know that $N_s \le_{\frak K} N$ and $N_t
\le_{\frak K} N$ and as $s \le^2_I t$ we have $f_s \subseteq f_t$
hence $N_s \subseteq N_t$.  By axiom V of a.e.c. it follows that $N_s
\le_{\frak K} N_t$.  Now $M_s \le_{\frak K} M_t$ as
$f_t$ is an isomorphism from $M_t$ onto $N_t$ mapping
$M_s$ onto $N_s$ (as it extends $f_s$ by the definition of $\le^2_I$)
and $\le_{\frak K}$ is preserved by any isomorphism.  So by the
definition of $\le^1_I$ we are done.]
\mr
\item "{$(*)_4$}"  if $s \in I$ \ub{then} for some $t \in I_1$ we have $s
\le^2_I t$ (hence $I_1 \ne \emptyset$).
\ermn
[Why?  First choose $N' \le_{\frak K} N$ of cardinality $\le \lambda$
such that $N_s \subseteq N'$, (possibly by the basic properties of
a.e.c. (see \sectioncite[\S1]{88r} or \chaptercite{300b})).  
Second we can find $t \in I$ such that
$N_t = N' \wedge f_s \subseteq f_t$ by the characterization of
$\equiv_{\Bbb L_{\infty,\kappa}}$ as in $(*)_2$.  
So $s \le^2_I t$ by the definition
of $\le^2_I$ and $N_t = N' \le_{\frak K} N$ hence $t \in I_1$ as
required.  Lastly, $I_1 \ne \emptyset$ as by
$(*)_0(\alpha)$ we know that $I \ne \emptyset$ and apply what we prove.]
\mr   
\item "{$(*)_5$}"  if $s \le^0_{I_1} t$ then $N_s \le_{\frak K} N_t$.
\ermn
[Why?  As in the proof of $(*)_3$ by AxV of a.e.c. we have $N_s
\le_{\frak K} N_t$ (not the part on the $M$'s!)]
\mr
\item "{$(*)_6$}"  if $s \in I_1,A \in [M]^{\le \lambda}$ and $B \in
[M]^{\le \lambda}$ \ub{then} for some $t$ we have $s \le^2_{I_1} t$
and $A \subseteq M_t,B \subseteq N_t$.
\ermn
[Why?  By $(*)_2$ there is $t_1$ such that $s \le^2_I t_1,A \subseteq
M_{t_1}$ and $B \subseteq N_{t_1}$.  By $(*)_4$ there is $t \in I_1$ such that
$t_1 \le^2_I t$ hence by $(*)_0(\alpha)$ we have $s \le^2_I t$.  As $s,t
\in I_1$ this implies $s \le^2_{I_1} t$.]

Note that it is unreasonable to have ``$(I_1,\le^2_{I_1})$-directed"
but
\mr
\item "{$(*)_7$}"  $(I_1,\le^1_{I_1})$ is directed.
\ermn
[Why?  Let $s_1,s_2 \in I_1$.  We now choose $t_n$ by induction on $n <
\omega$ such that
\mr
\item "{$(a)$}"  $t_n \in I_1$
\sn
\item "{$(b)$}"  $M_{t_n}$ includes $\cup\{M_{t_k}:k<n\} \cup M_{s_1}
\cup M_{s_2}$ if $n \ge 2$
\sn
\item "{$(c)$}"   $N_{t_n}$ includes $\cup\{N_{t_k}:k<n\} \cup N_{s_1}
\cup N_{s_2}$ if $n \ge 2$
\sn
\item "{$(d)$}"  $t_0 = s_1$
\sn
\item "{$(e)$}"   $t_1 = s_2$
\sn
\item "{$(f)$}"  if $n=m+1 \ge 2$ then $t_m \le^0_{I_1} t_n$
\sn
\item "{$(g)$}"  if $n=m+2$ then $t_m \le^2_I t_n$ hence
$t_m \le^2_{I_1} t_n$.
\ermn
For $n=0,1$ this is trivial.  For $n=m+2 \ge 2$, apply $(*)_6$
with $t_m,\cup\{M_{t_k}:k \le m+1\},\cup\{N_{t_k}:k \le m+1\}$ here
standing for $s,A,B$ there getting $t_n$, so we get 
$t_n \in I_1$ in particular $t_m \le^2_{I_1} t_n$, so clause (a)
is satisfied by $t_n$. By the choice of $t_n$ and as $s_1 = t_0,s_2
= t_1$, clauses (b) + (c) hold for $t_n$.  By the choice of $t_n$, 
obviously also clause (g).  Now why does clause (f) hold
(i.e. $t_{m+1} \le^0_I t_n$)?  It follows from clauses (a),(b),(c),
so $t_n$ is as required.  Hence we have carried the induction.  
Let $N^* = \cup\{N_{t_n}:2 \le n < \omega\}$, 
so clearly by $(*)_5$ and clause (f) we have $N_{t_n}
\le_{\frak K} N_{t_{n+1}}$ for $n \ge 1$, and clearly $M_{t_n}
\subseteq M_{t_{n+1}}$ for $n \ge 1$.
Let $M^* = \cup\{M_{t_n}:2 \le n < \omega\}$.  Note that by $(*)_3$ and clause
(g) we have $M_{t_n} \le_{\frak K} M_{t_{n+2}}$, so $\langle M_{t_{n+2}}:n <
\omega\rangle$ is $\subseteq$-increasing, and for $\ell=0,1$ the
sequence $\langle M_{t_{2n + \ell}}:n < \omega\rangle$ is 
$\le_{\frak K}$-increasing with union $M^*$, hence by the
basic properties of a.e.c. we have $M_{2n + \ell} \le_{\frak K} M^*$.
So $M_{s_1} = M_{t_0} \le_{\frak K} M^*,M_{s_2} = M_{t_1} 
\le_{\frak K} M^*$.  Now $M_{s_1},M_{s_2} \subseteq M_{t_2} 
\le_{\frak K} M^*$ hence $M_{s_1},M_{s_2} \le_{\frak K} M_{t_2}$,
Recall that 
$N_{s_1} = N_{t_0} \le_{\frak K} N_{t_2}$ was proved above and $N_{s_2}
= N_{t_1} \le_{\frak K} N_{t_2}$ was also proved above so $t_2$ is a common
$\le^1_I$-upper bound of $s_1,s_2$ as required.]
\mr
\item "{$(*)_8$}"  if $s \le^0_{I_1} t$ then $s \le^1_{I_1} t$.
\ermn
[Why?  By $(*)_7$ there is $t_1 \in I_1$ which is a common
$\le^1_{I_1}$-upper bound of $s,t$.  So $M_s \subseteq M_t$ (as $s
\le^0_{I_1} t$) and $M_s \le_{\frak K} M_{t_1}$ (as $s \le^1_{I_1}
t_1$) and $M_t \le_{\frak K} M_{t_1}$ (as $t \le^1_{I_1} t_1$).
Together by axiom V of a.e.c. we get $M_s \le_{\frak K} M_t$ and
by $(*)_5$ we have $N_s \le_{\frak K} N_t$.  Together $s \le^1_{I_1}
t$ as required.]
\mr
\item "{$(*)_9$}"  $\langle M_s:s \in (I_1,\le^1_{I_1})\rangle$ is
$\le_{\frak K}$-increasing, $(I_1,\le^1_{I_1})$ is directed and
$\cup\{M_s:s \in I_1\} = M$.
\ermn
[Why?  The first phrase by the definition of $\le^1_{I_1}$ in clause
$(c)(\beta)$ of $\boxdot$, the second
by $(*)_7$ and the third by $(*)_6 + (*)_4$.]

By the basic properties of a.e.c. (see \marginbf{!!}{\cprefix{88r}.\scite{88r-1.6}}) we deduce
\mr
\item "{$\odot$}"   $(a) \quad M \in K$
\sn
\item "{${{}}$}"  $(b) \quad t \in I_1 \Rightarrow M_t \le_{\frak K} M$.
\ermn
Now we strengthen the assumption $\boxtimes_1$ to
\mr
\item "{$\boxtimes_2$}"  the demands in $\boxtimes_1$ and
$M \prec_{\Bbb L_{\infty,\kappa}[\tau_{\frak K}]} N$.
\ermn
We note
\mr
\item "{$\circledast_1$}"   $(a) \quad$ if 
$\bar a \in {}^\alpha M,|\alpha| + \text{
LS}({\frak K}) \le \lambda < \kappa$ then for some $t \in
I_\lambda,f_t(\bar a) = \bar a$
\sn
\item "{${{}}$}"  $(b) \quad$ if $M' \subseteq M$ and $\|M\| \le
\lambda$ \ub{then} $(\text{id}_{M'},M',M') \in I_\lambda$
\sn
\item "{${{}}$}"  $(c) \quad$ if $M_1 \subseteq N_1 \subseteq N$ and
$M_1 \subseteq M$ and $\|N_1\| \le \lambda$ \ub{then} for some $ t \in
I$
\nl

\hskip25pt we have $N_t = N_1$ and id$_{M_1} \subseteq f_t$.
\ermn
[Why?  Clause (a) is a special case of clause (b) and clause (b) is a
special case of clause (c).  Lastly, clause (c) follows from the
assumption $M \prec_{\Bbb L_{\infty,\kappa}[\tau_{\frak K}]} N$ 
and \scite{734-0n.13}(2A),(2B).]

We next shall prove
\mr
\item "{$\circledast_2$}"  $M \le_{\frak K} N$.
\ermn
By \marginbf{!!}{\cprefix{88r}.\scite{88r-1.6}} and $(*)_9$ above for proving
$\circledast_2$ it suffices to prove:
\mr
\item "{$\circledast_3$}"  if $s \in I_1$ then $M_s \le_{\frak K} N$.
\ermn
[Why $\circledast_3$ holds?  As $M \subseteq N$ there is $N_*
\le_{\frak K} N$ of cardinality $\le \lambda$ such that $M_s \cup N_s
\subseteq N_*$.  By $\circledast_1(c)$ there is $t \in I$ such that $N_t
= N_*$ and id$_{M_s} \subseteq f_t$.  As $N_* \le_{\frak K} N$ it
follows that $t \in I_1$.  So by $\boxtimes_1 \Rightarrow \odot(b)$
applied to $s$ and to $t$ we can deduce $M_s \le_{\frak K} M$ and $M_t
\le_{\frak K} M$.  But as id$_{M_s} \subseteq f_t$ it follows that
$M_s \subseteq M_t$ hence by AxV of a.e.c. we know that $M_s
\le_{\frak K} M_t$.  But as $t \in I$ clearly $f_t$ is an isomorphism
from $N_t$ onto $M_t$ hence $f^{-1}_t(M_s) \le_{\frak K} N_t$, and as
id$_{M_s} \subseteq f_t$ this means that $M_s = f^{-1}_t(M_s)
\le_{\frak K} N_t$.  Recalling $N_t \le_{\frak K} N$ and $\le_{\frak
K}$ is transitive it follows that $M_s \le_{\frak K} N$ as required.]

Let us check parts (3) and (4) of the Fact.  Having proved $\boxtimes_1
\Rightarrow \odot(a)$, clearly in part (4) of the fact the first
conclusion there, $M \in K$, holds.  The second conclusion, $M
\equiv_{\Bbb L_{\infty,\kappa}[{\frak K}]} N$ holds by
\mr
\item "{$\circledast_4$}"  if $\varphi(\bar x) \in \Bbb
L_{\infty,\kappa}[{\frak K}]$ and $|\ell g(\bar x)| + \text{
LS}({\frak K}) \le \lambda <\kappa$ and $t \in I$ and $\bar a \in
{}^{\ell g(\bar x)}(M_t)$ then $M \models \varphi[\bar a]
\Leftrightarrow N \models \varphi[f_t(\bar a)]$.
\ermn
[Why?  Prove by induction on the depth of $\varphi$ for all $\lambda$
simultaneously.  For $\alpha = 0$,
first for the usual atomic formulas this should be clear.  Second,
by $(*)_4$ there is
$t_1$ such that $t \le^2_I t_1 \in I_1$ hence by $\circledast_3 +$
clause (d) of $\boxdot +$ clause (b) of $\odot$ we have 
$M_{t_1} \le_{\frak K} N \wedge N_{t_1} \le_{\frak K} N \wedge   
M_{t_1} \le_{\frak K} M$ respectively.  
So if $u \subseteq \ell g(\bar x)$ then $M
\restriction \text{ Rang}(\bar a \restriction u) \le_{\frak K} M
\Leftrightarrow M \restriction \text{ Rang}(\bar a \restriction u)
\le_{\frak K} M_{t_1} \Leftrightarrow N \restriction \text{
Rang}(f(\bar a) \restriction u) 
\le_{\frak K} N_{t_1} \Leftrightarrow N \restriction
\text{ Rang}(f(\bar a) \restriction u) \le_{\frak K} N$.  So we have
finished the case of atomic formulas, i.e. $\alpha = 0$.
For $\varphi(\bar x) = (\exists \bar y)\psi(\bar x,\bar y)$ use
$(*)_2$, the other cases are obvious.]
\sn
So part (4) holds.  As for part (3), the first statement, ``$M \in K$"
holds by part (4), the second statement, $M \le_{\frak K} N$, holds by
$\circledast_2$ and the third statement, 
$M \prec_{\Bbb L_{\infty,\kappa}[{\frak K}]} N$ follows by $\circledast_1(b) +
\circledast_4$.  As we have already
noted parts (1),(2),(5),(6) and part (7) is proved as
$\circledast_4$ is proved, we are done.  \hfill$\square_{\scite{734-11.1B}}$
\enddemo
\bigskip

\proclaim{\stag{734-11.1.21} Claim}  For a limit cardinal $\kappa >
\text{\rm LS}({\frak K})$:
\nl
1) $M \prec_{\Bbb L_{\infty,\kappa}[{\frak K}]} N$ provided that
\mr
\item "{$(a)$}"  if $\theta < \kappa$ and $\theta \in (\text{\rm LS}({\frak
K}),\kappa)$ then $M \prec_{\Bbb L_{\infty,\theta}[{\frak K}]} N$
\sn
\item "{$(b)$}"  for every $\partial < \kappa$ for some 
$\theta \in (\partial,\kappa)$ we have: if $\bar
a,\bar b \in {}^\partial M$ and $(M,\bar a) 
\equiv_{\Bbb L_{\infty,\theta}[{\frak K}]} (M,\bar b)$ \ub{then}
$(M,\bar a) \equiv_{\Bbb L_{\infty,\theta_1}[{\frak K}]} (M,\bar b)$
for every $\theta_1 \in [\theta,\kappa)$.
\ermn
1A) $M \equiv_{\Bbb L_{\infty,\kappa}[{\frak K}]} N$ provided that
\mr
\item "{$(a)$}"  if {\rm LS}$({\frak K}) < \theta < \kappa$ then $M
\equiv_{\Bbb L_{\infty,\theta}[{\frak K}]} N$
\sn
\item "{$(b)$}"   as in part (1).
\ermn
2) In parts (1) and (1A) we can conclude
\mr
\item "{$(b)^+$}"  for every $\partial < \kappa$ for some $\theta \in
(\partial,\kappa)$ we have: if $\bar a,\bar b \in {}^\partial M$ and
$(M,\bar a) \equiv_{\Bbb L_{\infty,\theta}[{\frak K}]} (M,\bar b)$ \ub{then}
$(M,\bar a) \equiv_{\Bbb L_{\infty,\kappa}[{\frak K}]}(M,\bar b)$.
\ermn
3) If {\rm cf}$(\kappa) = \aleph_0$ \ub{then} $M \cong N$ when
\mr
\item "{$(a)$}"  if $\theta < \kappa$ and $\theta \in
(\text{\rm LS}({\frak K}),\kappa)$ then $M 
\equiv_{\Bbb L_{\infty,\theta}[{\frak K}]} N$
\sn
\item "{$(b)$}"  as in part (1), i.e., for every $\partial \in
(\text{\rm LS}({\frak K}),\kappa)$ for some $\theta \in
(\partial,\kappa)$ we have: if $\bar a \in {}^\partial M$ and $\bar b
\in {}^\partial N$ and $(M,\bar a) \equiv_{\Bbb
L_{\infty,\theta}[{\frak K}]} (N,\bar b)$ \ub{then} $(M,\bar a)
\equiv_{\Bbb L_{\infty,\theta_1}[{\frak K}]} (N,\bar b)$ for every
$\theta_1 \in (\theta,\kappa)$
\sn
\item "{$(c)$}"  $M,N$ have cardinality $\kappa$.
\endroster
\endproclaim
\bigskip

\demo{Proof}  1) By \scite{734-11.1B}(3) it suffices to prove $M
\prec_{\Bbb L_{\infty,\kappa}} N$, for this it suffices to apply the
criterion from \scite{734-0n.13}(2A).

Let ${\Cal F}$ be the set of functions $f$ such that:
\mr
\item "{$\odot$}"  $(\alpha) \quad$ Dom$(f) \subseteq M$ has
cardinality $< \kappa$
\sn
\item "{${{}}$}"  $(\beta) \quad$ Rang$(f) \subseteq N$
\sn
\item "{${{}}$}"  $(\gamma) \quad$ if $\bar a$ lists Dom$(f)$ then for
every $\theta \in (\ell g(\bar a),\kappa)$ we have $\sftp_{\Bbb
L_{\infty,\theta}[{\frak K}]} (\bar a,\emptyset,M) =$
\nl

\hskip25pt $\sftp_{\Bbb L_{\infty,\theta}[{\frak K}]} (f(\bar a),\emptyset,N)$.
\ermn
1A)  Similarly. 
\nl
2) Similarly to part (1) using \scite{734-11.1B}(4) and \scite{734-0n.13}(2)
instead \scite{734-11.1B}(3),\scite{734-0n.13}(2A). 
\nl
3) Recall \scite{734-0n.13}(1).  \hfill$\square_{\scite{734-11.1.21}}$
\enddemo
\bigskip

\proclaim{\stag{734-11.2} Claim}  1) Assume \scite{734-11.0.3}(a) + (b),
i.e. ${\frak K}$ is categorical in $\mu > \text{\rm LS}({\frak K})$.
If $\mu = \mu^{< \kappa}$ and $\kappa > \text{\rm LS}({\frak K})$ 
\ub{then} for every $M \le_{\frak K} N$ from $K_\mu$ we have
$M \prec_{{\Bbb L}_{\infty,\kappa}[{\frak K}]} N$ (and there are
such $M <_{{\frak K}_\mu} N$).
\nl
2) Assume ${\frak K}$ is weakly or just pseudo 
$\mu$-solvable as witnessed by $\Phi$ 
(see Definition \scite{734-11.1} and Claim \scite{734-11.1.2}) and
$M^* = \text{\rm EM}_{\tau({\frak K})}(\mu,\Phi)$.  If 
$\mu = \mu^{< \kappa}$ and $\kappa > |\tau_\Phi|$ and 
$M \le_{\frak K} N$ are both isomorphic to $M^*$ \ub{then} $M \prec_{\Bbb
L_{\infty,\kappa}[{\frak K}]} N$.
\endproclaim
\bigskip

\demo{Proof}  1) We prove by induction on $\gamma$ that for any formula
$\varphi(\bar x)$ from $\Bbb L_{\infty,\kappa}
[{\frak K}]$ of quantifier depth
$\le \gamma$ (and necessarily $\ell g(\bar x) < \kappa$) we have
\mr
\item "{$(*)$}"  if $M \le_{\frak K} N$ are from $K_\mu$ and $\bar a
\in {}^{\ell g(\bar x)}M$ then $M \models \varphi[\bar a] \Leftrightarrow
N \models \varphi[\bar a]$. 
\ermn
If $\varphi(\bar x)$ is atomic this is clear
(for the ``$\{x_i:i < i^*\}$ is the universe of a $\le_{\frak K}$-submodel",
the implication $\Rightarrow$ holds as $\le_{\frak K}$ is transitive and
the implication $\Leftarrow$ as ${\frak K}$ satisfies AxV of a.e.c.). \nl
If $\varphi(\bar x)$ is a Boolean combination of formulas for which the 
assertion was proved, clearly it holds for $\varphi(\bar x)$.  So we
are left with the case 
$\varphi(\bar x) = (\exists \bar y)\psi(\bar y,\bar x)$, so
$\ell g(\bar y) < \kappa$.  The implication $\Rightarrow$ is trivial by
the induction hypothesis and so suppose that the other fails, say $N \models
\psi[\bar b,\bar a]$ and $M \models \neg(\exists \bar y)\psi
(\bar y,\bar a)$.  We choose by induction on $i < \mu^+$ a model $M_i \in
K_\mu,\le_{\frak K}$-increasing continuous, and for each $i$ in
addition we choose
an isomorphism $f_i$ from $M$ onto $M_i$ and if $i = j+1$ we shall choose an
isomorphism $g_j$ from $N$ onto $M_{j+1}$ extending $f_j$.  For $i=0$, let
$M_0 = M$, for $i$ limit let $M_i = 
\dbcu_{j < i} M_j$.  For any $i$, if $M_i$ was chosen, $f_i$ exists
as ${\frak K}$ is categorical in $\mu$.  Now if $i=j+1$ then $M_j,f_j$
are well defined and clearly we can choose $M_i = M_{j+1},g_j$
as required.  

By Fodor lemma, as $\mu = \mu^{< \kappa}$ and the set $\{\delta <
\mu^+:\text{\rm cf}(\delta) \ge \kappa\}$ is stationary, clearly for some 
$\alpha < \beta < \mu^+$ we have $f_\alpha
(\bar a) = f_\beta(\bar a)$, now (by the choice of 
$g_\alpha$) we  have $M_{\alpha +1} \models
\psi[g_\alpha(\bar b),g_\alpha(\bar a)]$, hence by the induction hypothesis 
applied to the pair 
$(M_{\alpha +1},M_\beta)$ we have $M_\beta \models \psi[g_\alpha
(\bar b),g_\alpha(\bar a)]$ so $M_\beta \models \varphi[g_\alpha(\bar a)]$.
But $g_\alpha(\bar a) = f_\alpha(\bar a) = f_\beta(\bar a)$, contradiction
to $M \models \neg \varphi[\bar a]$.
\nl
2) The same proof but we restrict ourselves to models in $K_{[M^*]}$ so,
e.g. in $(*)$ we have $M,N \in K_{[M^*]}$ recalling 
that ${\frak K}_{[M^*]}$ is a $\mu$-a.e.c., see Definition
\scite{734-0n.5}(3A) and Claim \scite{734-0n.6}(7).
\hfill$\square_{\scite{734-11.2}}$
\enddemo
\bn
\margintag{734-11n.2.7}\ub{\stag{734-11n.2.7} Exercise}:  1) For the proof (of \scite{734-11.2}(1)) it
suffices to assume ``$S \subseteq \{\delta < \mu^+:\text{cf}(\delta)
\ge \kappa\}$ is a stationary subset of $\mu^+$ and $M^* \in K_\mu$ is
locally $S$-weakly limit (see \marginbf{!!}{\cprefix{88r}.\scite{88r-3.1}}(5)).
\nl
2) Similarly we can weaken the demands ``$M^* = \text{ EM}_{\tau({\frak
K})}(\mu,\Phi)$ and $(K,\Phi)$ is pseudo solvable" to: for every $M
\le_{\frak K} N$ isomorphic to $M^*$ (which $\in K_\mu$) there is a
$\le_{\frak K}$-increasing sequence $\langle M_\alpha:\alpha <
\mu^+\rangle$ such that $\{\delta < \mu^+$: cf$(\delta) \ge \kappa$
and $(M_\delta,M_{\delta +1})$ is isomorphic to $(M,N)$ and $M_\delta
= \cup\{M_\alpha:\alpha < \delta\}\}$ is a stationary subset of $\mu^+$.
\bigskip

\proclaim{\stag{734-11.3} Claim}  Assume 
$\Phi \in \Upsilon^{\text{or}}_{<\kappa}[{\frak K}]$ satisfies
the conclusion of \scite{734-11.2}(2) for $(\mu,\kappa)$ 
and {\rm LS}$({\frak K}) < \kappa \le \mu$ and 
$J,I_1,I_2$ are linear orders and $I_1,I_2$ are $\kappa$-wide, see
Definition \scite{734-0n.2}(1).  \ub{Then}
\mr
\item "{$(a)$}"  If $I_1 \subseteq I_2$ \ub{then} 
{\rm EM}$_{\tau({\frak K})}(I_1,\Phi) 
\prec_{{\Bbb L}_{\infty,\kappa}[{\frak K}]}$  
{\rm EM}$_{\tau({\frak K})}(I_2,\Phi)$
\sn
\item "{$(b)$}"  Assume $J \subseteq I_1,J \subseteq I_2$; if 
$\varphi(\bar x) \in \Bbb L_{\infty,\kappa}[{\frak K}]$ 
so $\ell g(\bar x) < \kappa$ and $\bar a \in {}^{\ell g(\bar x)}
({\text{\rm EM\/}}(J,\Phi))$, 
\ub{then} {\rm EM}$_{\tau({\frak K})}(I_1,\Phi) \models
\varphi[\bar a] \Leftrightarrow { \text{\rm EM\/}}_{\tau({\frak K})}(I_2,\Phi) 
\models \varphi[\bar a]$
\sn
\item "{$(c)$}"  Assume $\bar\sigma = \langle
\sigma_i(\ldots,x_{\alpha(i,\ell)},\ldots)_{\ell < \ell(i)}:i <
i(*)\rangle$ where $i(*) < \kappa$, each $\sigma_i$ is a
$\tau(\Phi)$-term, $\alpha(i,\ell) < \alpha(*) < \kappa$.  If $\bar
t^\ell = \langle t^\ell_\alpha:\alpha < \alpha(*) \rangle$ is a sequence of
members of $I_\ell$ for $\ell=1,2$ and $\bar t^1,\bar t^2$ realizes the
same quantifier free type in $I_1,I_2$ respectively and $\bar a^\ell = 
\langle \sigma_i(\ldots,a_{t^\ell_{\alpha(i,j)}},\ldots)_{j<j(i)}:i <
i(*)\rangle$ for $\ell=1,2$ \ub{then} $\bar a^1,\bar a^2$ realize the same 
$\Bbb L_{\infty,\kappa}[{\frak K}]$ -type in 
{\rm EM}$_{\tau({\frak K})}(I_1,\Phi)$, 
{\rm EM}$_{\tau({\frak K})}(I_2,\Phi)$ respectively.
\endroster
\endproclaim
\bigskip

\demo{Proof of \scite{734-11.3}}  
\sn
\ub{Clause $(a)$}:  We prove that for $\varphi(\bar x) \in \Bbb
L_{\infty,\kappa}[{\frak K}]$ we have
\mr
\item "{$(*)_{\varphi(\bar x)}$}"  if $I_1 \subseteq I_2$ are
$\kappa$-wide linear orders of cardinality $\le \mu$ 
and $\bar a \in {}^{\ell g(\bar x)}(\text{EM}_{\tau({\frak K})}
(I,\Phi))$ \ub{then} EM$_{\tau({\frak K})}
(I_1,\Phi) \models \varphi[\bar a] \Leftrightarrow \text{
EM}_{\tau({\frak K})}(I_2,\Phi) \models \varphi[\bar a]$.
\ermn
This easily suffices as for any $I \in K^{\text{lin}}$, the model
EM$_{\tau({\frak K})}(I,\Phi)$ is the direct limit of $\langle
\text{EM}(I',\Phi):I' \subseteq I$ has cardinality $\le \mu\rangle$,
which is $\le_{\frak K}$-increasing and $\mu^+$-directed and as we have:
\mr
\item "{$\odot$}"  $M^1 \prec_{\Bbb L_{\infty,\kappa}[{\frak K}]} 
M^2$ \ub{when}:
{\roster
\itemitem{ $(a)$ }   $I$ is a $\kappa$-directed partial order
\sn
\itemitem{ $(b)$ }    $\bar M = \langle M_t:t \in I\rangle$
\sn
\itemitem{ $(c)$ }    $s <_I t \rightarrow M_s 
\prec_{\Bbb L_{\infty,\kappa}[{\frak K}]} M_t$
\sn
\itemitem{ $(d)$ }    $M^2 = \cup\{M_t:t \in I\}$
\sn
\itemitem{ $(e)$ }  $M^1 \in \{M_t:t \in I\}$ or for some
$\kappa$-directed $I' \subseteq I$ we have $M^1 = \cup\{M_t:t \in I'\}$.
\endroster}
\ermn
We prove $(*)_{\varphi(\bar x)}$ by induction on $\varphi$ (as in the
proof of \scite{734-11.2} above).  
The only non-obvious case is $\varphi(\bar x) = 
(\exists \bar y)\psi(\bar y,\bar x)$, so let $I_1 \subseteq I_2$ be
$\kappa$-wide linear orders of cardinality $\le \mu$ and
$\bar a \in {}^{\ell g(\bar x)}(\text{EM}_{\tau({\frak
K})} (I_1,\Phi))$.  Now if {\rm EM}$_{\tau({\frak K})}(I_1,\Phi) \models 
\varphi[\bar a]$ then for some $\bar b \in {}^{\ell g(\bar y)}
(\text{EM}_{\tau({\frak K})}(I_1,\Phi))$ we have 
EM$_{\tau({\frak K})}(I_1,\Phi) \models \psi[\bar b,\bar a]$ hence by
the induction hypothesis EM$_{\tau({\frak K})}(I_2,\Phi) \models
\psi[\bar b,\bar a]$ hence by the satisfaction definition
EM$_{\tau({\frak K})}(I_2,\Phi) \models \psi[\bar a]$, so we have
proved the implication $\Rightarrow$.

For the other implication assume that $\bar b \in {}^{\ell g(\bar y)}
(\text{EM}_{\tau({\frak K})}(I_2,\Phi))$ and EM$_{\tau({\frak
K})}(I_2,\Phi) \models \psi[\bar b,\bar a]$.  Let $\theta = |\ell
g(\bar a \char 94 \bar b)| + \aleph_0$, so $\theta < \kappa$ and
\wilog \, if $\kappa$ is singular then $\theta \ge \text{ cf}(\kappa)$.
Hence there is in $I_1$ a monotonic sequence 
$\bar c = \langle c_i:i < \theta^+ \rangle$, \wilog \, it is
increasing.  Clearly there is $I^*$
such that $\bar a \char 94 \bar b \in {}^{\ell g(\bar x \char 94
\bar y)}(\text{EM}(I^*,\Phi)),I^*
\subseteq I_2,|I^*| \le \theta$ and $\bar a \in {}^{\ell g(\bar x)}
(\text{ EM}(I^* \cap I_1,\Phi))$ and 
\wilog \, $i < \theta^+ \Rightarrow [c_0,c_i]_{I_2} 
\cap I^* = \emptyset$.  

Similarly without loss of generality
\mr
\item "{$(*)$}"  $I_1 \backslash \cup\{[c_0,c_i)_{I_1}:i <
\theta^+\}$ is $\kappa$-wide \ub{or} $\kappa = \theta^+$.
\ermn
Let $J_0 = I_2$; we can find $J_1$ such that $J_0 = I_2 \subseteq J_1$ and
$J_1 \backslash I_2 = \{d_\alpha:\alpha < \mu \times \theta^+\}$ 
with $d_\alpha$ being $<_{J_1}$-increasing with $\alpha$ and 
$(\forall x \in I_2)(x <_{J_1} d_\alpha \equiv
\dsize \bigvee_{i < \theta^+} x <_{J_1} c_i)$.

As EM$_{\tau({\frak K})}(I_2,\Phi) \models \psi[\bar b,\bar a]$ 
and $I_2 = J_0 \subseteq J_1,|J_1| \le \mu$ and $I_2$ is
$\kappa$-wide (and trivially $J_1$ is $\kappa$-wide), by 
the induction hypothesis EM$_{\tau({\frak K})}(J_1,\Phi) \models
\psi[\bar b,\bar a]$ hence EM$_{\tau({\frak K})}(J_1,\Phi) \models
\varphi[\bar a]$.  Let $J_2 = J_1 \restriction \{x:x \in J_1
\backslash J_0$ or $x \in I_1 \backslash \cup\{[c_0,c_i]_{I_1}:i <
\theta^+\}\}$.  So $J_1 \supseteq J_2$, both linear orders 
have cardinality $\mu$ and are $\kappa$-wide as witnessed by $\langle
d_\alpha:\alpha < \mu \times \theta^+\rangle$ for both hence
the conclusion of \scite{734-11.2} holds, i.e. EM$(J_2,\Phi) \prec_{\Bbb
L_{\infty,\kappa}[{\frak K}]} \text{ EM}(J_1,\Phi)$.  
Also $I^* \cap I_1 \subseteq
J_2$ and recall that $\bar a \in {}^{\ell g(\bar x)}(\text{EM}(I^*
\cap I_1,\Phi))$ hence $\bar a \in {}^{\ell g(\bar
x)}(\text{EM}(J_2,\Phi))$.  However, EM$_{\tau({\frak K})}(J_1,\Phi)
\models \varphi[\bar a]$, see above, hence by the last two sentences
EM$_{\tau({\frak K})}(J_2,\Phi) \models \varphi[\bar a]$.

So there is $\bar b^* \in {}^{\ell g(\bar y)}(\text{EM}_{\tau({\frak
K})}(J_2,\Phi))$ such that EM$_{\tau({\frak K})}(J_2,\Phi) \models
\psi[\bar b^*,\bar a]$. Let $J^* \subseteq J_2$ be of cardinality $\theta$
such that $\bar b^* \in {}^{\ell g(\bar y)}(\text{EM}_{\tau({\frak K})}
(J^*,\Phi))$ and $I^* \cap I_1 \subseteq J^*$ recalling $I^* \cap
[c_0,c_i)_{I_2} = \emptyset$ for $i < \theta^+$.  Now let 
$u \subseteq \mu \times \theta^+$ be such that $J^* \backslash
I_1 = \{d_\alpha:\alpha \in u\}$ so $|u| < \theta^+$.  Let $J_3
= J_2 \restriction \{t:t \in J_2 \cap I_1$ or $t = d_\alpha \wedge
\alpha  > \sup(u)$ or $t = d_\alpha \wedge \alpha \in u\}$; as cf$(\mu
\times \theta^+) = \theta^+ > |u|$, clearly sup$(u) < \mu \times
\theta^+$ hence $|J_3| = \mu$ and $J_3$ is $\kappa$-wide.  So by the
conclusion of \scite{734-11.2} (or by the induction hypothesis) also 
EM$_{\tau({\frak K})}(J_3,\Phi) \models \psi[\bar b^*,\bar a]$.  Let
$w = \{\alpha < \mu \times \theta^+:\alpha \in u$ or $\alpha > \sup(u)
\wedge (\alpha - \text{\rm sup}(u) < \theta^+)\}$, so otp$(w) = \theta^+$.

Let $J_4 = (J_3 \cap I_1) \cup \{d_\alpha:\alpha \in w\}$, so $J_4$ is
$\kappa$-wide as witnessed by $I_1  \backslash \cup\{[c_0,c_i):i <
\theta^+\}$ or by $\{d_\alpha:\alpha \in w\}$ recalling $(*)$ above
and $J_4 \subseteq J_3$ and $J^* \subseteq J_4$ hence
$\bar a,\bar b^* \subseteq {}^{\kappa >}(\text{EM}(J_4,\Phi))$
hence by the induction hypothesis EM$_{\tau({\frak K})}(J_4,\Phi)
\models \psi[\bar b^*,\bar a]$.

Let $J_5 = J_4 \cup \{c_i:i < \theta^+\} \backslash
\{d_\alpha:\alpha \in w\}$ equivalently $J_5 = (J_3 \cap I_1) \cup
\{c_\alpha:\alpha < \theta^+\} = (I_1 \backslash
\cup\{[c_0,c_i)_{I_1}:i < \theta^+\}) \cup \{c_i:i < \theta^+\}$ so
$J_5 \subseteq I_1$ and let $h:J_4 \rightarrow J_5$ be such that
$h(d_\alpha) = c_{\text{otp}(w \cap \alpha)}$ for $\alpha \in
w$ and $h(t) = t$ for others, i.e. for $t \in J_3 \cap I_1$. 
So $h$ is an isomorphism from $J_4$ onto $J_5$.  Recalling \scite{734-0n.em.7}
let $\hat h$ be the isomorphism from
EM$(J_4,\Phi)$ onto EM$(J_5,\Phi)$ which $h$ induces, so clearly $\hat
h(\bar a) = \bar a$.   Hence for some $\bar b^{**}$ we have $\bar b^{**} = \hat
h(\bar b^*) \in {}^{\ell g(\bar y)}(\text{EM}_{\tau({\frak K})}
(J_5,\Phi))$ and EM$_{\tau({\frak K})}(J_5,\Phi) \models \psi[\bar
b^{**},\bar a]$.  Note that by the choice of
$\langle c_i:i < \theta^+\rangle$, (see $(*)$ above), we know that
$J_5$ is $\kappa$-wide.  Also $J_5 \subseteq I_1$ so by 
the induction hypothesis applied to
$\psi(\bar y,\bar x),J_5,I_1$ we 
have EM$_{\tau({\frak K})}(I_1,\Phi) \models
\psi[\bar b^{**},\bar a]$ hence by the definition of satisfaction
EM$_{\tau({\frak K})}(I_1,\Phi) \models \varphi[\bar a]$, so we have
finished proving the implication $\Leftarrow$ hence clause $(a)$.
\mn
\ub{Clause (b)}:  Without loss of generality for some linear order $I$
we have $I_1 \subseteq I,I_2 \subseteq I$ and EM$(I_\ell,\Phi)
 \subseteq \text{ EM}(I,\Phi)$ for $\ell =1,2$ and use clause (a) twice. 
\mn
\ub{Clause (c)}:  Easy by now, e.g. using a linear order $I'$
extending $I_1,I_2$ which has an automorphism $h$ such that $h(t^1_\alpha)
= t^2_\alpha$ for $\alpha < \alpha(*)$.  \hfill$\square_{\scite{734-11.3}}$
\enddemo
\bigskip

\definition{\stag{734-11.4} Definition}  Fixing $\Phi \in 
\Upsilon^{\text{or}}_{\frak K}$.
\nl
1) For $\theta \ge \text{ LS}
({\frak K})$ let $K^*_\theta$, [let $K^{**}_\theta$] [let
$K^{*,*}_\theta$] be the family of
$M \in K_\theta$ isomorphic to some EM$_{\tau({\frak K})}(I,\Phi)$ 
where $I$ is a 
linear order of cardinality $\theta$ [which is $\theta$-wide][which
$\in K^{\text{flin}}_\theta$]. More
accurately we should write
$K^*_{\Phi,\theta},K^{**}_{\Phi,\theta},K^{*,*}_{\Phi,\theta}$; 
similarly below.
\nl
2) Let $K^*$ is the 
class $\cup\{K^*_\theta:\theta$ a cardinal 
$\ge \text{ LS}({\frak K})\}$, similarly $K^{\ast,\ast},K^*_{\ge
\lambda},K^{**}_{\ge \lambda}$, etc. 
\nl
3) Let ${\frak K}^* = {\frak K}^*_\Phi = 
(K^*,\le_{\frak K} \restriction K^*)$.
\nl
4) Let ${\frak K}^*_\lambda = K^*_{\Phi,\lambda}$ be
$(K^*_{\Phi,\lambda},\le_{\frak K} \restriction
K^*_{\Phi,\lambda})$.
\enddefinition
\bigskip

\proclaim{\stag{734-11.4A} Claim}  1) $K^{**}_\theta$ is categorical in $\theta$
if {\rm LS}$({\frak K}) < \theta \le \mu$,
{\rm cf}$(\theta) = \aleph_0$ and the conclusion of \scite{734-11.2}(2)
hence of \scite{734-11.3} holds for $\partial = \theta$ (and $\Phi$), e.g. 
${\frak K}$ is pseudo solvable in $\mu$ as witnessed by $\Phi$ and
$\mu = \mu^{< \theta}$.
\nl
2) $K^{*,*}_\theta,K^{**}_\theta \subseteq K^*_\theta$.
\nl
3) If $\theta$ is strong limit $> \text{\rm LS}({\frak K})$ \ub{then}
$K^{**}_\theta = K^*_\theta$.
\endproclaim
\bigskip

\demo{Proof}  1) By \scite{734-11.3} and \scite{734-0n.13}(1).
\nl
2) Read the definitions.
\nl
3) Recall \scite{734-0n.3}(2).  \hfill$\square_{\scite{734-11.4A}}$
\enddemo
\bigskip

\remark{\stag{734-11f.4B} Remark}  1) We will be specially interested in 
\scite{734-11.4A} in the case $(\mu,\lambda)$ is a ${\frak K}$-candidate
(see Definition \marginbf{!!}{\cprefix{734}.\scite{734-11.0.3}}) and $\theta = \lambda$.
\nl
2) Note that $K^*_\theta$ in general is not a $\theta$-a.e.c.
\nl
3) If we strengthen \scite{734-11.2A}(2) below, replacing
$(\mu,\lambda)$ by $(\mu,\lambda^+)$ then categoricity of
$K^*_\lambda$ and in fact Claim \scite{734-11.5}(4) follows immediately
from (or as in) Claim \scite{734-11.4A}(1).
\endremark
\bn
For the rest of this section we assume that the
triple $(\mu,\lambda,\Phi)$ is a pseudo ${\frak K}$-candidate 
(see Definition \scite{734-11.0.3})
and rather than $\mu = \mu^\lambda$ we assume just the conclusion of
\scite{734-11.2}, that is:
\demo{\stag{734-11.2A} Hypothesis}  1) The pair $(\mu,\lambda)$ is a pseudo 
${\frak K}$-candidate and $\Phi$ witnesses this,  so 
$|\tau_\Phi| \le \text{ LS}({\frak K}) <
\lambda = \beth_\lambda < \mu$ and $\Phi \in
\Upsilon^{\text{or}}_{\frak K}$ is as in Definition \scite{734-11.1} so
$I \in K^{\text{lin}}_\mu \Rightarrow 
\text{\rm EM}_{\tau({\frak K})}(I,\Phi) \in K^{\text{pl}}_\mu$.
\nl
2)  For every $\kappa \in (\text{LS}({\frak
K}),\lambda)$ the conclusion of \scite{734-11.2}(2) holds hence also of
\scite{734-11.3} (if $\mu = \mu^{< \lambda}$ this follows from (1) even 
for $\kappa = \lambda^+$ as $\mu^{< \kappa} = \mu^\lambda = \mu$ by cardinal
arithmetic).
\enddemo
\bigskip

\proclaim{\stag{734-11.5} Claim}  1) If $M_1 \le_{\frak K} M_2$ are from
$K^*_\lambda$ or just $K^*_{\ge \lambda}$ and {\rm LS}$({\frak K}) <
\theta < \lambda$ \ub{then} 
$M_1 \prec_{{\Bbb L}_{\infty,\theta}[{\frak K}]} M_2$; moreover $M_1
\prec_{\Bbb L_{\infty,\lambda}[{\frak K}]} M_2$.  
\nl
2) If $M_1 \le_{\frak K} M_2$ are from $K^*$ and 
$\|M_1\| \ge \kappa := \beth_{1,1}(\theta)$ (recall that this is
$\beth_{(2^\theta)^+}$) and $\mu > \theta \ge 
{ \text{\rm LS\/}}({\frak K})$ \ub{then} 
$M_1 \prec_{{\Bbb L}_{\infty,\theta^+}[{\frak K}]} M_2$.
\nl
3) Assume {\rm LS}$({\frak K}) < \theta < \kappa = \beth_{1,1}(\theta)
\le \chi < \mu,\chi_1 = \beth_{1,1}(\chi)$ and 
$M \in K^*_{\ge \chi_1}$ and $\bar a,\bar b \in {}^\gamma M$ where
$\gamma < \theta^+$ and $(M,\bar a) 
\equiv_{\Bbb  L_{\infty,\kappa}[{\frak K}]} (M,\bar b)$, i.e. $\varphi(\langle
 x_\beta:\beta < \gamma\rangle) \in \Bbb L_{\infty,\kappa^+}
[{\frak K}] \Rightarrow M \models \varphi[\bar a] \Leftrightarrow M \models
\varphi[\bar b]$.  \ub{Then} $(M,\bar a) 
\equiv_{\Bbb L_{\infty,\chi}[{\frak K}]} (M,\bar b)$.
\nl
4) $K^*_\lambda$ is categorical in $\lambda$ provided that 
{\rm cf}$(\lambda) = \aleph_0$.
\endproclaim
\bigskip

\remark{\stag{734-11.5.13} Remark}  1) What is the difference between say
\scite{734-11.5}(3) and clause (a) of \scite{734-11.3}?  
Here there is no connection between
the additional $\tau(\Phi)$-structures expanding $M_1,M_2$.
\nl
2) Note that $\Phi$ has the
$\kappa$-non-order property (see \scite{734-11.1.2}(2)(*)) when $\kappa
\ge \text{\rm LS}({\frak K}),\kappa^+ < \mu$ using \scite{734-11.5}(4).
\nl
3) Concerning \scite{734-11.5}(2), note that if $\|M_1\| \ge \mu$ it is 
easy to deduce this from
\scite{734-11.2A}(2), i.e, \scite{734-11.2}(2).  But the whole point in this
stage is to deduce something on cardinals $< \mu$.  
\nl
4) Note that the proof of \scite{734-11.5}(2) gives:
\mr
\item "{$\circledast$}"  assume LS$({\frak K}) \le \theta$ and 
$\delta(*) = \text{ Min}\{(2^\theta)^+,\delta(2^{\text{LS}({\frak K})} 
+ \theta)\}$ where on the function $\delta(-)$, see
\marginbf{!!}{\cprefix{300a}.\scite{300a-1.2.3}},\marginbf{!!}{\cprefix{300a}.\scite{300a-1.2}}, if $\beth_{\delta(*)} 
\le \mu$ then for some $\alpha(*) < \delta(*)$ we have:
{\roster
\itemitem{ $\odot$ }  if $M_1 \le_{\frak K} M_2$ are from $K^*$ and
$\|M_1\| \ge \beth_{\alpha(*)}$ then $M_1 \prec_{\Bbb
L_{\infty,\theta^+}[{\frak K}]} M_2$.
\endroster}
\ermn
5) Similarly for \scite{734-11.5}(3) so we can weaken the 
demand $M \in K^*_{\ge \chi_1}$
\nl
6) We use ``$\lambda$ has countable cofinality, i.e. cf$(\lambda) =
\aleph_0$" in the proof of part (4) of \scite{734-11.5}, but not in the
proof of the other parts. 
\nl
7) Recall that for notational simplicity we assume LS$({\frak K}) \ge
|\tau_{\frak K}|$ hence $\theta \ge |\tau_\Phi|$.
\nl
8) Note that for \scite{734-11.5}(2),(3) we can omit $\lambda$ from
Hypothesis \scite{734-11.2A}.
\nl
9) Note that we shall use not only \scite{734-11.5} but also its proof.
\endremark
\bigskip

\demo{Proof of \scite{734-11.5}}  1) The first phrase holds by part (2)
noting that $\kappa < \lambda$ if $\theta < \lambda$ as $\theta < 
\lambda = \beth_\lambda$.  The second phrase holds by \scite{734-11.1.21}
as its assumption holds by parts (1) and (3).
\nl
2) We prove by induction on the ordinal $\gamma$ that:
\mr
\item "{$(*)$}"  if $M_1 \le_{\frak K} M_2$ are from 
$K^*_{\ge \kappa}$ and the formula $\varphi(\bar x) \in 
\Bbb L_{\infty,\theta^+}[{\frak K}]$ has depth $\le \gamma$ 
(so necessarily $\ell g(\bar x) < \theta^+$) 
and $\bar a \in {}^{\ell g(\bar x)}(M_1)$ 
\ub{then} $M_1 \models \varphi[\bar a] \Leftrightarrow M_2 \models 
\varphi[\bar a]$.
\ermn
As in \scite{734-11.2}, the non-trivial case is to assume $\varphi(\bar x) =
(\exists \bar y)\psi(\bar y,\bar x)$ where $\bar a \in {}^{\ell g(\bar
x)}(M_1)$ and $M_2 \models \varphi[\bar a]$ and we shall prove $M_1
\models \varphi[\bar a]$, so necessarily 
$\ell g(\bar x) + \ell g(\bar y) < \theta^+$ and we can choose $\bar b
\in {}^{\ell g(\bar y)}(M_2)$ such that 
$M_2 \models \psi[\bar b,\bar a]$.
For $\ell=1,2$ as $M_\ell \in K^*_{\ge \kappa}$ there is an 
isomorphism $f_\ell$ from EM$_{\tau({\frak K})}(I_\ell,\Phi)$  onto
$M_\ell$ for some linear order $I_\ell$ of cardinality $\ge \kappa$.

So we can 
find $J_\ell \subseteq I_\ell$ of cardinality $\theta$ for $\ell =1,2$
such that $\bar a \subseteq M^-_1$ where $M^-_1 = 
f_1(\text{EM}_{\tau({\frak K})}(J_1,\Phi))$, 
and $\bar a \char 94 \bar b \subseteq M^-_2$ where $M^-_2 = f_2
(\text{EM}_{\tau({\frak K})}(J_2,\Phi))$ and
without loss of generality 
$M^-_1 = M^-_2 \cap M_1$.  
By \scite{734-11.2A}(1), i.e. \scite{734-0n.8}(1), clause (c) 
clearly $M^-_\ell \le_{\frak K} M_\ell$ and so by AxV of a.e.c. 
(see Definition \marginbf{!!}{\cprefix{600}.\scite{600-0.2}}), we have $M^-_1 \le_{\frak K} M^-_2$.  
First assume $\theta \ge 2^{\text{LS}({\frak K})}$; in fact it is not
a real loss to assume this.
By renaming 
\wilog \, there is a transitive
set $B$ (in the set theoretic sense) of 
cardinality $\le \theta$ such that the following objects belong to it:
\mr
\item "{$\oplus(a)$}"  $J_1,J_2$
\sn
\item "{$(b)$}"  $\Phi$ (i.e. 
$\tau_\Phi$ and $\langle (\text{EM}(n,\Phi),a_\ell)
_{\ell < n}:n < \omega \rangle)$
\sn
\item "{$(c)$}"  ${\frak K}$, i.e., $\tau_{\frak K}$ and
$\{(M,N):M \le_{\frak K} N$ have universe included in 
$\text{\rm LS}({\frak K})\}$
\sn
\item "{$(d)$}"  EM$(J_\ell,\Phi)$ and $\langle a_t:t \in
J_\ell\rangle$ for $\ell=1,2$.
\ermn
Let $\chi$ be large enough, ${\frak B} = ({\Cal H}(\chi),\in,<^*_\chi)$ and
${\frak B}^+$ be ${\frak B}$ expanded by the individual constants
$M^+_\ell =  \text{ EM}(I_\ell,\Phi),\langle
a^\ell_t:t \in I_\ell \rangle$ the skeleton, $M_\ell,M^-_\ell$ 
and $f_\ell$ (all for
$\ell = 1,2$), $\kappa,B$ and $x$ for each $x \in B$.  By the
assumption $\|M_1\| \ge \kappa = \beth_{1,1}(\theta)$, hence 
(see here \marginbf{!!}{\cprefix{300a}.\scite{300a-1.2}}) there is ${\frak C}$ such that
\mr
\item "{$\odot$}"  $(a) \quad {\frak C}$ is  
a $\tau({\frak B}^+)$-model elementarily equivalent to ${\frak B}^+$ 
\nl

\hskip25pt (that is, in first order logic)
\sn
\item "{${{}}$}"  $(b) \quad {\frak C}$ omits the
type $\{x \ne b \and x \in B:b \in B\}$ but
\sn
\item "{${{}}$}"  $(c) \quad |\{b:{\frak C} \models ``b \in
\kappa^{\frak C}"\}| = \mu = \|{\frak C}\|$.  
\ermn
Without loss of generality $b \in B \Rightarrow b^{\frak C} = b$.

Now
\mr
\item "{$\circledast_1$}"  if ${\frak C} \models
``M \in K"$, so $M$ is just a member of the model ${\frak C}$ \ub{then}
we can define a $\tau_{\frak K}$-model $M^{\frak C} = M[{\frak C}]$ as follows
{\roster
\itemitem{ $(a)$ }  the set of elements of $M^{\frak C}$ is
$\{a:{\frak C} \models ``a$ is a member of the model $M"\}$
\sn
\itemitem{ $(b)$ }  if $R \in \tau_K$ is an $n$-place predicate then
$R^{M[{\frak C}]} = \{\langle a_\ell:\ell < n\rangle:{\frak C}
\models ``\langle a_\ell:\ell < n \rangle \in R^M"\}$
\sn
\itemitem{ $(c)$ }   if $F \in \tau_K$ is an $n$-place function symbol,
$F^{M[{\frak C}]}$ is defined similarly.
\endroster}
\item "{$\circledast_2$}"  $(a) \quad$ if ${\frak C} \models
``I$ is a linear order" then we define $I^{\frak C}$ similarly
\sn
\item "{${{}}$}"  $(b) \quad$ similarly if ${\frak C} \models
``M$ is a $\tau(\Phi)$-model"
\sn
\item "{$\circledast_3$}"   if ${\frak C} \models
``I$ is a directed partial order, $\bar M = \langle M_s:s \in I\rangle$
satisfies $M_s \in K$ has cardinality LS$({\frak K})$ and $s \le_I t
\Rightarrow M_s \le_{\frak K} M_t$" \ub{then} also $\langle M^{\frak
C}_s:s \in I^{\frak C}\rangle$ satisfies this.
\ermn
By easy absoluteness (for clauses $(a)_1,(a)_2$ we use
\marginbf{!!}{\cprefix{88r}.\scite{88r-1.6}}, \marginbf{!!}{\cprefix{88r}.\scite{88r-1.7}} and $\circledast_3$):
\mr
\item "{$\boxtimes$}"  $(a)_1 \quad$ if ${\frak C} \models ``M \in K"$
then $M^{\frak C} \in K$
\sn
\item "{${{}}$}"  $(a)_2 \quad$ if ${\frak C} \models 
``M \le_{\frak K} N"$ then $M^{\frak C} \le_{\frak K} N^{\frak C}$
\sn
\item "{${{}}$}"  $(b)_1 \quad$ if ${\frak C} \models 
``I$ is a linear order" then $I^{\frak C} = I[{\frak C}]$ is a linear order
\sn
\item "{${{}}$}"  $(b)_2 \quad$ if ${\frak C} \models 
``I \subseteq J$ as linear orders" then $I^{\frak C} \subseteq J^{\frak C}$
\sn
\item "{${{}}$}"  $(c) \quad$ similarly for $\tau_\Phi$-models
\sn
\item "{${{}}$}"  $(d)_1 \quad$ if ${\frak C} \models 
``M = \text{ EM}(I,\Phi)"$ \ub{then} there is a canonical isomorphism
$f^{\frak C}_I$
\nl

\hskip25pt  from EM$(I^{\frak C},\Phi)$ onto $M^{\frak C}$ (hence
it is also an isomorphism from
\nl

\hskip25pt  EM$_{\tau({\frak K})}(I^{\frak C},\Phi)$
onto $M^{\frak C} \restriction \tau({\frak K})$)
\sn
\item "{${{}}$}"  $(d)_2 \quad$ if ${\frak C} \models 
``I \subseteq J$ as linear orders" then $f^{\frak C}_J$ extends
$f^{\frak C}_I$.
\ermn
Now clearly $J^{\frak C}_\ell = J_\ell$ and $I^{\frak C}
_\ell$ is a linear order of cardinality $\mu$ extending $J_\ell$ for
$\ell=1,2$.  Let $M^*_\ell = (M^-_\ell)^{\frak C}$ for $\ell = 1,2$.

So recalling clause (c) of $\odot$ we have: 
$M^{\frak C}_1,M^{\frak C}_2 \in K^*_\mu,M^{\frak C}_1 
\le_{\frak K} M^{\frak C}_2,M^*_\ell \le_{\frak K} M^{\frak C}_\ell,
M^*_1 \le_{\frak K} M^*_2$ and $f^{{\frak C}_0}_\ell,
f^{\frak C}_{I_\ell}$ are 
isomorphisms from EM$_{\tau({\frak K})}(I^{\frak C}_\ell,\Phi)$ onto 
$M^{\frak C}_\ell$, in fact, $f^{\frak C}_{I_\ell}$ is the identity on
EM$_{\tau({\frak K})}(J^{\frak C}_\ell,\Phi) = \text{ EM}_{\tau({\frak
K})}(J_\ell,\Phi)$ and $f^{\frak C}_\ell$ maps it onto $M^*_\ell$ 
for $\ell=1,2$.

Now $M_2 \models \psi[\bar a,\bar b]$, (why?  assumed above) hence
$M^{\frak C}_2 \models \psi[\bar a,\bar b]$
\nl
(why?  By \scite{734-11.3}, clause (b) or (c) and the situation recalling
\scite{734-11.2A}(2), of course noting that $I_2,I^{\frak C}_2$ are of
cardinality $\ge \kappa = \beth_{1,1}(\theta)$ hence are
$\theta^+$-wide), hence $M^{\frak C}_2 \models \varphi[\bar a]$ (by 
definition of satisfaction), hence $M^{\frak C}_1 \models \varphi[\bar a]$
(why?  as $M^{\frak C}_1,M^{\frak C}_2 \in K^*_\mu$ hence 
$M^{\frak C}_1 \prec_{{\Bbb L}_{\infty,\theta^+}
[{\frak K}]} M^{\frak C}_2$ by $\boxtimes$ and \scite{734-11.2A}(2) 
and recalling \scite{734-11.2}(2)) 
hence $M_1 \models \varphi[\bar a]$ 
(why?  by clause (b) of \scite{734-11.3} recalling \scite{734-11.2A}(2)) 
as required in \scite{734-11.5}(2).] 
\nl
So we are done except for a small debt: the case $\theta <
2^{\text{LS}({\frak K})}$ and $f^{\frak C}_\ell$ is an isomorphism
from EM$_{\tau({\frak K})}(I^{\frak C}_\ell,\Phi)$.

In this case choose two sets $B_1,B_2$ such that $|B_1| =
\theta,|B_2| = 2^{\text{LS}({\frak K})},B_1 \subseteq B_2$ and
concerning the demands in $\oplus$ above the objects from (a),(b),(d)
and $\tau_{\frak K}$ belong to $B_1$, the objects from (c) belong to $B_2$.

Again, \wilog \, $B_1,B_2$ are transitive sets and $B_1,B_2$ serve as
individual constants of ${\frak B}^+$ as well as each member of $B_1$.
Now concerning ${\frak C}$ we
demand that it is elementarily equivalent to ${\frak B}^+$; omit $\{x
\in B_1 \wedge x \ne b:b \in B_1\}$ and for some ${\frak B}^+_1 \prec
{\frak B}^+$ of cardinality $\theta$ we have ${\frak B}^+_1 \prec
{\frak C}$ and $\{b:{\frak C} \models b \in B_2\} \subseteq {\frak B}^+$.  
This influences just the proof of $\circledast_3$.
\nl
3) Without loss of generality $M = \text{\rm EM}_{\tau({\frak K})}
(I,\Phi)$ and $I \in K^{\text{lin}}_{\ge \chi_1}$.  As $\gamma <
\theta^+$ and $\bar a,\bar b \in {}^\gamma M$ there is $I_1 \subseteq
I$ of cardinality $\theta$ such that $\bar a,\bar b \in
{}^\gamma(M_1)$ where $M_1 = \text{ EM}_{\tau({\frak K})}(I_1,\Phi)$.
As $(M,\bar a) \equiv_{\Bbb L_{\infty,\kappa^+}[{\frak K}]} (M,\bar b)$  
necessarily there is $I_2 \subseteq I$ of cardinality $\kappa$ and
 automorphism $f$ of $M_2 = \text{ EM}_{\tau({\frak K})}(I_2,\Phi)$
 mapping $\bar a$ to $\bar b$ such that $I_1 \subseteq I_2$.  
Why?  Recalling \scite{734-0n.13}(2), by the hence
and forth argument as in the second part of the
 proof of \scite{734-11.1B}(3).
\nl
Now as in the proof of part (2) there is a linear order $I_3$
extending $I_1$ of cardinality $\chi_1$ and an automorphism $g$ of
$M_3 = \text{ EM}_{\tau({\frak K})}(I_3,\Phi)$ mapping $\bar a$ to
$\bar b$.  Without loss of generality for some linear order $I_4$ we
have $I \subseteq I_4$ and $I_3 \subseteq I_4$.

Let $M_4 = \text{ EM}_{\tau({\frak K})}(I_4,\Phi)$, now 
$M \prec_{\Bbb L_{\infty,\chi^+}[{\frak K}]} M_4$ by part (2), 
$M_3 \prec_{\Bbb L_{\infty,\chi^+}[{\frak K}]} M_4$ by part (3) and
$(M_3,\bar a) \equiv_{\Bbb L_{\infty,\chi^+}[{\frak K}]} (M_3,\bar b)$
by using the automorphism $g$ of $M_3$ so together we are done.
\nl
4) So let $M,N \in K^*_\lambda$ (in fact, hence $\in K^{**}_\lambda$ 
recalling $K^*_\lambda = K^{**}_\lambda$ by \scite{734-11.4A}(3) but not used).  
By parts (1),(3) the assumptions of
\scite{734-11.1.21}(3) holds with $\lambda$ here standing for $\kappa$
there, hence its conclusion, i.e. $M \cong N$.  
\nl
${{}}$  \hfill$\square_{\scite{734-11.5}}$
\enddemo
\bn
Note: here the types below are sets of formulas.
\definition{\stag{734-11.7} Definition}  Assume $M \in K,\bold I \subseteq
{}^\gamma M$ and ${\Cal L},{\Cal L}_1,{\Cal L}_2$ are languages
in the vocabulary $\tau_{\frak K}$.
\nl
1) We say that $\bold I$ is 
$({\Cal L},\partial,< \kappa)$-convergent in $M$,
if: $|\bold I| \ge \partial$ and for every $\bar b \in {}^{\kappa >} M$,
for some $\bold J \subseteq \bold I$ of cardinality $< \partial$ for 
some \footnote{We could have demanded it for every single formula,
here this distinction is not important} $p$ we have:
\mr
\item "{$(*)$}"   for every $\bar c \in \bold I 
\backslash \bold J$, the ${\Cal L}$-type of 
$\bar c \char 94 \bar b$ in $M$ is $p$. 
\ermn
2) Let Av$_{{\Cal L},\partial,< \kappa}(\bold I,M) = 
\{\varphi(\bar x,\bar b):\varphi(\bar x,\bar y)$ is an 
${\Cal L}$-formula, $\ell g(\bar y) < \kappa$ and $\bar a
\in \bold I \Rightarrow \ell g(\bar a) = \ell g(\bar x)$ and
$\bar b \in {}^{\ell g(\bar y)} M$ and for all but 
$< \partial$ of the sequences $\bar c \in \bold I$, the
sequence $\bar c$ satisfies $\varphi(\bar x,\bar b)$ in $M\}$.  If $\partial$
is missing, we mean $\partial = \kappa$.
In parts (1) and (2) we may write ``$\kappa$" instead of $<
\kappa^+$; similarly below.
\nl
3) We say that $\bold I$ is $({\Cal L}_1,
{\Cal L}_2,\partial,<\kappa)$-based on $A$ in $M$ (if ${\Cal L}_1 =
{\Cal L} = {\Cal L}_2$ we may write only ${\Cal L}$) \ub{when}:
\mr
\item "{$(a)$}"  $A \subseteq M$
\sn
\item "{$(b)$}"  $\bold I$ is $({\Cal L}_1,\partial,<\kappa)$-convergent,
\sn
\item "{$(c)$}"  Av$_{{\Cal L}_1,\partial,< \kappa}(\bold I,M)$ does not
$({\Cal L}_1,{\Cal L}_2,<\kappa)$-split over $A$, see below.
\ermn
4) We say that $p(\bar x) \in \Sfor^\alpha_{\Cal L}(B,M)$ does not
$({\Cal L}_1,{\Cal L}_2,<\kappa)$-split over $A$ \ub{when}: 
if $\varphi(\bar x,\bar y) \in {\Cal L}_1,\alpha = \ell g(\bar x) < 
\kappa,\ell g(\bar y) < \kappa$ 
and $\bar b,\bar c \in {}^{\ell g(\bar y)}B$ realize the 
same ${\Cal L}_2$-type in $M$ over $A$
\ub{then} $\varphi(\bar x,\bar b) \in p \Leftrightarrow
\varphi(\bar x,\bar c) \in p$; recalling that $\Sfor^\alpha_{\Cal L}(A,M)$ 
is defined in \scite{734-0n.4} and normally ${\Cal L}_1 = {\Cal L}_2$ or at
least ${\Cal L}_1 \subseteq {\Cal L}_2$.
\nl
5) Let Av$_{<\kappa}(\bold I,M)$ be Av$_{\Bbb L_{\infty,\kappa}[{\frak K}]}
(\bold I,M)$ and Av$_\kappa(\bold I,M)$ be 
Av$_{\Bbb L_{\infty,\kappa^+}[{\frak K}]}(\bold I,M)$.
\enddefinition
\bigskip

\remark{\stag{734-11.7.3} Remark}  1) See definition of $\Sav^\alpha(M)$ in
\scite{734-11.17}(2) below.
\nl
2) An alternative for clause (c) of \scite{734-11.7}(3) is:
\mr
\item "{$(c)'$}"  the set $\{\text{Av}_{{\Cal L},\partial,<\kappa}
(f(\bold I),M):f$ an automorphism of $M$ over $A\}$ has cardinality $\le
\beth_{1,1}(\text{\rm LS}({\frak K})+\theta+|A|) <\|M\|$.
\endroster
\endremark
\bigskip

\proclaim{\stag{734-11.8} Claim}  1) Assume 
that $M \in K,A \subseteq M,\bold I \subseteq {}^\theta M,|\bold I| 
\ge \partial =  { \text{\rm cf\/}}(\partial) > \kappa \ge \theta +
\text{\rm LS}({\frak K})$ and $\bold I$ 
is $({\Cal L},\partial,\kappa)$-convergent.  \ub{Then} the 
type $p = { \text{\rm Av\/}}_{{\Cal L},\partial,\kappa}(\bold I,M)$
belongs to $\Sfor^\theta_{\Cal L}(M)$, i.e., it is complete, recalling
Definition \scite{734-0n.4} (no
demand that it is realized in some $N,M \le_{\frak K} N!$).
\nl
2) Also $\bold I$ is $({\Cal L},\partial,\kappa)$-based on some set of 
cardinality $\le \partial$, even on $\cup \bold J$, for any $\bold J
\subseteq \bold I$ of cardinality $\ge \partial$.  
\endproclaim
\bigskip

\demo{Proof}  1) By the definition.
\nl
2) By the definitions: if $\bar b \in {}^{\kappa^+ >} M,\varphi =
\varphi(\bar x,\bar y) \in {\Cal L}$ and 
$\ell g(\bar b) = \ell g(\bar y),\ell g(\bar x) = \theta$, then by the
convergence

$$
\align
\varphi(\bar x,\bar b) \in p \Leftrightarrow &\text{ for all but } <
\partial \text{ members } \bar a \text{ of } \bold I,M \models
\varphi[\bar a,\bar b] \Leftrightarrow \\
   &\text{ for all but } < \partial \text{ members of } \bold J,M
\models \varphi[\bar a,\bar b].
\endalign
$$
\mn
So only \sftp$_{\Cal L}(\bar b,\cup \bold J,M)$ matters 
hence the non-splitting required in clause (c) of
Definition \scite{734-11.7}(3).  \hfill$\square_{\scite{734-11.8}}$
\enddemo
\bn
As in \marginbf{!!}{\cprefix{300a}.\scite{300a-1.7}}, we deduce non-splitting over a small set from
non-order.
\proclaim{\stag{734-11.7.7} Claim}  Assume $M = \text{\rm EM}_{\tau({\frak
K})}(I,\Phi),\theta + \text{\rm LS}({\frak K}) \le \kappa < \lambda$
and $\beth_{1,1}(\partial) \le |I|$ where $\partial = (2^{2^{\kappa}})^+$
or $I$ is well ordered and $\partial = (2^{\kappa})^+$.   If $M 
\prec_{\Bbb L_{\infty,\partial}[{\frak K}]} N$ \ub{then} for every $\bar a \in
{}^{\theta \ge} N$ there is $B \subseteq M$ of cardinality $<
\partial$ such that $\sftp_{\Bbb L_{\infty,\kappa^+}[{\frak K}]} (\bar
a,M,N)$ does not $(\Bbb L_{\infty,\kappa^+}[{\frak K}],
\Bbb L_{\infty,\kappa^+}[{\frak K}])$-split over $B$.
\endproclaim
\bigskip

\demo{Proof}  Let $\bar x = \langle x_i:i < \ell g(\bar a)\rangle$.

We try to choose $B_\alpha,\gamma_\alpha,\bar a_\alpha,\bar
b_\alpha,\bar c_\alpha,\varphi_\alpha(\bar x,\bar y_\alpha) \in
\Bbb L_{\infty,\kappa^+}[{\frak K}]$ by induction on $\alpha < \partial$
such that
\mr
\item "{$\circledast$}"  $(a) \quad B_\alpha = \cup\{\bar a_\beta:\beta <
\alpha\}$
\sn
\item "{${{}}$}"  $(b) \quad \bar b_\alpha,\bar c_\alpha \in 
{}^{\gamma_\alpha} M$ and $\gamma_\alpha < \kappa^+$
\sn
\item "{${{}}$}"  $(c) \quad \varphi_\alpha(\bar x,\bar y_\alpha) 
\in \Bbb L_{\infty,\kappa^+}[{\frak K}]$ such that $\ell g(\bar y_\alpha) =
\gamma_\alpha$ 
\sn
\item "{${{}}$}"  $(d) \quad N \models ``\varphi_\alpha[\bar a,\bar
b_\alpha] \equiv \neg \varphi_\alpha[\bar a,\bar c_\alpha]"$
\sn
\item "{${{}}$}"  $(e) \quad \bar a_\alpha \in {}^{\ell g(\bar a)}M$
realizes $\{\varphi_\beta(\bar x,\bar b_\beta) \equiv 
\neg \varphi_\beta(\bar x,\bar c_\beta):\beta < \alpha\}$ in $M$
\sn
\item "{${{}}$}"  $(f) \quad M \models ``\varphi_\alpha[\bar
a_\beta,\bar b_\alpha] \equiv \varphi_\alpha[\bar a_\beta,\bar
c_\alpha]"$ for $\beta \le \alpha$.
\ermn
If we are stuck at $\alpha(*) < \partial$ then we cannot choose
$\gamma_\alpha,\bar b_\alpha,\bar c_\alpha,\varphi_\alpha
(\bar x,\bar y_\alpha)$ clauses (b),(c),(d), 
because then $\bar a_\alpha$ as required
in clauses (e),(f) exists because $M \prec_{\Bbb
L_{\infty,\partial}[{\frak K}]} N$.  Hence $B := \cup\{\bar
a_\alpha:\alpha < \alpha(*)\}$ is as required.  So assume that we have
carried the induction.  As $\gamma_\alpha <\kappa^+ < \partial =
\text{cf}(\partial)$ without loss of generality 
$\gamma_\alpha= \gamma < \kappa^+$ for every $\alpha < \partial$.

Let $\partial_1 = (2^\kappa)^+$.

Now by \scite{734-11.9}(5) below when $I$ is not well ordered
and by \scite{734-11.9}(4) below when $I$ is well ordered
(and part (1) of
\scite{734-11.9}(1), recalling $I$ is $\kappa^+$-wide as $\kappa <
\partial$ and $\beth_{1,1}(\partial) \le |I|$) clearly
for some $S \subseteq \partial$ of order type
$\partial_1$, the sequence
$\langle \bar a_\alpha \char 94 \bar b_\alpha \char 94 \bar
c_\alpha:\alpha \in S\rangle$ is $(\Bbb L_{\infty,\kappa^+}[{\frak
K}],\kappa^+,\kappa)$-convergent and 
$(\Bbb L_{\infty,\kappa^+}[{\frak K}],< \omega)$-indiscernible in $M$
hence \wilog \, $\alpha \in S \Rightarrow
\varphi_\alpha = \varphi$.  But as $\partial_1 > \kappa^+$ this
contradicts (e) + (f) of $\circledast$ (if we use
$\partial_1 = \kappa^+$, we
can use a further conclusion of \scite{734-11.9}(1) stated in
\scite{734-11.9}(2), i.e., $\langle \bar a_\alpha \char 94 
\bar b_\alpha \char 94 \bar c_\alpha:\alpha
\in S\rangle$ is a $(\Bbb L_{\infty,\kappa}[{\frak K}],<
\omega)$-indiscernible set not just a sequence, contradiction to (e) +
(f) of $\circledast$).   \hfill$\square_{\scite{734-11.7.7}}$
\enddemo
\bigskip

\proclaim{\stag{734-11.9} Claim}    Assume $M = 
{ \text{\rm EM\/}}_{\tau({\frak K})}(I,\Phi),I$ is $\kappa^+$-wide,
$\kappa < \lambda$ and 
${\text{\rm LS\/}}({\frak K}) + \theta \le \kappa < \partial$. 
\nl
1) Assume that ${\Cal L} = \Bbb L_{\infty,\kappa^+}[{\frak K}]$ and
$\bar a_\alpha = \langle \sigma_i
(\ldots,a_{t(\alpha,i,\ell)},\ldots)_{\ell < n_i}:i < \theta \rangle$ for
$\alpha < \partial$ so $\sigma_i$ is a $\tau(\Phi)$-term, and
{\rm cf}$(\partial) > \kappa$.  Assume further that 
letting $\bar t_\alpha = \langle t(\alpha,i,\ell):
i < \theta,\ell < n_i \rangle$, 
the sequence $\langle \bar t_\alpha:\alpha <
\partial \rangle$ is indiscernible in $I$ for quantifier free formulas
(i.e. the truth values of 
$t(\alpha_1,i_1,\ell_1) < t(\alpha_2,i_2,\ell_2)$ depends only on
$i_1,\ell_1,i_2,\ell_2$ and the truth value of $\alpha_1 < \alpha_2,
\alpha_1 = \alpha_2,\alpha_1 > \alpha_2$).  \ub{Then} 
$\langle \bar a_\alpha:\alpha < \partial \rangle$ 
is $({\Cal L},\partial,\kappa)$-convergent in the model $M$.
\nl
2) In part (1), even dropping the assumption {\rm cf}$(\partial) > \kappa$,
moreover, the sequence $\langle \bar a_\alpha:\alpha <
 \partial\rangle$ is $({\Cal L},\kappa^+,\kappa)$-convergent
and $({\Cal L},< \omega)$-indiscernible in $M$.
\nl
3) In part (1) and in part (2), letting $J_0 = \{t(0,i,\ell):t(0,i,\ell) =
t(1,i,\ell)$ and $i < \theta,\ell < n_i\}$ assume 
$J_0 \subseteq J \subseteq I,J$ is $\kappa^+$-wide (e.g. $J =
\{t(\alpha,i,\ell):\alpha < \kappa^+,i < \theta,\ell < n_i\}$) and
$B$ is the universe of $\text{\rm EM}_{\tau({\frak K})}(J,\Phi)$ 
and $(i_1,i_2 < \theta,\ell_1 < n_{\ell_1},\ell_2 < n_{i_2}$
and $[\alpha,\beta < \partial \Rightarrow t(\alpha,i_1,\ell_1) <_I
t(\beta,i_2,\ell_2)] \Rightarrow \exists s \in J_0[\alpha,\beta
< \partial \Rightarrow t(\alpha,i_1,\ell_1) <_I t <_I t(\beta,i_2,\ell_2)]$
\ub{then} $B$ is a $(\partial,\kappa)$-base of
$\{\bar a_\alpha:\alpha < \partial\}$. 
\nl
4) If $I$ is well ordered (or just is {\rm EM}$_{\{<\}}(J,\Psi),\Psi
\in \Upsilon^{\text{or}},J$ well ordered),
{\rm LS}$({\frak K}) + \theta \le \kappa,2^\kappa < \partial,
(\forall \alpha < \partial)[|\alpha|^\theta < \partial 
= { \text{\rm cf\/}}(\partial)]$ and $\bar b_\alpha \in
{}^\theta M$ for $\alpha < \partial$, \ub{then} for some
stationary $S \subseteq \{\delta < \partial:\text{\rm cf}(\delta) 
\ge \theta^+\}$,
the sequence $\langle \bar b_\alpha:\alpha \in S \rangle$ is as in
part (1) hence is $(\kappa^+,\kappa)$-convergent in $M$.  
Moreover, if $S_0 \subseteq \{\delta < \partial:{\text{\rm cf\/}}(\delta) 
\ge \theta^+\}$ is stationary we can demand $S \subseteq S_0$.
\nl
5) If in (4) we omit the assumption ``$I$ is well ordered", and add
$\partial \rightarrow (\partial_1)^2_{2^\kappa}$, e.g. $\partial_1 =
(2^\kappa)^+,\partial = (2^{2^\kappa})^+$ \ub{then} we can find
$S \subseteq \partial,|S| = \partial_1$ such that $\langle \bar
a_\alpha:\alpha \in S \rangle$ is as in (1).
\endproclaim
\bigskip

\remark{Remark}  In fact the well order case always applies at least
if $\partial < \mu$.
\endremark
\bigskip

\demo{Proof}  1) Let $\bar b \in {}^\kappa M$, so $\bar b =  \langle
\sigma^*_j(\ldots,a_{s(j,\ell)},\ldots)_{\ell < m_j}:j < \kappa
\rangle$ where $\sigma^*_i$ is a $\tau(\Phi)$-term, $s(j,\ell) \in I$ 
and let $\bar s = \langle s(j,\ell):\ell < m_j,j < \kappa \rangle$.

Now for each $i_1 < \theta,\ell_1 < n_{i_1}$ and $j_1 < \kappa,k_1 <
m_{j_1}$ the sequence $\langle t(\alpha,i_1,\ell_1):\alpha < \partial
\rangle$ is monotonic (in $I$) hence there is
$\alpha(i_1,\ell_1,j_1,k_1) < \partial$ such that
\mr
\item "{$(*)_1$}"  if $\beta,\gamma \in \partial \backslash
\{\alpha(i_1,\ell_1,j_1,k_1)\}$ and $\beta <
\alpha(i_1,\ell_1,j_1,k_1) \equiv \gamma < \alpha(i_1,\ell_1,j_1,k_1)$
then $\bigl(t(\beta,i_1,\ell_1) <_I s(j_1,k_1)\bigr) 
\equiv \bigl(t(\gamma,i_1,\ell_1)
<_I s(j_1,k_1)\bigr)$ and $\bigl(t(\beta,i_1,\ell_1) >_I
s(j_1,k_1)\bigr) \equiv \bigl(t(\gamma,i_1,\ell_1) >_I s(j_1,k_1)\bigr)$.
\ermn
Let $u := \{\alpha(i_1,\ell_1,j_1,k_1):i_1 < \theta,\ell_1 < n_{i_1},j_1 <
\kappa,k_1 < m_{j_1}\}$, it is a subset of $\partial$ of cardinality $\le
\theta +\kappa = \kappa$.
\nl
Hence
\mr
\item "{$(*)_2$}"  if $\beta,\gamma \in \partial \backslash u$ and
$\beta {\Cal E}_u \gamma$ which is 
defined by $(\forall \alpha \in u)(\alpha < \beta
\equiv \alpha < \gamma)$ \ub{then} $\bar t_\beta \char 94 \bar s,\bar
t_\gamma {} \char 94 \bar s$ realizes the same quantifier free type in $I$
\ermn
Now by clause (c) of \scite{734-11.3} recalling $I$ is $\kappa^+$-wide we have
\mr
\item "{$(*)_3$}"  if $\beta,\gamma \in \partial \backslash u$ and
$\beta {\Cal E}_u \gamma$ then $\bar a_\beta \char 94 \bar b,\bar a_\gamma
\char 94 \bar b$ realizes the same $\Bbb L_{\infty,\kappa^+}[{\frak
K}]$-type in $M$.
\ermn
As $\bar b$ was any member of ${}^\kappa M$ we have gotten
\mr
\item "{$(*)_4$}"  if $\bar b \in {}^{\kappa \ge}M$, then for some $u
= u_{\bar b} \subseteq \partial$ of cardinality $\le \kappa$ we have:
\nl
if $\beta,\gamma \in \partial \backslash u$ and $\beta {\Cal E}_u \gamma$
then $\bar a_\beta \char 94 \bar b,\bar a_\gamma \char 94 \bar b$
realize the same $\Bbb L_{\infty,\kappa^+}[{\frak K}]$-type in $M$.
\ermn
As we are assuming cf$(\partial) > \kappa(\ge \theta + \text{
LS}({\frak K}) \ge |\tau_\Phi|)$ we can conclude that
\mr
\item "{$(*)_5$}"   $\langle \bar a_\alpha:\alpha < 
\partial \rangle$ is $({\Cal L},\partial,\kappa)$-convergent in $M$.
\ermn
So we have proved \scite{734-11.9}(1).
\nl
2) We start as in the proof of part (1). 
However, after $(*)_3$ above letting for simplicity 
$u^+ = \{\alpha < \partial$: for some $\beta \in u \cap \alpha$ we
have $\alpha + \kappa = \beta + \kappa\}$ we have
\mr
\item "{$(*)_6$}"  if $\beta,\gamma \in \partial \backslash u^+$ and
$\beta < \gamma,\neg(\beta E_{u^+} \gamma)$ \ub{then} we can find
$(\mu^+,I^+,\bar s',\bar b)$ such that
{\roster
\itemitem{ $(\alpha)$ }  $I \subseteq I^+ \in K^{\text{lin}}$
\sn
\itemitem{ $(\beta)$ }  $M^+ = \text{ EM}_{\tau({\frak K})}(I,\Phi)$
hence $M \prec_{\Bbb L_{\infty,\kappa^+}[{\frak K}]} N$
\sn
\itemitem{ $(\gamma)$ }  $\bar s = \langle
s'(j,k):k<m,j<\kappa\rangle$ a sequence of elements of $I^+$
\sn
\itemitem{ $(\delta)$ }  $\bar b' = \langle
\sigma^*_j(\ldots,a_{s'(j,\ell)},\ldots)_{\ell < m_j}:j <
\kappa\rangle \in {}^\kappa(M^+)$
\sn
\itemitem{ $(\varepsilon)$ }  $\bar b \char 94 \bar a_\gamma,\bar b'
\char 94 \bar a_\gamma$ realize the same $\Bbb
L_{\infty,\kappa^+}[{\frak K}]$-types in $M^+$ as $\bar b \char 94
\bar a_\gamma,\bar b \char 94 \bar a_\beta$ respectively
\sn
\itemitem{ $(\zeta)$ }  $\bar s \char 94 \bar t_\beta,\bar s' 
\char 94 \bar t_\beta$ form a $\Delta$-system pair, i.e. are as
in $\boxtimes$ from \scite{734-11.1.2}(2).
\endroster}
\ermn
[Why?  

Let $w^+ = \{(j,k):k < m_j$ and $j < \kappa$ and for some $\ell
< n_{i_1},i_1 < \theta$ we have $\alpha(i_1,\ell_1,j,k) \in
(\beta,\gamma)\}$

$$
w^- := \{(j,k):j < \kappa,k < m_j \text{ and } (j,\kappa) \notin
w^+\}.
$$
\mn
We choose $I^+$ extending $I$ and $\bar s_\varepsilon = \langle
s_i(j,k):k < m_j,j < \kappa\rangle$ for $\varepsilon < \kappa$ such
that
\mr
\item "{$(a)$}"  the set of elements of $I^+$ is the disjoint union of
$I$ and $\{s_\varepsilon(j,k):(j,k) \in w$ and $\varepsilon \in
(0,\kappa)\}$
\sn
\item "{$(b)$}"  $\bar s_\varepsilon,\bar s$ realize the same
quantifier-free type in $I^+$
\sn
\item "{$(c)$}"  if $\varepsilon,\zeta < \kappa$ then $\bar t_{\gamma
+ \varepsilon} \char 94 \bar s_\zeta$ realizes in $I^+$ the 
quantifier-free type $\sftp_{\text{qf}}(\bar t_\beta \char 94 \bar
s,\emptyset,I)$ if $\varepsilon < \zeta$ and $\sftp_q(\bar t_\gamma
\char 94 \bar s,\emptyset,I)$ if $\varepsilon \ge \zeta$
\sn
\item "{$(d)$}"  $\langle \bar t_{\gamma + \varepsilon} \char 94 \bar
s_\varepsilon:\varepsilon < \kappa\rangle$ is indiscernible for
quantifier-free formulas on $I^+$
\sn
\item "{$(e)$}"  $\bar s_0 = \bar s$.
\ermn
This is straight.  Using $\bar s' = \bar s_1$ we are done.]

Now as $\Phi$ has the $\kappa$-non-order property (by Claim \scite{734-11.1.2}(2)
which contains a definition, noting that the assumption of
\scite{734-11.1.2} holds by \scite{734-11.2A}(1) and also \scite{734-11.2A}(2)),
repeating $(*)_4,(*)_5$ we get
\mr
\item "{$(*)_7$}"  for every $\bar b \in {}^{\kappa \ge} M$, for
some $u = u^+_{\bar b} \in [\partial]^{\le \kappa}$ if $\beta,\gamma \in
\partial \backslash u^+$ then $\bar a_\beta \char 94 \bar
b,\bar a_\gamma \char 94 \bar b$ realizes the same $\Bbb
L_{\infty,\kappa^+}[{\frak K}]$-type in $M$.
\ermn
In other words
\mr
\item "{$(*)_8$}"  the sequence $\langle \bar a_\alpha:\alpha <
\partial\rangle$ is 
$(\Bbb L_{\infty,\kappa^+}[{\frak K}],\kappa^+)$-convergent.
\ermn
The proof that it is a $(\Bbb L_{\infty,\kappa^+}[{\frak K}],
< \omega)$-indiscernible set is similar.
\nl
3) Not used; easy by \scite{734-11.8}(2) and convergence.  [That is, note
that we can find $I^+$ and $\bar a'_\alpha = \langle
\sigma_i(\ldots,a_{t'(\alpha,i,\ell)},\ldots)_{\ell_i < n_i}:i <
\theta \rangle$ for $\alpha < \partial + \gamma$ such that:
\mr
\item "{$(a)$}"  $I^+ \in K^{\text{lin}}$ extend $I$
\sn
\item "{$(b)$}"  $t'(\alpha,i,\ell) \in I^+$
\sn
\item "{$(c)$}"  $\bar t'_\alpha = \langle t'(\alpha,i,\ell):i <
\theta,\ell < n_i\rangle$
\sn
\item "{$(d)$}"  $\langle \bar t'_\alpha:\alpha < \partial +
\gamma\rangle$ is indiscernible for quantifier-free formulas in $I^+$
\sn
\item "{$(e)$}"  $\langle \bar t_\alpha:\alpha < \partial\rangle \char
94 \langle \bar t'_\alpha:\alpha \in [\partial,\partial +
\partial)\rangle$ is indiscernible for quantifier-free formulas in $I'$
\sn
\item "{$(f)$}"  for each $i < \theta,\ell < n_i$ such that
$t(0,i,\ell) = (j,i,t)$ the convex hull $I_*$ of
$\{t'(\alpha,i,\ell):\alpha < \partial\}$ in $I^+$ is disjoint to $I$
and if $s_1 <_I s_2$ and $(s_1,s_2)_{I^*} \cap I_* = \emptyset$ then
$[s_1,s_2]_{I^*} \cap J_0 \ne \emptyset$.
\ermn
So we can average over $\langle \bar a'_\alpha:\alpha <
\partial\rangle$ instead averaging over $\langle \bar a_\alpha:\alpha <
\partial\rangle$, and this implies the result.  In fact we can weaken
the assumption.]
\nl
4) Should be clear.  [Still let $\bar t_\alpha = \langle t_{\alpha,i}:i <
\theta \rangle$ be such that \nl
$\bar b_\alpha = \langle
\sigma_{\alpha,j}(\ldots,a_{t_{\alpha,i(j,\alpha,\ell)}},\ldots)_{\ell
< n(\alpha,j)}:j < \theta \rangle$.  So as 
$(\text{LS}({\frak K})+|\tau_\Phi|)^\theta < \partial = 
\text{\rm cf}(\partial)$ for some stationary $S_1 \subseteq \{\delta <
\partial$: cf$(\delta) \ge \theta^+\}$ we have $\alpha \in S_1 \wedge
j < \theta \Rightarrow \sigma_{\alpha,j} = \sigma_j$ (hence $j <
\theta \Rightarrow n(\alpha,j) = n(j))$ and $\alpha \in S_1 \wedge j <
\theta \wedge \ell <n(j) \Rightarrow i(j,\alpha,\ell) = i(j,\ell)$ and
for every $i_1,i_2 < \theta$ we have $t_{\alpha,i_1} <_I
t_{\alpha,i_2} \equiv (i_1,i_2) \in W$ for some sequence $\bar \sigma
= \langle \sigma_j:j < \theta \rangle$ of $\tau_\Phi$-terms and $W
\subseteq \kappa \times \kappa$ and sequence $\left< \langle
i(j,\ell):\ell < n(j)\rangle:j < \theta\right>$.

If $I$ is well ordered, for $\delta \in S_1$ let $\gamma_\delta
= \text{\rm Min}\{\gamma$: if $i < \theta$ and there are $\beta <
\delta,j < \theta$ such that $t_{\delta,i} <_I t_{\beta,j}$ 
and \ub{then} letting $(\beta_{\delta,i},j_{\delta,i})$ be
such a pair with $t_{\beta_{\delta,i},j_{\delta,i}}$ being
$<_I$-minimals, we have $\beta_{\delta,i} < \gamma\}$;  clearly
$\gamma_\delta$ is well defined and $< \delta$ so by Fodor lemma for
some $\gamma_* < \partial$ the set $S_1 := \{\delta \in S_2:\gamma_\delta =
\gamma_*\}$ is stationary.  As $|\gamma_*|^\theta < \partial$, for
some $u \subseteq \theta $ and stationary $S_3 \subseteq S_2$ we have:
if $\delta \in S_3$ then $j \in u \Leftrightarrow
(\beta_{\delta,i},j_{\delta,i})$ well defined and $j \in u \wedge \alpha
\in S_3 \Rightarrow (\beta_{\delta,i},j_{\delta,i}) = (\beta_i,j_i)$
and for each $i \in u$ the truth value of ``$t_{\delta,i} =
t_{\beta_i,j_i}$" is the same for all $\delta \in S_3$.

Now apply part (1) to $\langle \bar b_\alpha:\alpha \in S_3\rangle$.]
\nl
5) By (1) and the definition of $\partial \rightarrow 
(\partial_1)^2_{2^\kappa}$.   \hfill$\square_{\scite{734-11.9}}$
\enddemo
\bigskip

\proclaim{\stag{734-11.10} Claim}  1) If 
$M \le_{\frak K} N$ are from $K^*_\lambda$
and $\kappa \in [{\text{\rm LS\/}}({\frak K}),\lambda),\kappa^+
< \partial = \text{\rm cf}(\partial) < \lambda$ and moreover $\theta \le
\kappa$ and $\bar a \in {}^\theta N$ \ub{then} there is a 
$(\kappa^+,\kappa)$-convergent set $\bold I \subseteq {}^\theta M$ of 
cardinality $\partial$ such that ${\text{\rm Av\/}}_\kappa(\bold I,M)$ is 
realized in $N$ by $\bar a$. 
\nl
2) In fact we can weaken $M,N \in K^*_\lambda$ to $M,N \in
K^*_{\ge \beth_{1,1}(\partial')}$ where, e.g. $\partial' = \beth_5(\kappa)^+$. 
\nl
3)  Assume $\theta \le \kappa,\kappa \in [\text{\rm LS}
({\frak K}),\lambda),\partial' = \beth_5(\kappa)^+$ and $M_1 \in K^*_{\ge
\beth_{1,1}(\partial')}$.  Assume further $M_1 \le_{\frak K} M_2 = 
\,\text{\rm EM}_{\tau({\frak K})}(I_2,\Phi),|\xi| = \theta$ and
$\bold I \subseteq {}^\xi (M_1)$ is a $(\kappa^+,\kappa)$-convergent set
(in $M_1$) of cardinality $\partial'$.  If $I_2 <^*_{K^{\text{flin}}}
I_3$ (or just $I_3$ is $\kappa^+$-wide over $I_2$, which follows as
$|I_2| \ge |\bold I| = \partial'$) and
$M_3 = \text{\rm EM}_{\tau({\frak K})}(I_3,M_3)$ \ub{then}
\mr
\item "{$(a)$}"  we can find $\bar d \in {}^\xi(M_3)$ realizing
{\rm Av}$_\kappa(M_2,\bold I)$ so well defined
\sn
\item "{$(b)$}"  if $M_1 \le_{\frak K} N \in K^*$ and $\bar d^* \in
{}^\xi N,|\xi| \le \theta$ \ub{then} we 
can find $\bar d \in {}^\xi(M_3)$
realizing $\sftp_{\Bbb L_{\infty,\kappa^+}[{\frak K}]}(\bar d^*,M_1,N)$
and $\sftp_{\Bbb L_{\infty,\kappa^+}[{\frak K}]}(\bar
d,M_2,M_3)$ is the average of some $(\kappa^+,\kappa)$-convergent
$\bold I' \subseteq {}^\alpha(M_1)$ of cardinality $\partial'$.
\endroster
\endproclaim
\bigskip

\remark{Remark}  The exact 
value of $\partial'$ have no influences for our purpose.
\endremark
\bigskip

\demo{Proof}  1) Without loss of generality $M = 
\text{\rm EM}_{\tau({\frak K})}(I,\Phi)$.  
Let $\partial_0 = \partial$ and $\partial_{\ell +1} =
\beth_2(\partial_\ell)^+$ for $\ell=0,1$ so $\partial_\ell <
\lambda$ and $\ell =1,2 \Rightarrow 
(\forall \alpha < \partial_\ell)(|\alpha|^{\kappa +\theta}
< \partial_\ell = \text{\rm cf}(\partial_\ell) < \lambda)$ (if $I$ is well
ordered (which is O.K. by \scite{734-11.5}(4)) and 
$(\forall \alpha < \partial)(|\alpha|^\kappa < \partial)$ then
we can use $\partial_\ell = \partial$).

By \scite{734-11.7.7} there is $B_* \subseteq M$ of cardinality $<
\partial_2$ (or just $\le 2^{2^\kappa} < \partial_2$)
such that $\sftp_{\Bbb L_{\infty,\kappa^+}[{\frak K}]}(\bar a,M,N)$ does
not $(\Bbb L_{\infty,\kappa^+}[{\frak K}],\Bbb
L_{\infty,\kappa^+}[{\frak K}])$-split over $B_*$.

Now by \scite{734-11.5}(1) for
every $B \subseteq M,|B| < \partial_2$ there is $\bar a' \in {}^\theta M$
realizing in $M$, equivalently in $N$ (with $\ell g(\bar x) = 
\theta$, of course), the type $\sftp_{\Bbb L_{\infty,\kappa^+}[{\frak
K}]}(\bar a,B,N) = \{\varphi(\bar x,\bar b):\bar b \in
{}^{\kappa \ge}B,\varphi(\bar x,\bar y) \in 
\Bbb L_{\infty,\kappa^+}[{\frak K}]$
and $N \models \varphi[\bar a,\bar b]\}$.

We can choose $J_\alpha,B_\alpha,\bar a_\alpha$
by induction on $\alpha < \partial_2$ such that $B_\alpha$ includes
$\cup\{\bar a_\beta:\beta < \alpha\} \cup B_*,B_\alpha$ is the universe of
EM$(J_\alpha,\Phi),J_\alpha \subseteq I,|J_\alpha| < \partial_2,J_\alpha$ 
increasing with $\alpha$ and 
$J_\alpha$ is quite closed (e.g. is ${\frak B}_\alpha \cap I$ where
${\frak B}_\alpha \prec_{\Bbb L_{\kappa^+,\kappa^+}}
({\Cal H}(\chi),\in,<^*_\chi)$ 
with $M,N,\text{EM}(I,\Phi),{\frak K},\langle \bar a_\beta:\beta <
\alpha \rangle,{\frak K},\kappa,\theta$ belonging to ${\frak
B}_\alpha$ and ${\frak B}_\alpha$ has cardinality $< \partial_2$ and ${\frak
B} \cap\partial_2 \in \partial_2$).
Then choose $\bar a' = \bar a_\alpha$ as above, i.e. 
$\bar a_\alpha \in {}^\theta M$ realizes the same
$\Bbb L_{\infty,\kappa^+} [{\frak K}]$-type as $\bar a$ 
over $B_\alpha = M \cap {\frak B}_\alpha = 
\text{ EM}_{\tau({\frak K})}(J_\alpha,\bar a)$ in $N$; 
such $\bar a_\alpha$ exists by \scite{734-11.5}(1).
So for some set $S_1 \subseteq \partial_2$ of order type $\partial_1$ 
the sequence $\bold I = \langle \bar a_\beta:\beta \in S_1\rangle$ is 
$(\kappa^+,\kappa)$-convergent (by \scite{734-11.9}(4),(5)).  

It is enough to show that $\bold I$ is as required, toward
contradiction assume that not.  Then
there is an appropriate formula $\varphi(\bar x,\bar y)$ with 
$\ell g(\bar x) = \theta,
\ell g(\bar y) = \kappa$ and $\bar b \in {}^\kappa M$ such that $N \models
\varphi[\bar a,\bar b]$ but $u := \{\alpha \in S_1:M \models \varphi
[\bar a_\alpha,\bar b]\}$ has cardinality $< \kappa^+$.  Now for 
$\alpha \in S_1$
as $J_\alpha$ was chosen ``closed enough", there is $\bar b_\alpha \in
{}^\kappa(\text{EM}_{\tau({\frak K})}(J_\alpha,\Phi)) \subseteq
{}^\kappa M$ 
realizing $\sftp_{\Bbb L_{\infty,\kappa^+}[{\frak K}]}(\bar b,B_*,M)$
such that $\beta \in S_1 \cap \alpha
\Rightarrow M \models ``\varphi[\bar a_\beta,\bar b] \equiv \varphi
[\bar a_\beta,\bar b_\alpha]"$ (possible, e.g. as $|B_\alpha|^{|S \cap
\alpha|} \le (2^{< \partial_1})^{< \partial_1} < \partial_2)$.

So, again by \scite{734-11.9}(4),(5), 
for some $S_0 \subseteq S_1$ of order type $\partial = \partial_0$,
the sequence $\langle \bar a_\alpha
\char 94 \bar b_\alpha:\alpha \in S_0 \rangle$ is 
$(\Bbb L_{\infty,\kappa^+},\kappa^+,\kappa)$-convergent in $M$ and 
$(\Bbb L_{\infty,\kappa^+},< \omega)$-indiscernible.
Let $\alpha \in S_0$ be such that $|S_0 \cap \alpha| > \kappa$,
possible as $|S_0| = \partial_0 > \kappa_0$.  So the set 
$\{\beta \in S_1 \cap \alpha:M \models \varphi
[\bar a_\beta,\bar b_\alpha]\}$ has cardinality $\le \kappa$ (being
equal to $\{\beta \in S_1 \cap \alpha:N \models 
\varphi[\bar a_\beta,\bar b]\})$ but $\alpha \in S_0 \subseteq S_1$
and $|S_0 \cap \alpha| > \kappa$, so for some $\beta < \alpha$ 
from $S_0,M \models \neg \varphi
[a_\beta,\bar b_\alpha]$ hence by the indiscernibility $M \models \neg
\varphi[\bar a_\beta,\bar b_\gamma]$ for every 
$\beta < \gamma$ from $S_0$.

On the other hand if $\alpha < \beta$ are from $S_0$ then by the
choice of $\bar b_\alpha$ the sequences $\bar b,\bar b_\alpha$
realizes the same $\Bbb L_{\infty,\kappa^+}[{\frak K}]$-type over
$B_*$.   Now $\sftp_{\Bbb L_{\infty,\kappa^+}[{\frak K}]}(\bar a,M,N)$
does not split over $B_*$ by the choice
of $B_*$ so we have $N \models ``\varphi[\bar a,\bar b] \equiv \varphi[\bar
a,\bar b_\alpha]"$ but by the choice of $\bar b$ we have 
$N \models \varphi[\bar a,\bar b]$ hence $N \models \varphi
[\bar a,\bar b_\alpha]$ hence $M \models \varphi[\bar a_\beta,\bar b_\alpha]$
by the choice of $\bar a_\beta$.  Together this contradicts \scite{734-11.1.2},
i.e., \scite{734-11.2A}(1).
\nl
2) Similarly (using \scite{734-11.5}(2) instead of \scite{734-11.5}(1).
\nl
3) \ub{Clause (a)}:

By \scite{734-11.3} and the LS argument (i.e. by \scite{734-0n.13}(4)) \wilog
\, $M_1 \in K^*_{< \lambda}$.  Let $\partial_\ell =
\beth_\ell(\kappa)^+$ for $\ell \le 5$ so $\partial' = \partial_5$
and for notational simplicity assume $\theta \ge
\aleph_0$.

Let $\{\bar a_\alpha:\alpha < \partial'\}$ list the members of $\bold
I$, so for each $\alpha < \partial'$ there is $I_{2,\alpha} \subseteq
I_2$ of cardinality $\theta$ such that $\bar a_\alpha$ is from
EM$_{\tau({\frak K})}(I_{2,\alpha},\Phi)$.

For each $\alpha < \partial'$ let $\bar t^\alpha = \langle
t^\alpha_i:i < \theta\rangle$ list $I_{2,\alpha}$ and so $\bar
a_\alpha = \langle \sigma_{\alpha,\zeta}(\bar t^\alpha):\zeta < \xi\rangle$
for some sequence $\langle \sigma_{\alpha,\zeta}(\bar x):\zeta <
\xi\rangle$ of $\tau_\Phi$-terms.  We can find $S \subseteq \partial'$
of order type $\partial_4$ such that $\zeta < \xi \wedge \alpha \in S
\Rightarrow \sigma_{\alpha,\zeta} = \sigma_\zeta$ and $\langle \bar
t^\alpha:\alpha \in S\rangle$ is an indiscernible sequence (for
quantifier free formulas, in $I_2$, of course).

By renaming $\kappa^+ \subseteq S$.  We define a partition $\langle u_{-1},
u_0,u_1\rangle$ of $\xi$ by

$$
u_0 = \{i < \theta:t^\alpha_i = t^\beta_i \text{ for } \alpha,\beta
\in S\}
$$

$$
u_1 = \{i < \theta:t^\alpha_i <_{I_2} t^\beta_i \text{ for } \alpha <\beta
\text{ from } S\}
$$

$$
u_{-1} = \{i < \theta:t^\beta_i <_{I_2} t^\alpha_i \text{ for } \alpha
< \beta \text{ from }  S\}.
$$
\mn
We define an equivalence relation $e$ on $u_{-1} \cup u_1$
\mr
\item "{$\odot$}"  $i_1 e i_2$ iff for some $\ell \in
\{1,-1\},i_1,i_2 \in u_\ell$ and $(t^\alpha_{i_1} <_I
t^\beta_{i_2}) \equiv (t^\alpha_{i_2} <_I t^\beta_{i_1})$ for every
(equivalently some) $\alpha < \beta$ from $S$.
\ermn
There is a natural set of representatives:  $W = \{\zeta < \theta:\zeta \in
u_{-1} \cup u_1$ and $\zeta = \text{ min}(\zeta/e)\}$.

We now define a linear order $I^+_2$; its set of elements is $\{t:t
\in I_2\} \cup\{t^*_i:i \in u_{-1} \cup u_1\}$ where, of course,
$t^*_i \in I^+_2$ are pairwise distinct and $\notin I_2$.  
The order is defined by (or see $\circledast_2$ and think)
\mr
\item "{$\circledast_1$}"  $s_1 <_{I^+_2} s_2$ iff
{\roster
\itemitem{ $(a)$ }  $s_1,s_2 \in I_2$ and $s_1 <_{I_2} s_2$
\sn
\itemitem{ $(b)$ }  $s_1 \in I_2,s_2 = t^*_i$ and $s_1 <_{I_2}
t^\alpha_i$ for every $\alpha < \kappa^+$ large enough
\sn
\itemitem{ $(c)$ }  $s_1 = t^*_i,s_2 \in I_2$ and $t^\alpha_i <_{I_2}
s_2$ for every $\alpha < \kappa^+$ large enough
\sn
\itemitem{ $(d)$ }  $s_1 = t^*_i,s_2 = t^*_j$ and $t^\alpha_i <_I
t^\alpha_j$ for every $\alpha < \kappa^+$.
\endroster}
\ermn
Let $t^*_i = t^\alpha_i$ for $i \in u_0$ and any $\alpha < \kappa^+$.
Let $M^+_2 = \text{ EM}_{\tau({\frak K})}(I^+_2,\Phi)$.

It is easy to check (by \scite{734-11.3}(a),(c)) that
\mr
\item "{$\circledast_2$}"  $(a) \quad I_2 \subseteq I^+_2$
\sn
\item "{${{}}$}"  $(b) \quad \bar t^* \in {}^\theta(I^+_2)$
\sn
\item "{${{}}$}"  $(c) \quad$ if $J \subseteq I_2$ has cardinality
$\le \kappa$ \ub{then} for every $\alpha < \kappa^+$ large enough,
\nl

\hskip23pt the sequences $\bar t^*,\bar t^\alpha$ realizes the same 
quantifier free type 
\nl

\hskip23pt over $J$ inside $I^+_2$.
\ermn
Let
\mr
\item "{$\circledast_3$}"  $\bar d := \langle \sigma_\zeta(\bar
t^*):\zeta < \xi \rangle \in {}^\xi(M^+_2)$.
\ermn
Recall that $\|M_2\| < \lambda$ hence $|I_2| <\lambda$ and $I_2$ is
$\kappa^+$-wide having cardinality $\ge \partial' > 2^\kappa$.

Note
\mr
\item "{$\circledast_4$}"  $\bar t^*$ realizes Av$_{\text{qf}}(\{\bar
t^\alpha:\alpha \in S\},I_2)$ in the linear order $I^+_2$.
\ermn
Without loss of generality $I^+_2 \cap I_3 = I_2$, so we can find a
linear order $I_4$ of cardinality $\lambda$ such that $I^+_2 \subseteq
I_4 \wedge I_3 \subseteq I_4$.  As $I_3$ is $\kappa^+$-wide over $I_2$
(see the assumption and Definition \scite{734-0n.2}(6)+(3)), 
there is a convex subset $I'_3$ of $I_3$
disjoint to $I_2$ which contains a monotonic sequence $\langle
s_\alpha:\alpha < \kappa^+\rangle$.  Without loss of generality there
are elements 
$s_\alpha \, (\alpha \in [\kappa^+,\lambda \times \kappa^+)$ in $I_4$
such that $\langle s_\alpha:\alpha < \lambda \times \kappa^+\rangle$
is monotonic (in $I_4$), and its
convex hull is disjoint to $I_2$.  Let $I^-_3 = I_2 \cup
\{s_\alpha:\alpha < \kappa^+\}$ and $I^\pm_3 = I_2 \cup \{s_\alpha:\alpha <
\lambda \times \kappa^+\}$.  
 
Now we use \scite{734-11.3} several times.  First, EM$_{\tau({\frak
K})}(I_2,\Phi) \prec_{\Bbb L_{\infty,\kappa^+}[{\frak K}]} \text{
EM}_{\tau({\frak K})}(I^+_2,\Phi) 
\prec_{\Bbb L_{\infty,\kappa^+}[{\frak K}]} \text{ EM}_{\tau({\frak K})}
(I_4,\Phi)$ as $I_2 \subseteq I^+_2 \subseteq I_4$ are $\kappa^+$-wide,
hence by $\circledast_4$ the sequence $\bar d$  
realizes $q:= \text{ Av}_\kappa(\{\langle \sigma_\zeta(\bar t^\alpha):
\zeta < \theta\rangle:\alpha < \kappa^+\},M_2) = \text{ Av}(\{\bar
a_\alpha:\alpha < \kappa^+\},M_2) = \text{ Av}_\kappa(\bold I,M_2)$ in
$M^+_2$ and also in EM$_{\tau({\frak K})}(I_4,\Phi)$.  Second, as $|I_2| <
\lambda,I_2 \subseteq I^\pm_3 \subseteq I_4$ and $|I^\pm_3| = |I_4| =
\lambda$, by \scite{734-11.5}(1) we have 
EM$_{\tau({\frak K})}(I^\pm_3,\Phi) \prec_{\Bbb
L_{\infty,\lambda}[{\frak K}]} \text{ EM}_{\tau({\frak K})}(I_4,\Phi)$
so some $\bar d' \in {}^\xi(\text{EM}_{\tau({\frak K})}(I^\pm_3,\Phi))$
realizes the type $q$ in EM$_{\tau({\frak K})}(I^\pm_3,\Phi)$.  Let
$w_1 \subseteq \lambda \times \kappa^+$ be of cardinality $\le \theta \le
\kappa$ such that $\bar d'$ belongs to EM$_{\tau({\frak K})}(I_2 \cup
\{s_\alpha:\alpha \in w_1\},\Phi)$.  Choose $w_2 \subseteq \lambda
\times \kappa^+$ of order type $\kappa^+$ including $w_1$, so
EM$_{\tau(({\frak K})}(I_2 \cup \{s_\alpha:\alpha \in w_2\},\Phi)
\prec_{\Bbb L_{\infty,\kappa^+}[{\frak K}]} \text{ EM}_{\tau({\frak
K})}(I^\pm_3,\Phi)$ and $\bar d'$ belongs to the former hence realizes $q$
in it.  But there is an isomorphism $h$ from $I_2 \cup
\{s_\alpha:\alpha \in w_2\}$ onto $I^-_3$ over $I_2$, hence it induces
an isomorphism $\hat h$ from EM$_{\tau({\frak K})}(I_2 \cup
\{s_\alpha:\alpha \in w_2\},\Phi)$ onto EM$_{\tau({\frak
K})}(I^-_3,\Phi)$ so $\hat h \, (\bar d')$ realizes $q$ in the
latter.  But $I^-_3 \subseteq I_3$ are both $\kappa^+$-wide 
hence by \scite{734-11.3} the sequence 
$\hat h(\bar d')$ realizes
$q$ in $M_3 = \text{ EM}_{\tau({\frak K})}(I_3,\Phi)$ as required.
\bn
\ub{Clause (b)}:

By part (2) we can find appropriate $\bold I$ and then apply clause (a).
  \hfill$\square_{\scite{734-11.10}}$
\enddemo
\bigskip

\remark{\stag{734-11.10D} Remark}  1) In fact in \scite{734-11.7.7}, we can
choose $B$ of cardinality $\kappa$, hence similarly in the proof of
\scite{734-11.10}(1).
\nl
2) Also using solvability to get well ordered $I$ we can prove : if $A
\subseteq M = \text{ EM}_{\tau({\frak K})}(\lambda,\Phi)$ and $|A|
< \lambda$ then the set of $\Bbb L_{\infty,\kappa^+}
[{\frak K}]$-types realized in $M$ over $A$ is $\le (|A| +2)^\kappa$.
\endremark
\bigskip

\proclaim{\stag{734-11.11} Claim}  1) If $M \in K^{**}_{\ge\kappa}$ 
and {\rm LS}$({\frak K}) \le \theta$ 
and $\partial = \beth_{1,1}(\theta) \le \kappa \le \lambda$, \ub{then} for 
$\bar a,\bar b \in {}^\theta M$ the following are equivalent: (the
difference is using $\partial$ or $\kappa$)
\mr
\item "{$(a)$}"   $\bar a,\bar b$ realize the same
$\Bbb L_{\infty,\partial}[{\frak K}]$-type in $M$
\sn
\item "{$(b)$}"  $\bar a,\bar b$ realize the same $\Bbb L_{\infty,\kappa}
[{\frak K}]$-type in $M$.
\ermn
2) For $M,\theta,\partial,\kappa$ as above, the number of 
$\Bbb L_{\infty,\partial}[{\frak K}]$-types of $\bar a \in {}^\theta
M$ where $M = { \text{\rm EM\/}}_{\tau({\frak K})}(I,\Phi),|I| \ge \partial$ 
is $\le 2^\theta$.
\endproclaim
\bigskip

\remark{Remark}  Part (1) improves \scite{734-11.5}(3).
\endremark
\bigskip

\demo{Proof}  1) Clearly  $(b) \Rightarrow (a)$, so assume
clause (a) holds.  As $M \in K^{**}_{\ge\kappa}$ without loss of
generality there is a $\kappa$-wide linear order $I$ 
such that $M = \text{ EM}_{\tau({\frak K})}(I,\Phi)$;
hence for some $J \subseteq I,|J| = \theta$ we have $\bar a,\bar b 
\in {}^\theta(\text{EM}_{\tau({\frak K})}(J,\Phi))$.  So
for every $\alpha < (2^\theta)^+$, by the hence and forth argument for
$\Bbb L_{\infty,\beth^+_\alpha}[{\frak K}]$ there are
$J_\alpha,f_\alpha$ such that $J \subseteq J_\alpha \subseteq I,
|J_\alpha| = \beth_\alpha$ and $f_\alpha$
is an automorphism of EM$_{\tau({\frak K})}
(J_\alpha,\Phi)$ which maps $\bar a$ to $\bar b$. 
Hence as in the proof of \scite{734-11.5} there is a linear order $J^+$ of
cardinality $\mu$ extending $J$ and an automorphism $f$ of $M^+ =
\text{ EM}_{\tau({\frak K})}(I^+,M)$ mapping $\bar a$ to $\bar b$.
By clause (b) of Claim \scite{734-11.3} we are done.
\nl
2) Easy by clause (c) of \scite{734-11.3}, i.e., by \scite{734-11.2A}.
  \hfill$\square_{\scite{734-11.11}}$
\enddemo
\bigskip

\proclaim{\stag{734-11.14} Claim}  Assume:
\mr
\item "{$(a)$}"  $I_1 \subseteq I_2,I_1 \ne I_2$, moreover $I_1
  <_{K^{\text{flin}}} I_2$, see Definition \scite{734-0n.2}(6)
\sn
\item "{$(b)$}"  $M_\ell = { \text{\rm EM\/}}_{\tau({\frak K})}
(I_\ell,\Phi)$ for $\ell=1,2$
\sn
\item "{$(c)$}"  $\bar b,\bar c \in {}^\alpha(M_2)$
\sn
\item "{$(d)$}"  $\theta \ge |\alpha| + { \text{\rm LS\/}}({\frak K})$
\sn
\item "{$(e)$}"   $\kappa = \beth_{1,1}(\theta_2) \le \lambda$ where
$\theta_1 = 2^\theta,\theta_2 = (2^{\theta_1})^+$
\sn
\item "{$(f)$}"  $|I_1| \ge \kappa$
\sn
\item "{$(g)$}"  $M_1 \le_{\frak K} M_2$, follows from (a) + (b)
\ermn
1) Assume that for every $\bar a \in {}^{\kappa >}(M_1)$ the sequences
$\bar a \char 94 \bar b,\bar a \char 94 \bar c$ realize the same 
$\Bbb L_{\infty,\kappa} [{\frak K}]$-type in $M_2$.  \ub{Then} there
are $I_3,M_3$ and $f$ such that $I_2 \le_{K^{\text{flin}}} I_3 \in
K^{\text{flin}}_\lambda,M_3 = \text{\rm
EM}_{\tau({\frak K})}(I_3,\Phi)$ and $f$ an automorphism of $M_3$ over
$M_1$ mapping $\bar b$ to $\bar c$. 
\nl
2) Assume that for every $\bar a \in {}^{\kappa >}(M_1)$ the sequences
$\bar a \char 94 \bar b,\bar a \char 94 \bar b$ realize the same 
$\Bbb L_{\infty,\kappa}[{\frak K}]$-type in $M_2$ (as in part (a)) and
$\beth_{1,1}(\partial) \le |I_1|$ and $\partial < \lambda$.  
\ub{Then} for every $\bar a \in
{}^{\kappa >}(M_1)$, the sequences $\bar a \char 94 \bar b,\bar a
\char 94 \bar c$ realize the same $\Bbb L_{\infty,\partial}[{\frak
K}]$-type in $M_2$.
\nl
3) Assume that {\rm cf}$(\lambda) = \aleph_0$ and $|I_1| = \lambda$
and recall $\lambda = \beth_\lambda > { \text{\rm LS\/}}({\frak K})$.
If $M_1 \le_{\frak K} M^*_2 \in K^*_\lambda$ \ub{then} for some $I_3$,
a linear order $\le_{K^{\text{flin}}_\lambda}$-extending $I_2$ the model
$M^*_2$ can be $\le_{\frak K}$-embedded into
$M_3 := \text{\rm EM}_{\tau({\frak K})}(I_3,\Phi)$ over $M_1$.
\endproclaim
\bigskip

\remark{Remark}  1) Under mild assumptions with somewhat more work 
in \scite{734-11.14}(1),(3) we can 
choose $I_3 = I_2$ (but for this has to be more careful with 
the linear orders).
Recall that for $I \in K^{\text{lin}}_\lambda$ like $I_2$ in
\scite{734-11.1.7}(c) we have $\alpha < \lambda^+ \Rightarrow I \times 
\alpha$ can be embedded into $I$ and \scite{734-11.1}(1)(d). 
\endremark
\bigskip

\demo{Proof}  1) There is $J_2 \subseteq I_2$ of cardinality 
$\le \theta$ such that $\bar b,\bar c \in {}^\alpha(\text{EM}_{\tau({\frak K})}
(J_2,\Phi))$; let $J_1 = I_1 \cap J_2$.

We define a two-place relation ${\Cal E}$ on $I_2 \backslash J_2:s
{\Cal E} t$ \ub{iff} $(\forall x \in J_2)(x <_{I_2} s \equiv x <_{I_2} t)$.  
Clearly ${\Cal E}$ is an equivalence relation.
As $I_1 <_{K^{\text{flin}}} I_2$ clearly
\mr
\item "{$\odot_1$}"  $(\alpha) \quad$ any interval of $I_1$ has
cardinality $|I_1| \ge \kappa$
\sn
\item "{${{}}$}"  $(\beta) \quad$ for every $t \in I_2 \backslash J_2$
the equivalence class $t/{\Cal E}$ is a singleton or has
\nl

\hskip25pt $|I_2| \ge \kappa$ members,
\sn
\item "{${{}}$}"  $(\gamma) \quad$ for every $t \in I_1 \backslash
J_1,(t/{\Cal E}) \cap I_1$ is a singleton or has $|I_1| \ge\kappa$
members
\sn
\item "{${{}}$}"  $(\delta) \quad I_1 \backslash J_2$
has at least $\kappa$ elements
\sn
\item "{${{}}$}"  $(\varepsilon) \quad {\Cal E}$ has 
$\le 2^{|J_2|} \le 2^\theta$ equivalence classes
\sn
\item "{${{}}$}"  $(\zeta) \quad$ we may
$\le_{K^{\text{flin}}}$-increase $I_2$ so \wilog 
{\roster
\itemitem{ ${{}}$ }  $(*)_1 \quad t \in I_2 \backslash J_2 
\Rightarrow |t/{\Cal E}| = |I_2|$
\sn
\itemitem{ ${{}}$ }  $(*)_2 \quad$ for every $t \in I_1$ for some
$s_1,s_2 \in I_2$ we have $s_1 <_{I_2} t <_{I_2} s_2$
\nl

\hskip45pt  and $(s_1,t_{I_2}),(t,s_2)_{I_2}$ are disjoint to $I_1$.
\endroster}
\ermn
Let $\langle {\Cal U}_i:i < i(*)\rangle$ 
list the equivalence classes of ${\Cal E}$, so without loss of
generality $i(*) \le 2^\theta$.
For $\ell = 0,1$ let $u_\ell = \{i < i(*):{\Cal U}_i \cap I_1$ has
exactly $\ell$ members$\}$ and let $u_2 = i(*) \backslash u_0
\backslash u_1$, so by clause $\odot_1(\gamma)$, i.e. the 
definition of $I_1 \in K^{\text{flin}}$ we
have $i \in u_2 \Rightarrow |{\Cal U}_i \cap I_1| = |I_1| \ge
\kappa$.  For $i \in u_1$ let $t^*_i$ be the unique member of ${\Cal
U}_i \cap I_1$.  

Without loss of generality $u_1 = \{i:i \in [j^*_0,j^*_1)\}$ for some
$j^*_0 \le j^*_1 \le i(*)$ and let $i'(*) = i(*) + (j^*_1 - j^*_0)$ 
and $u'_1 = [i(*),i'(*))$ and define ${\Cal U}'_i$ for $i < i'(*)$ by
\mr
\item "{$\odot_2$}"  $(a) \quad {\Cal U}'_i = {\Cal U}_i$ if $i \in
u_0 \cup u_2$
\sn
\item "{${{}}$}"  $(b) \quad {\Cal U}'_i = \{t \in {\Cal U}_i:t <
t^*_i\}$ if $i \in u_1$ and
\sn
\item "{${{}}$}"  $(c) \quad {\Cal U}'_i = \{t \in {\Cal U}_\iota:
t^*_\iota <_{I_2} t\}$ if $i \in [i(*),i'(*)],\iota \in
(j^*_0,j^*_1)$ and
\nl

\hskip25pt  $i-i(*) = \iota - j^*_0$.
\ermn
For $i < i'(*)$ let $\langle t_{i,\alpha}:\alpha <
\kappa \rangle$ be a sequence of pairwise distinct members of 
${\Cal U}'_i$ such that 
$i \in u_2 \Rightarrow t_{i,\alpha} \in I_1$ and $i \in u_0 \Rightarrow
t_{i,\alpha} \notin I_1$, this actually follows.  By $\odot_1(\zeta)$
and $\odot_1(\beta),(\gamma)$ we can find such $t_{i,\alpha}$'s.

For $\zeta < \theta_2$ (see clause (e) of the assumption so 
$\beth_\zeta < \kappa$) let
$J_{1,\zeta} = \{t_{i,\alpha}:i \in u_2,\alpha < \beth_\zeta\} \cup
J_1 \cup \{t^*_i:i \in u_1\}$.
Now by the hence and forth argument (or see \scite{734-0n.13}(2)) for 
each $\zeta < \theta_2$, 
there are $J_{2,\zeta}$ and $f_\zeta$ such that
$J_{2,\zeta} \subseteq I_2$ is of cardinality $\beth_\zeta$, it
includes $J_{1,\zeta} \cup J_2$ and also $\{t_{i,\alpha}:i < i'(*)$ and
$\alpha < \beth_\zeta\}$ and $f_\zeta$ is an automorphism of
EM$_{\tau({\frak K})}(J_{2,\zeta},\Phi)$ over EM$_{\tau({\frak
K})}(J_{1,\zeta},\Phi)$ mapping $\bar b$ to $\bar c$.

(Why?  Let $\bar a_0$ list EM$(J_{1,\zeta},\Phi)$ 
so $\bar a_0 \char 94 \bar b,\bar a_0
\char 94 \bar c$ realize the same $\Bbb L_{\infty,\beth^+_\zeta}[{\frak
K}]$-type in $M_2$, and $f$ be the mapping taking $\bar a_0 \char 94
\bar b$ to $\bar a_0 \char 94 \bar c$, etc.)
\nl 
Now we shall immitate the proof of \scite{734-11.5}.  By renaming
\wilog \, there is a transitive set $B$ (in the set theoretic sense)
of cardinality $\le \theta_1 = 2^\theta$ which includes
\mr
\item "{$\oplus(a)$}" $J_1,J_2$ 
\sn 
\item "{$(b)$}" $\Phi$ (i.e.
$\tau_\Phi$ and $\langle (\text{EM}(n,\Phi),a_\ell) _{\ell < n}:n <
\omega \rangle)$
\sn
\item "{$(c)$}" ${\frak K}$, i.e.,
$\tau_{\frak K}$ and $\{(M,N):M \le_{\frak K} N$ have universe
included in $\text{\rm LS}({\frak K})\}$
\sn
\item "{$(d)$}" $\langle t^*_i:i \in u_1\rangle$ 
so each $t^*_i$ for $i \in u_1$
\sn
\item "{$(e)$}" the ordinal $i(*)$.
\ermn 
Let $\chi$ be large enough, let ${\frak B} = ({\Cal H}(\chi),\in,
<^*_\chi)$ and let ${\frak B}^+_\zeta$ be ${\frak B}$ expanded by 
\mr 
\item "{$\circledast_1$}" $(a) \quad Q^{{\frak
B}_\zeta} = \{\alpha:\alpha < \beth_\zeta\}$ 
\sn 
\item "{${{}}$}" $(b) \quad
P^{{\frak B}_\zeta}_i = J_{2,\zeta} \cap {\Cal U}'_i$ for $i<i'(*)$ 
\sn 
\item "{${{}}$}" $(c) \quad F^{{\frak B}_\zeta}_2(t) = a_t$ for $t \in
I_2$ 
\sn 
\item "{${{}}$}" $(d) \quad H^{{\frak B}_\zeta} = f_\zeta$ and $Q^{{\frak
B}_\zeta}_1 = J_{1,\zeta},Q^{{\frak B}_\zeta}_2 = J_{2,\zeta}$ 
\sn 
\item "{${{}}$}" $(e) \quad$ for $i < i'(*),H^{{\frak B}_\zeta}_i$ is the
function mapping $\alpha < \beth_\zeta$ to $t_{i,\alpha}$
\sn
\item "{${{}}$}"  $(f) \quad$ individual constants for $B$ and
for each $x \in B$,  
\nl

\hskip25pt hence, e.g. for $t^*_i(i \in u_1),J_1,J_2,t$ for $t \in
J_2$
\sn
\item "{${{}}$}"  $(g) \quad$ individual constants $J_{1,*},J_{2,*}$
interpreted as the linear orders $J_{1,\zeta},J_{2,\zeta}$
\nl

\hskip25pt respectively and individual constants for 
$M^+_\ell = \text{\rm EM}(J_{0,\zeta},\Phi)$, and
\nl

\hskip25pt $\langle a_t:t \in I_\ell\rangle$ for $\ell=1,2$.
\ermn
As in the proof of \scite{734-11.5} there is a 
$\tau({\frak B}^+)$-model ${\frak C}$, such that
\mr
\item "{$\boxtimes$}"  $(a) \quad$ for some unbounded $S \subseteq
\theta_2$
{\roster
\itemitem{ ${{}}$ }  $(\alpha) \quad {\frak C}$ is a first order
elementarily equivalent to ${\frak B}^+_\zeta$ for every $\zeta \in S$
\sn
\itemitem{ ${{}}$ }  $(\beta) \quad {\frak C}$ omits every type 
omitted by ${\frak B}_\zeta$ for every $\zeta \in S$.
\nl

\hskip35pt In particular this gives
\sn
\itemitem{ ${{}}$ }  $(\gamma) \quad {\frak C}$ omits the type 
$\{x \ne b \wedge x \in B:b \in B\}$ so
\sn
\itemitem{ ${{}}$ }  $(\delta)$ \wilog \, $b \in 
B \Rightarrow b^{\frak C} = b$ 
\endroster}
\item "{${{}}$}"  $(b) \quad {\frak C}$
is the Skolem hull of some infinite indiscernible sequence 
\nl

\hskip25pt $\langle y_r:r \in I\rangle$, where $I$ an infinite linear 
order and $y_r \in Q^{\frak C}$ for $r \in I$.
\ermn
Without loss of generality $I \in K^{\text{flin}}$ and $I_2$ 
can be $\le_{K^{\text{flin}}}$-embedded into $I$ say by the
function $g$ such that $(\forall t \in I_2)(\exists s_1,s_2 \in I)[s_1
<_I g(t) <_I s_2 \wedge (\forall t' \in I_2)(t' <_{I_2} t \rightarrow
g(t') <_I s_1) \wedge (\forall t' \in I_2)(t <_{I_2} t' 
\rightarrow s_2 <_I g(t'))]$; and also
$\|{\frak C}\| = |I|$.   Hence for each $i < i'(*)$ there is an embedding
$h_i$ of the linear order ${\Cal U}'_i$, i.e., 
$I_2 \restriction {\Cal U}'_i$ into
$(P^{\frak C}_i,(<_{I_2})^{\frak C})$ such that $t \in {\Cal U}'_i
\Rightarrow (t \in I_1 \leftrightarrow h_i(t) \in Q_1^{\frak C})$.   
\nl
[Why?  \ub{Case 0}:  $i \in u_0$.

Trivial.
\nl
\ub{Case 1}:  $i \in u_1 \cup u'_1$.

Similar to Case 0 as ${\Cal U}'_i \cap I_1 = \emptyset$, of course, we
take care that $a = h_i(t) \wedge t \in {\Cal U}'_i \wedge i \in u_1
\Rightarrow {\frak C} \models ``a <_{I_2} t^*_i"$ and similarly for $u_{-1}$.
\nl
\ub{Case 2}: $i \in u_2$.

First approximation is $h'_i = (H^{\frak C}_i \circ (g \restriction
{\Cal U}_i))$, so $t \in {\Cal U}_i \Rightarrow h'_i(t) \in
Q_1^{\frak C}$.  However by the choice of $g$ we can find $\langle
(s^-_t,s^+_t):t \in {\Cal U}_i\rangle$ such that:
\mr
\item "{$(\alpha)$}"   $s^-_t,s^+_t \in Q_2^{\frak C}$
\sn
\item "{$(\beta)$}"   $(s^-_t,s^+_t)_{I^{\frak C}_2} \cap Q_2^{\frak C} =
\{h'_i(t)\}$.
\ermn
As $I_2$ is dense with no extremal members (being from
$K^{\text{flin}}$) clearly $t_1 <_{I_2 \restriction {\Cal U}'_i} t_2
\Rightarrow s^+_{t_1} <_{(I_2)^{\frak C}} s^-_{t_2}$.  Now choose
$h_i$ by: $h_i(t)$ is $h'_i(t)$ if $t \in I_1$ and is $s^+_{t_1}$ if
$t \in I_1 \backslash I_2$.] 

Hence there is an
embedding $h$ of the linear order $I_2$ into $J_{1,*}^{\frak C}$ such that:
\mr
\item "{$\circledast_2$}"  $h(t)$ is:
{\roster
\itemitem{ $(a)$ }  $t$ if $t \in J_2 \cup \{t^*_i:i \in u_1\}$
\sn
\itemitem{ $(b)$ }  $h_i(t)$ if $t \in {\Cal U}'_i$ and $i<i'(*)$.
\endroster}
\ermn
Note
\mr
\item "{$\circledast_3$}"  for every $t \in I_2 \backslash J_2$ for some
$i < i(*) \le \theta_1$ we have $(\forall s \in J_2)[s <_{I_2} t \equiv s 
<_{I_2} h_i(t_{i,0})]$
\ermn
hence by the omitting type demand in $\boxtimes(a)(\beta)$:
\mr
\item "{$\circledast'_3$}"  for $t \in I^{\frak C}_2 \backslash J_2$ for some
$i < i(*)$ we have $(\forall s \in J_2)[s <_{I^{\frak C}_2} t \equiv s 
<_{I^{\frak C}_2} (h_i(t_{i,0}))]$. 
\ermn
We can find a linear order $I_3,I_2 \subseteq I_3$ and an isomorphism
$h_*$ from $I_3$ onto $Q_2^{\frak C}$ extending $h$, so clearly
$I_3 \in K^{\text{flin}}$ and \wilog \, $h(I_2) <_{K^{\text{flin}}} I_3$.
Now let $\hat h_*$ be
the isomorphism which $h_*$ induces from EM$_{\tau({\frak K})}(I_3,\Phi)$ onto
$(\text{\rm EM}_{\tau({\frak K})}(J^{\frak C}_{2,*},\Phi))^{\frak C}$,
so e.g., it maps for each $t \in I_2$,
the member $a_t$ of the skeleton to $F^{\frak C}_2(h_*(t))$.

Note that $h_*$ maps ${\Cal U}_i \cap I_1$ into $Q^{\frak C}_1
\subseteq I^{\frak C}_1$ when ${\Cal U}_i \subseteq I_1$ and is the
identity on $J_1 \cup \{t^*_i:i \in u_1\}$ so recalling $Q^{{\frak
B}_\zeta} = J_{1,\zeta} = \{t_{i,\alpha}:i \in u_2$ and $\alpha <
\beth_\zeta\} \cup J_1 \cup\{t^*_i:i \in u_1\}$ 
hence it map $I_1$ into $Q^{\frak C}_1$ but ${\frak B}_\zeta \models
``H$ is a unary function, an automorphism of EM$_{\tau({\frak K})}
(J^{\frak C}_{2,*},\Phi)$ mapping $\bar b$ to $\bar c$ and is the identity on 
EM$_{\tau({\frak K})}(J^{\frak C}_{1,*},\Phi)"$.
 Now $(\hat h_*)^{-1}H^{\frak C}(\hat h_*)$ is an
automorphism of EM$_{\tau({\frak K})}(I_3,\Phi)$ as required.
\nl
2) By part (1), i.e. choose $I_3,M_3,f_3$ as there; so as $f$ is an
automorphism of $M_3$ over $M_1$ mapping $\bar b$ to $\bar c$,
clearly $\bar b,\bar c$ realize the same 
$\Bbb L_{\infty,\partial}[{\frak K}]$-type over $M_1$ inside $M_3$.  The
desired result (the type inside $M_2$ rather than inside $M_3$) follows
because $M_1 \prec_{\Bbb L_{\infty,\partial}[{\frak K}]} M_2
\prec_{\Bbb L_{\infty,\partial}[{\frak K}]} M_3$  by \scite{734-11.3}(a).
\nl 
3) Let $M^*_2 = \dbcu_{n < \omega} M^*_{2,n}$ be such that $n < \omega
\Rightarrow M^*_{2,n} \le_{\frak K}
M^*_{2,n+1}$ and $\|M^*_{2,n}\| < \lambda$.  Let 
$\bar c_n$ list $M^*_{2,n}$ for $n < \omega$ (with no repetitions) 
and be such that $\bar c_n
\triangleleft \bar c_{n+1}$.  Let $\theta_n = \|M^*_{2,n}\| +
\text{ LS}({\frak K})$ so \wilog \, $\theta_n = \ell g(\bar c_n)$ and
let $\theta'_n = \beth_3(\theta_n),\kappa_n = \beth_{1,1}(\theta'_n)$, 
\wilog \, $\kappa_n < \theta_{n+1}$ and we choose for each
$n < \omega$, a sequence $\bar b_n \in {}^{\ell g(\bar c_n)}(M_2)$ realizing 
\sftp$_{{\Bbb L}_{\infty,\kappa^+_n}[{\frak K}]}
(\bar c_n,M_1,M^*_2)$ in $M_2$.  
This is possible by \scite{734-11.10}(3) after possibly 
$<_{K^{\text{flin}}}$-increasing $I_2$.

Now we choose $(I_{3,n},f_n,M_{3,n},\bar b'_n)$ by induction on $n$ such that
\mr
\item "{$(*)$}"  $(a) \quad I_{3,0} = I_2$ and $I_{3,n} \in
K^{\text{lin}}_\lambda$
\sn
\item "{${{}}$}"  $(b) \quad n=m+1 \Rightarrow I_{3,m} <_{K^{\text{flin}}}
I_{3,n}$
\sn
\item "{${{}}$}"  $(c) \quad M_{3,n} = \text{\rm EM}_{\tau({\frak
K})}(I_{3,n},\Phi)$ (hence $n=m+1 \Rightarrow M_{3,m} \le_{{\frak
K}_\lambda} M_{3,n}$)
\sn
\item "{${{}}$}"  $(d) \quad f_n$ is an automorphism of $M_{3,n}$ over
$M_1$
\sn
\item "{${{}}$}"  $(e) \quad \bar b'_n \in {}^{\ell g(\bar
b_n)}(M_{3,n})$ realizes tp$_{\Bbb L_{\infty,\kappa^+_n}[{\frak K}]}
(\bar c_n,M_1,M^*_2)$
\sn
\item "{${{}}$}"  $(f) \quad$ if $n = m+1$ then $\bar  b'_m
\trianglelefteq \bar b'_n$
\sn
\item "{${{}}$}"  $(g) \quad$ if $n=m+1$ then $f_n$ maps $\bar b_{n+1}
\restriction \ell g(\bar b_n)$ to $\bar b'_n$ and $f_0$ 
maps $\bar b_0$ to $\bar b'_0$.
\ermn
For $n=0,I_{3,0},M_{3,0}$ are defined in clauses (a),(c) of $(*)$ and
we let $f_0 = \text{ id}_{M_2} = \text{ id}_{M_{3,n}},\bar b'_0
= \bar b_0$ this is trivially as required.  For $n=m+1$ we apply part
(1) with 
\mr
\item "{$\boxdot$}"  $I_1,I_{3,m},M_1,M_{3,m},\bar b_{n+1} \restriction 
\ell g(\bar c_m),\bar b'_m,\theta_m,\kappa_m$ here
\nl

standing for 
$I_1,I_2,M_1,M_2,\bar b,\bar c,\theta,\kappa$ there.
\ermn
Why its assumptions holds?  The main point is to check that for every
$\bar a \in {}^{\kappa_m>}(M_1)$ the sequences $\bar a \char 94 (\bar
b_{n+1} \restriction \theta_m),\bar a \char 94 \bar b'_m$ realize the
same $\Bbb L_{\infty,\kappa_m}[{\frak K}]$-type in $M_{3,m}$.  Now
$\bar a \char 94 (\bar b_{m+1} \restriction \theta_m),
\bar a \char 94 \bar b'_m$ realize the same 
$\Bbb L_{\infty,\kappa_n}[{\frak K}]$-type in $M_{3,m}$ by
the induction hypothesis.  Also the sequences $\bar b_{n+1} \restriction
\theta_m,\bar b_{m+1} \restriction \theta_m$ satisfy for any 
$\bar a \in {}^{\kappa_m}(M_1)$ the sequences $\bar a \char 94 (\bar
b_{n+1} \restriction \theta_m),\bar a \char 94 (\bar b_{m+1}
\restriction \theta_m)$ realize the same $\Bbb
L_{\infty,\kappa_m}[{\frak K}]$-type in $M_{3,m}$ because the 
$\Bbb L_{\infty,\kappa_m}[{\frak K}]$-type which $\bar a \char 94 (\bar
b_{n+1} \restriction \theta_m)$ realizes in $M_{3,m}$ is the same as the
$\Bbb L_{\infty,\kappa_m}[{\frak K}]$-type it realizes in $M_2 =
M_{3,0}$ which (by the choice of $\bar b_{n+1}$) is equal to the
$\Bbb L_{\infty,\kappa_m}[{\frak K}]$-type which $\bar a \char 94(\bar
c_{n+1} \restriction \theta_m)$ realizes in $M^*_2$ which is the same
as the $\Bbb L_{\infty,\kappa_m}[{\frak K}]$-type which $\bar a \char
94(\bar c_{m+1} \restriction \theta_m)$ realizes in $M^*_2$ which is
equal to the $\Bbb L_{\infty,\kappa_m}[{\frak K}]$-type which $\bar a
\char 94 (\bar b_{m+1} \restriction \theta_m)$ realizes in $M_{3,m}$.

By the last two sentences for every $\bar a \in {}^{\kappa_m >}(M_1)$
the sequences $\bar a \char 94 (\bar b_{n+1} \restriction
\theta_m),\bar a \char 94 \bar b'_m$ realizes the same $\Bbb
L_{\infty,\kappa_m}[{\frak K}]$-type in $M_{3,m}$, so indeed the
assumptions of part (1) holds for the case we are trying to apply it,
see $\boxdot$ above.

So we get the conclusion of part (1), i.e. we 
get $I_{3,n},f_n$ here standing for $I_3,f$ there so $I_{3,m}
<_{K^{\text{flin}}_\lambda} I_{3,n}$ and $f_n$ is an automorphism of
$M_{3,n} = \text{ EM}_{\tau({\frak K})}(I_{3,n},\Phi)$ over $M_1$
mapping $\bar b_{n+1} \restriction \theta_m$ to $\bar b'_m$.  Now 
we let $\bar b'_n = f_n(\bar b_{n+1} \restriction \theta_n)$ and can
check all the clauses in $(*)$.  Hence we have carried the induction.
So we can satisfy $(*)$.

So $\bar b'_n$ satisfies the requirements on
$\bar b_n$ and $\bar b'_n \triangleleft \bar b'_{n+1}$.
Let $I_3 = \cup\{I_{3,n}:n < \omega\}$ and let $M_3 = \text{
EM}_{\tau({\frak K})}(I_3,\Phi)$ and let
$g:M^*_2 \rightarrow M_3$ map $c_{n,i}$ to $b'_{n,i}$ for $i < \ell
g(\bar c_n),n < \omega$, easily it is as required. 
That is, $g(c_{n,i})$ is well defined as $c_{n,i} \mapsto b'_{n,i},(i
< \ell g(\bar c_n))$ is a well defined mapping for each $n$ and
$i < \ell g(\bar c_n)
\Rightarrow c_{n,i} = c_{n+1,i} \wedge b'_{n,i} = b'_{n+1,i}$.  Also
$g \restriction \{c_{n,i}:i < \ell g(\bar c_n)\}$ is a $\le_{\frak
K}$-embedding of $M^*_{2,n}$ into $M_3$ and is the identity on
$M^*_{2,n} \cap M_1$ as $\bar c_n$ list the
elements of $M_{2,i}$ and $\sftp_{\Bbb L_{\infty,\kappa^+_n}[{\frak
K}]} (\bar c_n,M_1,M^*_2) = \sftp_{\Bbb L_{\infty,\kappa^+_n}[{\frak
K}]} (\bar b'_n,M_1,M_3)$ by clause (e) of $(*)$.  
But $\langle g \restriction M^*_{2,n}:n < \omega\rangle$ is 
$\subseteq$-increasing with union $g$ so by Ax(V) of a.e.c.
$g$ is a $\le_{\frak K}$-embedding of $M^*_2$ into $M_3$.  
Lastly, obviously $g \supseteq \cup\{\text{id}_{M^*_{2,n} \cap M_1}:n
< \omega\} =  \text{ id}_{M_1}$, so we are done.
\hfill$\square_{\scite{734-11.14}}$
\enddemo
\bn
We arrive to the crucial advance:
\proclaim{\stag{734-11.15} The Amalgamation Theorem}  If 
{\rm cf}$(\lambda) = \aleph_0$, \ub{then}
${\frak K}^*_\lambda$, i.e., $(K^*_\lambda,\le_{\frak K} \restriction
K^*_\lambda)$ has amalgamation, even disjoint one.
\endproclaim
\bigskip

\demo{Proof}   So assume 
$M_0 \le_{{\frak K}^*_\lambda} M_\ell$ for $\ell=1,2$. 
Choose $I_0 \in K^{\text{flin}}_\lambda$ so $M'_0 := 
\text{ EM}_{\tau({\frak K})}(I_0,\Phi) \in K^*_\lambda$ but $K^*_\lambda$ is
categorical (see \scite{734-11.4A} or \scite{734-11.5}(4)) hence $M'_0 \cong
M_0$, so without loss of generality $M'_0 = M_0$.
Choose $I_1 \in K^{\text{flin}}_\lambda$ such that $I_0
<_{K^{\text{flin}}} I_1$ and let $M'_1 = \text{ EM}_{\tau({\frak
K})}(I_1,\Phi)$ so $M_0 \le_{\frak K} M'_1$.  
By applying \scite{734-11.14}(3) with $I_0,I_1,M_0,M'_1,M_1$ here
standing for $I_1,I_2,M_1,M_2,M^*_2$ there, we can find a 
pair $(I_2,f_1)$ such that $I_1 <_{K^{\text{flin}}_\lambda} I_2$
and $f_1$ is a $\le_{\frak K}$-embedding of $M_1$ into $M'_2 :=
\text{ EM}_{\tau({\frak K})}(I_2,\Phi)$ over $M_0$.  Apply
\scite{734-11.14}(3) again with $I_0,I_2,M_0,\text{EM}_{\tau({\frak
K})}(I_2,\Phi),M_2$ here standing for $I_1,I_2,M_1,M_2,M^*_2$ there.
So there is a pair $(I_3,f_2)$ such that
$I_2 <_{K^{\text{flin}}_\lambda} I_3$ and $f_2$ is $\le_{\frak
K}$-embedding $M_2$ into $M_3 := \text{ EM}_{\tau({\frak K})}
(I_3,\Phi)$ over $M_0 = \text{ EM}_{\tau({\frak
K})}(I_0,\Phi)$.  Of course, $M_3 \in K^*_\lambda$ and we are done
proving the ``has amalgamation".

Why disjoint?  Let $(I_4,h)$ be such that $I_3 
<_{K^{\text{flin}}_\lambda} I_4$ and $h$ is a
$\le_{K^{\text{flin}}}$-embedding of $I_3$ into $I_4$ over $I_0$ such
that $h(I_3) \cap I_3 = I_0$.  Now $h$ induces an
isomorphism $\hat h$ from EM$_{\tau({\frak K})}(I_3,\Phi)$ onto
EM$_{\tau({\frak K})}(h(I_3),\Phi) \le_{\frak K} M_3$.

Lastly, by our assumptions on $\Phi$ if $J_1,J_2 \subseteq J$ and 
$I_1 \cap I_2$ is a dense linear order
(in particular with neither first nor last member, e.g. are from
$K^{\text{flin}}_\lambda$ as in our case) then 
EM$_{\tau({\frak K})}(I_1,\Phi) \cap \text{ EM}_{\tau({\frak K})}
(I_2,\Phi) = \text{ EM}_{\tau({\frak K})}(I_1 \cap I_2,\Phi)$.
So in particular, above EM$_{\tau({\frak K})}(I_3,\Phi) \cap \text{
EM}_{\tau({\frak K})}(\hat h(I_3,\Phi) = \text{ EM}_{\tau({\frak
K})}(I_0,\Phi)$ and $f_1,\hat h \circ f_2$ are $\le_{\frak
K}$-embeddings of $M_1,M_2$ respectively over $M_0 = \text{
EM}_{\tau({\frak K})}(I_0,\Phi)$ into EM$_{\tau({\frak K})}(I_3,\Phi)
\le_{\frak K} \text{ EM}_{\tau({\frak K})}(I_4,\Phi)$ and
EM$_{\tau({\frak K})}(h(I_3),\Phi) \le_{\frak K} 
\text{ EM}_{\tau({\frak K})}(I_4,\Phi)$, respectively, so 
we are done.    \hfill$\square_{\scite{734-11.15}}$
\enddemo
\bigskip

\proclaim{\stag{734-11.16} Claim}  Assume {\rm cf}$(\lambda) = \aleph_0$.
If $\delta < \lambda^+$, the sequence
$\langle M_i:i < \delta \rangle$ 
is $\le_{\frak K}$-increasing continuous and $M_i \in K^*_\lambda$
for $i < \delta$, \ub{then} $M_\delta := \cup\{M_i:i < \delta\}$ can be
$\le_{\frak K}$-embedded into some member of $K^*_\lambda$.
\endproclaim
\bigskip

\demo{Proof}   We choose $I_i \in K^{\text{flin}}_\lambda$
by induction on $i \le \delta$, which is 
$<_{K^{\text{flin}}_\lambda}$-increasing continuous with $i$ and
 a $\le_{\frak K}$-embedding $f_i$ of $M_i$ into $N_i 
:= \text{\rm EM}_{\tau({\frak K})}(I_i,\Phi)$,
increasing continuous with $i$.  For $i=0$ choose $I_0 \in
K^{\text{flin}}_\lambda$, so $N_0 := \text{ EM}_{\tau({\frak
K})}(I_0,M)$ is isomorphic to $M_0$ hence $f_0$ exists; 
for $i$ limit use $I_i := \cup\{I_j:j<i\}$ and $f_i :=
\cup\{f_j:j < i\}$.  So assume $i =j+1$.
Now we can find $M'_i,f'_i$ satisfying:
$f'_i$ is an isomorphism from $M_i$ onto $M'_i$ extending $f_j$
such that $f_j(M_j) \le_{\frak K} M'_i$ (actually this trivially follows) 
and $M'_i \cap N_j = f_j(M_j)$; so
also $M'_i$ belongs to $K^*_\lambda$.  Now
$f_j(M_j),\text{EM}_{\tau({\frak K})}(I_j,\Phi),
M'_i$ can be disjointly amalgamated (by \scite{734-11.15}) in $(K^*_\lambda,\le_{\frak K})$, so there is $M^*_i \in
K^*_\lambda$ such that $N_j = \text{\rm EM}_{\tau({\frak K})}(I_j,\Phi) 
\le_{\frak K} M^*_i$ and $M'_i \le_{\frak K} M^*_i$.  
Now by \scite{734-11.14}(3) there are $I_i,g_i$ such that $I_j
<_{K^{\text{flin}}_\lambda} I_i$ and $g_i$ is a 
$\le_{\frak K}$-embedding of $M^*_i$ into $N_i := 
\text{\rm EM}_{\tau({\frak K})} (I_i,\Phi)$ over
EM$_\tau(I_j,\Phi)$.  
Let $f_i = g_i \circ f'_i$, clearly it is as required.
Having carried the induction, $f_\delta$ is a 
$\le_{\frak K}$-embedding of $M_\delta$ into EM$_{\tau({\frak K})}
(\dbcu_{j < \delta} I_j,\Phi)$, as promised.
\hfill$\square_{\scite{734-11.16}}$
\enddemo
\bigskip

\proclaim{\stag{734-11.16A} Claim}  1) Assume {\rm cf}$(\lambda) =
\aleph_0$.  For every $M_0 \in K^*_\lambda$ 
there is a $\le_{\frak K}$-extension $M_1 \in K^*_\lambda$ of $M_0$ 
such that: if $M_0 \le_{{\frak K}_\lambda} M_2 \in
K^*_\lambda$ and $\bar a \in {}^{\lambda >}(M_2)$ \ub{then} for some 
$(M_3,f)$ we have:
\sn
\block
$M_1 \le_{\frak K} M_3 \in K^*_\lambda,f$ is a $\le_{\frak K}$-embedding
of $M_2$ into $M_3$ over $M_0$ and $f(\bar a) \in {}^{\lambda >}(M_2)$.
\endblock
\bn
2) Assume {\rm cf}$(\lambda) = \aleph_0$. 
For every $M_0 \in K^*_\lambda$ there is a $\le_{\frak K}$-extension
$M_1 \in K^*_\lambda$ which is universal over $M_0$ for 
$\le_{{\frak K}_\lambda}$-extensions. 
\nl
3) If (a) then (b) where
\mr
\item "{$(a)$}"  $I_0 \le_{K^{\text{flin}}_\lambda} I'_1
<_{K^{\text{flin}}_\lambda} I_1$ 
\sn
\item "{$(b)$}"  if $I_0 \subseteq I_2 \in K^{\text{flin}}_\lambda$
and $\beta \le \gamma < \lambda,\bar b_1 \in
{}^\beta(\text{\rm EM}_{\tau({\frak K})}(I'_1,\Phi))$ and $\bar c_2 \in
{}^\gamma(\text{\rm EM}_{\tau({\frak K})}(I_2,\Phi))$ and $\bar b_2 = \bar
c_2 \restriction \beta$ and for every $\kappa < \lambda$ we have 
$$\align\text{\rm tp}_{\Bbb L_{\infty,\kappa}[{\frak K}]}(\bar b_1,
\text{\rm EM}_{\tau({\frak K})}(I_0,\Phi)), 
&\text{\rm EM}_{\tau({\frak K})}(I_1,\Phi))\  = \cr  
& = \text{\rm tp}_{\Bbb L_{\infty,\kappa}[{\frak K}]}(\bar b_2,
\text{\rm EM}_{\tau({\frak K})}(I_0,\Phi), \text {\rm EM}_{\tau({\frak
K})}(I_2,\Phi))
\endalign
$$ \ub{then} for some $(I^+_1,f)$ we have
$I_1 \le_{K^{\text{flin}}} I^+_1 \in K^{\text{flin}}_\lambda$ and $f$
is a $\le_{\frak K}$-embedding of $\text{\rm EM}_{\tau({\frak K})}
(I_2,\Phi)$ into {\rm EM}$_{\tau({\frak K})}(I^+_1,\Phi)$ over
{\rm EM}$_{\tau({\frak K})}(I_0,\Phi)$ mapping $\bar b_2$ to $\bar
b_1$ and $\bar c_2$ into {\rm EM}$_{\tau({\frak K})}(I_1,\Phi)$.
\ermn
4) Assume {\rm cf}$(\lambda) = \aleph_0$.
If (c) then (d) and moreover (d)$^+$ when
\mr
\item "{$(c)$}"  $\langle J_\alpha:\alpha \le \omega\rangle$ is 
$<_{K^{\text{flin}}_\lambda}$-increasing, $I_0 = J_0,I_1 = J_\omega$
\sn
\item "{$(d)$}"  if $I_0 \subseteq I_2 \in K^{\text{flin}}_\lambda$
then some $f$ is a $\le_{\frak K}$-embedding of
{\rm EM}$_{\tau({\frak K})}(I_2,\Phi)$ into 
{\rm EM}$_{\tau({\frak K})}(I_1,\Phi)$ over {\rm EM}$_{\tau({\frak
K})}(I_0,\Phi)$
\sn
\item "{$(d)^+$}"  {\rm EM}$_{\tau({\frak K})}(I_1,\Phi)$ is 
$\le_{{\frak K}^*_\lambda}$-universal over {\rm EM}$_{\tau({\frak
K})}(I_0,\Phi)$. 
\endroster
\endproclaim
\bigskip

\demo{Proof}  Note that by \scite{734-11.14}(3) clearly (3) $\Rightarrow$
(1) and (4) $\Rightarrow$ (2).  So we shall prove (3) and (4). 
\nl
3) First assume $\beta = 0,\gamma=1$ so $\bar c_2 = \langle c \rangle$.
Toward contradiction assume $I_0 \subseteq I_2 \in K^{\text{lin}}_\lambda,a \in
M_2 := \text{ EM}_{\tau({\frak K})}(I_2,\Phi)$ but there is no pair
$(I^+_1,f)$ as required in clause (b).  Without loss of
generality for some $I_3$ we have $I_0 \le_{K^{\text{flin}}_\lambda}
I_2 \le_{K^{\text{flin}}_\lambda} I_3$ and $I_0
\le_{K^{\text{flin}}_\lambda} I_1 \le_{K^{\text{flin}}_\lambda} I_3$.

Let EM$(I_2,\Phi) \models ``c_2 =
\sigma(a_{t^2_0},\dotsc,a_{t^2_{n-1}})"$ where
$\sigma(x_0,\dotsc,x_{n-1})$ a $\tau_\Phi$-term, $n < \omega$ and $I_2
\models ``t^2_0 < \ldots < t^2_{n-1}"$.  Let $u = \{\ell < n:t^2_\ell
\in I_0\}$.  As $I_0 <_{K^{\text{flin}}_\lambda} I_1$, we can find
$\langle t^1_0,\dotsc,t^1_{n_1}\rangle$ such that:
\mr
\item "{$\circledast$}"  $(a) \quad t^1_\ell \in I_1$ for $\ell < n$
\sn
\item "{${{}}$}"  $(b) \quad t^1_0 <_{I_1} \ldots <_{I_1} t^1_{n-1}$
\sn
\item "{${{}}$}"  $(c) \quad$ if $\ell \in u$ then $t^2_\ell =
t^1_\ell (\in I_0)$
\sn
\item "{${{}}$}"  $(d) \quad$ if $\ell < n \wedge \ell \notin u$ then
$t^1_\ell \in I_1 \backslash I_0$
\sn
\item "{${{}}$}"  $(e) \quad$ if $\ell_1 \le \ell_2 < n$ and $[\ell_1,\ell_2]
\cap u = \emptyset$ then $t^2_{\ell_2} <_{I_3} t^1_{\ell_1}$.
\ermn
Let $M_\ell = \text{ EM}_{\tau({\frak K})}(I_\ell,\Phi)$ for
$\ell=0,1,2,3$ and let $c_2 = c$ and $c_1 =
\sigma^{\text{EM}(I_1,\Phi)}(a_{t^1_0},\dotsc,a_{t^1_{n-1}})$.

Let $\kappa < \lambda$ be large enough such that 
$\sftp_{\Bbb L_{\infty,\kappa^+}[{\frak K}]}(c_\ell,M_0,M_\ell)$ for
$\ell=1,2$ be distinct (exists by \scite{734-11.14}(1) because its conclusion
fails by the ``toward contradiction").
We easily get contradiction to the non-order property (see $(*)$ of 
\scite{734-11.1.2}(2)).

Note that if in addition $\langle I_{1,\alpha}:\alpha \le
\lambda\rangle$ is $<_{K^{\text{flin}}_\lambda}$-increasing
continuous, $I_{1,0} = I'_1,I_{1,\lambda} = I_1$ then by what we have
just proved and the proof of \marginbf{!!}{\cprefix{600}.\scite{600-4a.2}} we can prove the
general case (and part (4)).  But we also give a direct proof.

In the general case, let $\theta = |\beta| + \aleph_0$,
so we assume clause (a) and the assumptions of clause (b) and \wilog \,
$I_1 \cap I_2 = I_0$ hence there is $I_3$ such that $I_\ell
<_{K^{\text{flin}}_\lambda} I_3$ for $\ell=1,2$.  Let $\kappa \in
(\theta,\lambda)$ be large enough.

Hence EM$_{\tau({\frak K})}(I_0,\Phi) \prec_{\Bbb L_{\infty,\lambda}[{\frak
K}]} \text{ EM}_{\tau({\frak K})}(I_\ell,\Phi) \prec_{\Bbb
L_{\infty,\lambda}[{\frak K}]} \text{ EM}_{\tau({\frak K})}(I_3,\Phi)$
for $\ell=1,2$.  Applying \scite{734-11.14}(1) with $I_1,I_2,\bar b,\bar
c$ there standing for $I_0,I_3,\bar b_1,\bar b_2$ here we can find a
pair $(I_4,f_4)$ such that $I_3 <_{K^{\text{flin}}_\lambda} I_4$ and
$f_4$ is an automorphism of $M_4 := \text{ EM}_{\tau({\frak
K})}(I_4,\Phi)$ over EM$_{\tau({\frak K})}(I_0,\Phi)$
mapping $\bar b_2$ to $\bar b_1$.  
\nl
Clearly 
$M_3 := \text{ EM}_{\tau({\frak K})}(I_3,\Phi) 
\prec_{\Bbb L_{\infty,\lambda}[{\frak K}]} 
\text{ EM}_{\tau({\frak K})}(I_4,\Phi)$.  So $f_4(\bar c_2)
\in{}^\gamma(M_4)$, hence we can apply clause (b) of Claim
\scite{734-11.10}(3) with $M_1,M_2,I_2,N,\xi,\bar d^*$ there standing for
$\text{EM}_{\tau({\frak K})}(I'_1,\Phi),\text{ EM}_{\tau({\frak
K})}(I_1,\Phi),I_1$,
\nl
$\text{ EM}_{\tau({\frak K})}
(I_4,\Phi),\gamma,f_4(\bar c_2)$ here.  Hence we can find $\bar c'_2 \in
{}^\gamma(M_1)$ realizing in $M_1$ the type 
$\sftp_{\Bbb L_{\infty,\kappa}[{\frak K}]}(f_4(\bar c_2),
\text{ EM}_{\tau({\frak K})}(I'_1,\Phi),\text{ EM}_{\tau({\frak
K})}(I_1,\Phi))$.

Lastly, applying Claim \scite{734-11.14}(1) with $I_1,I_2,\bar b,\bar c$
there standing for $I'_1,I_4,f_4(\bar c_2),\bar c'_2$ here, clearly  
there is a pair
$(I_5,f_5)$ such that $I_4 <_{K^{\text{flin}}_\lambda} I_5$ and $f_5$
is an automorphism of 
$\text{EM}_{\tau({\frak K})}(I_5,\Phi)$ over
$\text{ EM}(I'_1,\Phi)$ mapping to $f_4(\bar c_2)$ to $\bar c'_2$.
\nl
Let $I^+_1 := I_5,f = f'_5 \circ f'_4$ where $f'_5 = f_5 \restriction
\text{ EM}_{\tau({\frak K})}(I_4,\Phi)),f'_4 = f_4 \restriction 
\text{ EM}_{\tau({\frak K})}(I_2,\Phi)$; now $I^+_1,f$ are as required because
$f_4(\bar b_2) = \bar b_1$ while $f_5(\bar b_1) = \bar b_1$.
\nl
4) Easy by part (3).  First note that (d)$^+$ follows by (d) by
\scite{734-11.14}(3), so we shall ignore clause (d)$^+$.
Let EM$_{\tau({\frak K})}(I_2,\Phi)$ be
$\cup\{M_{2,n}:n < \omega\}$ where $M_{2,n} \in K_{< \lambda}$ and $n
< \omega \Rightarrow M_{2,n} \le_{\frak K} M_{2,n+1}$.

Let $\bar a_n$ list the elements of $M_{2,n}$ with no repetitions such
that $\bar a_n \triangleleft \bar a_{n+1}$ for $n < \omega$.  By
induction on $n$, we choose $\bar b_n$ such that
\mr
\item "{$\circledast$}" $(a) \quad \bar b_n \in {}^{\ell g(\bar
a_n)}(\text{EM}_{\tau({\frak K})}(J_{n+1},\Phi)$
\sn
\item "{${{}}$}"  $(b) \quad$ if $n=m+1$ then $\bar b_m \triangleleft
\bar b_n$
\sn
\item "{${{}}$}"  $(c) \quad$ for every $\kappa < \lambda$ the type
tp$_{\Bbb L_{\infty,\kappa}[{\frak K}]}(\bar
b_n,\text{EM}_{\tau({\frak K})}(I_0,\Phi),\text{EM}_{\tau({\frak
K})}(I_{n+1},\Phi)))$
\nl

\hskip25pt  is equal to the type 
tp$_{\Bbb L_{\infty,\kappa}[{\frak K}]}
(\bar a_n,\text{EM}_{\tau({\frak K})}(I_0,\Phi),\text{EM}_{\tau({\frak
K})}(I_2,\Phi))$. 
\ermn
The induction step is by part (3).  Let $f_n$ be the unique function
mapping $\bar a_n$ to $\bar b_n$ (with domain Rang$(\bar a_n))$.  So
$f_n \subseteq f_{n+1}$ and $f_n$ is a $\le_{\frak K}$-embedding of
$M_{2,n}$ into EM$_{\tau({\frak K})}(J_{n+1},\Phi)$ but $J_{n+1}
\subseteq I_1$ hence into EM$_{\tau({\frak K})}(I_1,\Phi)$.  So $f :=
\cup\{f_n:n < \omega\}$ is a $\le_{\frak K}$-embedding of
EM$_{\tau({\frak K})}(I_2,\Phi)$ into 
EM$_{\tau({\frak K})}(I_1,\Phi)$.  Also $f_n$ is the identity on
Rang$(\bar a_n) \cap \text{ EM}_{\tau({\frak K})}(I_0,\Phi)$ hence $f$
is the identity on $\dbcu_n(\text{Rang}(\bar a_n) \cap 
\text{ EM}_{\tau({\frak K})}(I_0,\Phi) = \text{ EM}_{\tau({\frak
K})}(I_0,\Phi)$ so $f$ is as required.   \hfill$\square_{\scite{734-11.16A}}$
\enddemo
\bn
\margintag{734-11.16E}\ub{\stag{734-11.16E} Exercise}:  1) Assume ${\frak K}_\lambda = 
(K_\lambda,\le_{{\frak K}_\lambda})$ satisfies axioms I,II (and 0, presented
below) and amalgmation.  Then
$\ortp(a,M,N)$ for $M \le_{{\frak K}_\lambda} N$ and $a \in N$ and ${\Cal
S}_{{\frak K}_\lambda}(M)$ are well defined and has the basic
properties of types from \sectioncite[\S1]{600}.
\nl
2) If in addition ${\frak K}_\lambda$ satisfies AxIII$^\odot$ below and ${\frak
K}_\lambda$ is stable (i.e. $|{\Cal S}_{{\frak K}_\lambda}(M)| \le
\lambda$ for $M \in K_\lambda$) \ub{then} every $M \in {\frak
K}_\lambda$ has a $\le_{\frak K}$-universal extension $N$ which means
$M \le_{{\frak K}_\lambda} N$ and $(\forall N')(M \le_{{\frak
K}_\lambda} N' \rightarrow (\exists f)[f$ is a $\le_{{\frak
K}_\lambda}$-embedding of $N'$ into $N$ over $M])$.
\nl 
3) AxIII (see \marginbf{!!}{\cprefix{600}.\scite{600-0.2}}) implies AxIII$^\odot$
\nl
where: 
\mn
Ax0: $K$ is a class of $\tau_{\frak K}$-models, $\le_{\frak K}$
a two place relation of $K_\lambda$, both preserved under isomorphisms
\mn
AxI: if $M \le_{{\frak K}_\lambda} N$ then $M \subseteq N$ (are
$\tau({\frak K}_\lambda)$-models of cardinality $\lambda$
\mn
AxII: $\le_{{\frak K}_\lambda}$ is a partial order (so $M \le_{{\frak
K}_\lambda} M$ for $M \in K_\lambda$)
\mn
AxIII$^\odot$: In following game the COM player has a winning
strategy.  A play last $\lambda$ moves, they construct a
$\le_{{\frak K}_\lambda}$-increasing continuous sequence $\langle
M_\alpha:\alpha \le \lambda\rangle$.  In the $\alpha$-th move
$M_\alpha$ is chosen, by INC if $\alpha$ is even by COM is $\alpha$ is
odd.  Now Com wins as long as INC has legal moves.
\mn
AxIV$^\odot$:  For each $M \in K_\lambda$, in the 
following game, INC has no winning strategy: a play
lasts $\lambda +1$ moves, in the $\alpha$-th move $f_\alpha,M_\alpha,
N_\alpha$ are chosen such that $f_\alpha$ is a $\le_{\frak
K}$-embedding of $M_\alpha$ into $N_\alpha$, both are $\le_{{\frak
K}_\lambda}$-increasing continuous, $f_\alpha$ is
$\subseteq$-increasing continuous, $M_0 = M$ and in the $\alpha$-th
move, $M_\alpha$ is chosen by INC, and the pair is chosen by the
player INC if $\alpha$ is even and by the player COM if $\alpha$ is odd.
 The player COM wins if INC has always a legal move (the player COM 
always has: he can choose $N_\alpha =
M_\alpha$)
\bigskip

\definition{\stag{734-11.17} Definition}  1) Let $<^*_\lambda =
<^*_{{\frak K}^*_\lambda}$ be the following
two-place relation on $K^*_\lambda$ (so $M \le^*_{{\frak K}^*_\lambda}
N$ mean $M = N \in {\frak K}^*_\lambda$ or $M <^*_{{\frak K}^*_\lambda} N$):
\block  $M_1 <^*_\lambda M_2$ iff $M_1 \le_{{\frak K}_\lambda} M_2$
are from $K^*_\lambda$ and $M_2$ is 
$\le_{{\frak K}_\lambda}$-universal over $M_1$.
\endblock
2) For $\alpha < \lambda,\kappa = \beth_{1,1}(|\alpha| + 
\text{LS}({\frak K}))$ and 
$M \in K^*_\lambda$ let $\Sav^{\text{bs},\alpha}(M)$ be the set
of $\{\text{Av}_\kappa(\bold I,M):\bold I$ is a
$((2^\kappa)^+,\kappa)$-convergent subset of ${}^\alpha M\}$.  We 
define \sftp$_*(\bar a,M,N)$ when $M \le_{\frak K} N$ are 
from $K^*_\lambda$ and $\bar a \in {}^\alpha N$, as
\sftp$_{{\Bbb L}_{\infty,\kappa}[{\frak K}]}(\bar a,M,N) 
\in \Sav^{\text{bs},\alpha}(M)$ naturally.
\nl
3) Let ${\frak K}^*_\lambda = (K^*_\lambda,\le_{\frak K} \restriction
{\frak K}^*_\lambda,\le^*_{{\frak K}^*_\lambda})$, see \scite{734-11.17.3}
below but if $(K^*_\lambda,\le_{\frak K} \restriction K^*_\lambda)$ is a
$\lambda$-a.e.c. then we omit $\le^*_{{\frak K}^*_\lambda}$.
\enddefinition
\bigskip

\remark{\stag{734-11.17.3} Remark}  1) Note that the relation $<^*_\lambda =
<^*_{{\frak K}_\lambda}$ seemingly depends on the choice of $\Phi$.
However, assuming $\mu$-solvability,
by \scite{734-11.19A}(2) below it does not depend.
\nl
2) The proof of \scite{734-11.19A} is like \marginbf{!!}{\cprefix{600}.\scite{600-0.22}}(3).
\nl
3) So ${\frak K}^*_\lambda$ is a semi-$\lambda$-a.e.c. (see
\chaptercite{E53}) but we do not use this notion here.
\endremark
\bigskip

\proclaim{\stag{734-11.18} Claim}  Assume {\rm cf}$(\lambda) = \aleph_0$.
\nl
0) If $M \in K^*_\lambda$ \ub{then} for some $N, 
M <^*_{{\frak K}^*_\lambda} N (\in K^*_\lambda)$.
\nl
1) If 
$M \le_{\frak K} N$ are from $K^*_\lambda,\alpha < \lambda$ and $\bar a \in
{}^\alpha N \backslash {}^\alpha M$ \ub{then} $\bar a$ realizes some 
$p \in \,\text{\rm Sav}^{\text{bs},\alpha}(M)$. 
\nl 
2) If $M_0 \le_{\frak K} M_1 <^*_{{\frak K}^*_\lambda} 
M_2 \le_{\frak K} M_3$ and 
$M_\ell \in K^*_\lambda$ for $\ell < 4$, \ub{then} $M_0 
<^*_{{\frak K}^*_\lambda} M_3$. 
\endproclaim
\bigskip

\demo{Proof}  0) As $K^*_\lambda$ is categorical (by \scite{734-11.4A}(1))
this follows by \scite{734-11.16A}(2).
\nl
1) A proof of this is included in the proof of \scite{734-11.14}(2),
i.e. by \scite{734-11.10}(1).
\nl
2) Easy recalling amalgamation.  \hfill$\square_{\scite{734-11.18}}$
\enddemo
\bigskip

\proclaim{\stag{734-11.19A} Claim}  Assume {\rm cf}$(\lambda) = \aleph_0$.
\nl
1) Assume $\langle M_i:i \le \delta
\rangle$ is $\le_{{\frak K}_\lambda}$-increasing continuous, $M_{2i+1}
<^*_{{\frak K}^*_\lambda} M_{2i+2}$ for $i < \delta$ \ub{then} $M_\delta \in
K^*_\lambda$.
\nl
2) Assume that $\langle M^\ell_i:i \le \delta \rangle$ is an
 $\le_{{\frak K}^*_\lambda}$-increasing continuous sequence such that
$M^\ell_{2i+1} <^*_{{\frak K}^*_\lambda} M^\ell_{2i+2}$ for 
$i < \delta$ all for
 $\ell=1,2$.  Any isomorphism $f$ from $M^1_0$ onto $M^2_0$ (or just a
 $\le_{{\frak K}_\lambda}$-embedding) can be extended to an
 isomorphism from $M^1_\delta$ onto $M^2_\delta$.
\endproclaim
\bigskip

\demo{Proof}  1) We prove this by induction on $\delta$, hence \wilog
 \, $i < \delta \Rightarrow M_i \in K^*_\lambda$.

Let $M^1_\alpha = M_\alpha$ for $\alpha \le \delta$
and let $\langle I_\alpha:\alpha \le \delta\rangle$ be
$<_{K^{\text{flin}}_\lambda}$-increasing.
Let $M^2_\alpha = \text{ EM}_{\tau({\frak K})}(I_\alpha,\Phi)$.  Now there
is an isomorphism $f$ from $M^1_0$ onto $M^2_0$ as $K^*_\lambda$ is
categorical, so by part (2) there is an isomorphism $g$ from
$M^1_\alpha$ onto $M^2_\alpha$, but $M^2_\alpha \in K^*_\lambda$ so we
are done.
\nl
2) Note 
\mr
\item "{$\boxtimes_2$}"   \wilog 
{\roster
\itemitem{ $\boxdot$ }   $M^2_i <^*_\lambda M^2_{i+1}$.
\endroster}
\ermn
[Why?  We can find $\langle M^3_i:i \le \delta \rangle$ which is
$\le^*_{{\frak K}^*_\lambda}$-increasing continuous and $M^3_0 = M^2_0$
and $M^3_i <^*_\lambda M^3_{i+1}$.  Now apply the restricted version
(i.e., with the assumption $\boxdot$) twice.]

By induction on $i \le \delta$ we choose $(f_i,N^1_i,N^2_i)$ such that
\mr
\item "{$\circledast$}"  $(a) \quad N^1_i,N^2_i$ belongs to
$K^*_\lambda$
\sn
\item "{${{}}$}"  $(b) \quad f_i$ is an isomorphism from $N^1_i$ onto
$N^2_i$
\sn
\item "{${{}}$}"  $(c) \quad N^1_i,N^2_i,f_i$ are increasing
continuous with $i$
\sn
\item "{${{}}$}"  $(d) \quad$ for $i=0,N^1_i = M^1_i,f_i=f$ and
$N^2_i$ is $f(M^1_i) = M^2_i$
\sn
\item "{${{}}$}"  $(e) \quad$ if $i>0$ is a limit ordinal then $N^1_i
= M^1_i$ and $N^2_i = M^2_i$
\sn
\item "{${{}}$}"  $(f) \quad$ when $i = \omega \alpha + 2n < \delta$ we have
{\roster
\itemitem{ ${{}}$ }  $(\alpha) \quad N^1_{\omega \alpha + 2n+1} = M^1_{\omega
\alpha + 2n+1}$
\sn
\itemitem{ ${{}}$ }  $(\beta) \quad N^2_{\omega \alpha + 2n+1} \le_{\frak K} 
M^2_{\omega \alpha + 2n+1}$
\sn
\itemitem{ ${{}}$ }  $(\gamma) \quad N^1_{\omega \alpha + 2n+2} \le_{\frak K} 
M^1_{\omega \alpha + 2n+2}$
\sn
\itemitem{ ${{}}$ }  $(\delta) \quad N^2_{\omega \alpha + 2n+2} = 
M^2_{\omega \alpha + 2n+2}$.
\endroster}
\ermn
\ub{Case 1}:  For $i=0$ this is trivial by clause (d) and the
assumption of the claim on $f$.
\mn
\ub{Case 2}:  $i = \omega \alpha + 2n +1$.

Note that $N^2_{\omega \alpha + 2n} = M^2_{\omega \alpha +2n}$.  (Why?
If $i = 0$ (i.e. $\alpha = 0 = n$) by $\circledast(d)$ and if $i$ is a
limit ordinal (i.e. $\alpha > 0 \wedge n=0$) by clause (e) of
$\circledast$ and if $n>0$ by clause ($(f)(\delta)$ of $\circledast$).

Now we let $N^1_i = N^1_{\omega \alpha + 2n+1} := M^1_{\omega \alpha +
2n+1}$ and hence satisfying clause $(f)(\alpha)$ of $\circledast$.
So $N^1_{i-1} = N^1_{\omega \alpha + 2n} \le_{\frak K} 
M^1_{\omega \alpha + 2n} \le_{\frak K} 
M^1_{\omega \alpha + 2n+1} = N^1_{\omega \alpha + 2n+1} = N^1_i$; 
and note that 
$N^2_{i-1} = N^2_{\omega \alpha + 2n} <^*_\lambda M^2_{\omega
\alpha+2n}$ by $\boxdot$ above hence we can apply Definition
\scite{734-11.17}(1) and find an
extension $f_i$ of $f_{i-1}$ to $\le_{\frak K}$-embedding of $N^1_i =
M^1_{\omega \alpha +2n+1}$ into $M^2_{\omega \alpha + 2n+1}$ and let
$N^2_i := f_i(N^1_i)$.
\mn
\ub{Case 3}:  $i = \omega \alpha + 2n+2$.

Note that $N^1_{\omega \alpha + 2n+1} = M^1_{\omega \alpha + 2n+1}$ by
clause $(f)(\alpha)$ of $\circledast$ hence by the assumption of the claim
$N^1_{\omega \alpha + 2n+1} <^*_{{\frak K}^*_\lambda} 
M^1_{\omega \alpha + 2n+2}$.  We choose 
$N^2_{\omega \alpha + 2n+2} := M^2_{\omega \alpha + 2n+2}$ hence
$N^2_{i-1} = N^2_{\omega \alpha + 2n+1} 
\le_{\frak K} M^2_{\omega \alpha + 2n+1}
\le_{\frak K} M^2_{\omega \alpha + 2n+2} = N^2_{\omega \alpha + 2n+2}
= N^2_i$.  Now we apply Definition \scite{734-11.17}(1) to find a $\le_{\frak
K}$-embedding $g_i$ of $N^2_{\omega \alpha + 2n+2}$ into 
$M^1_{\omega \alpha + 2n+2}$ extending $f^{-1}_{i-1}$.  

Lastly, let 
$f_i = g^{-1}_i$ and $N^1_i = M^1_i \restriction \text{ Dom}(f_i)$.
So we can carry the induction hence prove the claim.
\hfill$\square_{\scite{734-11.19A}}$ 
\enddemo
\bn
Note that now we use more than in Hypothesis \scite{734-11.2A}.
\proclaim{\stag{734-11.20} Claim}  Assume
\mr
\item "{$\boxtimes$}"   $(a) \quad \langle \lambda_n:n < \omega \rangle$ is 
increasing, $\lambda = \lambda_\omega = \dsize \sum_{n < \omega} \lambda_n$ 
satisfying
\nl

\hskip25pt $\lambda_n = \beth_{\lambda_n} > { \text{\rm LS\/}}({\frak K})$ 
and {\rm cf}$(\lambda_n) = \aleph_0$ for $n < \omega$
\sn
\item "{${{}}$}"  $(b) \quad \Phi \in \Upsilon^{\text{or}}_{\frak K}$
and for it each $\lambda_n$ and $\lambda = \lambda_\omega$ is 
as in Hypothesis \scite{734-11.2A} 
\nl

\hskip25pt or just satisfies all its conclusions so far.
\ermn
1) $K^*_\lambda$ is closed under unions $\le_{\frak K}$-increasing chains (of
length $< \lambda^+$). \nl
2) If $M_n \in K^*_{\lambda_n},M_n \le_{\frak K} M_{n+1}$ and $M 
= \dbcu_{n < \omega} M_n$ \ub{then} $M \in K^*_\lambda$. \nl
3) If $M \in K_\lambda$ and $\theta < \lambda \Rightarrow 
M \equiv_{{\Bbb L}_{\infty,\theta}[{\frak K}]} 
{\text{\rm EM\/}}_{\tau({\frak K})}
(\lambda,\Phi)$ \ub{then} $M \in K^*_\lambda$.
\nl
4) $K^*_\lambda$ is categorical.
\endproclaim
\bigskip

\demo{Proof of \scite{734-11.20}}  1) We rely on part (2) which is
proven below.

So let $\langle M_i:i < \delta \rangle$ be
$\le_{\frak K}$-increasing in $K^*_\lambda$ with 
$\delta < \lambda^+$.  Without loss
of generality $\delta = \text{ cf}(\delta)$ hence $\delta < \lambda$
so call it $\theta$ and we prove
this by induction on $\theta$, so without loss of generality 
$\langle M_i:i < \theta \rangle$ is $\le_{\frak K}$-increasing 
continuous such that $M_i \in K^*_\lambda$ for $i < \theta$, 
and let $M_\theta = \dbcu_{i < \theta} M_i$.  By renaming \wilog \,
$\theta < \lambda_0$.

Let $I_n,I'_n$ be such that:
\mr
\item "{$\odot_1$}"  $(a) \quad I_n$ is a linear order 
of cardinality $\lambda_n$ from $K^{\text{flin}}$
\sn
\item "{${{}}$}"  $(b) \quad I'_n$ is a linear order of cardinality
$2^{\lambda_n}$ from $K^{\text{flin}}$ 
\sn
\item "{${{}}$}"  $(c) \quad I'_n$ is $\lambda^+_n$-saturated (which means that
its cofinality is $> \lambda_n$, the
\nl

\hskip25pt  cofinality of its inverse is $>
\lambda_n$ and if $I'_n \models ``s_{\alpha_1} < s_{\beta_1} < t_{\beta_2} <
t_{\alpha_2}"$
\nl

\hskip25pt  where $\alpha_1 < \beta_1 < \gamma_1,\alpha_1 < \beta_2 <
\gamma_2$ and $|\gamma_1| + |\gamma_2| < \lambda^+_n$
then for 
\nl

\hskip25pt some $r$ we have $I'_n \models ``s_{\alpha_1} < r < t_{\alpha_2}"$ 
for $\alpha_1 < \gamma_1,\alpha_2 < \gamma_2$)
\sn
\item "{${{}}$}"  $(d) \quad I_n <_{K^{\text{flin}}} 
I'_n <_{K^{\text{flin}}} I_{n+1}$ for $n < \omega$.
\ermn
Let $I = \cup\{I_n:n < \omega\}$, so $I$ is a universal member of
$K^{\text{lin}}_\lambda$.  Let $M^* = \text{ EM}_{\tau({\frak K})}
(I,\Phi)$, so for every $i < \theta$ there is an isomorphism $f_i$ from
$M^*$ onto $M_i$, exists as $K^*_\lambda$ is categorical by
\scite{734-11.5}(4) as cf$(\lambda) = \aleph_0$.

Now
\mr
\item "{$\odot_2$}"  $(a) \quad$ every interval of $I$ is universal in
$K^{\text{lin}}_\lambda$
\sn
\item "{${{}}$}"  $(b) \quad$ if $n < \omega,J \subseteq I,\chi =
|J| < \lambda$ and ${\Cal E}_{J,I} = \{(t_1,t_2):t_1,t_2 \in I
\backslash J$ and
\nl

\hskip25pt  $s \in J \Rightarrow s <_I t_1 \equiv s <_J
t_2\}$ \ub{then} for at most $\chi$ elements of $t$ of $J \backslash
I$
\nl

\hskip25pt  the set $t/{\Cal E}_{J,I}$ is a singleton.
\ermn
[Why?  Clause (a) is obvious.  For clause (b) assume $\langle
t_\alpha:\alpha < \chi^+\rangle$ are pairwise distinct members of $J
\backslash I$ such that $t_\alpha/{\Cal E}_{J,I}$ is a singleton for
each $\alpha < \chi^+$.  Without
loss of generality for some $k< \omega$ we have $\alpha < \chi^+
\Rightarrow t_\alpha \in I_k$ hence $\chi \le \lambda_k$.
For each $\alpha < \chi^+$ we can choose $s_\alpha \in I'_k$ such that
$s_\alpha <_{I'_k} t_\alpha$ and $(s_\alpha,t_\alpha)_{I'_k} \cap J =
\emptyset$.  Clearly $\alpha < \beta < \chi^+ \Rightarrow (t_\alpha
<_I s_\beta \vee t_\beta <_I s_\alpha)$ hence $\langle
(s_\alpha,t_\alpha)_I:\alpha < \chi^+\rangle$ are pairwise disjoint
intervals of $I$, so for every $\alpha < \chi^+$ large enough,
$(s_\alpha,t_\alpha)_I \cap J = \emptyset$, but then
$(s_\alpha,t_\alpha)_I \subseteq t_\alpha/{\Cal E}_{J,I}$, contradiction.]
 
Now by induction on $n < \omega$ and
for each $n$ by induction on $\varepsilon \le \theta$ and for each 
$n < \omega$ and $\varepsilon \le \theta$ for $i \le \theta$, we choose
$J_{n,\varepsilon,i} \in K^{\text{flin}}_{\lambda_n}$ such that:
\mr
\item "{$\odot_3$}"  $(a) \quad J_{n,\varepsilon,i} \subseteq I$
\sn
\item "{${{}}$}"  $(b) \quad J_{n,\varepsilon,i}$ has cardinality $\lambda_n$
\sn
\item "{${{}}$}"  $(c) \quad I_n <_{K^{\text{flin}}} J_{n,0,i}$
\sn
\item "{${{}}$}"  $(d) \quad$ if $\zeta < \varepsilon \le \theta$ 
and $i \le \theta$  then $J_{n,\zeta,i} \subseteq
J_{n,\varepsilon,i}$, moreover if for
\nl

\hskip25pt some $\xi,\zeta = 2 \xi +1$ and $\varepsilon = 2 \xi+2$ then
there is a
\nl

\hskip25pt $<_{K^{\text{flin}}_{\lambda_n}}$-increasing continuous sequence of
length $\omega$ with first
\nl

\hskip25pt  member $J_{n,\zeta,i}$ and union $J_{n,\varepsilon,i}$ 
\sn
\item "{${{}}$}"  $(e) \quad$ for $\varepsilon$ limit, $J_{n,\varepsilon,i} =
\dbcu_{\zeta < \varepsilon} J_{n,\zeta,i}$
\sn
\item "{${{}}$}"  $(f) \quad$ if $\varepsilon$ is odd and 
$i < j < \theta$ then 
\nl

\hskip25pt $f_i(\text{EM}_{\tau({\frak K})}
(J_{n,\varepsilon,i},\Phi))
= M_i \cap f_j(\text{EM}_{\tau({\frak K})}(J_{n,\varepsilon,j},\Phi))$
\sn
\item "{${{}}$}"  $(g) \quad J_{n,\theta,i} \subseteq J_{n+1,0,i}$
\sn
\item "{${{}}$}"  $(h) \quad$ for every $k < \omega$ and $s <_I t$ from 
$J_{n,\varepsilon,i}$ if $[s,t]_I \cap I'_k \ne \emptyset$ then
\nl

\hskip25pt $[s,t]_I \cap I'_k \cap J_{n,\varepsilon,i} \ne \emptyset$
\sn
\item "{${{}}$}"  $(i) \quad$ if $\zeta$ is odd and
$\varepsilon = \zeta +1$ then EM$_{\tau({\frak K})}
(J_{n,\zeta,i},\Phi) <^*_{{\frak K}^*_{\lambda_n}} 
\text{ EM}_{\tau({\frak K})}(J_{n,\varepsilon,i},\Phi)$. 
\ermn
There is no problem to carry the definition, for $\varepsilon = 2 \xi
+2$ recalling $\odot_2$ above; 
the only non-trivial point is clause (i), which follows by
\scite{734-11.16A}(4) and clause (d) of $\odot_3$.  Clearly $\langle
J_{n,\varepsilon,i}:\varepsilon \le \theta\rangle$ is
$\subseteq$-increasing continuous by $\odot_3(d) + (e)$.
\nl
Let $M^*_{n,\varepsilon,i} =
f_i(\text{EM}_{\tau({\frak K})}(J_{n,\varepsilon,i},\Phi))$ and
$M^*_{n,\varepsilon} = M^*_{n,2\varepsilon,\varepsilon}$.  So clearly
$M^*_{n,\varepsilon,i} \in K^*_{\lambda_n}$ by $\odot_3(b)$ and the
choice of $M^*_{n,\varepsilon,i}$ the sequence
$\langle M^*_{n,\varepsilon}:\varepsilon
< \theta \rangle$ is $\le_{\frak K}$-increasing 
continuous, all members in $K^*_{\lambda_n}$.

Now
\mr
\item "{$\odot_4$}"  $\langle M^*_{n,\varepsilon}:\varepsilon <
\theta\rangle$ is $<^*_{{\frak K}^*_{\lambda_n}}$-increasing.
\ermn
[Why? As $\zeta < \varepsilon < \theta \Rightarrow M^*_{n,\zeta} =
M_{n,2 \zeta,\zeta} \le_{{\frak K}^*_{\lambda_n}} M_{n,2 \zeta
+1,\zeta} \le_{{\frak K}^*_{\lambda_n}} M_{n,2 \zeta +1,\varepsilon}
<^*_{{\frak K}^*_{\lambda_n}} M_{n,2 \zeta +2,\varepsilon}
\le_{{\frak K}^*_{\lambda_n}} M_{n,2 \varepsilon,\varepsilon} =
M^*_{n,\varepsilon}$ by the choice of $M^*_{n,\zeta}$, by
$\odot_3(d)$ and Ax(V) of a.e.c., by $\odot_3(f)$ and Ax(V) of a.e.c., by
$\odot_3(i)$, by $\odot_3(d)$ + Ax(V) of a.e.c.(e),
by the choice of $M^*_{n,\varepsilon}$ respectively).  
Now by \scite{734-11.18}(2) this
argument shows that $\zeta < \varepsilon < \theta \Rightarrow
M^*_{n,\zeta} <^*_{{\frak K}^*_{\lambda_n}} M^*_{n,\varepsilon}$.]

We can conclude by using \scite{734-11.19A}(1) for ${\frak K}^*_{\lambda_n}$, that 
$M^*_n := \dbcu_{\varepsilon < \theta} M^*_{n,\varepsilon}$
belongs to $K^*_{\lambda_n}$.  Also as $M^*_{n,\varepsilon} \le_{\frak K} 
M_\varepsilon \le_{\frak K} M_\delta$ for 
$\varepsilon < \theta = \delta$ by AxIV of a.e.c. we have
$M^*_n \le_{\frak K} M_\delta$ and similarly 
$M^*_n \le_{\frak K} M^*_{n+1}$, and 
obviously for each $i < \theta$ we have $\dbcu_{n < \omega} M^*_n$
includes $\cup\{M^*_{n,\varepsilon}:n < \omega,\varepsilon < \theta\}
= \cup\{M^*_{n,2,\varepsilon,\varepsilon}:n < \omega,\varepsilon <
\theta\} = \cup\{M^*_{n,2\varepsilon,i}:n < \omega,i <
\theta,\varepsilon <\theta\} = 
\dbcu_{n < \omega} M^*_{n,0,i}$ which recalling
the choice of $M^*_{n,0,i}$ includes
$\dbcu_n f_i(\text{EM}_{\tau({\frak K})}(J_{n,0,i},\Phi)) \supseteq 
\dbcu_{n < \omega} f_i(\text{EM}_{\tau({\frak K})}(I_n,\Phi)) = 
f_i(\text{EM}_{\tau({\frak K})}(I,\Phi)) = M_i$.
As this holds for every $i < \theta$ we get $\dbcu_{n < \omega} M^*_n 
=  M_\delta$.  So by part (2) we are done.
\nl
2) We choose $I_n$ by induction on $n$ such that:
\mr 
\item "{$\odot_5$}"  $(a) \quad I_n \in K^{\text{flin}}_{\lambda_n}$
\sn
\item "{${{}}$}"  $(b) \quad I_m <_{K^{\text{flin}}} I_n$ if $n = m+1$.
\ermn
Let $N_n = \text{ EM}_{\tau({\frak K})}(I_n,\Phi)$.

We now choose $(g_n,I'_n,I''_n,M'_n,M''_n,N'_n,N''_n)$ by 
induction on $n < \omega$ such that:
\mr 
\item "{$\odot_6$}"  $(a) \quad g_n$ is an isomorphism from $N''_n$
onto $M''_n$
\sn
\item "{${{}}$}"  $(b) \quad I_n \subseteq I'_n \subseteq I''_n 
\subseteq I_{n+2}$ and
$|I'_n| = \lambda_n,|I''_n| = \lambda_{n+1}$ and $I_{n+1} \subseteq I''_n$
\sn
\item "{${{}}$}"  $(c) \quad N'_n = \text{ EM}_{\tau({\frak
K})}(I'_n,\Phi)$ and $N''_n = \text{ EM}_{\tau({\frak
K})}(I''_n,\Phi)$
\sn
\item "{${{}}$}"  $(d) \quad M_n \le_{{\frak K}^*_{\lambda_n}} M'_n
\le_{{\frak K}^*} M''_n \le_{{\frak K}^*} M_{n+2}$
and $M_{n+1} \le_{{\frak K}^*_{\lambda_{n+1}}} M''_n$
\sn
\item "{${{}}$}"  $(e) \quad g_n$ maps $N'_n = 
\text{ EM}_{\tau({\frak K})}(I'_n,\Phi)$ onto $M'_n$
\sn
\item "{${{}}$}"  $(f) \quad g_n$ extends $g_m \restriction N'_m$ if
$n = m+1$
\sn
\item "{${{}}$}"  $(g) \quad I'_n \subseteq I'_{n+1}$.
\ermn
\ub{Case 1}:  For $n=0$.

First, let $M''_n = M_1,I''_n = I_1$ so also $N''_n$ is defined.  
Second, choose $g_n$ satisfying (a) of $\odot_6$ by
\scite{734-11.4A}(1), i.e. \scite{734-11.5}(4), categoricity in
$K^*_{\lambda_n}$.  Third, choose $I^*_n \subseteq I''_n = I_1$ of
cardinality $\lambda_n$ such that $g_n(\text{EM}_{\tau({\frak
K})}(I^*_n,\Phi))$ includes $M_0$.  Fourth, let $I'_n = I^*_n \cup
I_n$ and $N'_n = \text{ EM}_{\tau({\frak K})}(I'_n,\Phi)$ and let
$M'_n = g_n(N'_n)$.  
\mn
\ub{Case 2}:  For $n=m+1$.

Let $k=n+2$, let $\bar a \in {}^{\lambda_m}(M'_m)$ list $M'_m$ (with
no repetitions).

Now
\mr
\item "{$(*)_1$}"  If $\theta < \lambda_n$ then $\sftp_{\Bbb
L_{\infty,\theta}[{\frak K}]}(\bar a,\emptyset,N_k) = \sftp_{\Bbb
L_{\infty,\theta}[{\frak K}]}(\bar a,\emptyset,N''_m)$.
\ermn
[Why?  As EM$_{\tau({\frak K})}(I''_m,\Phi) \prec_{\Bbb
L_{\infty,\theta}[{\frak K}]} \text{ EM}_{\tau({\frak K})}(I_k,\Phi)$
by \scite{734-11.3}(a) as $I''_m \subseteq I_k$.]
\mr
\item "{$(*)_2$}"  if $\theta < \lambda_n = \lambda_{m+1}$ \ub{then}
\nl
$\sftp_{\Bbb L_{\infty,\theta}}(\bar a,\emptyset,N''_m) = 
\sftp_{\Bbb L_{\infty,\theta}}(g_m(\bar a),\emptyset,M''_m)$.
\ermn
[Why?  As $g_m$ is an isomorphism from $N''_m$ onto $M''_m$ by
$\odot_6(a)$, i.e. the induction hypothesis.]
\mr
\item "{$(*)_3$}"   if $\theta < \lambda_n$ then
$\sftp_{\Bbb L_{\infty,\theta}[{\frak K}]}(g_m(\bar a),\emptyset,M''_m) 
= \sftp_{\Bbb L_{\infty,\theta}[{\frak K}]} (g_m(\bar a),\emptyset,M_k)$.
\ermn
[Why?  This follows from $M'_m \prec_{\Bbb L_{\infty,\theta}[{\frak
K}]} M_k$ which we can deduce by \scite{734-11.5}(1) as $M''_m \in
K^*_{\lambda_{m+1}} = K^*_{\lambda_n}$ by clause (d) of $\odot_6$,
$M_k \in K^*_k$ by an assumption of the claim, $M''_m \le_{{\frak
K}_\lambda} M_k$ by clause (d) of $\odot_6$.]
\mr
\item "{$(*)_4$}"   if $\theta < \lambda_n$ then
$\sftp_{\Bbb L_{\infty,\theta}[{\frak K}]}(\bar a,\emptyset,N_k) 
= \sftp_{\Bbb L_{\infty,\theta}[{\frak K}]}(g_m(\bar a),\emptyset,M_k)$.
\ermn
[Why?  By $(*)_1 + (*)_2 + (*)_3$.]
\mr
\item "{$(*)_5$}"   $\sftp_{\Bbb L_{\infty,\lambda^+_{n+1}}[{\frak K}]}
(\bar a,\emptyset,N_k) = 
\sftp_{\Bbb L_{\infty,\lambda^+_{n+1}}[{\frak K}]}(g_m(\bar a),\emptyset,M_k)$.
\ermn
[Why?  Clearly $N_k,M_k \in K^*_{\lambda_k}$ hence by \scite{734-11.5}(4)
there is an isomorphism $f_n$ from $N_k$ onto $M_k$, so obviously
$\sftp_{\Bbb L_{\infty,\theta}[{\frak K}]}(\bar a,\emptyset,N_k) =
\sftp_{\Bbb L_{\infty,\theta}[{\frak K}]}(f_n(\bar a),\emptyset,N_k)$ so by
$(*)_4$ we have $\sftp_{\Bbb L_{\infty,\theta}[{\frak K}]}(g_m(\bar
a),\emptyset,M_k) = \sftp_{\Bbb L_{\infty,\theta}[{\frak K}]}(\bar
a,\emptyset,N_k) = \sftp_{\Bbb L_{\infty,\theta}[{\frak K}]}(f_n(\bar
a),\emptyset,M_k)$ so by \scite{734-11.5}(3) we have $\sftp_{\Bbb
L_{\infty,\lambda^+_{n+1}}[{\frak K}]}(g_n(\bar a),\emptyset,M_k) =
\sftp_{\Bbb L_{\infty,\lambda^+_{n+1}}[{\frak K}]}(f_n(\bar a),\emptyset,M_k)$.
But as $f_n$ is an isomorphism from $N_k$ onto $M_k$ and the previous
sentence we get $\sftp_{\Bbb L_{\infty,\lambda_{n+1}}[{\frak K}]}(\bar
a,\emptyset,N_k) = \sftp_{\Bbb L_{\infty,\lambda^+_{n+1}}[{\frak
K}]}(f_n(\bar a),\emptyset,M_k) = 
\sftp_{\Bbb L_{\infty,\lambda}}(g_n(\bar a),\emptyset,M_k)$ as required.] 
\mr
\item "{$(*)_6$}"  there are $g_n,I''_n,N''_n,M''_n$ as required in
the relevant parts of $\odot_6$ (ignoring $I'_n,N'_n,M'_n$),
i.e. clauses (a),(f) and the relevant parts of (b),(c),(d):
{\roster
\itemitem{ $(b)'$ }  $I_n \subseteq I''_n \subseteq I_{n+2} = I_k$ and
$|I''_n| = \lambda_{n+1}$ and $I_{n+1} \subseteq I'_n$
\sn
\itemitem{ $(c)'$ }  $N''_n = \text{ EM}_{\tau({\frak
K})}(I''_n,\Phi)$
\sn
\itemitem{ $(d)'$ }  $M_n \le_{{\frak K}^*} M''_n \le_{{\frak K}^*}
M_{n+2}$ and $M_{n+1} \le_{{\frak K}^*_{\lambda_{n+2}}} M''_n$.
\endroster}
\ermn
[Why?  By the hence and forth argument, but let us elaborate.

First, let $\bar a'$ be a sequence of length $\lambda_{n+1}$ listing 
(without repetitions) the
set of elements of $M_{n+1}$ and \wilog \, $g(\bar a) \triangleleft
\bar a'$.  Note that Rang$(g_m) \subseteq M_{m+2} = M_{n+1}$.

Second, let $g'$ be a function from Rang$(\bar a')$ into $N_k$ extending
$(g_m \restriction N'_m)^{-1} = (g_m \restriction \text{ Rang}(\bar
a))^{-1}$ such that $\sftp_{\Bbb L_{\infty,\lambda^+_{n+1}}[{\frak K}]}
(g'(\bar a'),\emptyset,N_k) = \sftp_{\Bbb
L_{\infty,\lambda^+_{n+1}}[{\frak K}]}(\bar a',\emptyset,M_k)$, exists
by $(*)_5$.  Let
$I''_n \subseteq I_k$ of cardinality $\lambda_{n+1}$ be
such that Rang$(g') \subseteq \text{ EM}(I''_n,\Phi)$ and $I_{n+1}
\subseteq I''_n$.  Let $\bar a''$ list the elements of
EM$_{\tau({\frak K})}(I''_n,\Phi) \subseteq N_k$ and without
loss of generality 
$g'(\bar a') \triangleleft \bar a''$
and let $g_n$ be a function from EM$_{\tau({\frak K})}(I''_n,\Phi)$ to
$M_k$ extending
$(g')^{-1}$ such that $\sftp_{\Bbb L_{\infty,\lambda^+_{n+1}}[{\frak
K}]} (\bar a'',\emptyset,N_k) = \sftp_{\Bbb
L_{\infty,\lambda^+_{n+1}}[{\frak K}]}(g_n(\bar a''),\emptyset,M_k)$.

Lastly, let $N''_n = \text{ EM}_{\tau({\frak K})}(I''_n,\Phi)$ and
$M''_n = g_n(N'_n)$ so we are done.]
\mr
\item "{$(*)_7$}"  there are $I'_n,N'_n,M'_n$ as required.
\ermn
[Why?  By the LS argument we can choose $I'_n$ and define $N'_n,M'_n$
accordingly.] 
 
So we can carry the induction.  Now $N'_n \le_{\frak K} N'_{n+1}$ 
(by clauses (g),(c) of $\odot_6$) and $g_n \restriction
N'_n \subseteq g_{n+1} \restriction N'_{n+1}$ 
(by clause (f) + the previous statement).  Hence
$g = \cup\{g_n \restriction N'_n:n < \omega\}$ is an 
isomorphism from $\cup\{N'_n:n < \omega\}$ onto $\cup\{M'_n:
n < \omega\}$.  But $N = \cup\{N_n:n < \omega\} 
\subseteq \cup\{N'_n:n < \omega\} \subseteq \text{ Dom}(g)
\subseteq N$ and $M = \cup\{M_n:n < \omega\} \subseteq \cup\{M'_n:n <
\omega\} \subseteq \text{ Rang}(g) \subseteq M$.  Together $g$ is an
isomorphism from $N$ onto $M$ but obviously $N \in K^*_\lambda$ hence
$M \in K^*_\lambda$ is as required.
\nl
3),4)  Should be clear and depends just on \scite{734-11.5}(4).  
\hfill$\square_{\scite{734-11.20}}$
\enddemo
\bigskip

\demo{\stag{734-11.21} Conclusion}   Let $\lambda$ be as in $\boxtimes$ of
\scite{734-11.20}. \nl
1)  ${\frak K}^*_\lambda$ is a 
$\lambda$-a.e.c. (with $\le_{\frak K} \restriction K^*_\lambda$) and
it has amalgamation and is categorical. \nl
2) ${\frak K}^\oplus_{\ge \lambda}$ is an a.e.c., 
LS$({\frak K}^\oplus_{\ge \lambda}) = \lambda$ and $({\frak
K}^*_\lambda)^{\text{up}} = K^\oplus_{\ge \lambda}$ and
$({\frak K}^\oplus_{\ge \lambda})_\lambda = 
{\frak K}^*_\lambda$, see Definition below.
\enddemo
\bigskip

\definition{\stag{734-11.22} Definition}  Let 
${\frak K}^\oplus_{\ge \lambda} =
{\frak K} \restriction K^\oplus_{\ge \lambda}$ where
$K^\oplus_{\ge \lambda} = \{M \in K_\lambda:M 
\equiv_{{\Bbb L}_{\infty,\lambda}
[{\frak K}]} \text{ EM}_{\tau({\frak K})}(\lambda,\Phi)\}$.
\enddefinition

\demo{Proof}  1) It was clear defining $(K^*_\lambda,\le_{\frak K}
\restriction K^*_\lambda)$ that it is of the right form and ``$M \in
K^*_\lambda$", ``$M \le_{{\frak K}^*_\lambda} N$" are preserved by 
isomorphisms.  Obviously ``$\le_{\frak K} \restriction
K^*_\lambda$ is a partial order", so AxI, AxII hold and obviously
 AxV holds (see \marginbf{!!}{\cprefix{600}.\scite{600-0.2}}).  The missing point was AxIII,
about $\le_{\frak K}$-increasing union and it holds by
\scite{734-11.20}(1).  Then AxIV becomes easy by the definition of
$\le_{{\frak K}^*_\lambda} = \le_{\frak K} \restriction K^*_\lambda$
and lastly the amalgamation holds by \scite{734-11.15}.
\nl
2) By \sectioncite[\S1]{600} we can ``lift ${\frak K}^*_\lambda$ up", the
result is ${\frak K}^\oplus_{\ge \lambda}$ (see
\marginbf{!!}{\cprefix{600}.\scite{600-0.31}},\marginbf{!!}{\cprefix{600}.\scite{600-0.32}}).  \hfill$\square_{\scite{734-11.21}}$
\enddemo
\bn
Let us formulate a major conclusion in ways less buried
inside our notation.
\demo{\stag{734-11.28} Conclusion}  Assume $({\frak K},\Phi)$ is pseudo
solvable in $\mu$, \ub{then} $({\frak K},\Phi)$ is pseudo solvable in
$\lambda$ provided that LS$({\frak K}) < \lambda,\mu = \mu^{<
\lambda}$ (or just the hypothesis \scite{734-11.2A} holds), cf$(\lambda) =
\aleph_0$ and $\lambda$ is an accumulation point of the class of the 
fix point of the sequence of the $\beth$'s.
\enddemo
\bigskip

\demo{Proof}  By \scite{734-11.21}(1).  \hfill$\square_{\scite{734-11.28}}$
\enddemo
\bigskip

\remark{Remark}  About [weak] solvability, see \cite{Sh:F782}.
\endremark
\goodbreak

\head {\S2 Trying to Eliminate $\mu = \mu^{< \lambda}$} \endhead  \resetall \sectno=2
 \spuriousreset
\bigskip

There was one point in \S1 where we use $\mu = \mu^\lambda$ (i.e.
in \scite{734-11.2}, more accurately in justifying hypothesis \scite{734-11.2A}(1)).
In this section we try to eliminate it.  So we try to prove $M_1
\le_{{\frak K}_\mu} M_2 \Rightarrow M_1 
\prec_{\Bbb L_{\infty,\theta}[{\frak K}]} M_2$ for $\theta < \lambda$, 
hence we fix ${\frak K},\mu,\theta$.  We succeed to do it with ``few
exceptions".  
\bigskip

\demo{\stag{734-eq.0} Hypothesis}  (We shall mention $(b)_\mu$ or $(b)^-_\mu,
(c),(d)$ when used! but not clause (a))
\mr
\widestnumber\item{$(b)_\mu$}
\item "{$(a)$}"  ${\frak K}$ is an a.e.c. and $\Phi \in
\Upsilon^{\text{or}}_{\frak K}$
\sn
\item "{$(b)_\mu$}"  ${\frak K}$ categorical in $\mu$ and $\Phi \in
\Upsilon^{\text{or}}_{\frak K}$, or at least
\sn
\item "{$(b)^-_\mu$}"   ${\frak K}$ is pseudo $\mu$-solvable as witnessed by
$\Phi \in \Upsilon^{\text{or}}_{\frak K}$, see Definition \scite{734-11.1} 
in particular EM$_{\tau({\frak K})}(I,\mu)$ is pseudo superlimit for $I \in
K^{\text{lin}}_\lambda$, 
\sn
\item "{$(c)$}"  $\mu \ge \beth_{1,1}(\text{LS}({\frak K}))$
\sn
\item "{$(d)$}"   $\mu > \text{ LS}({\frak K})$.
\endroster
\enddemo
\bn
\margintag{734-eq.0.8}\ub{\stag{734-eq.0.8} Convention}:  $K^*_\lambda = K^*_{\Phi,\lambda}$,
etc., see Definition \scite{734-11.4}.
\bigskip

\definition{\stag{734-eq.1} Definition}  Assume
\mr
\item "{$\boxdot$}"   $\mu \ge \chi \ge \theta > \text{ LS}({\frak K})$ 
\ermn
1) We let $K^1_{\mu,\chi} = 
\{(M,N):N \le_{\frak K} M,N \in K_\chi,M \in K_\mu$ and $\mu = \chi
\Rightarrow M=N\}$ and let $\le_{\frak K} = \le_{{\frak K},\mu,\chi}$
be the following partial order on $K_{\mu,\chi},(M_0,N_0)
\le_{\frak K} (M_1,N_1)$ iff $M_0 \le_{\frak K} M_1,N_0 \le_{\frak K} N_1$
(formally we should have written $\le_{{\frak K},\mu,\chi}$).
Note that each pair $(M,N) \in K_{\mu,\chi}$ determine $\mu,\chi$.
So if $\chi = \mu,K_{\mu,\chi}$ is essentially ${\frak K}_\mu$.  Let
$K^1_\mu = K_\mu$ and let $\cup\{(M_i,N_i):i < \delta\} = (\cup\{M_i:i
< \delta\}, \cup\{N_i:i < \delta\})$ for any 
$\le_{\frak K}$-increasing sequence $\langle (M_i,N_i):i < \delta\rangle$. 
\nl
1A) Let $K_{\mu,\chi} = K^2_{\mu,\chi} = \{(M,N) \in K^1_{\mu,\chi}:M
\in K^*_\mu\}$ and $K^2_\mu = K^*_\mu$ but we use them 
only when $\Phi$ witnesses ${\frak K}$ is pseudo
$\mu$-solvable, i.e. $(b)^-_\mu$ from Hypothesis \scite{734-eq.0} holds.
\nl 
2) For $k \in \{1,2\}$ a formula 
$\varphi(\bar x) \in \Bbb L_{\infty,\theta}[{\frak K}]$ (so
$\ell g(\bar x) < \theta$), cardinal $\kappa \ge \theta$ the main case
being $\kappa = \mu$; we may omit $k$ if $k=2$, and 
$M \in K^k_\kappa,\bar a \in {}^{\ell g(\bar x)}M$ \ub{we define} 
when $M \Vdash_k \varphi[\bar a]$ by
induction on the depth of $\varphi(\bar x) \in \Bbb
L_{\infty,\theta}[{\frak K}]$, so the least obvious case is:
\mr
\item  "{$(*)$}"  $M \Vdash_k (\exists \bar y)\psi(\bar y,\bar a)$
\ub{when} for
every $M_1 \in K^k_\kappa$ such that $M \le_{\frak K} M_1$ there is $M_2 \in
K^k_\kappa$ satisfying $M_1 \le_{\frak K} M_2$ and $\bar b \in 
{}^{\ell g(\bar y)} M_2$ such that $M_2 \Vdash_k \psi[\bar b,\bar a]$. 
\ermn
Of course
\mr
\item "{$(\alpha)$}"  for $\varphi$ atomic, $M \Vdash_k \varphi[\bar a]$
\ub{iff} $M \models \varphi[\bar a]$
\sn
\item "{$(\beta)$}"  for $\varphi(\bar x) = \dsize \bigwedge_{i < \alpha}
\varphi_i(\bar x)$ let $M \Vdash_k \varphi[\bar a]$ iff 
$M \Vdash_k \varphi_i[\bar a]$
for each $i < \alpha$
\sn
\item "{$(\gamma)$}"  $M \Vdash_k \neg \varphi[\bar a]$ iff for no $N$ do we
have $M \le_{\frak K} N \in K^k_\kappa$ and $N \Vdash_k \varphi[\bar a]$.
\ermn
3) Let $k \in \{1,2\},\Lambda \subseteq \Bbb L_{\infty,\theta}[{\frak K}]$ 
(each formula with $< \theta$ free variables, of course):
\mr
\item "{$(a)$}"  $\Lambda$ is downward closed if it is closed under
subformulas
\sn
\item "{$(b)$}"  $\Lambda$ is $(\mu,\chi)$-model$^k$ complete 
(when $\mu$ is clear from the context we may write 
$\chi$-model$^k$ complete) \ub{if} 
$|\Lambda| < \mu$, and for every $(M_0,N_0) \in 
K^k_{\mu,\chi}$ we can find $(M,N) \in K^2_{\mu,\chi}$ above $(M_0,N_0)$ which 
is $\Lambda$-generic, where:
\sn
\item "{$(c)$}"  $(M,N) \in K^k_{\mu,\chi}$ is $\Lambda$-generic$^k$ when:
\nl
if $\varphi(\bar x) \in \Lambda$ and 
$\bar a \in {}^{\ell g(\bar x)}N$ then 
\nl
$M \Vdash_k \varphi[\bar a] \Leftrightarrow N \models \varphi[\bar a]$
(yes! neither $(M,N) \Vdash_k \varphi[\bar a]$ which was not defined, nor
``$M \models \varphi[\bar a]$")
\sn
\item "{$(d)$}"  $\Lambda$ is called $(\mu,< \mu)$-model$^k$ complete 
when $|\Lambda| + \theta_\Lambda < \mu$ and for every $\chi$: 
if $|\Lambda| + \theta_\Lambda \le \chi < \mu$ 
then $\Lambda$ is $\chi$-model$^k$ complete where $\theta_\Lambda :=
\text{ min}\{\partial:\partial > \text{\rm LS}({\frak K})$ and
$\Lambda \subseteq \Bbb L_{\infty,\partial}[{\frak K}]\}$.
We say $\Lambda$ is model$^k$ complete if it is $(\mu,< \mu)$-model$^k$
complete and $\mu$ is understood from the context
\sn
\item "{$(e)$}"  above if $\Phi$ or $({\frak K},\Phi)$ is not clear
from the context we may replace $\Lambda$ by $(\Lambda,\Phi)$ or by
$(\Lambda,\Phi,{\frak K})$.
\ermn
4) For $M \in K^k_\kappa,\bar a \in {}^{\theta >}M$ and 
$\Lambda \subseteq \Bbb L_{\infty,\theta}
[{\frak K}]$ let gtp$^k_\Lambda(\bar a,\emptyset,M) = 
\{\varphi[\bar a]:M \Vdash_k \varphi[\bar a]\}$; if we write $\theta$ 
instead of $\Lambda$ we mean
$\Bbb L_{\infty,\theta}[{\frak K}]$ (note: this type is not a priori
complete) and we say that $\bar a$ materializes this type in $M$.  To
stress $\kappa$ we may write gtp$^{\kappa,k}_\Lambda(\bar
a,\emptyset,M)$ or gtp$^{\kappa,k}_\theta(\bar a,\emptyset,M)$ though
$M$ determines $\kappa$.
\nl
5)  We say $M \in K_\kappa$ is $\Lambda$-generic$^k$ \ub{when} for every
$\varphi(\bar x) \in \Lambda$ and $\bar a \in {}^{\ell g(\bar x)}M$ we have
$M \Vdash_k \varphi[\bar a] \Leftrightarrow  M \models 
\varphi[\bar a]$.  So $M \in K^k_\mu$ is $\Lambda$-generic$^k$ 
iff $(M,M) \in K^k_{\mu,\mu}$ is $\Lambda$-generic$^k$.  We say
$\Lambda$ is $\kappa$-model$^k$ complete \ub{when} every $M \in K^k_\kappa$
has a $\Lambda$-generic $\le_{\frak K}$-extension in $K^k_\kappa$
(so depend on ${\frak K}$ and if $k=2$ also on $\Phi$).
\nl
6) In all cases above, if $k=2$ we may omit it.
\enddefinition
\bigskip

\proclaim{\stag{734-eq.2} Claim}  Assume that {\rm LS}$({\frak K}) < 
\theta \le \chi < \mu$ and
$\kappa > \theta$ and $k \in \{1,2\}$ so if $k=2$ then \scite{734-eq.0}(b)$^-_\mu$
holds, see \scite{734-eq.1}(1A).
\nl
1) $(K^k_{\mu,\chi},\le_{\frak K})$ is a partial
order and chains of length $\delta < \chi^+$ of members has a
$\le_{\frak K}$-{\rm lub}, this is the union, see \scite{734-eq.1}(1).  
If {\rm EM}$_{\tau({\frak K})}(\mu,\Phi)$ is superlimit (not just pseudo
superlimit) \ub{then} $K^2_{\mu,\chi}$ is a
dense subclass of $K^1_{\mu,\chi}$ under $\le_{\frak K}$. 
\nl
2) If $M_1 \Vdash_k \varphi(\bar a)$ and 
$M_1 \le_{\frak K} M_2$ are from $K^k_\kappa$ \ub{then} 
$M_2 \Vdash_k \varphi[\bar a]$. 
\nl
3) If $(M_\ell,N_\ell) \in K^k_{\mu,\chi}$ are $\Lambda$-generic$^k$ for
$\ell=1,2$ and $(M_1,N_1)
\le_{\frak K} (M_2,N_2)$ \ub{then} $N_1 \prec_\Lambda N_2$.
\nl
4) If $M_i \in K^k_\kappa$ for $i < \delta$ is $\le_{\frak K}$-increasing, 
$\delta < \kappa^+$, {\rm cf}$(\delta) \ge \theta,\Lambda
\subseteq \Bbb L_{\infty,\theta}
[{\frak K}]$ and each $M_i$ is $\Lambda$-generic$^k$, \ub{then}
$M_\delta := \dbcu_{i < \delta} M_i$ is $\Lambda$-generic$^k$ and $i < \delta
\Rightarrow M_i \prec_\Lambda M_\delta$. 
\nl
5) If $(M_i,N_i) \in K^k_{\mu,\chi}$ for $i < \delta$ is 
$\le_{\frak K}$-increasing, $\delta < \chi^+$, 
{\rm cf}$(\delta) \ge \theta,\Lambda \subseteq 
\Bbb L_{\infty,\theta}[{\frak K}]$ and each $(M_i,N_i)$ is
$\Lambda$-generic$^k$, 
\ub{then} $(\dbcu_{i < \delta} M_i,\dbcu_{i < \delta} N_i)$
is $\Lambda$-generic$^k$ and $N_j \prec_\Lambda \dbcu_{i < \delta} N_i$
for each $j < \delta$. 
\endproclaim
\bigskip

\demo{Proof}  Should be clear; in part (1) for $k=2$ we use clause $(b)^-_\mu$
of \scite{734-eq.0}.    In part (5) note that $\cup\{M_i:I < \delta\} \in
K^*_\mu$ by Clause (b)$^-_\mu$ of \scite{734-eq.0}.
\hfill$\square_{\scite{734-eq.2}}$
\enddemo
\bn
\margintag{734-eq.3n.21}\ub{\stag{734-eq.3n.21} Exercise}:  If $(M,N)$ is $\Lambda$-generic$^k$ and
$(M,N) \le_{\frak K} (M',N) \in K^k_{\mu,\chi}$ \ub{then} $(M',N)$ is
$\Lambda$-generic$^k$.
\bigskip

\proclaim{\stag{734-eq.3} Claim}  Assume that $\mu \ge \chi \ge \theta >
\text{\rm LS}({\frak K})$ and $k \in \{1,2\}$.
\nl
1) The set of quantifier free formulas in
$\Bbb L_{\infty,\theta}[{\frak K}]$ is $(\mu,\chi)$-model$^k$ complete. 
\nl
2) If $\Lambda_\varepsilon \subseteq \Bbb
L_{\infty,\theta}(\tau_{\frak K})$ is 
downward closed, $(\mu,\chi)$-model$^k$ complete
for $\varepsilon < \varepsilon^*$, and $\Lambda := \dbcu_{\varepsilon <
\varepsilon^*} \Lambda_\varepsilon,\theta = \text{\rm cf}(\theta) \le \chi \vee
\theta < \chi,\varepsilon^* < \chi^+$ (and 
$\mu > \theta \vee \mu = \theta = { \text{\rm cf\/}}(\theta))$
\ub{then} $\Lambda$ is $(\mu,\chi)$-model$^k$ complete. 
\endproclaim
\bigskip

\demo{Proof}  1) Easy.  
\nl
2) Given $(M,N) \in K^k_{\mu,\chi}$ let $\theta_r$ be
   min$\{\partial:\partial \ge \theta$ is regular$\}$.  Clearly
$\theta_r \le \chi$ and we choose $(M_i,N_i) \in K^k_{\mu,\chi}$ 
for $i \le \varepsilon^* \times \theta_r$ such that
\mr
\item "{$\circledast$}"  $(a) \quad \langle M_i:i \le \varepsilon^*
\times \theta_r \rangle$ is $\le_{\frak K}$-increasing continuous
\sn
\item "{${{}}$}"  $(b) \quad \langle N_i:i \le \varepsilon^* \times \theta_r
\rangle$ is $\le_{\frak K}$-increasing continuous
\sn
\item "{${{}}$}"  $(c) \quad$ if $i = \varepsilon^* \times \gamma +
\varepsilon$ and $\varepsilon < \varepsilon^*$ then $(M_{i+1},N_{i+1})$ is
$\Lambda_\varepsilon$-generic$^k$
\sn
\item "{${{}}$}"  $(d) \quad (M_0,N_0) = (M,N)$.
\ermn
There is no problem to do this.

Now for each $\varepsilon < \varepsilon^*$ the sequence
$\langle(M_{\varepsilon^* \times \gamma + \varepsilon +1},
N_{\varepsilon^* \times \gamma + \varepsilon +1}):\gamma <
\theta_r\rangle$ is $\le_{{\frak K},\mu,\chi}$-increasing with union
$(M_{\varepsilon^* \times \theta_r},N_{\varepsilon^* \times
\theta_r})$, and each member of the sequence is
$\Lambda_\varepsilon$-generic$^k$ hence by \scite{734-eq.2}(5) we know
that the pair $(M_{\varepsilon^* \times \theta_r},N_{\varepsilon^*
\times \theta_r})$ is $\Lambda_\varepsilon$-generic$^k$.  As this holds
for each $\Lambda_\varepsilon$ it holds for $\Lambda$ so
$(M_{\varepsilon^* \times \theta_r},N_{\varepsilon^* \times
\theta_r})$ is as required.  \hfill$\square_{\scite{734-eq.3}}$
\enddemo
\bn
\relax From now on in this section
\demo{\stag{734-eq.3.1} Hypothesis}  We assume (a) + (b)$^-_\mu$ of
\scite{734-eq.0} and we omit $k$ using Definition \scite{734-eq.1} meaning $k=2$.
\enddemo
\bigskip

\proclaim{\stag{734-eq.3.3} Claim}  1) For $M \in K^*_\mu$ and 
{\rm LS}$({\frak K}) < \theta < \mu$ the number of complete $\Bbb
L_{\infty,\theta}[{\frak K}]$-types realized by sequences from
${}^{\theta >}M$ is $\le 2^{< \theta}$, moreover, the relation 
${\Cal E}^{< \theta}_M := 
\{(\bar a,\bar b):\bar a,\bar b \in {}^{\theta >} M$ and
some automorphism of $M$ maps $\bar a$ to $\bar b\}$ is an equivalence
relation with $\le 2^{< \theta}$ equivalence classes.
\nl
2) Hence there is a set $\Lambda_* = \Lambda^*_\theta =
 \Lambda^*_{{\frak K},\Phi,\mu,\theta} \subseteq 
\Bbb L_{\infty,\theta}[{\frak K}]$ such that:
\mr
\item "{$(a)$}"  $|\Lambda_*| \le 2^{< \theta}$ and $\Lambda_*
\subseteq \Bbb L_{(2^{< \theta})^+,\theta}[{\frak K}]$
\sn
\item "{$(b)$}"  $\Lambda_*$ is closed under sub-formulas and finitary
operations 
\sn
\item "{$(c)$}"  each $\varphi(\bar x) \in \Lambda_*$ has quantifier
depth $< \gamma^*$ for some $\gamma^* < (2^{< \theta})^+$
\sn
\item "{$(d)$}"  for $\alpha < \theta,M \in K^*_\mu$ and $\bar a
\in {}^\alpha M$, the $\Lambda_*$-type which $\bar a$ realizes in $M$
determines the $\Bbb L_{\infty,\theta}[{\frak K}]$-type which $\bar a$
realizes in $M$, moreover one formula in the type determine it
\sn
\item "{$(e)$}"  similarly for materialize in $M \in K^*_\mu$, see Definition
\scite{734-eq.1}(4)
\sn
\item "{$(f)$}"  if {\rm LS}$({\frak K}) \le \chi < \mu$ and
$(M,N) \in K_{\mu,\chi}$ is $\Lambda_*$-generic \ub{then} it is
$\Bbb L_{\infty,\theta}[{\frak K}]$-generic
\sn
\item "{$(g)$}"  if $M \in K^2_\mu$ is $\Lambda_*$-generic \ub{then} it is
 $\Bbb L_{\infty,\theta}[{\frak K}]$-generic.
\endroster
\endproclaim
\bigskip

\remark{Remark}  Part (1) can also be proved using just $(\lambda +1)
\times I_*$ with $I_*$ a $\theta$-saturated dense linear order with
neither first nor last element, but this is not clear for \scite{734-eq.4}(1). 
\endremark
\bigskip

\demo{Proof}  1) By \scite{734-am3.2.3}(1) and categoricity of $K^*_\lambda$.
\nl
2) Follows but we elaborate.

Let $\{\bar a_\alpha:\alpha < \alpha^* \le 2^{< \theta}\}$ be a set of
representatives of the ${\Cal E}^{< \theta}_M$-equivalence classes.
For each $\alpha \ne \beta$ such that $\ell g(\bar a_n) = \ell g(\bar
a_\beta)$, let $\bar x_\alpha = \langle x_i:i < \ell g(\bar
a_\alpha)\rangle$ and choose $\varphi_{\alpha,\beta}(\bar x_\alpha),
\psi_{\alpha,\beta}(\bar x_\alpha) \in \Bbb L_{(2^{< \theta})^+,\theta}
[{\frak K}]$ such that, if possible we have $M \models
\varphi_{\alpha,\beta}[\bar a_\alpha] \wedge \neg
\varphi_{\alpha,\beta}[\bar a_\beta]$ and, under this, if possible $M \Vdash
``\psi_{\alpha,\beta}(\bar a_\alpha) \wedge \neg \psi_{\alpha,\beta}
(\bar a_\beta)$ but in any case $M \models \varphi_{\alpha,\beta}[\bar
a_\alpha]$ and $M \Vdash \psi_{\alpha,\beta}[\bar a_\alpha]$.  Let
$\varphi_\alpha(\bar x) = \wedge\{\varphi_{\alpha,\beta}(\bar
x_\alpha):\beta < \alpha^*,\beta \ne \alpha$ and $\ell g(\bar a_\beta)
=\ell g(\bar a_\alpha)\}$ and similarly $\psi_\alpha(\bar x_\alpha)$.  
Let $\Lambda_*$ be the closure of
$\{\varphi_{\alpha,\beta},\psi_{\alpha,\beta},\varphi_\alpha,\psi_\alpha:
\alpha \ne \beta < \alpha^*\}$ under subformulas 
and finitary operations.  Obviously,
clauses (a),(b) hold hence the existence of $\gamma^* < (2^{<
\theta})^+$ as required in clause (c) follows.  Clause (d) holds as
$\bar a {\Cal E}^{< \theta}_M \bar b \Rightarrow \sftp_{\Bbb
L_{\infty,\theta}[{\frak K}]} (\bar a,\emptyset,M) = \sftp_{\Bbb
L_{\infty,\theta}[{\frak K}]} (\bar b,\emptyset,M)$ using the
automorphisms and for $\alpha,\beta < \alpha_*$ such that $\ell
g(\bar a_\alpha) = \ell g(\bar a_\beta)$ we have $M \models (\forall
\bar x_\alpha)(\varphi_\alpha(\bar x_\alpha) = 
\varphi_\beta(\bar x_\beta)$ 
implies $\sftp_{\Bbb L_{(2^{< \theta})^+,\theta}[{\frak K}]}
(\bar a_\alpha,\emptyset,M) =
\sftp_{\Bbb L_{(2^{< \theta})^+,\theta}[{\frak K}]}(\bar a_\beta,
\emptyset,M)$ and even $\sftp_{\Bbb L_{\infty,\theta}[{\frak K}]}
(\bar a_\alpha,\emptyset,M) = \sftp_{\Bbb L_{\infty,\theta}[{\frak
K}]} (\bar a_\beta,\emptyset,M)$ recalling the 
choice of the $\varphi_{\alpha,\beta}$'s.

Clause (e) holds similarly by the choice of the
$\psi_{\alpha,\beta}$'s.  Clauses (f),(g) should also be clear.  (The
proof is similar to the proof of the classical \scite{734-0n.13}(3).)
  \hfill$\square_{\scite{734-eq.3.3}}$
\enddemo
\bigskip

\demo{\stag{734-eq.6} Observation}  Assume (\scite{734-eq.0}(b)$^-_\mu$ of
course and) $\Lambda \subseteq
\Bbb L_{\infty,\theta}[{\frak K}]$ and $\mu > 2^{< \theta}$ and
$\theta > \text{ LS}({\frak K})$.
\nl
1) The number of complete $\Bbb L_{\infty,\theta}[{\frak K}]$-types
realized in some $M \in K^*_\mu$, by a sequence of length $< \theta$ of
course, is $\le 2^{< \theta}$.  Hence every formula in $\Bbb
L_{\infty,\theta}[{\frak K}]$ is equivalent, for models from $K^*_\mu$
to a formula of quantifier depth $<(2^{<\theta})^+$, even from 
$\Lambda_* \subseteq \Bbb L_{(2^{< \theta})^+,\theta}[{\frak K}]$
where $\Lambda_*$ is in \scite{734-eq.3.3}(2).
\nl
2) Assume that $I_1 \subseteq I_2$ are well ordered, cf$(I_1)$, cf$(I_2) >
2^{<\theta}$ and $t \in I_2 \backslash I_1 \Rightarrow 2^{<\theta} <
\text{ cf}(I_1 \restriction \{s \in I_1:s <_{I_2} t\})$ and
$t \in I_2 \backslash I_1 \Rightarrow 2^{<\theta} < \text{\rm cf}(I_2
\restriction \{s \in I_2:(\forall r \in I_1)(r <_{I_2} t \equiv r
<_{I_2} s)\})$.  \ub{Then} EM$_{\tau({\frak K})}(I_1,\Phi) 
\prec_{\Bbb L_{\infty,\theta}[{\frak K}]} \text{ EM}_{\tau({\frak K})}
(I_2,\Phi)$.
\nl
3)  If $M = \text{ EM}_{\tau({\frak K})}(I,\Phi),
|I| = \mu,I$ well ordered of cofinality $> 2^{< \theta},
\bar a \in {}^\alpha M$ where $\alpha < \theta$ 
and $a_i = \sigma_i(\ldots,a_{t_{i,\ell}},
\ldots)_{\ell < n(i)}$ for $i < \alpha$ \ub{then} \sftp$_{\Lambda_*}
(\bar a,\emptyset,M)$ is determined by $\langle \sigma_i
(x_0,\dots,x_{n(\ell)-1}):i < \ell g(\bar a) \rangle$ and
the essential $\theta$-type of $\langle t_{i,\ell}:i < \ell
g(\bar a),\ell < n(i)\rangle$, see Definition \scite{734-eq.6.6} below.
\enddemo
\bn
Before proving \scite{734-eq.6}
\definition{\stag{734-eq.6.6} Definition}  1) For 
$\bar t = \langle t_i:i < \alpha \rangle \in {}^\alpha I,I$
well ordered, let the essential $\theta$-type of $\bar t$ in $I$ be
the essential $(\theta,(2^{< \theta})^+)$-type where for an ordinal
$\gamma$ we let the essential $(\theta,\gamma)$-type 
of $\bar t$ in $I$, {\rm estp}$_{\theta,\gamma}(\bar t,\emptyset,I)$ be the
following information stipulating $t_\alpha = \infty$:
\mr
\item "{$(a)$}"   the truth value of $t_i < t_j$ (for $i,j < \alpha$)
\sn
\item "{$(b)$}"   otp$([r_i,t_i)_I)$ for $i < \alpha$ where for $i \le
\alpha$ we let $r_i$ be the minimal member $r$ of $I$ such that 
otp$([r,t_i)_I) < \theta \times \gamma$ and $r \le_I t_i$ and $j <
\alpha \wedge t_j < t_i \Rightarrow t_j \le r$
\sn
\item "{$(c)$}"   Min$\{\theta \times \gamma,\text{otp}[s_i,r_i)_I\}$ 
for $i \le \alpha$ where we let $s_i$ be the minimal member of $I$ 
such that $(\forall j < \alpha)[t_j <_I t_i \Rightarrow t_j <_I s_i]$
\sn
\item "{$(d)$}"    Min$\{\theta,\text{cf}(I \restriction 
\{s:s <_I r_i\})\}$ for $i \le \alpha$ which may be zero.
\ermn
2) Let the function implicit in \scite{734-eq.6}(3) be called
$\bold t^\mu_\Lambda = \bold t^\mu_{{\frak K},\Lambda} = \bold
t^\mu_{{\frak K},\Phi,\Lambda}$, i.e., $\bold t^\mu_\Lambda(\bold
s,\bar \sigma) = \sftp_\Lambda(\bar a,\emptyset,M)$ when $\bar a = \langle
\sigma_i(\ldots,a_{t_{\beta(i,\ell)}},\ldots)_{\ell < n_i}:i < \ell g(\bar
a)\rangle,\bar \sigma = \langle 
\sigma_i(\ldots,x_{\beta(i,\ell)},\ldots)_{\ell < n};i < \ell
g(\bar a)\rangle$ and $\bold s$ is the essential $\theta$-type of $\langle
t_{i,\ell}:i < \ell g(\bar a),\ell < n_i\rangle$ in $I$. 
\nl
If $\Lambda = \Bbb L_{\infty,\theta}[{\frak K}]$ we may write just $\theta$.
\enddefinition
\bigskip

\demo{Proof of \scite{734-eq.6}}  1) By \scite{734-eq.3.3}(1) this 
holds for each $M \in K^*_\mu$.
\nl
2) It is known by Kino \cite{Kin66} that $I_1 \prec_{\Cal L} I_2$ if ${\Cal L}
\subseteq \{\varphi \in \Bbb L_{\infty,\theta}(\{<\}):\varphi$ has
quantifier depth $< (2^{< \theta})^+\}$.  From this the result
follows by part (1).

More fully let $\theta_r$ be the first regular cardinal $\ge \theta$, and
we say that the pair $(I_1,I_2)$ is $\gamma$-suitable when
we replace in the assumptions ``of cofinality $> 2^{< \theta}$" by
``of cofinality $\ge \theta$ and of order type divisible by $\theta
\times \gamma$".  Now we prove by induction on $\gamma$ that 
\mr
\item "{$\odot_1$}"  assume that for $\alpha < \theta$ and for
$\ell=1,2$ we have: $I_\ell$ is a well ordering, $\bar t^\ell =
\langle t^\ell_i:i < \alpha\rangle$ is $<_{I_\ell}$-increasing,
$t^\ell_0$ is the first element of $I_\ell$, we stipulate
$t^\ell_\alpha = \infty$ and otp$([t^\ell_i,t^\ell_{i+1})_{I_0}) = 
\theta_r \gamma \alpha^\ell_i + \beta_i$ where $\beta_i < \theta
\gamma$ and $(\text{cf}(\alpha^1_i) = \text{ cf}(\alpha^1_i)) \vee
(\text{ cf}(\alpha^1_i) \ge \theta \wedge \text{ cf}(\alpha^2_i) \ge
\theta)$.   
\nl
\ub{Then} for any formula $\varphi(\langle x_i:i < \alpha\rangle) \in
\Bbb L_{\infty,\theta}(\{<\})$ of quantifier depth $\le \gamma$
we have $I_1 \models \varphi[\bar t^1] \Leftrightarrow I_2 \models
\varphi[\bar t^2]$.
\ermn
Hence
\mr
\item "{$\odot_2$}"  if
$\vartheta(\bar x) \in \Bbb L_{\infty,\theta}(\{<\})$ has quantifier
depth $< \gamma$ and $(I_1,I_2)$ is $\gamma$-suitable and $\bar t \in
{}^{\ell g(\bar x)}(I_1)$ then $I_1 \models \varphi[\bar t]
\Leftrightarrow I_2 \models \theta[\bar t]$.
\ermn
3) Follows by part (2).  \hfill$\square_{\scite{734-eq.6}}$ 
\enddemo
\bigskip



\proclaim{\stag{734-eq.4} Claim}   Assume
\mr
\item "{$\boxdot$}"   $(a) \quad M \in K^*_\mu$
\sn
\item "{${{}}$}"  $(b) \quad \Lambda \subseteq 
\Bbb L_{\infty,\theta}[{\frak K}]$ is downward
closed, $|\Lambda| \le \chi$, {\rm LS}$({\frak K}) < \theta \le 
\chi < \mu$ and $2^{< \theta} \le \chi$ and 
$\theta = { \text{\rm cf\/}}(\theta) \vee \theta < \chi$ so $\Lambda =
\Lambda_*$ from \scite{734-eq.3.3} is O.K.
\sn
\item "{${{}}$}"  $(c) \quad$ in part (3),(4),(5) 
we assume $(\chi^{< \theta} \le \mu) \vee (\text{\rm cf}(\mu) \ge \theta)$
\sn
\item "{${{}}$}"  $(d) \quad$ for part (6) we assume 
{\rm cf}$(\mu) \ge \theta$ (hence the demand in clause (c)
\nl

\hskip35pt  holds).
\ermn
1) If $M \in K^*_\mu$ then 
$\{{\text{\rm gtp\/}}_\Lambda(\bar a,\emptyset,M):
\bar a \in {}^{\theta >} M\}$ has cardinality $\le 2^{<\theta}$.
\nl
2) If $(M,N) \in 
K_{\mu,\chi}$ \ub{then} we can find $N',(M,N) \le_{\frak K}
(M,N') \in K_{\mu,\chi}$ such that
\mr
\item "{$(*)$}"  if $\alpha < \theta$ and 
$\bar b \in {}^\alpha M$ and $\Lambda \subseteq \Bbb
L_{\infty,\theta}[{\frak K}]$ then for some
$\bar b' \in {}^\alpha(N')$ we have: for every $\bar a \in 
{}^{\theta >} N$, {\rm gtp}$_\Lambda(\bar a \char 94 \bar b,
\emptyset,M) = { \text{\rm gtp\/}}_\Lambda(\bar a \char 94 
\bar b',\emptyset,M)$.
\ermn
3) If $(M,N) \in K_{\mu,\chi}$, \ub{then} we 
can find $(M_1,N_1)$ such that $(M,N) \le_{\frak K}
(M_1,N_1) \in K_{\mu,\chi}$ and (note that $\bar y$ may be the empty sequence)
\mr
\item "{$(*)$}"  if $\exists \bar y \varphi(\bar y,\bar x) 
\in \Lambda$ and $\bar a \in 
{}^{\ell g(\bar x)}N$ then $M_1 \Vdash \neg \exists \bar y \varphi(\bar y,
\bar x)$ or for some \nl
$\bar b \in {}^{\ell g(\bar y)}(N_1)$ we have 
$M_1 \Vdash \varphi[\bar b,\bar a]$.
\ermn
4) In part (3) we can demand
\mr
\item "{$(*)^+$}"  if $\exists \bar y \varphi(\bar y,\bar x) \in
\Lambda$ and $\bar a \in {}^{\ell g(\bar x)}(N_1)$ 
then $M_1 \Vdash \neg(\exists \bar y) \varphi
(\bar y,\bar x)$ or for some $\bar b \in {}^{\ell g(\bar y)}
(N_1)$ we have $M_1 \models \varphi[\bar b,\bar a]$.
\ermn
5) In part (4) it follows that the pair 
$(M_1,N_1)$ is $\Lambda$-generic (most interesting for $\Lambda_*$, 
see \scite{734-eq.3.3}).
\nl
6) If $M_1 \in K^*_\mu$ then it is $\Lambda$-generic.
\endproclaim
\bigskip

\demo{Proof}  1) Proved just like \scite{734-eq.3.3}(1).
\nl
2) First assume $\theta$ is a successor cardinal.
As $M \in K^*_\mu$ \wilog \, $M = \text{ EM}_{\tau({\frak K})}
(I,\Phi)$ for some linear order
$I$ of cardinality $\mu$ as in \scite{734-am3.2.3}(1),(4) with
$\theta^-,\theta,\chi^+,\mu$ here standing for
$\mu,\theta_1,\theta_2,\lambda$ there.  It follows that for 
some $J \subseteq I$ of cardinality $\chi$ we have $N \subseteq
\text{ EM}_{\tau({\frak K})}(J,\Phi)$, and let $J^+ \subseteq I$ 
be such that $J \subseteq J^+,|J^+ \backslash J| = \chi$ and for every
$\bar t \in {}^{\theta >} I$ there is an automorphism $f$ of $I$ over
$J$ which maps $\bar t$ to some member of ${}^{\ell g(\bar t)}(J^+)$.

Lastly, let $N' = \text{ EM}_{\tau({\frak K})}(J^+,\Phi)$, it is easy
to check (see \scite{734-11.1}) that $(*)$ holds.  If $\theta$ is a limit
ordinal it is enough to prove for each $\partial < \theta$, a version of
$(*)$ with $\alpha < \partial$; and this gives $N'_\partial$.  Now we
choose $N'$ such that $\partial < \theta \Rightarrow N'_\partial
\le_{\frak K} N'$ and $(M,N') \in K_{\mu,\chi}$.
\nl
3),4),5),6)  We prove by induction on $\gamma$ that if we
let $\Lambda_\gamma$ be $\{\varphi(\bar x):\varphi(\bar x) \in \Lambda$ 
has quantifier depth $< 1+\gamma\}$ then parts (3),(4),(5),(6) holds for
$\Lambda_\gamma$.   For all four parts, $|\Lambda| \le \chi$ hence 
$|\Lambda_\gamma| \le \chi$ and it suffices to consider $\gamma < \chi^+$.  
For $\gamma = 0$ they are trivial and for $\gamma$ limit
also easy (let $\theta_r$ be the first regular $\ge \theta$ and extend 
$|\gamma|^+ \times \theta_r$ times taking care of $\Lambda_\beta$ in stage
$\gamma \times \zeta + \beta$ for each 
$\beta < \gamma$).  So let $\gamma = \beta +1$.

We first prove (3), but we have two cases (see clause (c)) of the
assumption.  If $\chi^{< \theta} \le \mu$ this is straight by
bookkeeping.  So assume {\rm cf}$(\mu) \ge \theta$.
Given $(M,N) \in K_{\mu,\chi}$ we try to 
choose by induction on $i < \chi^+$ a pair $(M_i,N_i)$ and for $i$ odd also
$\psi_i(\bar y_i,\bar x_i),\bar a_i,\bar b_i$ such that
\mr
\item "{$\circledast_1$}"  $(a) \quad (M_0,N_0) = (M,N)$
\sn
\item "{${{}}$}"  $(b) \quad (M_i,N_i) \in 
K_{\mu,\chi}$ is $\le_{\frak K}$-increasing continuous
\sn
\item "{${{}}$}"  $(c) \quad M_{i+1}$ is
$\Lambda_\beta$-generic for $i$ even
\sn
\item "{${{}}$}"   $(d) \quad$ for $i$ odd $\psi_i(\bar y_i,\bar x_i) 
\in \Lambda_\beta$ and
$\bar a_i \in {}^{\theta >} N$ and $\bar b_i \in {}^{\theta
>}(N_{i+1})$ are such \nl

\hskip25pt that $\ell g(\bar a_i) = 
\ell g(\bar x_i),\ell g(\bar b_i) = \ell g(\bar y_i)$ and
{\roster
\itemitem{ ${{}}$ }  $(\alpha) \quad \bar b \in {}^{\ell g(\bar y_i)}(M_i) 
\Rightarrow M_i \nVdash \psi_i[\bar b_i,\bar a]$ but
\sn
\itemitem{ ${{}}$ }   $(\beta) \quad M_{i+1} 
\Vdash \psi_i[\bar b_i,\bar a_i]$
\sn
\itemitem{ ${{}}$ }   $(\gamma) \quad$ for every $\bar b \in
{}^{\theta >}(M_{i+1})$ there is an automorphism of $M_{i+1}$ over
\nl

\hskip35pt $N_i$ mapping $\bar b$ into $N_{i+1}$. 
\endroster}
\ermn
If we succeed, by part (2) applied to the pair of models 
$(\dbcu_{i < \chi^+} M_i,N)$ as $\chi^+ \le \mu$ this pair belongs to 
$K_{\mu,\chi}$ we get
$N'$ as there, hence for some odd $i < \chi^+,N' \subseteq M_i$, let
$\zeta = i+2$ and this gives
a contradiction to the choice of $(\psi_\zeta,\bar a_\zeta,\bar b_\zeta)$.
\nl
[Why?  There is an automorphism $f$ of $M := \cup\{M_j:j < \chi^+\}$
over $N$ mapping $\bar b_\zeta$ into $N'$ hence into $M_i$ hence
$f(\bar b_\zeta) \in {}^{\theta >}(M_\zeta)$.  We know 
(by clause $(d)(\beta)$ above) that $M_{\zeta +1} \Vdash
\psi_\zeta[\bar b_\zeta,\bar a_\zeta]$ but $M_{\zeta +1} 
\le_{{\frak K}_\mu} M$ hence
$M \Vdash \psi_\zeta[\bar a_\zeta,\bar b_\zeta]$.  Recall that 
$f$ is an automorphism of $M$ over $N$
hence $M \Vdash \psi_\zeta[f(\bar b_\zeta),f(\bar a_\zeta)]$, but 
$\bar a_\zeta \in {}^{\theta >} N$ so $f(\bar a_\zeta) = 
\bar a_\zeta$ hence $M \Vdash \psi_\zeta[\bar b_\zeta,f(\bar
a_\zeta)]$ but $M_\zeta \le_{{\frak K}_\mu} M$ and $\bar a,f(\bar
b_\zeta)$ are from $M_\zeta$  hence
$M_\zeta \nVdash \neg\psi_\zeta[f(\bar b_\zeta),(\bar a_\zeta)]$.  However 
by clause (d)$(\alpha)$ of $\circledast_1$ we have $M_\zeta \nVdash 
\psi_\zeta[f(\bar b_\zeta),\bar a_\zeta]$.  
But as $i$ is an odd ordinal the last
two sentences contradicts clause $(c)$ of $\circledast_1$ applied to $i+1$.]
\nl
Hence we are stuck for some $i < \chi^+$.  Now for $i=0$ clause
$\circledast(a)$ gives a permissible value and for $i$ limit take
unions noting that clauses (c),(d) required nothing.  So $i=j+1$; if
$j$ is even we apply the induction hypothesis to part (6) for the pair
$(M_i,N_i)$.  Hence $j$ is odd so we cannot choose $\psi_j(\bar y,\bar x),
\bar a_j,\bar b_j$, recalling part (2) so the 
pair $(M_j,N_j)$ is as required thus proving (3) (for $\Lambda_\gamma$).

Second, we prove part (4).  We can now again 
try to choose by induction on $i < \chi^+$ a pair $(M_i,N_i)$ satisfying
\mr
\item "{$\circledast_2$}"  $(a) \quad (M_0,N_0) = (M,N)$
\sn
\item "{${{}}$}"  $(b) \quad (M_i,N_i) 
\in K_{\mu,\chi}$ is $\le_{\frak K}$-increasing continuous
\sn
\item "{${{}}$}"  $(c) \quad$ if $i = 2j+1$, \ub{then}
$(M_{i+1},N_{i+1})$ is as in part (3) for $\Lambda_\gamma$ 
with 
\nl

\hskip25pt $(M_i,N_i),(M_{i+1},N_{i+1})$ here standing for $(M,N),(M_1,N_1)$
there
\sn
\item "{${{}}$}"  $(d) \quad$ if $i=2j$ \ub{then} for some 
$\psi_i(\bar y_i,\bar x_i) \in \Lambda_\beta$ and
$\bar a_i \in {}^{(\ell g(\bar x_i))}(N_i)$ and
\nl

\hskip25pt  $\bar b_i \in {}^{(\ell g(\bar y_i))}(N_{i+1})$ we 
 have $M_{i+1} \Vdash \psi_i(\bar b_i,\bar a_i)$ but 
\nl

\hskip25pt $\bar b \in {}^{\ell g(\bar y_i)}(M_i) \Rightarrow M_i
\nVdash \psi_i[\bar b,\bar a_i]$.
\ermn
If we succeed, let $S_0 = \{\delta < \chi^+:\text{cf}(\delta) \ge \theta\}$,
so by an assumption $S$ is a stationary subset of $\chi^+$, i.e. as by
clause $\boxdot(b)$ we have $\theta = \text{\rm cf}(\theta) \le \chi
\vee \theta < \chi$; also for 
$\delta \in S_0$, as $\langle N_i:i < \delta \rangle$ is increasing 
with union $N_\delta$,
and $\delta = 2 \delta$ clearly $\bar a_\delta$ is well defined, so
for some $i(\delta) < \delta$ we have $\bar a_\delta \in 
{}^{\theta >}(N_{i(\delta)})$ and \wilog \, $i(\delta) = 2j(\delta)
+1$ for some $j(\delta)$ hence by clause (c) of $\circledast_2$ the pair
$(M_{i(\delta)+1},N_{i(\delta)+1})$ is as required there
contradiction as in the proof for part (3).  
Hence for some $i$ we cannot choose $(M_i,N_i)$.

For $i=0$ let $(M_i,N_i) = (M,N)$ so only clauses (a) + (b) of
$\circledast_2$ apply and are satisfied.  For $i$ limit take unions.
So $i=j+1$.
If $j=1$ mod 2, clause (d) of $\circledast_2$ is
relevant and we use part (3) for $\Lambda_\beta$ which holds as we have
just proved it.

Lastly, if $j=2$ mod 2 and we are stuck then 
the pair $(M_j,N_j)$ is as required. 
\nl
Third, Part (5) should be clear but we elaborate.

We prove by induction on $\gamma'$ that if $\varphi(\bar x) \in
\Lambda_\gamma$ has quantifier depth $< 1 + \gamma'$ then
for every $\bar a \in {}^{\ell g(\bar x)}(N_1)$ we have $M_1 \models
\varphi[\bar a] \Leftrightarrow N_1 \models \varphi[\bar a]$.  For
atomic $\varphi$ this is obvious and for $\varphi = \dsize
\bigwedge_{i < \alpha} \varphi_i$ should be clear.  If $\varphi(\bar
x) = \neg \psi(\bar x)$ note that in $(*)^+$ of part (4) we can use empty
$\bar y$ so $\neg(\exists \bar y)\psi(\bar x) = \neg \psi(\bar x)$.
Also for $\varphi(\bar x) = (\exists \bar y)\varphi'(\bar y,\bar x)$
we apply part (4).    

Fourth, we deal with part (6), so (see clause (d) of the assumption)
we have cf$(\mu) \ge \theta$.
Let $\chi = \langle \chi_i:
i < \text{ cf}(\mu)\rangle$ be constantly $\mu^-$ (so $\mu =
\chi^+_i$) if $\mu$ is a successor cardinal, and 
be increasing continuous with limit $\mu,2^{< \theta} <
\chi_i < \mu$ if $\mu$ is a limit cardinal recalling $2^{< \theta} <
\mu$ by $\boxdot(b)$.  Consider 
$K_{\mu,\bar \chi} = \{\bar M:\bar M = \langle M_i:
i \le \text{ cf}(\mu) \rangle$ is $\le_{\frak K}$-increasing continuous,
$M_{\text{cf}(\mu)} \in K^*_\mu$ and
$M_i \in K_{\chi_i}$ for $i < \text{ cf}(\mu)\}$ ordered by $\bar M^1
\le_{\frak K} \bar M^2$ \ub{iff} $i \le \text{ cf}(\mu) \Rightarrow M^1_i
\le_{\frak K} M^2_i$.

By \scite{734-eq.4} and part (5) for $\Lambda_\gamma$ which we proved
 we can easily find $\bar M \in K_{\mu,\bar \chi}$ 
such that $i < \text{ cf}(\mu) \Rightarrow (M_{\text{cf}(\mu)},
M_{i+1})$ is $\Lambda_\gamma$-generic; such
$\bar M$ we call $\Lambda_*$-generic.
Next
\mr
\item "{$\boxtimes$}"  if 
$\varphi(\bar x) \in \Lambda_\gamma$ and $\bar M$ is 
$\Lambda_\gamma$-generic, $\bar a \in {}^{\theta >}(M_i),i$
successor, $\varphi(\bar x) \in \Bbb L_{\infty,\theta}[{\frak K}]$ and
$\ell g(\bar x) = \ell g(\bar a)$ \ub{then} $M_{\text{cf}(\mu)} \models
\varphi[\bar a] \Leftrightarrow M_{\text{cf}(\mu)} \Vdash \varphi[\bar a]$.
\ermn
[Why?  Recalling cf$(\mu) \ge \theta$, we prove this by 
induction on the quantifier depth of $\varphi$.]

By the definition of ``$M$ is $\Lambda$-generic" and categoricity of
$K^*_\mu$ we are done.  \hfill$\square_{\scite{734-eq.4}}$
\enddemo
\bigskip

\demo{\stag{734-eq.5} Conclusion}  If $\mu \ge (2^{< \theta})^+,\theta
> \text{ LS}({\frak K})$ and cf$(\mu) \ge \theta > \text{ LS}({\frak K})$ 
\ub{then} every $M \in K^*_\mu$ is $\Bbb L_{\infty,\theta}[{\frak
K}]$-generic, hence if $M_1 \le_{\frak K} M_2$ are from $K^*_\mu$ 
\ub{then} $M_1 \prec_{{\Bbb L}_{\infty,\theta}[{\frak K}]} M_2$.
\enddemo
\bigskip

\remark{Remark}  1) With a little more care, if $\mu = \mu^+_0$ also
$\theta = \mu$ is O.K. but here this is prepheral.
\nl
2) $\theta \le \text{ LS}({\frak K})$ 
is not problematic, we just ignore it. 
\nl
3)  So \scite{734-eq.5} improve \scite{734-11.2}, i.e. we need cf$(\mu) \ge
\lambda (> \text{ LS}({\frak K}))$ instead $\mu = \mu^{< \lambda}$ 
but still there is a class of $\mu$ which are not covered.
\endremark
\bigskip

\demo{Proof}   Let $\Lambda_*$ be as in \scite{734-eq.3.3}(2) so in
particular $|\Lambda_*| \le 2^{< \theta}$.   Now \scite{734-eq.4}(6) and 
clause (g) of \scite{734-eq.3.3} proves the first assertion 
in \scite{734-eq.5}.  For the second
assume that $M_1 \le_{{\frak K}_\mu} M_2$ and we shall prove that $M_1
\prec_{\Bbb L_{\infty,\theta}[{\frak K}]} M_2$.

By the categoricity of ${\frak K}$ in $\mu$ or clause $(b)^-_\mu$ of
Hypothesis \scite{734-eq.0}, $K^*$ is categorical in $\mu$ hence 
$M_1,M_2 \in K^*_\mu$ are $\Lambda_*$-generic.
Suppose $\bar a \in {}^{(\ell g(\bar x))}(M_1),
\varphi(\bar x) \in \Lambda_*$, so by $M'_1$ being $\Lambda_*$-generic
(or $\boxtimes$ from the end of the proof of \scite{734-eq.4}
applied to $\bar M^2$) we have
\mr
\item "{$(*)_1$}"  $M_1 \models \varphi[\bar a] \Rightarrow M_1 \Vdash
\varphi[\bar a] \Rightarrow M_1 \models \varphi[\bar a]$
\ermn
and by $M_2$ being $\Lambda_*$-generic (or $\boxtimes$ from the end of
the proof of \scite{734-eq.4} applied to $\bar M^2$) we have
\mr
\item "{$(*)_2$}"   $M_2 \models \varphi[\bar a] \Rightarrow M_2 \Vdash
\varphi[\bar a] \Rightarrow M_2 \models \varphi[\bar a]$
\ermn
and by the definition of ``$M \Vdash \varphi[\bar a]$" recalling $M_1
\le_{{\frak K}_\mu} M_2$,
\mr
\item "{$(*)_3$}"  if $M_1 \Vdash \varphi'[\bar a]$ then $M_2 \Vdash
\varphi'[\bar a]$ for $\varphi'(\bar x) \in \{\varphi(\bar x),\neg
\varphi(\bar x)\}$.
\ermn
So both $M_1$ and $M_2$ satisfy $\varphi[\bar a]$ if $M_1$ satisfy it,
but this applies to $\neg \varphi[\bar a]$ too; so we are done.
\hfill$\square_{\scite{734-eq.5}}$
\enddemo
\bigskip

\proclaim{\stag{734-eq.7} Claim}  If $K$ is categorical also in $\mu^*$ 
or just Hypothesis \scite{734-eq.3.1} apply also to $\mu^*$, too, (with the
same $\Phi$) and $\mu^* \ge \mu^{< \theta} > \mu > \theta 
> { \text{\rm LS\/}}({\frak K})$ and $(*)$ below,
\ub{then} every $M \in K^*_\mu$ is 
$\Bbb L_{\infty,\theta}[{\frak K}]$-generic and $M_1 \in K^*_\mu
\wedge M_2 \in K^*_\mu \wedge M_1 \le_{{\frak K}_\mu} M_2 
\Rightarrow M_1 \prec_{\Bbb L_{\infty,\theta}[{\frak K}]} M_2$,
i.e. the conclusions of \scite{734-11.2}, \scite{734-eq.5} hold where
\mr
\item "{$(*)$}"  if $M \in K^*_{\mu^*}$ and 
$A \in [M]^\mu$ \ub{then} we can find $N \le_{\frak K} M$ 
such that $A \subseteq N \in K^*_\mu$ and for every $\varphi(\bar x) \in \Bbb
L_{\infty,\theta}[{\frak K}]$ and $\bar a \in {}^{\ell g(\bar x)}N$ we
have $M \Vdash \varphi[\bar a] \Leftrightarrow N \Vdash \varphi[\bar a]$.
\endroster
\endproclaim
\bigskip

\demo{Proof}  We shall choose $(M_i,N_i) \in K_{\mu^*,\mu}$ 
by induction on $i \le \theta^+$ such that not only $M_i \in
K^*_{\mu^*}$ (see the definition of $K_{\mu^*,\mu}$) but 
also $N_i \in K^*_\mu$ and this sequence of pairs is
$\le_{\frak K}$-increasing continuous.  For $i=0$ use any pair,
e.g. $M_0 = \text{ EM}_{\tau({\frak K})}(\mu^*,\Phi)$ and $N_0 =
\text{ EM}_{\tau({\frak K})}(\mu,\Phi)$.

For $i$ limit take unions, recalling $M_j,N_j$ are pseudo superlimit
for $j<i$.  

For $i = j+1$, let $N^+_j \le_{\frak K}
M_j$ be such that $N_j \subseteq N^+_j \in K_\mu$ and $(M_j,N^+_j)$
satisfies $(*)$ of the claim (standing for $(M,N))$.  Let $\Lambda_*$
be as in \scite{734-eq.3.3} for $\mu^*$.  Then by \scite{734-eq.4}(5) with
$(\mu^*,\mu,\theta)$ here standing for $(\mu,\chi,\theta)$ there
noting that in $\boxdot(c)$ there we use the case $\chi^{< \theta} \le
\mu$ which here means $\mu = \mu^{< \theta}$, we can 
choose a $\Lambda_*$-generic 
pair $(M_i,N_i) \in K_{\mu^*,\mu}$ above $(M_j,N^+_j)$ hence by 
\scite{734-eq.3.3}(2)(g) also it is a $\Bbb L_{\infty,\theta}[{\frak
K}]$-generic pair.  Now for
$j < \theta^+$, for $\bar a \in {}^{\theta >}(N_j)$,
we can read gtp$^{\mu^*}_\theta(\bar a,\emptyset,M_{j+1})$ 
and it is complete, but as by our use of $(*)$ it is
 the same as gtp$^\mu_\theta(\bar a,\emptyset,
N^+_{j+1})$.  So gtp$^\mu_\theta(\bar a,\emptyset,N^+_{j+1})$ is 
complete for every $\bar a \in
{}^{\theta >}(N_j)$, so also gtp$^\mu(\bar a,\emptyset,N_{\theta^+})$
is complete by monotonicity.  

Now if $\bar a \in {}^{\theta >}(N_{\theta^+})$ then for some $j <
\theta^+$ we have $\bar a \in {}^{\theta >}(N_j)$, so by the above
$p_{\bar a} := \text{\rm gpt}^{\mu^*}_\theta \negthinspace 
(\bar a,\emptyset,M_{j+1}) = \text{ gtp}^\mu_\theta
(\bar a,\emptyset,N^+_{j+1}) 
= \text{\rm gtp}^\mu_\theta(\bar a,\emptyset,N_{\theta^+})$ is 
complete and does not depend on $j$ as long as $j$ is large enough.

Now we prove that if $\bar a \in {}^{\theta >}(N_{\theta^+})$ then
$\varphi(\bar x) \in p_{\bar a} \Rightarrow
N_{\theta^+} \models \varphi[\bar a]$; and we prove this by induction
on the quantifier depth of $\varphi(\bar x)$; as usual the real case
is $\varphi(\bar x) = (\exists \bar y)\varphi(\bar y,\bar x)$.  Let $j
< \theta^+$ be such that $\bar a \in {}^{\ell g(\bar x)}(N_j)$, so
$p_{\bar a} = \text{ gtp}^{\mu^*}_\theta(\bar a,M_{j+1})$ so $M_{j+1}
\Vdash \varphi[\bar a]$ and by the choice of $(M_{j+1},N_{j+1})$ 
it follows that $N_{j+1} \models \varphi[\bar a]$ hence for some $\bar b \in
{}^{\ell g(\bar y)}(N_{j+1})$ we have 
$N_{j+1} \models \psi[\bar b,\bar a]$ 
hence $M_{j+1} \Vdash \psi(\bar b,\bar a)$, hence 
$\psi(\bar y,\bar x) \in p_{\bar b \char 94 \bar a}$ 
hence by the induction hypothesis
$N_{\theta^+} \models \psi[\bar b,\bar a]$ hence $N_{\theta^+} \models
\varphi[\bar a]$.   
\nl
${{}}$   \hfill$\square_{\scite{734-eq.7}}$
\enddemo
\bigskip

\demo{\stag{734-eq.8} Conclusion}  1) For each $\theta \ge 
\text{ LS}({\frak K})$ the family of $\mu > 2^{<\theta}$ in which $K$ is 
categorical but some (equivalent every) $M \in K_\mu$ is not
$\Bbb L_{\infty,\theta}[{\frak K}]$-generic is $\subseteq
\{[\mu_i,\mu^{< \theta}_i]:i < 2^{2^\theta}\}$ for some sequence
$\langle \mu_i:i < 2^{2^\theta}\rangle$ of cardinals.
\nl
2) Similarly for pseudo solvable, i.e. for each 
$\theta \ge \text{ LS}({\frak K})$ and $\Phi \in
\Upsilon^{\text{or}}_\theta$ for at most 
$\beth_2(\theta)$ cardinals $\mu > 2^{<\theta}$ we have $(\forall \alpha <
\mu)(|\alpha|)^{< \theta} < \mu)$ and for some $\mu^* \in [\mu,\mu^{<
\theta}]$ the pair $({\frak K},\Phi)$ is pseudo $\mu^*$-solvable but some
   $\equiv$ every $M \in K^*_{\Phi,\mu^*}$ is not 
$\Bbb L_{\infty,\theta^+}[{\frak K}]$-generic.
\enddemo
\bigskip

\demo{Proof}  Straight.
Note that it is enough to prove this for each $\Phi$ separately.

Toward contradiction assume $\langle \mu_\varepsilon:\varepsilon <
(\beth_2(\theta))^+\rangle$ is an increasing sequence of such
cardinals, satisfying $(\mu_\varepsilon)^{< \theta} < \mu_{\varepsilon
+1}$ and choose $I_\varepsilon \times \mu_\varepsilon \times (2^{<
\theta})^+$, hence 
$\langle I_\varepsilon:\varepsilon < (\beth_2(\theta))^+\rangle$ is an
increasing sequence  of linear orders as in \scite{734-eq.6}, in
particular, well ordered.  Let $\Lambda_\varepsilon =
\Lambda^*_{{\frak K},\Phi,\mu_\varepsilon,\theta}$ be from
\scite{734-eq.3.3}(2) applied to $I = I_\varepsilon$ hence to any ordinal
$< \mu^+_\varepsilon$ of cofinality $> 2^{< \theta}$.
Now the number of functions 
$\bold t^{\mu_\varepsilon}_{{\frak K},\Lambda_\varepsilon}$ (see Observation
\scite{734-eq.6}(3) and Definition \scite{734-eq.6.6}(2)) is at most
$\beth_2(\theta)$, so for some $\varepsilon < \zeta <
(\beth_2(\theta))^+$ we have 
$\bold t^{\mu_\varepsilon}_{{\frak K},\Lambda_\varepsilon} =
\bold t^{\mu_\zeta}_{{\frak K},\Lambda_\zeta}$.

Now apply \scite{734-eq.7} with $(\mu_\zeta,\mu_\varepsilon)$ here
standing for $(\mu^*,\mu)$ there, $(*)$ there holds easily by
\scite{734-eq.6}(3) so we get a contradiction. \hfill$\square_{\scite{734-eq.8}}$
\enddemo
\bn
\centerline {$* \qquad * \qquad *$}
\bn
For the rest of this section we note some 
basic facts on the dependency on $\Phi$ (not used here).
\definition{\stag{734-ep.1} Definition}  1) We define a two-place relation
${\Cal E}_\kappa = {\Cal E}^{\text{or}}_\kappa[{\frak K}]$ on
$\Upsilon^{\text{or}}_\kappa[{\frak K}]$, so $\kappa \ge \text{
LS}({\frak K}):\Phi_1 {\Cal E}_\kappa \Phi_2$ \ub{iff} for every linear orders
$I_1,I_2$ there are linear orders $J_1,J_2$ extending $I_1,I_2$
respectively such that EM$_{\tau({\frak K})}(J_1,\Phi)$,
EM$_{\tau({\frak K})}(J_2,\Phi)$ are isomorphic.
\nl
2) We define $\le^{\text{or}}_\kappa = \le^{\text{or}}_\kappa[{\frak K}]$, a
two-place relation on $\Upsilon^{\text{or}}_\kappa[{\frak K}]$ as
in part (1) only in the end EM$_{\tau({\frak K})}(J_1,\Phi_1)$ can be
$\le_{\frak K}$-embedded into $\text{EM}_{\tau({\frak K})}(J_2,\Phi_2)$.
\enddefinition
\bigskip

\proclaim{\stag{734-ep.2} Claim}  1) The following conditions on
$\Phi_1,\Phi_2 \in \Upsilon^{\text{or}}_\kappa[{\frak K}]$ are
equivalent
\mr
\item "{$(a)$}"  $\Phi_1 {\Cal E}_\kappa \Phi_2$
\sn
\item "{$(b)$}"  there are $I_1,I_2 \in K^{\text{lin}}$ of cardinality
$\ge \beth_{1,1}(\kappa)$ such that {\rm EM}$_{\tau({\frak
K})}(I_1,\Phi_1)$, {\rm EM}$_{\tau({\frak K})}(I_2,\Phi)$ are
isomorphic
\sn
\item "{$(c)$}"  there are $\Phi'_1,\Phi'_2$ satisfying $\Phi_\ell
\le^\otimes \Phi'_\ell \in \Upsilon^{\text{or}}_\kappa[{\frak K}]$ for
$\ell=1,2$ such that $\Phi'_1,\Phi'_2$ are essentially equal (see
Definition \scite{734-ep.3} below).
\ermn
2) The following conditions are equivalent
\mr
\item "{$(a)$}"  $\Phi_1 \le^{\text{or}}_\kappa \Phi_2$ recall
$\le_\kappa = \le^{\text{or}}_\kappa [{\frak K}]$
\sn
\item "{$(b)$}"   there are $I_1,I_2 \in K^{\text{lin}}$ of
cardinality $\ge \beth_{1,1}(\kappa)$ such that {\rm EM}$_{\tau({\frak
K})}(I_1,\Phi_1)$ can be $\le_{\frak K}$-embedded into
{\rm EM}$_{\tau({\frak K})}(I_2,\Phi_2)$
\sn
\item "{$(c)$}"   for every $I_1 \in K^{\text{lin}}$ there is $I_2 \in
K^{\text{lin}}$ such that {\rm EM}$_{\tau({\frak K})}(I_1,\Phi_1)$ can
be $\le_{\frak K}$-embedded into {\rm EM}$_{\tau({\frak K})}(I_2,\Phi_2)$.
\endroster
\endproclaim
\bigskip

\definition{\stag{734-ep.3} Definition}  $\Phi_1,\Phi_2 \in
\Upsilon^{\text{or}}_\kappa[{\frak K}]$ are essentially equal when
for every linear order $I$ there is an isomorphism $f$ from
EM$_{\tau({\frak K})}(I,\Phi_1)$ onto EM$_{\tau({\frak K})}(I,\Phi_2)$
such that for any $\tau_{\Phi_1}$-term $\sigma_1(x_0,\dotsc,x_{n-1})$
there is a $\tau_{\Phi_2}$-term $\sigma_2(x_0,\dotsc,x_{n-1})$ such
that: $t_0 <_I \ldots <_I t_{n-1} \Rightarrow f(a_1) = a_2$, where
$a_\ell$ is $\sigma_\ell(a_{t_0},\dotsc,a_{t_{n-1}})$ as computed in
EM$(I,\Phi_\ell)$ for $\ell=1,2$.
\enddefinition
\bigskip

\demo{Proof of \scite{734-ep.2}}  Straight (particularly recalling such
proof in \scite{734-11.14}(1)).  \hfill$\square_{\scite{734-ep.3}}$
\enddemo
\bigskip

\proclaim{\stag{734-ep.4} Claim}  1) ${\Cal E}_\kappa =
{\Cal E}^{\text{or}}_\kappa[{\frak K}]$ is an equivalence relation and
$\Phi_1 {\Cal E}^{\text{or}}_\kappa[{\frak K}] \Phi_2 \Rightarrow
\Phi_1 \le^{\text{or}}_\kappa [{\frak K}] \Phi_2$.
\nl
1A) In fact if $\langle \Phi_\varepsilon:\varepsilon <
\varepsilon(*)\rangle$ are pairwise ${\Cal E}_\kappa$-equivalent and 
$\varepsilon(*) \le \kappa$ \ub{then} we can find $\langle
\Phi'_\varepsilon:\varepsilon < \kappa \rangle$ satisfying
$\Phi'_\varepsilon \le^\otimes \Phi'_\varepsilon$ for $\varepsilon <
\varepsilon(*)$ such that the $\Phi'_\varepsilon$ for $\varepsilon <
\varepsilon(*)$ are pairwise essentially equal.
\nl
2) $\le^{\text{or}}_\kappa$ is a partial order.
\nl
3) If $\Phi_1,\Phi_2 \in \Upsilon^{\text{or}}_\kappa[{\frak K}]$ are
essentially equal \ub{then} $({\frak K},\Phi_1)$ is psuedo/weakly/strongly 
$(\mu,\kappa)$-solvable iff $({\frak K},\Phi_2)$ is
pseudo/weakly/strongly $(\mu,\kappa)$-solvable.
\nl
4) If $\Phi_1 \in \Upsilon^{\text{or}}_\kappa[{\frak K}]$ is strongly
$(\mu,\kappa)$-solvable and $\Phi_2$ exemplifies ${\frak K}$ is
$(\mu,\kappa)$-solvable \ub{then} $\Phi_1 {\Cal E}_\kappa \Phi_2$.
\nl
5) If ${\frak K}$ is categorical in $\mu$ and $\mu > \kappa \ge 
\text{\rm LS}({\frak K})$ \ub{then} every $\Phi \in
\Upsilon^{\text{or}}_\kappa[{\frak K}]$ is strongly
$(\mu,\kappa)$-solvable.
\nl
6) Assume $({\frak K},\Phi_\ell)$ is pseudo $(\mu,\kappa)$-solvable
and $\mu \ge \beth_{1,1}(\kappa)$ for $\ell=1,2$.  \ub{Then} 
$\Phi_1 {\Cal E}_\kappa \Phi_2$ iff $\Phi_1
\le^{\text{or}}_\kappa[{\frak K}] \Phi_2 \wedge \Phi_2
\le^{\text{or}}_\kappa [{\frak K}] \Phi_1$.
\nl
7) If $\Phi_1 \le^{\text{or}}_\kappa \Phi_2$ and $\Phi_1$
 is strongly $(\mu,\kappa)$-solvable or just pseudo $(\mu,\kappa)$-solvable
\ub{then} $\Phi_1,\Phi_2$ are
 ${\Cal E}^{\text{or}}_\kappa[{\frak K}]$-equivalent.
\endproclaim
\bigskip

\demo{Proof}  Easy, use \scite{734-11.14}(1) and its proof.
 \hfill$\square_{\scite{734-ep.4}}$ 
\enddemo
\goodbreak

\head {\S3 Categoricity for cardinals on a club} \endhead  \resetall \sectno=3
 \spuriousreset
\bigskip

We draw here an easy conclusion from \S2, getting that on a closed
unbounded class of cardinals which is $\aleph_0$-closed we get a
constant answer to being categorical.  This is, of course,
considerably weaker than conjecture \scite{734-0n.0} but still is a
progress, e.g. it shows that the categoricity spectrum is not totally chaotic.

We concentrate on the case the results of \S1 holds (e.g. $\mu =
\mu^\lambda$) for the $\lambda$'s with which we deal.  To eliminate
this extra assumption we need \S2.  This section is not used later.
Note that \scite{734-e.1.7} is continued (and improved) in \cite{Sh:F820}
and Exercise \scite{734-e.2.35}, \cite{Sh:F782} improve 
\scite{734-e.2.1}; similarly \scite{734-e.2.28}.
\bn
In the claims below we concentrate on fix points of the sequence of
$\beth_\alpha$'s. 
\demo{\stag{734-e.0} Hypothesis}  As in Hypothesis \scite{734-11.0},
(i.e. ${\frak K}$ is an a.e.c. with models of arbitrarily large cardinality).
\enddemo
\bigskip

\definition{\stag{734-e.1} Definition}   1) Let Cat$_{\frak K}$ be the class of
cardinals in which ${\frak K}$ is categorical. 
\nl
1A) Let Sol = Sol$_{{\frak K},\Phi} = \text{ Sol}^1_{{\frak K},\Phi}$
be the class of $\mu > 
\text{ LS}[{\frak K}]$ such that $({\frak K},\Phi)$ is pseudo
$\mu$-solvable.  Let Sol$^2_{{\frak K},\Phi}[\text{Sol}^3_{{\frak
K},\Phi}]$ be the class of $\mu >
\text{ LS}({\frak K})$ such that $({\frak K},\Phi)$ is weakly [strongly]
$\mu$-solvable. 
\nl
2) Let mod-com$_{{\frak K},\Phi}$ be the class of pairs $(\mu,\theta)$
such that: $\mu > \theta \ge \text{ LS}({\frak K})$ and 
$\Bbb L_{\infty,\theta^+}[{\frak K}]$ is $\mu$-model complete (on
$K^*_{\Phi,\mu}$, see Definition \scite{734-eq.1}(3)(b), \scite{734-eq.1}(5)).
\nl
3) Let Cat$'_{\frak K}$ be the class of $\mu \in 
\text{ Cat}_{\frak K}$ such that: $\mu \ge
\beth_{1,1}(\text{LS}({\frak K}))$ and if LS$({\frak K}) \le \theta$ and
$\beth_{1,1}(\theta) \le \mu$ \ub{then} $\Bbb L_{\infty,\theta^+}
[{\frak K}]$ is $\mu$-model complete.
\nl
3A) For $\Phi \in \Upsilon^{\text{or}}_{\frak K}$ let 
Sol$^{k,*}_{{\frak K},\Phi}$ be the class of $\mu \in 
\text{ Sol}^k_{{\frak K},\Phi}$ 
such that $\mu \ge \beth_{1,1}(\text{LS}({\frak K}))$ and: 
if LS$({\frak K}) \le \theta$ and
$\beth_{1,1}(\theta) \le \mu$ then the pair $(\Bbb L_{\infty,\theta^+}[{\frak
K}],\Phi)$ is $\mu$-model complete. 

Let Sol$^{\ell,< \theta}_{{\frak K},\Phi}$ be the class of $\lambda
\in \text{ Sol}^\ell_{{\frak K},\Phi}$ such that $\Bbb
L_{\infty,\theta}[{\frak K}]$ is $\mu$-model complete (see \sectioncite[\S2]{734}).

Let Sol$'_{{\frak K},\Phi} = \text{ Sol}^{1,*}_{{\frak K},\Phi}$.
 Instead $k,*$ we may write $3+k$.
\nl
4) Let $\bold C = \{\lambda:\lambda = \beth_\lambda$ and cf$(\lambda)
= \aleph_0\}$. 
\enddefinition
\bn
\margintag{734-e.1.3}\ub{\stag{734-e.1.3} Exercise}:  1) The conclusion of \scite{734-11.2}(1)
equivalently \scite{734-11.2}(2) means that $\theta \le \lambda \Rightarrow
(\mu,\theta) \in \text{ mod-com}_{{\frak K},\Phi}$.
\nl
2) Write down the obvious implications.
\bigskip

\proclaim{\stag{734-e.1.7} Claim}  If $\mu > \lambda = \beth_\lambda >
\kappa \ge \text{\rm LS}({\frak K})$ and $\Phi \in
\Upsilon^{\text{or}}_\kappa[{\frak K}]$, 
{\rm cf}$(\lambda) = \aleph_0$ \ub{then}
$\mu = \mu^{< \lambda} \Rightarrow \mu \in \text{\rm Sol}'_{{\frak
K},\Phi} \Rightarrow \lambda \in \text{\rm Sol}'_{{\frak K},\Phi}$.
\endproclaim
\bigskip

\demo{Proof}  The first implication holds by \scite{734-11.2}(2) and
\scite{734-e.1.3}.  The second implication, its assumption implies
Hypothesis \scite{734-11.2A}, see \scite{734-e.1.3}(1) hence its conclusion holds
by \scite{734-11.28}.  
\nl
${{}}$  \hfill$\square_{\scite{734-e.1.7}}$
\enddemo
\bigskip

\demo{\stag{734-e.2} Observation}   $K_\lambda$ is categorical in $\lambda$
(hence Hypothesis \scite{734-11.2A} holds), \ub{if}:
\mr
\item "{$\circledast_\lambda$}"  $\lambda = \beth_\lambda =
\sup(\lambda \cap \text{\rm Cat}'_{\frak K}) >  
{ \text{\rm LS\/}}({\frak K})$ and $\aleph_0 = { \text{\rm
cf\/}}(\lambda)$. 
\endroster
\enddemo
\bigskip

\demo{Proof}  Fix $\Phi \in \Upsilon^{\text{or}}_{\frak K}$, now clearly
Sol$'_{{\frak K},\Phi} \supseteq \text{ Cat}'_{\frak K}$ by their definitions.

By the assumptions we can find $\langle \mu_n:n < \omega \rangle$ such that 
$\lambda = \Sigma\{\mu_n:n < \omega\}$, LS$({\frak K}) <
\mu_n \in \text{ Cat}'_{\frak K}$ and $\beth_{1,1}(\mu'_n) <
\mu_{n+1}$ where $\mu'_n = \beth_{1,1}(\mu_n)$.
As every $M \in K_{\mu_{n+1}}$ is $\Bbb L_{\infty,\mu'_n}
[{\frak K}]$-generic (as $K_{\mu_{n+1}} \subseteq K_{\Phi,\mu_{n+1}}$
and $\mu_{n+1} \in \text{ Cat}'_{\frak K}$) easily
\mr
\item "{$(*)_0$}"  if $M \le_{\frak K} N$ are from $K^*_{\Phi,\ge
\mu_{n+1}}$ then $M \prec_{\Bbb L_{\infty,\mu'_n}[{\frak K}]} N$.
\ermn
Let $M^\ell \in K_\lambda$, for $\ell \in \{1,2\}$; so we can find a
$\le_{\frak K}$-increasing sequence $\langle M^\ell_n:n < \omega
\rangle$ such that $M^\ell_n \in K_{\mu_n},M^\ell_n \le_{\frak K}
M^\ell_{n+1} \le_{\frak K} M^\ell$ and $M^\ell = 
\cup\{M^\ell_n:n < \omega\}$.  Now
\mr
\item "{$(*)_1$}"  $M^\ell_n \in K^*_{\Phi,\mu_n}$. 
\ermn
[Why?  As ${\frak K}$ is categorical in $\mu_n = \|M^\ell_n\|$.]
\mr
\item "{$(*)_2$}"  if $\alpha \le \mu_n,n < m < k$ and $\bar a,\bar b
\in {}^\alpha(M^\ell_m)$ \ub{then}: 
{\roster
\itemitem{ $(a)$ }   $\sftp_{\Bbb L_{\infty,\mu'_n}[{\frak K}]}
(\bar a,\emptyset,M^\ell_m) = \sftp_{\Bbb L_{\infty,\mu'_n}}
(\bar b,\emptyset,M^\ell_m)$ 
\ub{iff} 
\sftp$_{{\Bbb L}_{\infty,\mu'_n}[{\frak K}]}(\bar a,\emptyset,M^\ell_k) =
\sftp_{\Bbb L_{\infty,\mu'_n}}(\bar b,\emptyset,M^\ell_k)$. 
\sn
\itemitem{ $(b)$ }   if $\sftp_{\Bbb L_{\infty,\mu'_n}[{\frak K}]}
(\bar a,\emptyset,M^\ell_k) = \sftp_{\Bbb L_{\infty,\mu'_n}[{\frak K}]}
(\bar b,\emptyset,M^\ell_k)$ \ub{then} 
$\sftp_{\Bbb L_{\infty,\mu'_m}[{\frak K}]}(\bar a,\emptyset,M^\ell_k) =
\sftp_{\Bbb L_{\infty,\mu'_m}[{\frak K}]}(\bar b,\emptyset,M^\ell_k)$.
\endroster}
\ermn
[Why?  Clause (a) by $(*)_0$, clause (b) by \scite{734-11.5}(3).]
\mr
\item "{$(*)_3$}"  $M^1_n \cong M^2_n$. 
\ermn
[Why?  As ${\frak K}$ is categorical in $\mu_n$.]

We now proceed as in the proof of \scite{734-11.20}.
Let ${\Cal F}_n = \{f$: for some $\bar a_1,\bar a_2$ and $\alpha
<\mu_n$ we have
$\bar a_\ell \in {}^\alpha(M^\ell_{n+2})$ for $\ell=1,2$,
\sftp$_{{\Bbb L}_\infty,\mu_{n+1}[{\frak K}]}(\bar a_1,\emptyset,M^1_{n+2}) =
\text{ \sftp }_{{\Bbb L}_\infty,\mu_{n+1}[{\frak K}]}
(\bar a_2,\emptyset,M^2_{n+1})$
and $f$ is the function which maps $\bar a_1$ into $\bar a_2\}$,
(actually can use $\alpha = \mu_n$).
\sn
By the hence and forth argument we can find $f_n \in {\Cal F}_n$
by induction on $n < \omega$ such that $M^1_n \subseteq \text{
Dom}(f_{2n+2}),M^2_n \subseteq \text{ Rang}(f_{2n+2})$ and $f_n
\subseteq f_{n+1}$; hence $\cup\{f_n:n < \omega\}$ is an isomorphism from $M^1$
onto $M^1$.  \hfill$\square_{\scite{734-e.1.7}}$
\enddemo
\bigskip

\proclaim{\stag{734-e.2.1} Claim}  ${\frak K}$ is categorical in $\lambda$
\ub{when}: 
\mr
\item "{$\circledast^+_\lambda$}"  $\lambda = \beth_\lambda >
\text{\rm LS}({\frak K})$ and $\lambda = 
\text{\rm otp}(\text{\rm Cat}_{\frak K} \cap \lambda \cap \bold C)$ 
and {\rm cf}$(\lambda) = \aleph_0$.
\endroster  
\endproclaim
\bigskip

\demo{Proof}  Fix $\Phi$ as in the proof of \scite{734-e.1.7}.  
Let $\langle \theta_n:n < \omega \rangle$ be increasing
such that $\lambda = \Sigma\{\theta_n:n < \omega\}$ and LS$({\frak K})
< \theta_0$.  For each $n$, by \scite{734-eq.8} we know $\{\mu \in
\text{\rm Cat}_{\frak K}:\mu > \theta_n$ and the $M \in K_\mu$ is not
$\Bbb L_{\infty,\theta^+_n}$-generic$\}$ is ``not too large", i.e. is
included in the union of at most $\beth_2(\theta_n)$ intervals of the
form $[\chi,\chi^{\theta_n}]$.  Now we choose $(n(\ell),\mu_\ell)$ by
induction on $\ell < \omega$ such that
\mr
\item "{$\circledast$}"  $(a) \quad n(\ell) < \omega$ and $\mu_\ell \in
\text{\rm Cat}_{\frak K} \cap\lambda$
\sn
\item "{${{}}$}"  $(b) \quad$ if $\ell = k+1$ then $n(\ell) >
n(k),\theta_{n(\ell)} > \mu_k,\mu_\ell \in \text{\rm Cat}_{\frak K} \cap
\lambda \backslash \theta^+_{n(\ell)}$ and the 
\nl

\hskip25pt $M \in K_{\mu_\ell}$ is
$\Bbb L_{\infty,\theta_{n(\ell)}}[{\frak K}]$-generic (hence
$\Bbb L_{\infty,\mu^+_k}[{\frak K}]$-generic). 
\ermn
This is easy and then continue as in \scite{734-e.2}. 
\hfill$\square_{\scite{734-e.2.1}}$ 
\enddemo
\bn
We have essentially proved
\proclaim{\stag{734-e.2.28} Theorem}  In \scite{734-e.2}, \scite{734-e.2.1} we can
use {\rm Sol}$_{{\frak K},\Phi}$, {\rm Sol}$'_{{\frak K},\Phi}$
instead of {\rm Cat}$_{\frak K}$, {\rm Cat}$'_{\frak K}$.
\endproclaim
\bn
\margintag{734-e.2.35}\ub{\stag{734-e.2.35} Exercise}:  For 
Claim \scite{734-11.20}(2), Hypothesis \scite{734-11.2A} suffice.
\nl
[Hint: The proof is similar to the existing one using \scite{734-11.5}.]
\goodbreak

\head {\S4 Good Frames} \endhead  \resetall \sectno=4
 \spuriousreset
\bigskip

Here comes the main result of \chaptercite{734}: from categoricity (or
solvability) assumptions we derive the existence of good $\lambda$-frames.

Our assumption is such that we can apply \S1.
\demo{\stag{734-f.0} Hypothesis}  1) 
\mr
\item "{$(a)$}"   ${\frak K}$ is an a.e.c.
\sn
\item "{$(b)$}"   $\mu > \lambda = \beth_\lambda > \text{ LS}({\frak
K})$ and cf$(\lambda) =\aleph_0$; 
\sn
\item "{$(c)$}"  $\Phi \in \Upsilon^{\text{or}}_{\frak K}$
\sn
\item "{$(d)$}"  ${\frak K}$ is categorical in $\mu$ or just
\sn
\item "{$(d)^-$}"  $({\frak K},\Phi)$ is pseudo superlimit in $\mu$
(this means $\Phi \in \text{ Sol}^1_{{\frak K},\Phi}$; so
\scite{734-11.2A}(1) holds)
\sn
\item "{$(e)$}"  also \scite{734-11.2A}(2)(a) holds, i.e. the 
conclusion of \scite{734-11.2}(2) holds.
\ermn
2) In addition we may use some of the following but then we mention them and 
(we add superscript $*$ when used; note that $(g) \Rightarrow (f)$ by
\scite{734-11.21})
\mr
\item "{$(f)$}"  $K^*_\lambda$ is closed under $\le_{\frak K}$-increasing
unions (justified by \scite{734-11.20})
\sn
\item "{$(g)$}"  $\langle \lambda_n:n < \omega\rangle$ is increasing,
$\lambda_0 > \text{\rm LS}({\frak K}),\lambda = \Sigma\{\lambda_n:n <
\omega\}$ and the assumptions of \scite{734-11.20} holds.
\endroster
\enddemo
\bigskip

\demo{\stag{734-f.0.7} Observation}  1) ${\frak K}^*_\lambda$ is
categorical.
\nl
2) ${\frak K}^*_\lambda$ has amalgamation.
\nl
3)$^*$ (We assume (f) of \scite{734-f.0}(2)).  ${\frak K}_\lambda$ is a
$\lambda$-a.e.c. 
\enddemo
\bigskip

\demo{Proof}  1) By \scite{734-11.4A}(1) or \scite{734-11.5}(4) as
cf$(\lambda) = \aleph_0$.
\nl
2) By \scite{734-11.15}(1).
\nl
3) As in \scite{734-11.21}, (i.e. as
 $\le_{{\frak K}^*_\lambda} = \le_{\frak K} \restriction 
{\frak K}$, closure under unions of $\le_{\frak K}$-increasing chains
is the only problematic point and it holds by (f) of \scite{734-f.0}(2)).
\hfill$\square_{\scite{734-f.0.7}}$
\enddemo
\bigskip

\remark{\stag{734-f.0.5} Remark}  1) Why do we not assume \scite{734-f.0}(1),(2) all
the time?  The main reason is that for proving some of the results
assuming \scite{734-f.0}(1),(2) we use some such results on smaller
cardinals on which we use \scite{734-f.0}(1) only.
\nl
2) Note that it is not clear whether improvement by using \scite{734-f.0}(1)
only will have any affect when (or should we say if) we succeed to
have the parallel of \sectioncite[\S12]{705}.
\endremark
\bigskip


\proclaim{\stag{734-f.1} Claim}  1) Assume $M_0 \le_{{\frak K}^*_\lambda}
M_\ell,\alpha < \lambda$ and $\bar a_\ell \in {}^\alpha(M_\ell)$ for
$\ell = 1,2$ and $\kappa := \beth_{1,1}(\beth_2(\theta)^+)$ where
$\theta := |\alpha| + \text{\rm LS}({\frak K})$ so $\kappa < \lambda$.  If 
$\sftp_{{\Bbb L}_{\infty,\kappa}[{\frak K}]}
(\bar a_1,M_0,M_1) = \sftp_{{\Bbb L}_{\infty,\kappa}[{\frak
K}]}(\bar a_2,M_0,M_2)$ \ub{then} {\rm \ortp}$_{{\frak K}^*_\lambda}
(\bar a_1,M_0,M_1) = { \text{\rm \ortp\/}}_{{\frak K}^*_\lambda}
(\bar a_2,M_0,M_2)$.
\nl
2) If $M_1 \le_{{\frak K}^*_\lambda} M_2$ \ub{then} $M_1 
\prec_{\Bbb L_{\infty,\theta}[{\frak K}]} M_2$ for every $\theta <
\lambda$, and moreover $M_1 \prec_{\Bbb L_{\infty,\lambda}[{\frak K}]} M_2$. 
\nl
2A) If $M_0 \le_{{\frak K}^*_\lambda} M_\ell$ for $\ell=1,2$ and
$\ortp_{{\frak K}^*_\lambda}(\bar a_1,M_0,M_1) = 
\ortp_{{\frak K}^*_\lambda}(\bar a_2,M_0,M_2)$ and $\bar a_\ell \in
{}^\alpha(M_0),\alpha < \kappa \le \lambda$ \ub{then} 
$\sftp_{\Bbb L_{\infty,\kappa}[{\frak K}]}(\bar a_1),M_0,M_1) = 
 \sftp_{\Bbb L_{\infty,\kappa}[{\frak K}]}(\bar a_2,M_0,M_2)$.
\nl
2B) In part (1), if $M_\ell \le_{{\frak K}^*_\lambda} M'_\ell$ for
$\ell=1,2$ \ub{then} $\sftp_{\Bbb L_{\infty,\kappa}[{\frak K}]}(\bar
a_1,M,M'_1) = \sftp_{\Bbb L_{\infty,\kappa}[{\frak K}]}(\bar a_2,M,M'_2)$.
\nl
3) Assume that $M_0 \le_{{\frak K}^*_\lambda} M_1 
\le_{{\frak K}^*_\lambda} M_2 \le_{{\frak K}^*_\lambda} M_3,\bar a \in
{}^\alpha(M_2),\alpha < \lambda$ and $\kappa = \beth_{1,1}(|\alpha| +
\text{\rm LS}({\frak K})) < \theta < \lambda$.  \ub{Then} 
\mr
\item "{$(a)$}"  from $\sftp_{\Bbb L_{\infty,\kappa}[{\frak K}]}(\bar
a,M_1,M_2)$ we can compute
$\sftp_{\Bbb L_{\infty,\theta}[{\frak K}]}(\bar a,M_1,M_2)$ and
$\sftp_{\Bbb L_{\infty,\lambda}[{\frak K}]}(\bar a,M_0,M_3)$
\sn
\item "{$(b)$}"  from $\sftp_{\Bbb L_{\infty,\kappa}[{\frak K}]}(\bar
a,\emptyset,M_2)$ we can compute 
$\sftp_{\Bbb L_{\infty,\theta}[{\frak K}]}(\bar a,\emptyset,M_2)$ and
even $\sftp_{\Bbb L_{\infty,\lambda}[{\frak K}]}(\bar
a,\emptyset,M_2)$
\sn
\item "{$(c)$}"  from $\ortp_{{\frak K}^*_\lambda}(\bar a,M_1,M_2)$ 
we can compute $\sftp_{\Bbb L_{\infty,\lambda}[{\frak K}]}(\bar a,M_1,M_2)$ and
$\ortp_{{\frak K}^*_\lambda}(\bar a,M_0,M_3)$.
\ermn
4) If $M_1 \le_{{\frak K}^*_\lambda} M_2$ and $\alpha < \kappa^* <
\lambda,\bold I_\ell \subseteq {}^\alpha(M_1),|\bold I_\ell| >
\kappa,\bold I_\ell$ is $(\Bbb L_{\infty,\theta}[{\frak K}],\kappa^*)$-
convergent in $M_1$ for $\ell=1,2$ and {\rm Av}$_{< \kappa}(\bold I_1,M_1) =
\text{\rm Av}_{< \kappa}(\bold I_1,M_1)$ \ub{then} $\bold I_\ell$ is
$(\Bbb L_{\infty,\kappa}[{\frak K}],\kappa^*)$-convergent in $M_\ell$ for
$\ell=1,2$ and {\rm Av}$_{< \kappa}(\bold I_1,M_\ell) = 
\text{\rm Av}_{<\kappa}(\bold I_1,M_2)$.
\endproclaim
\bigskip 

\demo{Proof}  1) Without loss of generality $M_0 = \text{
EM}_{\tau({\frak K})}(I_0,\Phi)$ and $I_0 \in
K^{\text{flin}}_\lambda$.  By \scite{734-11.14}(3) for $\ell=1,2$ there is a
pair $(I_\ell,f_\ell)$ such that $I_0 \le_{K^{\text{flin}}} I_\ell \in
K^{\text{flin}}_\lambda$ and $f_\ell$ is a $\le_{\frak K}$-embedding of
$M_\ell$ into $M'_\ell = \text{ EM}_{\tau({\frak K})}(I_\ell,\Phi)$
over $M_0$.  By renaming without loss of generality $f_\ell$ is the identity on
$M_\ell$ hence $M_\ell \le_{\frak K} M'_\ell$.  By \scite{734-11.5}(1) we
know that $M_\ell \prec_{\Bbb L_{\infty,\kappa}[{\frak K}]} M'_\ell$
hence $\sftp_{\Bbb L_{\infty,\kappa}[{\frak K}]}(\bar a_1,M_0,M'_1) =
\sftp_{\Bbb L_{\infty,\kappa}[{\frak K}]}(\bar a_1,M_0,M_1) =
\sftp_{\Bbb L_{\infty,\kappa}[{\frak K}]}(\bar a_2,M_0,M_2) =
\sftp_{\Bbb L_{\infty,\kappa}[{\frak K}]}(\bar a_2,M_0,M'_2)$.

By \scite{734-11.14}(1) we can find $(I_3,g_1,g_2,h)$ such that $I_0
\le_{K^{\text{flin}}} I_3 \in K^{\text{flin}}_\lambda,g_\ell$ is a
$\le_{\frak K}$-embedding of $M'_\ell$ into $M_4 := \text{
EM}_{\tau({\frak K})}(I_3,\Phi)$ over $M_0$ for $\ell=1,2$ and $h$ is
an automorphism of $M_4$ over $M_0$ mapping $g_1(\bar a_1)$ to
$g_2(\bar a_2)$.  By the definition of orbital types, this 
gives $\ortp_{{\frak K}^*_\lambda}
(\bar a_1,M_0,M_1) = \ortp_{{\frak K}^*_\lambda}(\bar a_2,M_0,M_2)$ as 
required.
\nl
2) This holds by \scite{734-11.5}(1) for $\theta \in (\text{LS}({\frak
K}),\lambda)$, hence by \scite{734-11.1.21}(1) also for $\theta = \lambda$ (the 
assumptions of \scite{734-11.1.21} hold
as clause (a) there holds by the case above $\theta < \lambda$ and clause
(b) there holds by \scite{734-11.11}(1)).
\nl
2A) Should be clear:
\mr
\item "{$(a)$}"  by part (2) this holds if $\bar a_1 = \bar a_2$ and
$M_1 \le_{\frak K} M_2$
\sn
\item "{$(b)$}"  trivially it holds if there is an isomorphism from
$M_1$ onto $M_2$ over $M_0$ mapping $\bar a_1$ to $\bar a_2$
\sn
\item "{$(c)$}"  by the definition of $\ortp$ we are done.
\ermn
2B) Should be clear by part (2).
\nl
3) \ub{Clause (a)}: 

By parts (1) + (2).
\mn
\ub{Clause (b)}:  By \scite{734-11.11}(1). 
\mn
\ub{Clause (c)}:  By part (2A) and the definition of $\ortp$.  
\nl
4) Easy, too.  \hfill$\square_{\scite{734-f.1}}$
\enddemo
\bigskip

\definition{\stag{734-f.2} Definition}   Assume $M_0 \le_{{\frak K}^*_\lambda}
M_1 \le_{{\frak K}^*_\lambda} M_2,\alpha < \lambda$ and 
$\bar a \in {}^\alpha(M_2)$ and $p = \ortp_{{\frak K}^*_\lambda}
(\bar a,M_1,M_2)$.  We say that $p$ does not fork over $M_0$ (for
${\frak K}^*_\lambda$) \ub{when}, letting $\theta_0 =  
|\alpha| + \text{ LS}({\frak K})$, $\theta_1 =
\beth_{1,1}(\beth_2(\theta_0)^+),\theta_2 = 2^{\theta_1},\theta_2 =
\beth_2(\theta_1)$ we have: 
\mr 
\item "{$(*)$}"  for some $N \le_{{\frak K}^*} M_0$ satisfying
$\|N\| \le \theta_2$ we have $\sftp_{{\Bbb L}_{\infty,\theta_1}[{\frak K}]}
(\bar a,M_1,M_2)$ does not split over $N$.
\ermn
We now would like to show that there is ${\frak s}_\lambda$ which
fits \chaptercite{600} and \chaptercite{705} and ${\frak K}_{{\frak
s}_\lambda} = {\frak K}^*_\lambda$.
\enddefinition
\bigskip

\demo{\stag{734-f.2.7} Observation}  Assume that $M_0 \le_{{\frak
K}^*_\lambda} M_1 \le_{{\frak K}^*_\lambda} M_2,\bar a \in
{}^\alpha(M_2),\alpha < \lambda,\lambda > \kappa_0 \ge |\alpha| + 
\text{\rm LS}({\frak K}),\kappa_1 = \beth_{1,1}(\beth_2(\kappa_0)^+)$ 
and $\kappa_2 = \beth_2(\kappa_1)$.  
\ub{Then} the following conditions are equivalent
\mr
\item "{$(a)$}"  \ortp$_{{\frak K}^*_\lambda}(\bar a,M_1,M_2)$ does not
fork over $M_0$
\sn
\item "{$(b)$}"  for some $(\kappa^+_1,\kappa_1)$-convergent
$\bold I \subseteq {}^\alpha(M_0)$ of cardinality $> \kappa_2$ we have
\nl
$\sftp_{\Bbb L_{\infty,\kappa_1}[{\frak K}]}(\bar a,M_1,M_2) 
= \text{\rm Av}_{<\kappa_1}(\bold I,M_1)$ hence this type does not
split over $\cup \bold I'$ for any $\bold I' \subseteq \bold I$ of
cardinality $> \kappa_1$
\sn
\item "{$(c)$}"  for every $N \le_{\frak K} M_0$ of cardinality 
$\le \kappa_2$, if $\sftp_{\Bbb L_{\infty,\kappa_1}[{\frak K}]}
(\bar a,M_0,M_2)$ does not split over $N$ then 
the type $\sftp_{\Bbb L_{\infty,\kappa_1}[{\frak K}]}
(\bar a,M_1,M_2)$ does not split over $N$.
\endroster
\enddemo
\bigskip

\remark{\stag{734-f.2.9} Remark}  1) See verification 
of axiom (E)(c) in the proof of Theorem \scite{734-f.20}.
\nl
2) Note that have we used $\beth_7(\kappa_1)^+$ instead of $\kappa_1$
in \scite{734-f.2}, \scite{734-f.2.7}, the difference would be small.
\nl
3) We could in clause (c) of \scite{734-f.2.7} use ``for some $N
\le_{\frak K} M_0$ of cardinality $< \kappa_1,
\sftp_{\Bbb L_{\infty,\kappa_1}[{\frak K}]}$ ..."
The proof is the same.
\nl
4) We can allow below $M_0 \le_{\frak K} M_1$ if $M_0 \in K_{\ge \kappa_2}$.  
\endremark
\bigskip

\demo{Proof}  \ub{$(a) \Rightarrow (b)$}

Let $\theta_0,\theta_1,\theta_2$ be as in Definition \scite{734-f.2}.
By Definition \scite{734-f.2} there is $N \le_{\frak K} M_0$ of
cardinality $\le \theta_2$ such that
\mr
\item "{$(*)_1$}"  the type 
$\sftp_{\Bbb L_{\infty,\theta_1}[{\frak K}]}(\bar a,M_1,M_2)$ does 
not split over $N$.
\ermn
By Claim \scite{734-11.10}(1) there is a
$(\kappa^+_1,\kappa_1)$-convergent set $\bold I \subseteq
{}^\alpha(M_0)$ of cardinality $\kappa^+_2$ (convergence in $M_0$, of
course) such that $\sftp_{\Bbb L_{\infty,\kappa_1}[{\frak K}]}(\bar
a,M_0,M_2) = \text{ Av}_{< \kappa_1}(\bold I,M_0)$.  So as $M_0
\prec_{\Bbb L_{\infty,\lambda}[{\frak K}]} M_1 \prec_{\Bbb
L_{\infty,\lambda}[{\frak K}]} M_2$, by Claim \scite{734-f.1}(2), clearly
$\bold I$ is $(\kappa^+_1,\kappa_1)$-convergent also in $M_1$ and in $M_2$
hence Av$_{< \kappa_1}(\bold I,M_1)$ is well defined.  Hence, by Claims
\scite{734-11.8}(2), \scite{734-11.7}(3) the type
Av$_{<\kappa_1}(\bold I,M_1)$ does not split over $\cup \bold I$ but
$\theta_2 \le \kappa_2$ and $\cup \bold I \subseteq \cup \bold I \cup N$
hence
\mr
\item "{$(*)_2$}"   Av$_{<\theta_1}(\bold I,M_1)$ does not split 
 over $\cup \bold I \cup N$.
\ermn
But also
\mr
\item "{$(*)_3$}"  $\sftp_{\Bbb L_{\infty,\theta_1}
[{\frak K}]}(\bar a,M_1,M_2)$ does not split over $N$ (by the choice
of $N$) hence over $\cup \bold I \cup N$.
\ermn
As $M_0 \prec_{\Bbb L_{\infty,\lambda}[{\frak K}]} M_1$ 
and $|\cup \bold I \cup N| <
\lambda$ and $\sftp_{\Bbb L_{\infty,\theta_1}[{\frak K}]} (\bar
a,M_0,M_2) = \text{ Av}_{< \theta_1}(\bold I,M_0)$ clearly, by $(*)_2
+ (*)_3$ we have $\sftp_{\Bbb L_{\infty,\theta_1}[{\frak K}]} (\bar
a,M_1,M_2) = \text{ Av}_{< \theta_1}(\bold I,M_1)$. \nl
 Now there is a pair
$(M'_2,\bar a')$ satisfying that $M_1 \le_{\frak K} M'_2 \in K^*_\lambda$
and $\bar a' \in {}^\alpha(M'_2)$ such that 
$\sftp_{\Bbb L_{\infty,\theta_1}[{\frak K}]} (\bar a',M_1,M'_2) =
\text{ Av}_{< \theta_1}(\bold I,M_1)$ hence by the previous sentence
$\sftp_{\Bbb L_{\infty,\theta_1}[{\frak K}]} (\bar a',M_1,M'_2) = 
\sftp_{\Bbb L_{\infty,\theta_1}[{\frak K}]} (\bar a,M_1,M_2)$.  Now
 by \scite{734-f.1}(1) and then \scite{734-f.1}(2A) it follows that
$\sftp_{\Bbb L_{\infty,\kappa_1}[{\frak K}]} (\bar a,M_1,M_0) =
\text{ Av}_{< \kappa_1}(\bold I,M_1)$ as required.
\mn
\ub{$(b) \Rightarrow (c)$}

Let $\bold I$ be as in clause (b), so $\bold I$ is
$(\kappa^+_1,\kappa_1)$-convergence in $M_0$ and is of cardinality $>
\kappa_1$.  We know that $M_0 \prec_{\Bbb L_{\infty,\lambda}[{\frak
K}]} M_1$, so by the previous sentence, $\bold I$ is
$(\kappa^+_1,\kappa_1)$-convergent in $M_1$.  To prove clause (c) assume
that $N \le_{\frak K} M_0$ is of cardinality $\kappa_2$ and
$\sftp_{\Bbb L_{\infty,\kappa_1}[{\frak K}]} (\bar a,M_0,M_2)$ does
not split over $N$.  Hence Av$_{<\kappa_1}(\bold I,M_0) = \sftp_{\Bbb
L_{\infty,\kappa_1}[{\frak K}]} (\bar a,M_0,M_2)$ does not split over
$N$.  Again as $M_0 \prec_{\Bbb L_{\infty,\lambda}[{\frak K}]} M_1$ we
can deduce that Av$_{<\kappa_1}(\bold I,M_1)$ does not split over $N$
but by the choice of $\bold I$ it is equal to $\sftp_{\Bbb
L_{\infty,\kappa_1}[{\frak K}]} (\bar a,M_1,M_2)$, so we are done.
\mn
\ub{$(c) \Rightarrow (a)$}

By Claim \scite{734-11.7.7} there is $B \subseteq M_0$ of cardinality $\le
\kappa_2$ such that
$\sftp_{\Bbb L_{\infty,\kappa_1}[{\frak K}]}(\bar a,M_0,M_2)$ does not
split over $B$.

As we can increase $B$ as long as we preserve ``of cardinality $\le
\kappa_2$", without loss of generality $B = |N|$ where $N \le_{\frak
K} M_0$. So the antecedent of clause (c) holds, but we are assuming clause
(c) so the conclusion of clause (c) holds, that is
$\sftp_{\Bbb L_{\infty,\kappa_1}[{\frak K}]}(\bar a,M_1,M_2)$ does not
split over $N$.

Also by \scite{734-11.10}(1) there is $\bold I_1 \subseteq {}^\alpha(M_0)$
of cardinality $\kappa^+_2$ which is $(\kappa^+_1,\kappa_1)$-convergent 
and Av$_{< \kappa_1}(\bold I_1,M_0)
= \sftp_{\Bbb L_{\infty,\kappa_1}[{\frak K}]}(\bar a,M_0,M_1)$.
Clearly $\kappa_1 \ge \theta_1$ hence $\kappa_2 = (\kappa_2)^{\theta_1}$.
Now as $K^*_\lambda$ is categorical clearly $M_0 \cong \text{
EM}_{\tau({\frak K})}(\lambda,\Phi)$ hence applying 
\scite{734-11.9}(4) we can find $\bold I_2 \subseteq \bold I_1$
of cardinality $\kappa^+_2$ which is
$(\theta^+_1,\theta_1)$-convergent.  As above $M_0 
\prec_{\Bbb L_{\infty,\kappa_1}[{\frak K}]} M_1$ so we 
deduce that $\bold I_2$ is
$(\theta^+_1,\theta_1)$-convergent and
$(\kappa^+_1,\kappa_1)$-convergent also in $M_1$.

As above we have 
$M_0 \prec_{\Bbb L_{\infty,\kappa_1}[{\frak K}]} M_1$ by
\scite{734-11.5}(1) hence 
Av$_{< \kappa_1}(\bold I_2,M_1)$ is well defined and does not split
over $N$ hence is equal to 
$\sftp_{\Bbb L_{\infty,\kappa_1}[{\frak K}]}(\bar a,M_1,M_2)$.  This
implies that Av$_{< \theta_1}(\bold I_2,M_1) = 
\sftp_{\Bbb L_{\infty,\theta_1}[{\frak K}]}(\bar a,M_1,M_2)$.

Now choose $\bold I_3 \subseteq \bold I_2 \subseteq M_0$ of
 cardinality $\theta_2$ and $N_3 \le_{\frak K} M_0$ of cardinality
 $\theta_2$ such that $\bold I_3 \subseteq {}^\alpha(N_3)$.  Now by
 \scite{734-11.8}(2) we know that 
$\sftp_{\Bbb L_{\infty,\theta_1}[{\frak K}]}(\bar a,M_1,M_2)$ does not 
 split over $\bold I_3$ hence it does not split over $N_3$, so $N_3$
 witnesses clause (a).  \hfill$\square_{\scite{734-f.2.7}}$
\enddemo
\bigskip

\definition{\stag{734-f.3} Definition}  We define a pre-frame 
${\frak s}_\lambda = ({\frak K}_{{\frak s}_\lambda},
\nonfork{}{}_{{\frak s}_\lambda},{\Cal S}^{\text{bs}}_{{\frak
s}_\lambda})$ as follows:  
\mr
\item "{$(a)$}"  ${\frak K}_{{\frak s}_\lambda} = {\frak K}^*_\lambda$
\sn
\item "{$(b)$}"  ${\Cal S}^{\text{bs}}_{{\frak s}_\lambda}$ is defined by
${\Cal S}^{\text{bs}}_{{\frak s},\lambda}(M) := 
\{\ortp_{{\frak K}^*_\lambda}(a,M,N):M \le_{{\frak K}^*_\lambda}
N,a \in N \backslash M\}$,
\sn
\item "{$(c)$}"  $\nonfork{}{}_{{\frak s}_\lambda} =
\{(M_0,M_1,a,M_3):M_0 \le_{{\frak K}^*_\lambda} M_1 
\le_{{\frak K}^*_\lambda} M_2$ and $\text{\ortp}_{{\frak K}^*_\lambda}
(a,M_1,M_3)$ does not fork over $M_0\}$, see Definition \scite{734-f.2}.
\endroster
\enddefinition
\bigskip

\remark{\stag{734-f.3.3} Remark}  1) Recall $\le_{{\frak s}_\lambda} =
\le_{\frak K} \restriction K_{{\frak s}_\lambda} = \le_{{\frak
K}^*_\lambda}$.
\nl
2) Concerning the proof of \scite{734-f.20} below we mention a variant which the 
reader may ignore.   This variant, from weaker assumptions  
gets weaker conclusions.
In detail, define the weak versions $(f)^-$ of $(f)$ of
\scite{734-f.0}(2); see Definition \scite{734-11.17} and Claim \scite{734-11.19A}(1)
\mr
\item "{$(f)^-$}"  if $\langle M_\alpha:\alpha \le \delta\rangle$ is
$\le_{\frak K}$-increasing continuous and $\alpha < \delta \Rightarrow
M_{2 \alpha+1} <^*_{{\frak K}^*_\lambda} M_{2 \alpha+2}$ (e.g. 
$M_{2 \alpha +2}$ is $\le_{{\frak K}^*_\lambda}$-universal over $M_{2 \alpha
+1}$) hence both are from $K^*_\lambda$ \ub{then} $M_\delta \in K^*_\lambda$.
\ermn
Assuming only \scite{734-f.0}(1) + (f)$^-$ we do not know whether ${\frak
K}^*_\lambda$ is a $\lambda$-a.e.c. but still 
$(K^*_\lambda,\le_{\frak K} \restriction K^*_\lambda,<^*_{{\frak
K}^*_\lambda})$, see Definition \scite{734-11.17}, is a so called semi
$\lambda$-a.e.c., see \chaptercite{E53}.

If clause $(f)$ from \scite{734-f.0}(2) holds 
(i.e., $K_{{\frak s}_\lambda}$ is closed under unions), we can
omit ``$<^*_{{\frak s}_\lambda}$".
\nl
3)  It will be less good but not a disaster if we have assumed below
$\lambda = \sup(\text{Cat}'_{\frak K} \cap \lambda)$.
\nl
4) It will be better to have ${\frak K}_{{\frak s}_\lambda} = K_\lambda$; of
courses, this follows from categoricity so by \S3 is not unreasonable for
conjecture \scite{734-0n.0}.
\nl
5) But we can ask only for $M \in K_{{\frak s}_\lambda}$ to be
universal in ${\frak K}_\lambda$,
\nl
6) We can ask that for every $\mu > \lambda$ large enough, for every $M \in
K_\mu$ for a club of $N \in K_\lambda$ satisfying 
$N \le_{\frak K} M$ we have $N \in K_{{\frak s}_\lambda}$.
\endremark
\bigskip

\proclaim{\stag{734-f.20} Theorem$^*$}   (Assume \scite{734-f.0}(2),(g) hence (f)).
\nl
${\frak s}_\lambda$ is a good $\lambda$-frame categorical in $\lambda$
and is full.
\endproclaim
\bigskip

\demo{Proof}  We check the clauses in the definition
\marginbf{!!}{\cprefix{600}.\scite{600-1.1}}.
\mn
\ub{Clause $(A)$}:

By observation \scite{734-f.0.7}(3), [in the weak version using
(f)$^-$ from \scite{734-f.3.3}(1)].
\mn
\ub{Clause $(B)$}:

Categoricity holds by \scite{734-11.4A} (or \scite{734-f.0.7}(1))
and this implies ``there is a superlimit
model", the non-maximality by $\le_{{\frak K}^*_\lambda}$ holds by the 
choice of $\Phi$.
\mn
\ub{Clause $(C)$}:

Observation \scite{734-f.0.7}(2) guarantee amalgamation, categoricity (of
${\frak K}^*_\lambda$ by \scite{734-f.0.7}(1))
implies the JEP and ``no-maximal model" holds by clause (B).
\mn
\ub{Clause $(D)(a),(b)$}:

Obvious by the definition.  
\mr
\item "{$(D)$}"  $(c) \quad$ (density).
\ermn
Assume $M <_{{\frak K}^*_\lambda} N$, then there are $a \in N
\backslash M$ and for any such $a$ the type $\ortp_{{\frak
K}^*_\lambda}(a,M,N)$ belongs to ${\Cal S}^{\text{bs}}_{{\frak
s}_\lambda}(M)$.  In fact
\mr
\item "{$\circledast$}"  ${\frak s}_\lambda$ is type-full
\sn
\item "{$(D)$}"  $(d) \quad$ (bs-stability).  
\ermn
The demand means $M \in K^*_\lambda \Rightarrow 
|{\Cal S}^1_{{\frak K}^*_\lambda}(M)| \le \lambda$.
 
This holds by \scite{734-11.16A}(2) (and amalgamation).
\enddemo
\bn
$(E)(a),(b)$.  By the definition.
\mn
$(E)(c)$ (local character)

This says that if $\langle M_i:i \le \delta +1 \rangle$ is $\le_{{\frak
s}_\lambda}$-increasing continuous and $p = \ortp_{{\frak
s}_\lambda}(a,M_\delta,M_{\delta +1}) \in 
{\Cal S}^{\text{bs}}_{{\frak s}_\lambda}(M_\delta)$ \ub{then} for some 
$i < \delta$ the type $p$ does not fork over $M_i$ (for ${\frak s}_\lambda$). 

\relax From now on (in the proof of \scite{734-f.20}) we 
use \scite{734-f.2.7} freely and let (noting cf$(\delta) <
\lambda$ as $\lambda$ is singular)
\mr
\item "{$\odot$}"  $\kappa_0 = \text{ LS}({\frak K}) + \text{ cf}(\delta),
\kappa_1 = \beth_{1,1}(\beth_2(\kappa_0))^+,\kappa_2 = \beth_2(\kappa_1)$.
\ermn
Now by \scite{734-f.2.7} there is a $(\kappa^+_1,\kappa_1)$-convergent 
$\bold I \subseteq M_\delta$ with Av$_{<\kappa_1}(\bold I,
M_\delta) = \sftp_{{\Bbb L}_{\infty,\kappa_1}[{\frak
K}]}(a,M_\delta,M_{\delta +1})$ such that $\bold I$ is
of cardinality $> \kappa_2$.  
For some $i(*) < \delta,|\bold I \cap M_{i(*)}| > \kappa_2$, so without loss
of generality $\bold I \subseteq M_{i(*)}$, so by \scite{734-f.2.7} we are
done.
\mn
$(E)(d)$ Transitivity of non-forking

We are given $M_0 \le_{{\frak s}_\lambda} M_1 \le_{{\frak s}_\lambda}
M_2 \le_{{\frak K}_{\frak s}} M_3$ and $a \in M_3$ such that
\ortp$_{{\frak s}_\lambda}(a,M_{\ell +1},M_3)$ does not fork over $M_\ell$
for $\ell=0,1$.  So for $\ell=0,1$ 
there is $\bold I_\ell \subseteq M_\ell$ which is
$(\kappa^+_1,\kappa_1)$-convergent in $M_{\ell+1}$ of cardinality 
$\kappa^+_2$ such that Av$_{< \kappa_1}(\bold I_\ell,M_{\ell +1}) = 
\sftp_{\Bbb L_{\infty,\kappa_1}[{\frak K}]}(a,M_{\ell +1},M_3)$.  
As Av$_{< \kappa_1}(\bold I_0,M_1) = \text{ Av}_{< \kappa_1}
(\bold I_1,M_1)$ (being both realized by $a$) because $M_1 
\prec_{\Bbb L_{\infty,\lambda}[{\frak K}]} M_2$ by \scite{734-f.1}(4)
 clearly we
have Av$_{< \kappa_1}(\bold I_0,M_2) = \text{ Av}_{< \kappa_1}
(\bold I_1,M_2) = \sftp_{\Bbb L_{\infty,\kappa_1}
[{\frak K}]}(a,M_2,M_3)$ all well defined.   
So $\bold I_0$ witness by \scite{734-f.2.7} 
that $\sftp_{\Bbb L_{\infty,\kappa_1}[{\frak K}]}
(a,M_2,M_3)$ does not fork over $M_0$, which means that $\ortp_{{\frak
K}^*_\lambda}(a,M_2,M_3)$ does not fork over $M_0$ as required.  
\mn
$(E)(e)$ Uniqueness.

Recalling \scite{734-f.1}(1), the proof is similar to $(E)(d)$; 
the two witnesses are now in $M_0$.
\mn
$(E)(f)$ Symmetry

Toward contradiction, recalling \marginbf{!!}{\cprefix{600}.\scite{600-1.16E}}
 assume $M_0 \le_{{\frak K}^*_\lambda} 
M_1 \le_{{\frak K}^*_\lambda} M_2 \le_{{\frak K}^*_\lambda} M_3$ and
$a_\ell \in M_{\ell +1} \backslash M_\ell$ for $\ell = 0,1,2$ 
are such that $p_\ell = \ortp_{{\frak K}^*_\lambda}
(a_\ell,M_\ell,M_{\ell +1})$ does not fork over $M_0$ for $\ell=0,1,2$ and
$\ortp_{{\frak K}^*_\lambda}(a_0,M_0,M_1) = \ortp_{{\frak
K}^*_\lambda}(a_2,M_0,M_3)$ but $\ortp_{{\frak K}^*_\lambda}(\langle
a_0,a_1\rangle,M_0,M_3) \ne \ortp_{{\frak K}^*_\lambda}(\langle
a_2,a_1\rangle,M_0,M_3)$. 

By \scite{734-f.2.7} we can deal with $p_\ell = 
\sftp_{\Bbb L_{\infty,\kappa_1}[{\frak K}]}
(a_\ell,M_\ell,M_{\ell +1})$ for $\ell=0,1,2$.  For each $\ell \le 2$, we can
find convergent $\bold I_\ell = \{a^\ell_\alpha:\alpha < \kappa^+_2\} 
\subseteq M_0$ which is $(\kappa^+_1,\kappa_1)$-convergent such 
that Av$_{< \kappa_1}(\bold I_\ell,M_\ell) = p_\ell$.  

So as $M_0 \prec_{\Bbb L_{\infty,\kappa_1}[{\frak K}]} M_k$ we deduce
the set $\bold I_\ell$ is
$(\kappa^+_1,\kappa_1)$-convergent in $M_k$ for $\ell,k = 0,1,2$, also
Av$_{< \kappa_1}(\bold I_0,M_0) = \text{ Av}_{< \kappa_1}
(\bold I_2,M_0)$ hence Av$_{< \kappa_1}(\bold I_0,M_2) = 
\text{ Av}_{< \kappa_1}(\bold I_2,M_2)$ so without loss of generality 
$\bold I_0 = \bold I_2$.
\nl
Now use the non-order property to get symmetry.
\mn
$(E)(g)$ Existence

So assume $M \le_{{\frak s}_\lambda} N$ and $p \in {\Cal
S}^{\text{bs}}_{{\frak s}_\lambda}(M)$.  So we can find a pair $(M',a)$
such that $M \le_{{\frak s}_\lambda} M',a \in M_1$ and $p =
\ortp_{{\frak s}_\lambda}(a,M,M')$.  By \scite{734-11.10}(1) there is a
$(\kappa^+_1,\kappa_1)$- convergent $\bold I \subseteq M$ of
cardinality $\kappa^+_2$ such that
Av$_{< \kappa_1}(M,\bold I) = \sftp_{\Bbb L_{\infty,\kappa_1}[{\frak
K}]}(a,M,M')$.  By \scite{734-11.10}(3) + \scite{734-f.2.7} there is a pair
$(N',a')$ such that $N \le_{{\frak s}_\lambda} N',a' \in N'$ and
$\sftp_{\Bbb L_{\infty,\kappa_1}}(a',N,N') = \text{ Av}_{< \kappa_1}
(\bold I,N)$.  So by \scite{734-f.2.7} the type $\ortp_{{\frak
s}_\lambda}(a',N,N')$ easily $\in {\Cal S}^{\text{bs}}_{{\frak
s}_\lambda}(N)$, does not fork over $N$ and extend $p$, as required.
\mn
$(E)(h)$ Continuity 

Follow by \marginbf{!!}{\cprefix{600}.\scite{600-1.16A}}.  Alternatively assume 
$\langle M_i:i \le \delta +1 \rangle$ is $\le_{{\frak
s}_\lambda}$-increasing continuous, and $a \in M_{\delta +1}
\backslash M_\delta$ and $\ortp_{{\frak s}_\lambda}(a,M_i,M_{\delta +1})$
does not fork over $M_0$ for $i < \delta$.  So there is a convergent
$\bold I_i \subseteq M_0$ such that $i < \delta \Rightarrow 
\sftp_{\Bbb L_{\infty,\kappa}[{\frak K}]}(a,M_i,M_{\delta +1}) = \text{
Av}_\kappa(\bold I,M_i)$. 
\nl
As above, \wilog \, $\bold I_i = \bold I_0$.  We can
find a convergent $\bold I \subseteq M_\delta$ of cardinality $>
\text{ cf}(\delta) + \kappa$ (recall cf$(\delta) < \lambda$!) such
that $\sftp_{\Bbb L_{\infty,\kappa}[{\frak K}]}(a,M_0,M_{\delta +1}) =
\text{ Av}_\kappa(\bold I,M_\delta)$.  So for some $i(*) <
\delta,|\bold I \cap M_{i(*)}| > \kappa$ so \wilog \, (by equivalence)
$\bold I \subseteq M_{i(*)}$.  We finish as in $(E)(f)$.
\bn
\ub{Axiom (E)(i)}:

Follows by \marginbf{!!}{\cprefix{600}.\scite{600-1.15}}.  \hfill$\square_{\scite{734-f.20}}$
\bn
\margintag{734-f.27}\ub{\stag{734-f.27} Exercise}:  Replace above Av$_{<\kappa_1}(\bold I,M)$ by
$\cup\{\text{Av}_{\beth_\zeta(\kappa_0)}(\bold I,M):\zeta <
(2^{\kappa_0})^+\}$. 
\goodbreak

\head {\S5 Homogeneous enough linear orders} \endhead  \resetall \sectno=5
 \spuriousreset
\bigskip

\proclaim{\stag{734-am3.2.3} Claim}   Assume $\mu^+ = \theta_1 = 
\text{\rm cf}(\theta_1) < \theta_2 = \text{\rm cf}(\theta_2) <\lambda$.
\nl
1) \ub{Then} there is a linear order $I$ of cardinality $\lambda$
such that:   the following equivalence
relation ${\Cal E} = {\Cal E}^{\text{aut}}_{I,\mu}$ 
on ${}^\mu I$ has $\le 2^\mu$ equivalence classes, where

$\eta_1 {\Cal E} \eta_2$ \ub{iff} there is an automorphism of $I$ mapping
$\eta_1$ to $\eta_2$.
\nl
2) Moreover if $I' \subseteq I$ has cardinality $< \theta_2$ and $n < \omega$
\ub{then} the following equivalence relation ${\Cal E}$ 
on ${}^n I$ has $\le \mu + |I'|$ equivalence classes:
\nl
$\bar s {\Cal E} \bar t$ iff there is an automorphism $h$ of $I$
over $I'$ mapping $\bar s$ to $\bar t$.
\nl
3) Moreover, there is $\Psi$ proper for $K^{\text{lin}}_{\tau^*_2}$ 
(i.e. $\Psi \in \Upsilon^{\text{lin}}_{\aleph_0}[2]$, see Definition
\scite{734-X1.2}(5),\scite{734-0n.2}(9)) with $\tau(\Psi)$ countable such that $I =
\, \text{\rm EM}_{\{<\}}(I^{\text{lin}}_{\theta_2,\lambda \times
\theta_2},\Phi)$ where $I^{\text{lin}}_{\theta_2,\zeta} = (\zeta,
<,P_0,P_1),P_\ell = \{\alpha < \zeta$:({\rm cf}$(\alpha) < \theta_2) 
\equiv (\ell = 0)\}$.
\nl
4) If $I^*_0 \subseteq I$ has cardinality $< \theta_2$ \ub{then} for some
$I^*_1 \subseteq I$ of cardinality $\le \mu^+ + |I^*_0|$ 
for every $J \subseteq I$
of cardinality $\le \mu$ there is an automorphism of $I$ over
$I^*_0$ mapping $J$ into $I^*_1$.
\nl
5) If $I^*_1,I^*_2 \subseteq I^{\text{lin}}_{\mu,\lambda \times \mu^+}$ has
cardinality $\le \mu$ and $h$ is an isomorphism from $I^*_1$
onto $I^*_2$ \ub{then} there is an automorphism $\hat h$ of the linear
order $I = \,\text{\rm EM}_{\{<\}}(I^{\text{lin}}_{\theta,\lambda},\Psi)$
extending the natural isomorphism $\check h$ from
{\rm EM}$_{\{<\}}(I^*_1,\Psi)$ onto {\rm EM}$_{\{<\}}(I^*_2,\Psi)$. 
\endproclaim
\bigskip

\remark{Remark}  1) Of course, if $\lambda = \lambda^{< \theta_2}$ and
$I$ is a dense linear order of cardinality $\lambda$ which is
$\theta$-strongly saturated (hence $\theta$-homogeneous) then the demand in
\scite{734-am3.2.3}(1) is satisfied (and in part (2) of \scite{734-am3.2.3} 
the number (of ${\Cal E}$ equivalence classes) is $\le 2^\chi$
for every $\chi \in [\aleph_0,\theta_2))$.  Also if $\lambda = \dsize
\sum_{i < \delta} \lambda _i,\delta < \theta_2$ and $i < \delta
\Rightarrow \lambda^{<\theta_2}_i = \lambda$ we have such order.
\nl 
2) Laver \cite[\S2]{Lv71} deals with related linear orders but for his
aims $I_1,I_2$ are equivalent if each is embeddable into the other; see more in
\cite[AP,\S2]{Sh:e}.  For a cardinal $\partial$ and linear order $I$
let $\Theta_{I,\partial} = \{\text{cf}(J)$: for some $<_I$-decreasing
sequence $\langle t_i:i < \partial\rangle$ we have $J = I \restriction
\{t \in I:t <_I t_i$ for every $i < \partial\}\}$.  So if $\partial \le
\mu$ then $({}^\mu I)/E^{\text{aut}}_{I,\mu}$ has $\ge
|\Theta_{I,\partial}|$.  So we have to be careful to make
$\Theta_{I,\partial}$ small.  We choose a very concrete construction
which leads quickly to defining $I$ and the checking is straight so we
thought it would be easy but a posteriori the checking is lengthy;
\cite[AP,\S2]{Sh:e} is an anti-thetical approach.
\nl
3) We can replace $\theta_1 = \mu^+$ by $\theta_1 = \text{\rm cf}(\theta_1) >
\aleph_0$ and ``of cardinality $\le \mu$" by ``of cardinality $<
\theta_1$". 
\nl
4) In \scite{734-eq.3.3}(1), \scite{734-eq.4}(2) we use parts (1),(1)+(4)
respectively.   Also we use \scite{734-am3.2.3} in the proof of 
\scite{734-am2.18}.
\nl
5) The case $2^\mu \ge \lambda$ in \scite{734-am3.2.3}(1) says nothing, 
in fact if $2^\mu \ge \lambda$ then $2^\mu = \lambda^\mu = ({}^\mu
M)/{\Cal E}^{\text{aut}}_{I,\mu}$ for any model $M$ of cardinality
$\le 2^\mu$ but $\ge 2$,  for any vocabulary $\tau_M$.
\nl
6) Claim \scite{734-am3.2.3}(1),(2) holds also if we replace $\mu$ by
$\chi \in [\mu,\theta_2)$. 
\endremark
\bigskip

\demo{Proof}   1) Fix an ordinal 
$\zeta,\lambda \le \zeta < \lambda^+$ such that
cf$(\zeta) = \theta_2$, e.g., $\zeta = \lambda \times \theta_2$
(almost always cf$(\zeta) \ge \theta_2$ suffice).

Let $I_1$ be the following linear order, its set of
elements is $\{(\ell,\alpha):\ell \in \{-2,-1,1,2\}$, $\alpha 
< \zeta + \omega\}$ ordered by 
$(\ell_1,\alpha_1) <_{I_1} (\ell_2,\alpha_2)$ iff $\ell_1 <
\ell_2$ or $\ell_1 = \ell_2 \in \{-1,2\} \wedge \alpha_1 < \alpha_2$
or $\ell_1 = \ell_2 \in \{-2,1\} \wedge \alpha_1 > \alpha_2$.

For $t \in I_1$ let $t = (\ell^t,\alpha^t)$.

Let $I^*_2$ be the set $\{\eta:\eta$ is a finite sequence of members
of $I_1\}$ ordered by $\eta_1 <_{I_2} \eta_2$ \ub{iff} $(\exists n)(n <
\ell g(\eta_1) \wedge n < \ell g(\eta_2) \wedge \eta_1 \restriction n
= \eta_1 \restriction n \and \eta_1(n) <_{I_1} \eta_2(n))$ or $\eta_1
\triangleleft \eta_2 \wedge \ell^{\eta_2(\ell g(\eta_1))} \in \{1,2\}$ or
$\eta_2 \triangleleft \eta_1 \wedge \ell^{\eta_1(\ell g(\eta_2))} \in 
\{-2,-1\}$.

Let $I_2$ be $I^*_2$ restricted to the set of $\eta \in I^*_2$ satisfying
$\circledast$ where
\mr
\item "{$\circledast$}"  for no $n < \omega$ do we have:
{\roster
\itemitem{ $(a)$ }  $\ell g(\eta) > n+1$
\sn
\itemitem{ $(b)$ }  $\alpha^{\eta(n)}$ is a limit ordinal of
cofinality $\ge \theta_1$
\sn
\itemitem{ $(c)$ }  $\alpha^{\eta(n+1)} \ge \zeta$
\sn
\itemitem{ $(d)$ }  $\ell^{\eta(n)} \in \{-1,2\},\ell^{\eta(n+1)} =
-2$ \ub{or} $\ell^{\eta(n)} \in \{-2,1\},\ell^{\eta(n+1)} = 2$.
\endroster}
\ermn
Let $M_0$ be the following ordered field:
\mr
\item "{$(*)_1$}"  $(a) \quad M_0$ as a field, is $\Bbb Q(a_t:t \in
I_2)$, the field of rational functions with
\nl

\hskip25pt $\{a_t:t \in I_2\}$ algebraically independent
\sn
\item "{${{}}$}"  $(b) \quad$ the order of $M_0$ is determined by
{\roster
\itemitem{ ${{}}$ }  $(\alpha) \quad$ if $t \in I_2,n <\omega$ then $M_0
\models n < a_t$
\sn
\itemitem{ ${{}}$ }   $(\beta) \quad$ if $s <_{I_2} t$ and $n<\omega$
then $M_0 \models ``(a_s)^n < a_t"$.
\endroster}
\item "{${{}}$}"  $(c) \quad$ let $M$ be the real \footnote{in fact,
we could just use $M_0$} (algebraic) closure of 
$M_0$ (i.e. the elements 
\nl

\hskip25pt algebraic over $M_0$ in the closure by
adding elements realizing any
\nl

\hskip25pt  Dedekind cut of $M_0$).
\ermn
Now we shall prove that $I$, which is $M$ as a linear order, is as requested.
\mr
\item "{$\boxtimes_1$}"  each of $I_1,I^*_2$ and
$I_2$ is anti-isomorphic to itself.
\ermn
[Why?  Let $g:I_1 \rightarrow I_1$ be $g(t) = (-\ell^t,\alpha^t)$, clearly it
is an anti-isomorphism of $I_1$.  Let $\hat g:I^*_2 \rightarrow
I^*_2$ be defined by 
$\hat g(\eta) = \langle g(\eta(m)):m < \ell g(\eta) \rangle$, it is
an anti-isomorphism of $I^*_2$.  Lastly $\hat g$ maps $I_2$ onto
itself, in particular by the character of clause (d) of $\circledast$,
i.e. the two cases are interchanged by $\hat g$]
\mr
\item "{$\boxtimes_2$}"   $(a) \quad I_1,I^*_2,I_2$ have cofinality $\aleph_0$.
\sn
\item "{${{}}$}"  $(b) \quad$ if $t \in I_2$ then $I_{2,< t} :=
I_2 \restriction \{s:s <_{I_2} t\}$ has cofinality $\aleph_0$.
\ermn
[Why?  For clause (a), $\{(2,\lambda + n):n <
\omega\}$ is a cofinal subset of $I_1$ of order type $\omega$
 and
$\{<t>:t \in I_1\}$ is a cofinal subset of $I^*_2$ and of $I_2$ of
order type being the same as $I_1$.  For
clause (b) for $\eta \in I_2$ the set $\{\eta \char 94 \langle (-1,\lambda +
n)\rangle:n <\omega\}$ is a cofinal subset of $I_{2,< \eta}$
of order type $\omega$ by $\boxdot$ below.]

Now 
\mr
\item "{$\boxdot$}"   if $\eta$ satisfies $\circledast$ and
$\ell \in \{1,-1\}$ then also $\eta \char 94 \langle
(\ell,\alpha)\rangle$ satisfies $\circledast$ for any 
$\alpha < \lambda + \omega$. 
\ermn
[Why?  By clause (d) of $\circledast$ as the only value of $n$ there which
is not obvious is $n = \ell g(\eta)-1$, but to be problematic we
should have $\ell^{(\eta \char 94 <(\ell,\alpha)>)(n+1)} \in\{-2,2\}$ whereas
$\ell = -1$.]
\mr
\item "{$\boxtimes_3$}"   if $\partial = \text{ cf}(\partial)$ so $\partial$
is $0,1$ or an infinite regular cardinal and
$\bar \eta = \langle \eta_i:i < \partial \rangle$ is a $<_{I_2}$-decreasing
sequence and we let $J_{\bar\eta} = 
\{s \in I_2:s <_{I_2} \eta_i$ for every $i < \partial\}$ \ub{then} (clearly
exactly one of the following clauses applies)
{\roster
\itemitem{ $(a)$ }  if $J_{\bar \eta} = \emptyset$ then $\partial =
\aleph_0$
\sn
\itemitem{ $(b)$ }  if cf$(J_{\bar\eta}) =1$ then $\partial = \aleph_0$
\sn
\itemitem{ $(c)$ }  if cf$(J_{\bar\eta}) = \aleph_0$ then $\partial < \theta_1$
\sn
\itemitem{ $(d)$ }  if $\aleph_1 \le \text{ cf}(J_{\bar\eta}) < \theta_1$ then
$\partial = \aleph_0$ and for some $\ell \in \{-1,2\},\nu \in I_2$ and ordinal
$\delta < \zeta$ of cofinality cf$(J_{\bar \eta})$ the
set $\langle \nu \char 94 \langle (\ell,\alpha)\rangle:\alpha <\delta
\rangle$ is an unbounded subset of $J_{\bar\eta}$
\sn
\itemitem{ $(e)$ }  if $\theta_1 \le \text{\rm cf}(J_{\bar\eta})$ 
then $\partial \ge \theta_1$ and moreover $\partial = \theta_2 \vee 
\text{ cf}(J_{\bar \eta}) = \theta_2$.
\endroster}
\ermn
[Why does $\boxtimes_3$ hold?  The proof is split into cases and
finishing a case we can then assume it does not occur.

Clearly we can replace $\bar \eta$ by $\langle \eta_i:i \in u \rangle$
for any unbounded subset $u$ of $\partial$ and by $\langle \nu_i:i \in
u\rangle$ if $\eta_{\zeta_{2i+1}} \le_{I_2} \nu_i \le_{I_2} \eta_{\zeta_{2i}}
$ and $\langle \zeta_i:i < \partial\rangle$ an
increasing sequence of ordinals $< \partial$.
We shall use this freely.
\enddemo
\bn
\ub{Case 0}:  $\partial = 0$ or $\partial = 1$.

By $\boxtimes_2$ clearly clause (c) of $\boxtimes_3$ holds.
\bn
\ub{Case 1}:  $\partial = \aleph_0$ and there is $\nu \in
{}^\omega(I_1)$ such that $(\forall n < \omega)(\exists i <
\partial)(\eta_i \restriction n \triangleleft \nu)$.

Let $n_i = \ell g(\eta_i \cap \nu)$, it is impossible that
$\{i:n_i=k\}$ is infinite for some $k$, so  \wilog \, 
$\langle n_i:i < \omega \rangle$ is an increasing sequence and $n_0 >0$.

For every $i < \omega$ we have $\nu \restriction (n_i +1)
\trianglelefteq \eta_{i+1}$ and $\eta_{i+1} <_{I_2} \eta_i$, so by the
definition of $<_{I_2}$ also $\nu \restriction (n_i +1) <_{I_2}
\eta_i$, and we choose $\beta_{n_i} < \zeta + \omega$ so that
$(-2,\beta_{n_i}) <_{I_1} \nu(n_i)$ hence letting $\rho_i = \nu
\restriction n_i \char 94 \langle (-2,\beta_{n_i})\rangle$ we
have $\rho_i \in I_2$.  This can be done, e.g. because we can choose
$\beta_{n_i}$ such that $\beta_{n_i} = \alpha^{\nu(n_i)}+1$ if
$\ell^{\nu(n_i)} = -2$ and $\beta_{n_i} = 0$ otherwise.

For every $i,j < \omega$ we
have $\rho_i <_{I_2} \rho_{i+1} <_{I_2} \eta_{i+1} <_{I_2} \eta_i$, so
if $i \le j$ then $\rho_i <_{I_2} \rho_j <_{I_2} < \eta_j$, and if
$i>j$ then $\rho_i <_{I_2} \eta_i <_{I_2} \eta_j$, so $\rho_i \in
J_{\bar\eta}$.

Now $\langle \rho_i:i < \omega \rangle$ is $<_{I_2}$-increasing also it
is cofinal in $J_{\bar\eta}$,
for if $\rho \in J_{\bar\eta}$ let $n = \ell g(\rho \cap \nu)$, so for
$i < \omega$ such that $n_i \le n < n_{i+1}$ we have $\rho <_{I_2}
\eta_{i+1}$ so $\rho(n) <_{I_1} \eta_{i+1}(n) = \rho_{i+1}(n)$ and as $\rho
\restriction n = \nu \restriction n = \rho_{i+1} \restriction n$ we
have $\rho <_{I_2} \rho_{i+1}$.

As $\langle \rho_i:i < \omega\rangle$ is of order type $\omega$
clearly cf$(J_{\bar \eta}) = \aleph_0 = \partial$ hence clause (c) of 
$\boxtimes_3$ applies, and we are done. 

So from now on assume that case 1 fails.

As $\ell g(\eta_i) <
\omega$ and as not Case 1 \wilog \, for some $n$, we have $i < \partial
\Rightarrow \ell g(\eta_i)=n$.  Similarly 
\wilog\, for some $m$ and $\nu \in I_2$ we have $i
< \partial \Rightarrow \eta_i \restriction m = \nu$ and $\langle
\eta_i(m):i < \partial\rangle$ with no repetitions so $m<n$.  Without loss
of generality $i < \partial \Rightarrow \ell^{\eta_i(m)} =\ell^*$ and
so $\langle \alpha^{\eta_i(m)}:i < \partial \rangle$ is with no
repetitions; and \wilog \, is monotonic hence, as $\partial \ge \aleph_0$ 
is an increasing sequence of ordinals.  As $\bar \eta$ is
$<_{I_2}$-decreasing necessarily $\ell^* \in \{-2,1\}$ and let $\delta
= \cup\{\alpha^{\eta_i(m)}:i < \partial\}$, so clearly cf$(\delta) =
\partial$ and $\delta$ is a limit ordinal $\le \zeta + \omega$.  
Now those $\ell^*,\delta$ will be used till the end of the
proof of $\boxtimes_3$.  So for the rest of the proof we are assuming
\mr
\item "{$\odot$}"  $(a) \quad i < \partial \Rightarrow \eta_i
\restriction m = \nu$
\sn
\item "{${{}}$}"  $(b) \quad \langle \eta_i(m):i < \partial\rangle$ is
(strictly) increasing with limit $\delta$
\sn
\item "{${{}}$}"  $(c) \quad \ell^{\eta_i(m)} = \ell^* \in \{-2,1\}$
\sn
\item "{${{}}$}"  $(d) \quad$ cf$(\delta) = \partial,
\delta \le \zeta + \omega$.
\ermn
Also note by $\circledast$ that 
$\nu \char 94 \,\langle(\ell^*,\delta)\rangle \notin I_2
\Rightarrow \delta \in \{\zeta + \omega,\zeta\}$ and
if $\delta = \zeta \wedge \nu \char 94 \langle
(\ell^*,\delta)\rangle \notin I_2$ then $\ell g(\nu) > 0$ and 
the ordinal $\alpha^{\nu(\ell g(\nu)-1)}$ is limit of cofinality 
$\ge \theta_1$
(and more).
\bn
\ub{Case 2}:  $J_{\bar \eta}=\emptyset$.

Clearly $m=0 \wedge\ell^* = -2 \wedge \delta = \zeta + \omega$ hence 
$\partial = \aleph_0$ so clause (a) of $\boxtimes_3$ holds.  
\bn
\ub{Case 3}:  $\ell^*=1$ and 
$\nu \char 94 \langle (\ell^*,\delta)\rangle \notin
I_2$.

As $\ell^*=1$ clearly we cannot have $\delta = \zeta$ by clause
(d) of $\circledast$ so $\delta = \zeta + \omega$ and
recalling $\partial = \text{ cf}(\delta)$ we have
$\partial = \aleph_0$.  Now clearly $J_{\bar \eta}$ has a last element,
$\nu$, so case (b) of $\boxtimes_3$ applies.
\bn
\ub{Case 4}:  $\ell^* = -2,\partial = \aleph_0$ and $\nu \char 94
\langle (\ell^*,\delta)\rangle \notin I_2$.  

Again $\delta = \zeta + \omega$ as $\aleph_0 = \partial = \text{\rm
cf}(\delta)$ and cf$(\zeta) = \theta_2 > \mu \ge 
\aleph_0$ making $\delta = \zeta$ 
impossible; now $\ell g(\nu) > 0$ (as we have discarded the
case $J_{\bar \eta} = \emptyset$, i.e. Case 2); and let $k= \ell
g(\nu)-1$.  Now we prove case 4 by splitting to several subcases.
\bn
\ub{Subcase 4A}:  $\ell^{\nu(k)} \in \{-2,1\}$.

Let $\nu_1 = (\nu \restriction k)
\char 94 \langle (\ell^{\nu(k)},\alpha^{\nu(k)}+1)\rangle$, note that
$\nu_1 \in I_2$ as $\nu \in I_2 \wedge (\alpha^{\nu(k)} < \zeta 
\equiv \alpha^{\nu(k)} + 1 < \zeta)$ 
and (as $\ell^{\nu(k)} \in \{-2,1\}$) clearly
$\{\rho:\nu_1 \trianglelefteq \rho \in I_2\}$ is a cofinal subset of
$J_{\bar \eta}$ even an end segment.  Now for $n <\omega$ we have
$\nu_1 \char 94 \langle(2,\zeta + n)\rangle \in I^*_2$ and  
it satisfies $\circledast$.  (Why?  As $\nu_1
\in I_2$, only $n=k$ may be problematic, but $\alpha^{\nu(k)} +1 =
\alpha^{\nu_1(k)}$ here stands for $\alpha^{\eta(n)}$ there hence 
clause (b) of $\circledast$ does not apply), so by
the definition of $I_2$, clearly
$\{\nu_1 \char 94 \langle (2,\zeta + n)\rangle:n <\omega\}$
is $\subseteq I_2$ and is a cofinal subset of
$J_{\bar\eta}$ so $\partial = \aleph_0 = \text{ cf}(J_{\bar \eta})$ and
clause (c) of $\boxtimes_3$ holds.
\bn
\ub{Subcase 4B}:  $\ell^{\nu(k)} \in \{-1,2\}$ and 
$\alpha^{\nu(k)}$ is a successor ordinal.

Let $\nu_1 = (\nu \restriction k) \char 94 \langle(\ell^{\nu(k)},
\alpha^{\nu(k)}-1)\rangle$, of course $\nu_1 \in I^*_2$ and as $\nu
\in I_2$ clearly $\nu_1 \in I_2$ so the set $\{\rho:\nu_1
\trianglelefteq \rho \in I_2\}$ is an end segment of $J_{\bar \eta}$ and
has cofinality $\aleph_0$ because $n < \omega \Rightarrow 
\nu_1 \char 94 \langle(2,\zeta +n)\rangle \in I_2$.  
(Why?  It $\in I^*_2$ and as $\nu_1 \in
I_2$ checking $\circledast$ only 
$n=k$ may be problematic, but $(\ell^{\nu(k)},2)$ here stand
for $(\ell^{\eta(n)},\ell^{\eta(n+1)})$ there but presently
$\ell^{\nu(k)} \in \{-1,2\}$ contradicting clause (d) of
$\circledast$).  So clause (c) of $\boxtimes_3$.  
\bn
\ub{Subcase 4C}:  $\ell^{\nu(k)} \in \{-1,2\}$ and
$\alpha^{\nu(k)}=0$.

Then let $\nu_1 = (\nu \restriction k) \char 94 \langle(\ell^{\nu(k)}
-1,0)\rangle$.  Now $\nu_1 \in I_2$ as $\nu \restriction k \in I_2$
and for $n=k-1$ clause (c) of $\circledast$ fails
and $\nu_1 \char 94 \langle(2,\zeta +n)\rangle \in I_2$
because of $\nu_1 \in I_2$ and for $n=k$ the failure of clause 
(b) of $\circledast$ so continue as in Subcase 4B above.  

Lastly,
\nl
\ub{Subcase 4D}:  $\ell^{\nu(k)} \in \{-1,2\}$ and $\alpha^{\nu(k)}$ is a
limit ordinal.

Then $\{(\nu \restriction k) \char 94
\langle(\ell^{\nu(k)},\alpha)\rangle:\alpha < \alpha^{\nu(k)}\}$ 
is $\subseteq I_2$ and is an unbounded subset of $J_{\bar \eta}$ hence
cf$(J_{\bar \eta}) = \text{ cf}(\alpha^{\nu(k)})$.  If
cf$(\alpha^{\nu(k)}) = \aleph_0$, then clause (c) in $\boxtimes_3$
holds, and if cf$(\alpha^{\nu(k)}) \in [\aleph_1,\theta_1)$ then
necessarily $\alpha^{\nu(k)} \ne \zeta$ so being 
a limit ordinal $< \zeta + \omega$ clearly 
$\alpha^{\nu(k)} < \zeta$ so clause (d) from
$\boxtimes_3$ holds. To finish this subcase note that 
cf$(\alpha^{\nu(k)}) \ge \theta_1$ is impossible.
\nl
[Why ``impossible"?   Clearly for large enough $i < \partial$ we have 
$\eta_i(m) \ge \zeta$ (because $\delta = \zeta + \omega$ as 
said in the beginning of the
case) and recall $\nu \triangleleft \eta_i \in I_2$.  We now show that
clauses (a)-(d) of
$\circledast$ hold with $\eta_i,k$ here standing for $\eta,n$ there.
For clause (a) recall $\ell g(\eta_i) \ge \ell g(\nu) + 1$ and
$m = \ell g(\nu) = k+1$.  Now
$\ell^{\eta_i(k+1)} = \ell^{\eta_i(m)} =
\ell^* = -2$ as $\ell^* = -2$ is part of the case, 
$\ell^{\eta_i(k)} = \ell^{\nu(k)} \in \{-1,2\}$ in this subcase, 
so clause (d) of $\circledast$
holds.  Also $\alpha^{\eta_i(k+1)} = \alpha^{\eta_i(m)} \ge \zeta$ as 
said above so
clause (c) of $\circledast$ holds and cf$(\alpha^{\eta_i(k)}) 
= \text{ cf}(\alpha^{\nu(k)}) \ge \theta_1$ (as
we are trying to prove ``impossible"), so clause (b) of $\circledast$
holds.  Together we have proved 
(a)-(d) of $\circledast$.  But $\eta_i \in I_2$, contradiction.] 

Now subcases 4A,4B,4C,4D cover all the possibilities hence we are done
with case 4.
\bn
\ub{Case 5}:  $\ell^* = -2,\partial > \aleph_0$ and $\nu \char 94
\langle(\ell^*,\delta)\rangle \notin I_2$.

Recalling $\delta$ is the limit of the increasing sequence $\langle
\alpha^{\eta_i(m)}:i < \partial \rangle$ hence
$\text{cf}(\delta) = \partial > \aleph_0$ and $\nu \char 94 
\langle(-2,\delta)\rangle \notin I_2$, necessarily
$\delta = \zeta$ so $\partial = \theta_2$.  As $\nu \char 94
\langle(-2,\delta)\rangle \notin I_2$ necessarily clauses (a) - (d) of
$\circledast$ hold for some $n$ and as $\nu \in I_2$, clearly $n = \ell
g(\nu)-1$ (see clause (a) of $\circledast$) so we have $\ell g(\nu)
>0$, and letting
$k = \ell g(\nu)-1$, by clause (d) of $\circledast$ the
$\ell^{\eta(n+1)}$ there stands for $\ell^* = -2$ here so we have 
$\ell^{\nu(k)} \in \{-1,2\}$ and by clause (b) of $\circledast$ we have
cf$(\alpha^{\nu(k)}) \ge \theta_1$.  Hence $\{(\nu \restriction k) \char 94
\langle(\ell^{\nu(k)},\beta)\rangle:\beta < \alpha^{\nu(k)}\}$ is
cofinal in $J_{\bar \eta}$ and its cofinality is cf$(\alpha^{\nu(k)})$ as
$(\nu \restriction k) \char 94 \langle(\ell^{\nu(k)},\beta)\rangle$
increase (by $\le_{I_2}$) 
with $\beta$ as $\ell^{\nu(k)} \in \{-1,2\}$.  But
cf$(\alpha^{\nu(k)}) \ge \theta_1$ and $\partial = \theta_2$ (see first
sentence of the present case), so clause (e) of $\boxtimes_3$ holds.
\bn
\ub{Case 6}:  $\nu \char 94 \langle(\ell^*,\delta)\rangle \in I_2$.
\mn
\ub{Subcase 6A}:   $\nu \char 94 \langle(\ell^*,\delta),
(2,\zeta)\rangle \in I_2$.

Note that for $m = \ell g(\nu)$ and the pair $(\nu \char 94
\langle(\ell^*,\delta),(2,\zeta)\rangle,m)$ standing for $(\eta,n)$ in
$\circledast$, clauses (a),(c),(d) of
$\circledast$ hold (recall $\ell^* \in \{-2,1\}$, see the discussion
after case 1) so necessarily clause $(b)$ of $\circledast$ fails hence 
cf$(\delta) < \theta_1$ but $\partial = \text{ cf}(\delta)$ so
$\partial < \theta_1$.  Now as $\nu \char 94 \langle
(\ell^*,\delta),(2,\zeta)\rangle \in I_2$ clearly if $\ell < \omega$,
then $\nu \char 94 \langle (\ell^*,\delta),(2,\zeta + \ell)\rangle$
belongs to $I_2$ hence $\{\nu \char 94
\langle(\ell^*,\delta),(2,\zeta + \ell)\rangle:\ell < \omega\}$ is
a cofinal subset of $J_{\bar \eta}$ by the choice of
$I_2$ hence cf$(J_{\bar\eta}) = \aleph_0$ so clause (c) of $\boxtimes_3$
applies.
\bn
\ub{Subcase 6B}:  $\nu \char 94 \langle (\ell^*,\delta),
(2,\zeta)\rangle \notin I_2$.

As $\nu \char 94 \langle(\ell^*,\delta)\rangle \in I_2$, necessarily
clauses (a)-(d) of $\circledast$ hold with $(\nu \char 94 \langle
(\ell^*,\delta),(2,\zeta)\rangle,m)$ here standing for $(\eta,n)$
there, recalling $m = \ell g(\nu)$ so by clause (b) of $\circledast$ we
know that cf$(\delta) \ge \theta_1$ but $\partial = \text{ cf}(\delta)$
hence $\partial \ge \theta_1$.  Also 
$\{\nu \char 94 \langle (\ell^*,\delta),(2,\alpha)\rangle:\alpha <
\zeta\}$ is a subset of $I_2$ and cofinal in 
$J_{\bar\eta}$ and is increasing
with $\alpha$ so cf$(J_{\bar \eta}) = \theta_2$ so clause 
(e) of $\boxtimes_3$ applies.  

As the two subcases 6A,6B are complimentary case 6 is done.
\bn
\ub{Finishing the proof of $\boxtimes_3$}:

It is easy to check that our cases cover all the possibilities 
(as after discarding cases 0,1, if not case (6) then 
$\nu \char 94 \langle (\ell^*,\delta)\rangle
\notin I_2$, as not case (3), $\ell^* \ne 1$ but (see clause
$\odot(c)$ before case 2),
$\ell^* \in \{-2,1\}$ so necessarily $\ell^* = -2$, so case (4),(5) cover the
rest).  Together we have proved $\boxtimes_3$.]
\mr
\item "{$\boxtimes_4$}"  recall $\aleph_0 \le \mu < \theta_1 < \theta_2$; if
$X \subseteq I_2,|X| < \theta_2$ then we can
find $Y$ such that $X \subseteq Y \subseteq I_2,|Y| = \mu + |X|,Y$ is unbounded
in $I_2$ from below and from above and for every $\nu \in I_2 
\backslash Y$ the following linear orders have cofinality $\aleph_0$:
{\roster
\itemitem{ $(a)$ }  $J^2_{Y,\nu} = I_2 \restriction \{\eta \in I_2
\backslash Y:(\forall \rho \in Y)(\rho <_{I_2} \nu \equiv \rho <_{I_2} \eta)\}$
\sn
\itemitem{ $(b)$ }  the inverse of $J^2_{Y,\nu}$ 
\sn
\itemitem{ $(c)$ }   $J^-_{Y,\nu} = I_2 \restriction 
\{\eta \in I_2:(\forall \rho \in J^2_{Y,\nu})(\eta <_{I_2} \rho)\}$
\sn
\itemitem{ $(d)$ }  the inverse of $J^+_{Y,\nu}:= I_2 \restriction 
\{\eta \in I_2:(\forall \rho \in J^2_{Y,\nu})(\rho <_{I_2} \eta)\}$.
\endroster}
\ermn
[Why?  Let ${\Cal U} = \{\alpha^{\eta(\ell)}:\eta \in X$ and $\ell <
\ell g(\eta)\}$.
\nl
We choose $W_n$ by induction on $n < \omega$ such that 
\mr
\item "{$\boxdot_1$}"  $(a) \quad {\Cal U} \subseteq W_n \subseteq 
\zeta + \omega$
\sn
\item "{${{}}$}"  $(b) \quad W_n$ has cardinality $\mu + |{\Cal U}| =
\mu + |X|$ and $m<n \Rightarrow W_m \subseteq W_n$
\sn
\item "{${{}}$}"  $(c) \quad \mu \subseteq W_0$ and $\zeta + n \in
W_0$ for $n < \omega$
\sn
\item "{${{}}$}"  $(d) \quad \alpha \in W_n \Rightarrow \alpha +1
\in W_{n+1}$
\sn
\item "{${{}}$}"  $(e) \quad \alpha +1 \in W_n \Rightarrow \alpha \in
W_{n+1}$
\sn
\item "{${{}}$}"  $(f) \quad$ if $\delta \in W_n$ is a limit ordinal
of cofinality $< \theta_1$ then $\delta = \sup(\delta \cap W_{n+1})$
\sn
\item "{${{}}$}"  $(g) \quad$ if $\delta \in W_n$ and cf$(\delta) \ge
\theta_1$ (or just cf$(\delta) \le \mu + |X|$)
 then 
\nl

\hskip35pt sup$(\delta \cap W_n)+1 \in W_{n+1}$.
\ermn
This is straight.  Let $W = \cup\{W_n:n < \omega\}$, so
\mr
\item "{$\boxdot_2$}"  ${\Cal U} \subseteq W$ and $|W| = \mu + |X|$
and $W$ satisfies
{\roster
\itemitem{ $(a)$ }  $W \subseteq \zeta + \omega$
\sn
\itemitem{ $(b)$ }  $|W| < \theta_2$
\sn
\itemitem{ $(c)$ }  $0 \in W$ and $\{\zeta + m:m < \omega\} \subseteq W$
\sn
\itemitem{ $(d)$ }   $\alpha \in W \Leftrightarrow \alpha +1 \in W$
\sn
\itemitem{ $(e)$ }   if $\delta \in W$ and $\aleph_0 < \text{ cf}(\delta) \le
\mu$ \ub{then} $\delta = \text{\rm sup}(W \cap \delta)$
\sn
\itemitem{ $(f)$ }    if $\delta \in W$ and
cf$(\delta) \ge \theta_1$ or $\text{cf}(\delta) = \aleph_0$ then 
$\text{\rm cf(otp}(W \cap \delta))) = \aleph_0$.
\endroster}
\ermn
Let $Y = \{\eta \in I_2:\alpha^{\eta(\ell)} \in W$ for every $\ell <
\ell g(\eta)\}$.  Clearly $X \subseteq Y$ and $|Y| = \aleph_0 + |W| = 
\mu + |{\Cal U}| < \theta_2$.
It suffices to check that $Y$ is as required in $\boxtimes_4$.  From
now on we shall use only the choice of $Y$ and clauses (a)-(f) of $\boxdot_2$.
By $\boxdot_2(c)$ and the choice of $Y$ clearly $Y$ is unbounded in
$I_2$ from above and from below.  

So let $\nu \in I_2 \backslash Y$, as $\nu 
\restriction 0 \in Y$ there is $n < \ell
g(\nu)$ such that $\nu \restriction n \in Y,\nu \restriction (n+1)
\notin Y$, so $\alpha^{\nu(n)} < \zeta + \omega$, and
$\alpha^{\nu(n)} \notin W$, but by clause (c) of $\boxdot_2$ we have
$\{\zeta + m:m < \omega\} \subseteq W$ hence $\alpha^{\nu(n)} < \zeta$ and
so $\alpha_1 := \text{\rm Min}(W \backslash
\alpha^{\nu(n)})$ is well defined, is $\le \zeta$ 
and $> \alpha^{\nu(n)}$.  As clearly $0
\in W,\beta \in W \Leftrightarrow \beta +1 \in W$ by the choice of $W$,
obviously $\alpha_1$ is a limit ordinal.  By clause (e) of $\boxdot_2$
 clearly $\alpha_1$ is of cofinality $\aleph_0$ or $\ge \theta_1 =
 \mu^+$.  So clearly $\alpha_0 := \text{\rm sup}
(W \cap \alpha^{\nu(n)}) = \sup(W \cap \alpha_1) = 
\text{ min}\{\alpha:W \cap \alpha = W \cap \alpha^{\nu(n)}\}$ is a 
limit ordinal 
$\le \alpha^{\nu(n)}$ and $\alpha_0 \notin W$ so cf$(\alpha_0) \le |W|
< \theta_2$ but by the assumption on $W$, 
(see clause (f) of $\boxdot_2$) we have
cf$(\alpha_0) = \aleph_0$.  So $(\nu
\restriction n) \char 94 \langle(\ell^{\nu(n)},\alpha_0)\rangle
\in J^2_{Y,\nu}$; moreover
\mr
\item "{$\boxdot_3$}"  $\rho \in J^2_{Y,\nu}$ \ub{iff} $\rho \in I_2$
satisfies one of the following:
{\roster
\itemitem{ $(a)$ }  $(i) \quad \nu \rest n=\rho\rest n$, and
$\ell^{\nu(n)}=\ell^{\rho(n)}$,
\sn
\itemitem{ ${{}}$ }  $(ii) \quad \alpha^{\rho(n)} \in 
[\alpha_0,\alpha_1)$
\sn
\itemitem{ $(b)$ }  $(i) \quad \nu \rest n = \rho \rest n$, and
$\ell^{\nu(n)}=\ell^{\rho(n)}$,
\sn
\itemitem{ ${{}}$ }  $(ii) \quad \alpha^{\rho(n)} = \alpha_1$ and 
$\alpha^{\rho(n+1)} \in [\sup(W \cap \zeta),\zeta)$
\sn
\itemitem{ ${{}}$ }  $(iii) \quad (\ell^{\rho(n+1)},\ell^{\rho(n)}) =
(\ell^{\rho(n_1)},\ell^{\nu(n)}) \in \{(2,-2),(2,1),(-2,-1),(-2,2)\}$
\sn
\itemitem{ $(c)$ }  $(i) \quad \alpha_1 = \zeta$ and $n > \theta$ and 
$(\nu \rest n) \char 94 (\ell^{\nu(n)},\alpha_1) \notin I_2$
\sn
\itemitem{ ${{}}$ }  $(ii) \quad (\ell^{\nu(n)},\ell^{\nu(n-1)}) \in
\{(2,-2),(2,1),(-2,2),(-2,-1)\}$
\sn
\itemitem{ ${{}}$ }  $(iii) \quad$ cf$(\nu(n)) \ge \theta_1$ and 
$\nu(n) > \sup(W \cap \nu(n))$
\sn
\itemitem{ ${{}}$ }  $(iv) \quad \rho \rest (n-1) = \nu
\rest(n-1),\ell^{\rho(n-1)} = \ell^{\nu(n-1)}$ 
\sn
\itemitem{ ${{}}$ }  $(v) \quad \alpha^{\rho(n-1)} \in [\sup(\nu(n-1)
\cap W),\nu(n-1))$.
\endroster}
\ermn
[Why?  First note that if $\rho\in J^2_{Y,\nu}$ and $\rho\rest
k=\nu\rest k$, $\rho(k) \ne \nu(k)$, and $k \le n$ then necessarily $k
= n \wedge \ell^{\rho(k)} = \ell^{\nu(k)}$.  We now proceed to check ``if".
Let $f:\{-2,-1,1,2\} \rightarrow \{2,-2\}$ so that
$f^{-1}(\{2\}) = \{-2,1\}$ and $f^{-1}(\{-2\}) = \{-1,2\}$.
Case (a) is obvious.  In case (b) in order for $\eta\in Y$ to
separate between $\nu$ and $\rho$ it is necessary that
$\eta\rest(n+1)=\rho\rest(n+1)$,
$\ell^{\eta(n+1)}=\ell^{\rho(n+1)}=f(\ell^{\rho(n)})$ and that
$\alpha^{\eta(n+1)} \ge \zeta$, but then $\eta \notin I_2$.
In case (c) in order to separate between $\rho$ and $\nu$ by
$\eta\in Y$ there are two possibilities.  Either $\eta\rest n =
\nu\rest n$ and then
$\ell^{\eta(n)} = \ell^{\nu(n)} = f(\ell^{\nu(n-1)})$ (recall that
$\nu \rest n \char 94 \langle(\ell^{\nu(n)},\alpha_1)\rangle \notin I_2$), and
$\alpha^{\eta(n)} \ge \zeta$, but then also $\eta\notin
I_2$. The other possibility is that $\eta \rest (n-1) = \nu \rest
(n-1),\ell^{\eta(n-1)} = \ell^{\nu(n-1)}$ and
$\alpha=\alpha^{\eta(n-1)}$ is such that $\alpha\in W$  and
$\alpha^{\rho(n-1)}<\alpha<\alpha^{\nu(n-1)}$ which is also
impossible by the choice of $\alpha^{\rho(n-1)}$.  Showing that
these are the only cases (the ``only if" direction) is similar and is
actually done below.]

Now we proceed to check that clauses of $\boxtimes_4$ hold.
\bn
\ub{Clause (a)}: 

First assume $\ell^{\nu(n)} \in \{-2,1\}$, and let $J = \{\nu \rest
n \char 94 \langle(\ell^{\nu(n)},\alpha_0),(2,\zeta +m):m<\omega\}$.  Now
$J \subseteq I_2$ [why? clearly if $\rho\in J$ then
$\rho\rest(n+1)\in I_2$ so we only need to check $\circledast$ for
$n$, recall that cf$(\alpha_0) = \aleph_0<\theta_1$, hence clause
$(b)$ of $\circledast$ fails].  Now by clause $(a)$ of $\boxdot_3$
we have that $J \subseteq J^2_{Y,\nu}$, and we claim that it is
also cofinal in it. [Why? Note that as
$\ell^{\nu(n)} \in \{-2,1\}$ then $\nu \rest n \char 94
\langle(\ell^{\nu(n)},\alpha_0) <_{I_2} \nu \rest (n+1)$, and if
$\rho \in J^2_{Y,\nu}$ is as in clauses $(a)$ or $(b)$ of
$\boxdot_3$ then for every $m$ large enough $\rho <_{I_2} \nu \rest
n \char 94 \langle(\ell^{\nu(n)},\alpha_0),(2,\zeta +m)\rangle$. If
$\rho \in J^2_{Y,\nu}$ is as in clause $(c)$  of $\boxdot_3$ then
$\ell^{\nu(n)} \in \{-2,2\}$ by (ii) there, and as in this case
$\ell^{\nu(n)} \in \{-2,1\}$,  necessarily  $\ell^{\nu(n)}=-2$ and so
by (ii) of (c) of $\boxdot_3$ we have
$\ell^{\nu(n-1)} \in \{-1,2\}$, but then $\rho <_{I_2} \nu$ and so it is
below every element in $J$.]

Second, assume $\ell^{\nu(n)} \in \{-1,2\}$ and 
$\nu \rest n \char 94 \langle(\ell^{\nu(n)},\alpha_1)\rangle \in I_2$;
let $\delta^* = \sup(W \cap \zeta)$, so as above $\delta^* \notin W$,
and has cofinality $\aleph_0$ (which is less than $\theta_1$),
 recall also that cf$(\alpha_1) \ge \theta_1$.  So (for $\ell \in
 \{-2,-1,1,2\}$) by $\circledast$ we have
$(\nu \restriction n) \char 94
\langle(\ell^{\nu(n)},\alpha_1),(\ell,\beta)\rangle \in I_2$ iff
$\beta < \zeta \wedge \ell \in \{-2,-1,1,2\}$ or $(\zeta \le \beta <
\zeta + \omega \wedge \ell \ne -2)$.  Hence we have
$(\nu \restriction n) \char 94 \langle
(\ell^{\nu(n)},\alpha_1),(-2,\beta)) \in I_2 
\Leftrightarrow \beta < \zeta$.  Also $(\nu \restriction n) \char 94 \langle
(\ell^{\nu(n)},\alpha_1),(-2,\beta)\rangle \in Y \Leftrightarrow
\beta \in W$,  and as $\nu(n) < \alpha_1 \wedge \ell^{\nu(n)} \in
\{-1,2\}$ clearly $\nu <_{I_2} (\nu \rest n) \char 94 \langle
(\ell^{\nu(n)},\alpha_1),(-2,\beta)\rangle$.  Easily
$\{(\nu \restriction n) \char 94 \langle(\ell^{\nu(n)},\alpha_1),
(-2,\varepsilon)\rangle:\varepsilon \in W \cap \zeta)\}$ is a subset of
$\{\eta \in Y:\nu <_{I_2} \eta\}$ unbounded from below in it.

So $\{(\nu\restriction n) \char 94 \langle
(\ell^{\nu(n)},\alpha_1),(-2,\delta^*),(2,\alpha)\rangle:\zeta 
< \alpha < \zeta + \omega\}$ is included in $I_2$
(recalling clause (b) of $\circledast$ as cf$(\delta^*) = \aleph_0$)
and moreover
is a cofinal subset of $J^2_{Y,\nu}$ of order type $\omega$,
so cf$(J^2_{Y,\nu}) = \aleph_0$ as required.

Third, assume $\rho^{\nu(n)} \in \{-1,2\}$ and $(\nu \restriction n)
\char 94 \langle(\ell^{\nu(n)},\alpha_1)\rangle \in I_2$ and
cf$(\alpha_1) < \theta_1$, equivalently cf$(\alpha_1) = \aleph_0$ by
clause (e) of $\boxdot_2$.  In this case $\{(\nu \restriction n) \char 94
\langle(\ell^{\nu(n)},\alpha)(-2,\beta)\rangle:\zeta \le
\beta < \zeta + \omega\}$ is included in $I_2$ (recalling
clause (b) of $\circledast$) and in $Y$, hence recalling
$\boxdot_3(a)$ the set $\{(\nu \restriction n) \char 94
\langle(\ell^{\nu(n)},\alpha)\rangle:\alpha \in [\alpha_0,\alpha_1)\}$
is a cofinal subset of $J^2_{Y,\nu}$ hence its cofinality is
cf$(\alpha_1) = \aleph_0$ as required.

Fourth, we are left with the case $\ell^{\nu(n)} \in \{-1,2\}$ and
$(\nu \restriction n) \char 94 \langle
(\ell^{\nu(n)},\alpha_1)\rangle \notin I_2$ so necessarily $n > 0$ and 
clauses (a)-(d) of
$\circledast$ hold for it for $n-1$; then by clause (c) of $\circledast$
(recalling $\alpha_1 \le \zeta$ as shown before
$\boxdot_3$) necessarily $\alpha_1 = \zeta$.  
Clearly $k := n-1 \ge 0$ and as clause (d) of $\circledast$
holds and it says there ``$\ell^{\eta(n+1)} \in \{2,-2\}$" which means here
$\ell^{\nu(n)} \in \{2,-2\}$ but we are assuming presently
$\ell^{\nu(n)} \in \{-1,2\}$  
hence $\ell^{\nu(n)} = \ell^{\nu(k+1)} = 2$ so using clause (d) of
$\circledast$, see above, it follows that $\ell^{\nu(k)} \in
\{-2,1\}$ and by clause (b) of $\circledast$ we have
cf$(\alpha^{\nu(k)}) \ge \theta_1$.  Let $\delta_* = \sup(W
\cap \alpha^{\nu(k)})$.  Now if $\delta_* < \alpha^{\nu(k)}$ then
 by clause (f) of $\boxdot_2$ we know
cf$(\delta_*) = \aleph_0$ and $\{(\nu \restriction k) \char 94 
\langle(\ell^{\nu(k)},\delta_*)(2,\zeta +m)\rangle:
m < \omega\}$ is included in $I_2$ (as $\nu \in I_2$
and $\delta_* \le \alpha^{\nu(k)}$ we have to check in $\circledast$
only with $k+1$ here standing for $n$ there, but cf$(\delta_*) =
\aleph_0$ so clause (b) there fails) and so recalling $\boxdot_3(c)$ this set
is a cofinal subset of $J^2_{Y,\nu}$
exemplifying that its cofinality is $\aleph_0$.

Lastly, if $\delta_* = \alpha^{\nu(k)}$ then $\langle(\nu \rest n)
\char 94 \langle (\ell^{\nu(n)},\alpha)\rangle:\alpha \in W \cap
\zeta\rangle$ is $<_{I_2}$-increasing with $\alpha$, all members in
$Y$, and in $J^2_{Y,\nu}$, cofinal in it and has order type otp$(W
\cap \zeta)$ which has cofinality $\aleph_0$ so also $J^2_{Y,\nu}$ has
cofinality $\aleph_0$ as required.
\bn
\ub{Clause (b)}:  What about the cofinality of the inverse?  Recall
that $I_2$ is isomorphic to its inverse by the mapping $(\ell,\beta)
\mapsto (-\ell,\beta)$, but this isomorphism maps $Y$ onto itself
hence it maps $J^2_{Y,\nu}$ onto $J^2_{Y,\nu'}$ for some $\nu' \in I_2
\backslash Y$, but clause (a) was proved also for $\nu'$, so this follows.
\bn
\ub{Clause (c)}:  As $Y$ is unbounded from below in $I_2$ (containing
$\{\langle(-2,\zeta + n)\rangle:n < \omega\})$ it follows that
$J^-_{Y,\nu}$ is non-empty, hence cf$(J^-_{Y,\nu}) \ne 0$, 
but what is cf$(J^-_{Y,\nu})$?  

First, if $\ell^{\nu(n)} \in \{-1,2\}$ then $\{(\nu \restriction n) \char 94
\langle (\ell^{\nu(n)},\alpha)\rangle:\alpha < \alpha_0\}$ is an unbounded
subset of $J^-_{Y,\nu}$ of order type $\alpha_0$ hence cf$(J^-_{Y,\nu}) =
\text{\rm cf}(\alpha_0) = \aleph_0$ (see the assumption on $W$ 
and the choice of $\alpha_0$).

Second, if $\ell^{\nu(n)} = \{-2,1\}$ and $(\nu \restriction n) \char 94
\langle (\ell^{\nu(n)},\alpha_1)\rangle \in I_2$ and cf$(\alpha_1) \ge
\theta_1$ then as in the proof
of clause (a) we have $\{(\nu \restriction n) \char 94
\langle(\ell^{\nu(n)},\alpha_1),(2,\zeta + m)\rangle \notin
I_2$ for $m < \omega$ and again letting $\delta^* = \sup(W \cap
\zeta)$ we have $\{(\nu \restriction n) \char 94
\langle(\ell^{\nu(n)},\alpha_1),(2,\beta)\rangle:\beta \in W \cap
\zeta\}$ is included in $I_2$ and in $J^-_{Y,\nu}$ and even
is an unbounded subset of $J^-_{Y,\nu}$ of order type
otp$(W \cap \delta^*)$ which has the same cofinality as $\delta^*$
which is $\aleph_0$.

Third, if $\ell^{\nu(n)} \in \{-2,1\}$ and $(\nu \restriction n) \char
94 \langle(\ell^{\nu(n)},\alpha_1)\rangle \in I_2$ and cf$(\alpha_1) <
\theta_1$, equivalently cf$(\alpha_1) = \aleph_0$, \ub{then} $\{(\nu
\restriction n) \char 94 \langle(\ell^{\nu(n)},\alpha_1),(2,\zeta + m)
\rangle:m < \omega\}$ is a subset of $I_2$ (as
cf$(\alpha_1) = \aleph_0$) is included in $J^-_{Y,\nu}$, unbounded in
it and has cofinality $\aleph_0$, so we are done.

Fourth and lastly, if $\ell^{\nu(n)} \in \{-2,1\}$ and $(\nu \restriction n)
\char 94 \langle (\ell^{\nu(n)},\alpha_1)\rangle \notin I_2$ then as
in the proof of clause (a) we have
$\alpha_1 = \zeta$ and again letting $\delta^* = \sup(W \cap \zeta)$
we have cf$(\delta^*) = \aleph_0$ and $(\nu \restriction n)
\char 94 \langle (\ell^{\nu(n)},\delta^*)\rangle \in I_2$ and $\{(\nu
\restriction n) \char 94 \langle (\ell^{\nu(n)},\delta^*),(2,\zeta +m)
\rangle:m < \omega\}$ is a subset of $I_2$, moreover a
subset of $J^-_{Y,\nu}$ unbounded in it and $(\nu \restriction n)
\char 94 \langle (\ell^{\nu(n)},\delta^*),(2,\zeta + m)\rangle$ is
$<_{I_2}$--increasing with $m$.  So indeed $J^-_{Y,\nu}$ has
cofinality $\aleph_0$.
\bn
\ub{Clause (d)}:  As in clause (b) we use the anti-isomorphism.
\nl
So $\boxtimes_4$ holds.]
\mr
\item "{$\boxtimes_5$}"  if $I' \subseteq I_2$ then the number of cuts
of $I'$ induced by members of $I_2 \backslash I'$, that is $\{\{s \in
I':s <_{I_2} t\}:t \in I_2 \backslash I'\}$ is $\le |I'| +1$.
\nl
[Why?  Let ${\Cal U} := \{\alpha^{\eta(\ell)}:\ell < \ell g(\eta)$ and $\eta
\in I'\}$, it belongs to $[\zeta + \omega]^{\le \mu}$.
Now (by inspection) $\eta_1,\eta_2 \in I_2 \backslash I'$ 
realizes the same cut of $I'$ \ub{when}:
{\roster
\itemitem{ $(a)$ } $\ell g(\eta_1) = \ell g(\eta_2)$
\sn
\itemitem{ $(b)$ }  $\ell^{\eta_1(n)} = \ell^{\eta_2(n)}$ for $n <
\ell g(\eta_1)$
\sn
\itemitem{ $(c)$ }  $\alpha^{\eta_1(n)} \in {\Cal U} \Leftrightarrow
\alpha^{\eta_2(n)} \in {\Cal U} \Rightarrow \alpha^{\eta_1(n)} =
\alpha^{\eta_2(n)}$ for $n < \omega$
\sn
\itemitem{ $(d)$ }   $\beta < \alpha^{\eta_1(n)} \equiv \beta <
\alpha^{\eta_2(n)}$ for $\beta \in {\Cal U}$ and $n < \omega$
\endroster}
\ermn
[Why?  Now clauses (a)-(d) define an equivalence relation on $I_2 \backslash
I'$ which refines ``inducing the same cut" and has $\le |{\Cal U}| +
\aleph_0 = |I'| + \aleph_0$ equivalence classes.  
As the case $I'$ is finite is trivial, we
are done proving $\boxtimes_5$.]
\mr
\item "{$\boxtimes_6$}"  if $\partial$ is regular uncountable, $n^* <
\omega$ and $t_{\varepsilon,\ell} \in I_2$ for $\varepsilon < 
\partial,\ell < n^*$ and
$t_{\varepsilon,0} <_{I_2} \ldots <_{I_2} t_{\varepsilon,n^*-1}$ for 
$\varepsilon <
\partial$ \ub{then} for some unbounded (and even stationary) set $S
\subseteq \partial,m\le n^*$ and $0=k_0 < k_1 < \ldots < k_m =n^*$
stipulating $t_{\varepsilon,k_m} = \infty$ and letting $\varepsilon(*) = 
\text{\rm Min}(S)$ we have:
{\roster
\itemitem{ $(a)$ }   for each $i < m$:
\sn
\itemitem{ ${{}}$ }   $(\alpha) \quad$ if $\varepsilon < \xi$ are from $S$
and $\ell_1,\ell_2 \in [k_i,k_{i+1})$ then $t_{\varepsilon,\ell_1} <_{I_2}
t_{\xi,\ell_2}$ \ub{or}
\sn
\itemitem{ ${{}}$ }   $(\beta) \quad$ if $\varepsilon < \xi$ are from $S$ and
$\ell_1,\ell_2 \in [k_i,k_{i+1})$ then $t_{\xi,\ell_2} <_{I_2}
t_{\varepsilon,\ell_1}$ \ub{or}
\sn 
\itemitem{ ${{}}$ }   $(\gamma) \quad k_{i+1} = k_i+1$ and for every
$\varepsilon \in S$ we have $t_{\varepsilon,k_i} = t_{\varepsilon(*),k_i}$
\sn
\itemitem{ $(b)$ }   there is a sequence $\langle s^-_i,s^+_i:i < m
\rangle$ such that
\sn
\itemitem{ ${{}}$ }   $(\alpha) \quad i < m 
\Rightarrow s^-_i <_{I_2} s^+_i$ 
\sn
\itemitem{ ${{}}$ }   $(\beta) \quad$ if 
$i < m-1$ then $s^+_i < s^-_{i+1}$ except
possibly when $\langle t_{\varepsilon,k_i}:\varepsilon <
\partial \rangle$ is $<_{I_2}$-decreasing and there is no $t \in I_2$
such that $\varepsilon < \partial \Rightarrow t_{\varepsilon,k_i}
<_{I_2} t <_{I_2} t_{\varepsilon,k_{i+1}}$, hence (by $\boxtimes_3$)
we have $\partial \ge \theta_2$
\sn
\itemitem{ ${{}}$ }   $(\gamma) \quad$ for each $i < m$ the set
$\{t_{\varepsilon,\ell}:\varepsilon \in S$ and 
$\ell \in [k_i,k_{i+1})\}$ is included in the interval
$(s^-_i,s^+_i)_{I_2}$.
\endroster}
\ermn
[Why?  Straight.   For some stationary $S_1 \subseteq \partial$ and
$\langle n_k:k < n^*\rangle$ we have $\varepsilon \in S_1 \wedge k <
n^* \Rightarrow \ell g(t_{\varepsilon,k})=n_k$.  
Without loss of generality also $\langle
\ell^{t_{\varepsilon,k}(i)}:i < n_k\rangle$ does not depend on $\varepsilon \in
S_1$.  By $\dsize \sum_{k<n^*} n_k$ application of $\partial \rightarrow
(\partial,\omega)^2$, \wilog \, for each $k<n^*$ and $i < n_k$ the
sequence $\langle \alpha^{t_{\varepsilon,k}(i)}:\varepsilon \in S_1\rangle$ is
constant or increasing.  Cleaning a little more we are done.
\nl
So $\boxtimes_6$ holds.]
\bn
Lastly, recall that we chose $I$ to be 
$(|M|,<^M)$, where $M$ was the real closure of $M_0$ and (see
$(*)_1$), $M_0$ the ordered field generated over $\Bbb Q$ by $\{a_t:t
\in I_2\}$ as described in $(*)_1$ above and for every $u \subseteq \zeta$ let:
\mr
\item "{$(*)_2$}"  $(a) \quad I^1_u = \{(\ell,\beta) \in I_1:\beta \in
u$ or $\beta \in [\zeta,\zeta + \omega)\}$
\sn
\item "{${{}}$}"  $(b) \quad I^{*,2}_u = \{\eta \in
I^*_2:\alpha^{\eta(\ell)} \in I^1_u$ for every $\ell < \ell g(\eta)\}$
\sn
\item "{${{}}$}"  $(c) \quad I^2_u = \{\eta \in I_2:\alpha^{\eta(\ell)} \in
I^1_u$ for every $\ell <\ell g(\eta)\}$
\sn 
\item "{${{}}$}"  $(d) \quad I_u =$ the real closure of
$\Bbb Q(a_t:t \in I^2_u)$ in $M$
\sn
\item "{${{}}$}"  $(e) \quad$ for $t \in I_2 \backslash I^2_u$ let \nl

\hskip20pt $I^2_{u,t} = I_2 \restriction \{s \in I_2:s \notin
I^2_u$ and for every $r \in I^2_u$ 
\nl

\hskip100pt we have $r <_{I_2} t \equiv r <_{I_2} s\}$
\sn
\item "{${{}}$}"  $(f) \quad$ for $x \in I \backslash I_u$ let \nl

\hskip25pt $I_{u,x} = I \restriction \{y \in I:y \notin I_u$ and 
$(\forall a \in I_u)(a <_I y \equiv a <_I x)\}$
\sn
\item "{${{}}$}"  $(g) \quad$ let $\hat I_u$ be the set $I_u \cup
\{I_{u,a}:a \in I \backslash I_u\}$ ordered by: 
$x <_{\hat I_u} y$ \ub{iff} \nl

\hskip25pt one of the following holds:
{\roster
\itemitem{ ${{}}$ }  $(\alpha) \quad x,y \in I_u$ and $x <_{I_u} y$ 
\sn
\itemitem{ ${{}}$ }  $(\beta) \quad x \in I_u,y = I_{u,b}$ and $x <_{I_u} b$ 
\sn
\itemitem{ ${{}}$ }  $(\gamma) \quad x = I_{u,a},y \in I_u$ and $a <_{I_u} y$
\sn
\itemitem{ ${{}}$ }  $(\delta) \quad x = I_{u,a},y = I_{u,b}$ and $a <_{I_u} b$
(can use it more!)
\nl

\hskip35pt (note that by 
$\boxtimes_5$, $|u| \le \mu \Rightarrow |\hat I_u| \le \mu$).
\endroster}
\ermn
Now observe
\mr
\item "{$(*)_3$}"  for $u \subseteq \zeta,I^2_u$ is unbounded in $I_2$
from below and from above.
\ermn
We define
\mr
\item "{$(*)_4$}"  we say\footnote{we may in clauses (e) + (c) replace
$\mu$ by $\mu + |{\Cal U}|$, no harm and it makes $(c)(\beta)$ of
$(*)_1$, redundant} that $u$ is $\mu$-reasonable if:
{\roster
\itemitem{ $(a)$ }  $u \subseteq \zeta,|u| < \theta_2$ and $\mu \subseteq u$
\sn
\itemitem{ $(b)$ }  $\alpha \in u \equiv \alpha +1 \in u$ for every $\alpha$
\sn
\itemitem{ $(c)$ }  if $\delta \in u$ and $\aleph_0 \le
\text{ cf}(\delta) \le \mu$ then $\delta = \text{\rm sup}(u \cap \delta)$
\sn
\itemitem{ $(d)$ }  if $\delta \le \zeta$ and cf$(\delta) > \mu$ then
cf$(\text{otp}(\delta \cap u)) = \aleph_0$.
\endroster}
\ermn
Now we note
\mr
\item "{$(*)_5$}"  if $X \subseteq I$ has cardinality $< \theta_2$ and
$u_* \subseteq \zeta$ has cardinality $< \theta_2$ then we can find a
$\mu$-reasonable $u$ such that $X \subseteq I_u$ and $u_* \subseteq u$
and $|u| = \mu + |X| + |u_*|$.
\ermn
[Why?  By the proof of $\boxtimes_4$.]
\mr
\item "{$(*)_6$}"  if $u$ is $\mu$-reasonable then $Y := I^2_u$
satisfies the conclusions of $\boxtimes_4$.
\ermn
[Why?  By the proof of $\boxtimes_4$, that is if $u^+ := u \cup
\{\zeta +n:n < \omega\}$ then $Y$ as defined in the proof there using
$u^+$ for $W$, is $I^2_u$ from
$(*)_2(c)$, and it satisfies demands (a)-(f) from $\boxdot_2$ so the
proof there applies.]
\mr
\item "{$(*)_7$}"  if $u$ is $\mu$-reasonable and $x \in I
\backslash I_u$ then cf$(I_{u,x}) \le \aleph_0$.
\ermn
Why?  The proof takes awhile.  Toward contradiction assume 
$\partial = \text{\rm cf}(I_{u,x})$ is $> \aleph_0$ and let
$\langle b_\varepsilon:\varepsilon < \partial \rangle$ 
be an increasing sequence of
members of $I_{u,x}$ unbounded in it.  So for each $\varepsilon < \partial$
there is a definable function $f_\varepsilon(x_0,\dotsc,x_{n(\varepsilon)-1})$
where definable of course means
in the theory of real closed fields and $t_{\varepsilon,0} <_{I_2}
t_{\varepsilon,1} <_{I_2} \ldots <_{I_2}
t_{\varepsilon,n(\varepsilon)-1}$ from $I_2$
such that $M \models ``b_\varepsilon = f_\varepsilon(a_{t_{\varepsilon,0}},
\dotsc,a_{t_{\varepsilon,n(\varepsilon)-1}})"$ and $n(\varepsilon)$ is
minimal. As Th$(\Bbb R)$ is countable and $\aleph_0 < \partial = \text{\rm
cf}(\partial)$, \wilog \, $\varepsilon < \partial
\Rightarrow f_\varepsilon = f_*$ so
$\varepsilon < \partial \Rightarrow n(\varepsilon) = n(*)$.

Apply $\boxtimes_6$ to $\langle \bar t^\varepsilon = \langle 
t_{\varepsilon,\ell}:\ell < n(*)\rangle:\varepsilon < \partial \rangle$
and get $S \subseteq \partial$ 
and $0 = k_0 < k_1 <
\ldots < k_m = n(*)$ and 
$\langle (s^-_i,s^+_i):i < m \rangle$ and $\varepsilon(*)
= \text{\rm Min}(S)$ as there.  Without loss of generality the truth
value of $``t_{\varepsilon,\ell} \in I^2_u"$ for $\varepsilon \in 
S$, depends just on $\ell$.  Let $w_1 = \{i < m:(\forall \varepsilon \in
S)(t_{\varepsilon,k_i} = t_{\varepsilon(*),k_i})\},
w_2 = \{\ell < n(*):t_{\varepsilon(*),\ell} \in I^2_u\}$; clearly for
every $\ell < n(*)$ we have $(\forall \varepsilon \in
S)(t_{\varepsilon,\ell} = t_{\varepsilon(*),\ell}) \Leftrightarrow \ell
\in \{k_i:i \in w_1\}$ and $i \in w_1 \Rightarrow k_i +1 = k_{i+1}$.  

Let $t^*_{k_i} = t_{\varepsilon,k_i}$ for ($\varepsilon < \partial$ and
$i \in w_1$).
Renaming \wilog \, $S = \partial$ and $\varepsilon(*)=0$.

We have some free choice in choosing $\langle b_\varepsilon:\varepsilon 
< \partial \rangle$ (as long as it is cofinal in $I_{u,x}$), so \wilog \, we
choose it such that $n(*)$ is minimal and then $|w_1|$ is maximal and
then $|w_2|$ is maximal.

Now does the exceptional cae in $(b)(\beta)$ of $\boxtimes_6$ occurs?
This is an easier case and we delay it to the end.

As $I_2$ and $I_{2,<t}$ for $t \in I_2$ have cofinality $\aleph_0$ (see
$\boxtimes_2(a),(b))$ and $\boxtimes_3$ and this holds for 
the inverse of $I_2$, too, while
$\partial = \text{ cf}(\partial) > \aleph_0$ and we can replace $\langle
b_\varepsilon:\varepsilon < \partial\rangle$ by $\langle
b_{n(*) + \varepsilon}:\varepsilon < \partial\rangle$ we can find
$t_{\partial,\ell} \, (\ell < n(*))$ such that
\mr
\item "{$\odot$}"  $(a) \quad t_{\partial,0} <_{I_2} t_{\partial,1}
<_{I_2} \ldots <_{I_2} t_{\partial,n(*)-1}$
\sn
\item "{${{}}$}"  $(b) \quad$ if $\varepsilon < \xi < \partial$ and
$\ell_1,\ell_2 < n(*)$ then $(t_{\varepsilon,\ell_1} <_{I_2}
t_{\partial,\ell_2}) \equiv (t_{\varepsilon,\ell_1} <_{I_2} t_{\xi,\ell_2})$
\nl

\hskip25pt  and $(t_{\partial,\ell_1}  <_{I_2} t_{\varepsilon,\ell_2}) \equiv
(t_{\xi,\ell_1} < t_{\varepsilon,\ell_2})$
\sn
\item "{${{}}$}"  $(c) \quad$ if $\ell \in [k_i,k_{i+1})$ then
$t_{\partial,\ell} \in (s^-_i,s^+_i)_{I_2}$.
\endroster 
\bn
\ub{Case 0}: $\{0,\dotsc,m-1\} = w_1$.

This implies $i < m \Rightarrow k_i +1 = k_{i+1}$ hence $m=n$ hence
$\ell < n \Rightarrow t_{\xi,\ell} = t^*_\ell$ and so 
contradicts ``$\langle b_\varepsilon:\varepsilon < \partial\rangle$ 
is increasing" (as it becomes constant).
\bn
\ub{Case 1}:  $[0,m) \backslash w_1$ is not a singleton.

It cannot be empty by ``not case 1".  
Choose $i(*) \in \{0,\dotsc,m-1\} \backslash w_1$
and for $\varepsilon,\xi < \partial$ let $\bar t^{\varepsilon,\xi} = \langle
t^{\varepsilon,\xi}_\ell:\ell < n(*)\rangle$ be defined by:
$t^{\varepsilon,\xi}_\ell$ is $t_{\varepsilon,\ell}$ if $\ell \in
[k_{i(*)},k_{i(*)+1})$ and $t_{\xi,\ell}$ otherwise.  Let
$b_{\varepsilon,\xi} = f_*(a_{t^{\varepsilon,\xi}_0},
\dotsc,a_{t^{\varepsilon,\xi}_{n(*) -1}}) \in M$.
\sn
Clearly
\mr
\item "{$\circledast_0$}"  for any
$\varepsilon_1,\varepsilon_2,\xi_1,\xi_2 \le \partial$ 
the truth value of $b_{\varepsilon_1,\xi_1} < b_{\varepsilon_2,\xi_2}$
depend just on the inequalities which $\langle \varepsilon_1,
\varepsilon_2,\xi_1,\xi_2\rangle$ satisfies and even just on the
inequalities which the
$t_{\varepsilon_1,\ell},t_{\varepsilon_2,\ell},t_{\xi_1,\ell},t_{\xi_2,\ell}
\, (\ell < n(*))$ satisfy.
\ermn
[Why?  Recall $\left< \langle t_{\varepsilon,\ell}:
\ell < n(*)\rangle:\varepsilon \in S \right>$ is an indiscernible 
sequence in the linear order $I_2$ (for
quantifier free formulas) and $M$ has elimination of quantifiers.]
\mr
\item "{$\circledast_1$}"  $\dsize \bigwedge_{\ell=1,2}
\varepsilon(0) < \varepsilon_\ell < \varepsilon(1) < \partial 
\Rightarrow b_{\varepsilon(0)} <_I b_{\varepsilon_1,\varepsilon_2} <_I
b_{\varepsilon(1)}$.
\ermn
[Why?  By $\circledast_0$ the desire statement, $b_{\varepsilon(0)} <_I
b_{\varepsilon_1,\varepsilon_2} <_I b_{\varepsilon(1)}$ is equivalent to 
$b_{\varepsilon(0)} < b_{\varepsilon_1,\varepsilon_1} <
b_{\varepsilon(1)}$ which means $b_{\varepsilon(0)} < b_{\varepsilon_1} <
b_{\varepsilon(1)}$ which holds.]
\mr
\item "{$\circledast_2$}"  $b_{0,2} <_I b_1$.
\ermn
[Why?  Otherwise $b_1 \le_I b_{0,2}$ hence $\varepsilon \in (0,\partial)
\Rightarrow b_\varepsilon <_I b_{0,\varepsilon +1} <_I b_{\varepsilon +2}$ (by 
$\circledast_0 + \circledast_1)$ so $\langle b_{0,\varepsilon}:
\varepsilon \in (1,\partial)\rangle$ is also 
an increasing sequence unbounded in 
$I_{u,x}$ contradiction to ``$w_1$ maximal".]
\mr
\item "{$\circledast_3$}"  $b_{0,2} < b_{1,2}$.
\ermn
[Why?  By $\circledast_0 + \circledast_2$ we have $b_{0,4} < b_1$ 
and by $\circledast_1$ we have $b_1 < b_{2,4}$ together $b_{0,4} <
b_{2,4}$ so by $\circledast_0$ we have $b_{0,2} < b_{1,2}$.]

But then $\langle b_{\varepsilon,\partial}:\varepsilon < \partial 
\rangle$ increases (by $\circledast_3 + \circledast_0$) 
and $\varepsilon < \partial \Rightarrow b_\varepsilon =
b_{\varepsilon,\varepsilon} <  b_{\varepsilon +1,\partial} <
 b_{\varepsilon +2}$ (by $\circledast_1$ and $\circledast_2$
 respectively) hence is an unbounded subset of $I_{u,x}$ contradiction to 
the maximality of $|w_1|$.
\bn
\ub{Case 2}:  $m \backslash w_1 = \{0,\dotsc,m-1\}
 \backslash w_1$ is $\{i(*)\}$.
\bn
\ub{Subcase 2A}:  For some $i < m,i \ne i(*)$ and $j := k_i \notin w_2$.

Choose such $i$ with $|i-i(*)|$ maximal.  For any $s$ let
$t_{\varepsilon,\ell,s}$ be $t_{\varepsilon,\ell}$ if 
$\ell \ne j$ and be $s$ if $\ell=j$.

Let $I' = \{s \in I^2_{u,t_{\varepsilon(*),j}}:s,t_{\varepsilon(*),j}$
realize the same cut of $\{t_{\varepsilon,\ell}:\varepsilon < 
\partial,\ell \ne j\}\}$, note that $k_{j+1} = k_j +1$.  Recalling
$\boxtimes_2(b)$, the cofinality of $I_{2,< t_{\varepsilon(*),j}}$
is $\aleph_0$ and also the cofinality of the inverse of 
$I_{2,> t_{\varepsilon(*),j}}$ is $\aleph_0$ recalling the choice
of $\langle(s^-_\iota,s^+_\iota):\iota < m\rangle$ there is an open
interval\footnote{if we allow $+
\infty,-\infty$ as end points} of $I_2$ around
$t_{\varepsilon(*),j}$ which is $\subseteq I'$.  Note that $I'$
is dense in itself and has neither first nor last member by
$\boxtimes_2 + \boxtimes_4(a),(b)$.

As $f_*$ is definable, by the choice of $M_0$ and $M$ and of $I'
\subseteq I^2_{u,t_{\varepsilon(*),j}}$ we have: 
if $\varepsilon < \partial \wedge s \in I'$ then $t_{\varepsilon(*),j}$
and $s$ realize the same cut of $I^2_u \cup
\{t_{\varepsilon,\ell}:\varepsilon < \partial,j \ne \ell\}$ hence
$f^M_*(\ldots,a_{t_{\varepsilon,\ell,s}},\ldots)_{\ell < n},b_\varepsilon$ 
realize the same cut of $I_u$ which means that
$f_*(\ldots,a_{t_{\varepsilon,\ell,s}},\ldots)_{\ell <
n} \in I_{u,x}$ hence by the choice of $\langle
b_\varepsilon:\varepsilon < \partial\rangle$ we have $(\exists \xi <
\partial)(f_*(\ldots,a_{t_{\varepsilon,\ell,s}},\ldots) < b_\xi)$.

So again by the definability (and indiscernibility)
\mr
\item "{$\circledast_4$}"  $\varepsilon < \partial \wedge s \in I' \Rightarrow
f^M_*(\ldots,a_{t_{\varepsilon,\ell,s}},\ldots) < b_{\varepsilon +1}$.
\ermn
As $I'$ is dense in itself, what we say on the pair 
$(s,t_{\varepsilon(*),j})$ when $s \in I' \wedge s <_{I_2}
t_{\varepsilon(*),j}$ holds for the pair $(t_{\varepsilon(*),j},s)$
when $s \in I' \wedge t_{\varepsilon(*),j} <_I s$ so
\mr
\item "{$\circledast_5$}"   $\varepsilon < \partial 
\wedge s \in I' \Rightarrow b_\varepsilon < 
f^M_*(\ldots,a_{t_{\varepsilon +1,\ell,s}},\ldots)$
\ermn
(more fully let $s_1 <_{I_2} t_{\varepsilon(*),j} <_{I_2} s_2$ and $s_1,s_2
\in I'$ then the sequences $\langle t_{\varepsilon,\ell}:\ell \ne
j,\ell < n(*)\rangle \char 94 \langle s_1\rangle \char 94 \langle
t_{\varepsilon +1,\ell}:\ell \ne j,\ell <n(*)\rangle \char 94 \langle
t_{\varepsilon(*),j}\rangle$ and $\langle t_{\varepsilon,\ell}:\ell
\ne j,\ell < n(*)\rangle \char 94 \langle t_{\varepsilon(*),j}\rangle
\char 94 \langle t_{\varepsilon +1,\ell}:\ell \ne j,\ell < n(*)\rangle
\char 94 \langle s_2\rangle$ realizes the same quantifier free type in
$I_2$, (recalling $t_{\varepsilon,j} = t_{\varepsilon(*),j}$).

By $\circledast_4 + \circledast_5$ and indiscernibiity
we can replace $t_{\varepsilon(*),j}$ by any
$t' \in I'$ which realizes the same cut as $t_{\varepsilon(*),j}$ of
$\{t_{\varepsilon,\ell}:\varepsilon < \partial,\ell \ne j\}$.  
But if $j>i(*)$ then $\{t^*_{j+1},\dotsc,t^*_{n(*)-1}\} \subseteq I^2_u$
by the choice of $j$, and the set $I'' = \{t \in I_2$: 
if $\varepsilon < \partial,\ell \ne j$ then
$t \ne t_{\varepsilon,\ell}$ and $t_{\varepsilon,\ell} 
<_{I_2} t \equiv t_{\varepsilon,\ell} <_{I_2} t^*_j\}$ include an
initial segment of $J^+_{I^2_u,t_{\varepsilon(*),j}}$, see
$\boxtimes_4(d)$, i.e.  $(*)_6$ so its inverse
has cofinality $\aleph_0$, say $\langle s^*_n:n < \omega\rangle$
exemplifies this, so $n < \omega \Rightarrow s^*_{n+1} <_{I_2} s^*_n$.
So for every $\varepsilon < \partial$ for some $n <
\omega,f^M_*(\ldots,a_{t_{\varepsilon +1,\ell},s^*_n},\ldots) \in
(b_\varepsilon,b_{\varepsilon +1})_I$.  So for some $n_* < \omega$
this holds for unboundedly many $\varepsilon < \partial$, 
contradictory to ``$|w_2|$ is maximal".  Similarly if $j < i(*)$.
\mn
\ub{Subcase 2B}:  For every $\varepsilon < \partial$ for some $\xi \in
(\varepsilon,\partial)$, the interval of $I_2$ which is defined by
$t_{\varepsilon,k_{i(*)}},t_{\xi,k_{i(*)}}$ is not disjoint to
$I^2_u$ [so \wilog \, has $\ge k_{i(*)+1} - k_{i(*)}$ members of $I^2_u$].

In this case as in case 1, \wilog \, $\{k_{i(*)},\dotsc,k_{i(*)+1}\}
\subseteq w_2$ so as $|w_2|$ is maximal this
holds.  So as not subcase 2A,  
$\{t_{\varepsilon,\ell}:\varepsilon < \partial,\ell < n\} \subseteq
I^2_u$ hence $\{b_\varepsilon:\varepsilon < \partial\} \subseteq I_u$,
contradiction.
\mn
\ub{Subcase 2C}:  None of the above.

As not subcase(2B), \wilog \, 
$\{t_{\varepsilon,\ell}:\varepsilon < \partial$ and $\ell \in
[k_{i(*)},k_{i(*)+1})\} \subseteq I^2_{u,t_{\varepsilon(*),k_{i(*)}}}$. Then as
in subcase(2A) the sequence 
$\langle t_{\varepsilon,k_{i(*)}}:\varepsilon < \partial\rangle$ is
increasing/decreasing and is unbounded from above/below in
$I^2_{u,t_{\varepsilon(*),k_{i(*)}}}$ contradiction to $(*)_6$.  

In more detail, so $I' := I^2_{u,t_{0,k_{i(*)}}}$ includes all
$\{t_{\varepsilon,\ell}:\varepsilon < \partial$ and $\ell \in
[k_{i(*)},k_{i(*)+1})\}$.  Also $I'$ and its inverse are of 
cofinality $\aleph_0$ by $(*)_6$ hence \wilog \,  we can find (new) $\langle
t_{\partial,\ell}:\ell \in [k_{i(*)},k_{i(*)+1})\rangle$ such that
$t_{\partial,\ell} <_{I_2} t_{\partial,\ell+1},t_{\partial,\ell} \in
(s^-_{i(*)},s^+_{i(*)})_{I_2}$ and $\varepsilon < \partial \Rightarrow
t_{\varepsilon,\ell_1} <_{I_2}
t_{\partial,\ell} \equiv t_{\varepsilon,\ell_1} < 
t_{\varepsilon +1,\ell_2}$ and
the convex hull in $I_2$ of $\{t_{\zeta,\ell}:\zeta \le \partial$ and
$\ell \in [k_{i(*)},k_{i(*)+1}]\}$ is disjoint to $I^2_u$.  Let
$t_{\partial,\ell} = t_{\partial,\ell}$ for $\ell \notin
[k_{i(*)},k_{i(*)+1}],\ell < m,b_\partial =
f_*(a_{t_{\partial,0}},\dotsc,a_{t_{\partial,n-1}})$.  

Easily $\varepsilon < \partial \Rightarrow b_\varepsilon 
<_I b_\partial$.  As $\varepsilon < \xi < \partial 
\Rightarrow (b_\varepsilon,b_\xi)_{I_2} \cap u = \emptyset$
easily $\varepsilon < \partial \Rightarrow
(b_\varepsilon,b_\partial)_{I_2} \cap u = 0$, 
contradiction to $\langle b_\varepsilon:\varepsilon 
< \partial \rangle$ being cofinal in $I_{u,x}$.  

To finish proving $(*)_7$, we have to consider the possibility that
applying $\boxtimes_6$, the exceptional case in $(b)(\beta)$ of
$\boxtimes_6$ occurs for some $i<m$ say for $i(*)$; see before $\odot$.

Also \wilog \, as $\partial \ge \theta_2$ then \wilog \, $\ell \in w_2
\Rightarrow t_{\varepsilon,\ell} = t_{\varepsilon(*),\ell}$ and for
each $\ell < n(*)$ we have $(\forall \varepsilon,\zeta <
\partial)(\forall s \in I^2_u)(s <_{I_2} t_{\varepsilon,\ell} \equiv s
<_{I_2} t_{\zeta,\ell})$.

Now we can define $\bar t^{\varepsilon,\xi} = \langle
t^{\varepsilon,\xi}_\ell:\ell < n(*)\rangle$ as in case 1 and prove
$\circledast_0 - \circledast_3$ there.

Clearly all members of
 $\{t_{\varepsilon,\ell}:\varepsilon < \partial,\ell \in
[k_{i(*)},k_{i(*)+2})\}$ realize the same cut of $I^2_u$ and we get
easy contradiction.

As we can use only $\langle t_{n(*),\varepsilon}:\varepsilon <
\partial\rangle$ and add to $f_*$ dummy variables, without loss of
generality $k_{i(*)+1} - k_{i(*)} = 
k_{i(*)+2} - k_{i(*)+1}$.  Let $J$ be $\{1,-1\} \times
\partial$ ordered by $(\ell_1,\varepsilon_1) <_J (\ell_2,\varepsilon_2)$
iff $\ell_1 = 1 \wedge \ell_2 = -1$ or $\ell_1= 1 = \ell_2 \wedge
\varepsilon_1 < \varepsilon_2$ or $\ell_1 = -1 = \ell_2 \wedge
\varepsilon_1 > \varepsilon_2$.

For $\iota \in J$ let $\iota = (\ell^\iota,\varepsilon^\iota) =
(\ell[\iota],\varepsilon[\iota])$.  For $\zeta < \partial$ and
$\iota_1,\iota_2 \in J$ we define $\bar t_{\zeta,\iota_1,\iota_2} =
\langle t_{\zeta,\iota_1,\iota_2,n}:n < n(*)\rangle$ by
$t_{\zeta,\iota_1,\iota_2,n}$ is $t_{\varepsilon[\iota_1],n}$ if $n \in
[k_{i(*)},k_{i(*)+1}),t_{\varepsilon[\iota_2],n}$ if $n \in
[k_{i(*)+1},k_{i(*)+2})$ and $t_{\zeta,n}$ otherwise.  Now letting
$b_{\zeta,\iota_1,\iota_2} := f_*(\bar t_{\zeta,\iota_1,\iota_2})$
\mr
\item "{$\circledast_6$}"  all $b_{\zeta,\iota_1,\iota_2}$ realize
the ame cut of $I^2_u$.
\ermn
Now
\mr
\item "{$\circledast_7$}"  indiscernibility as in $\circledast_0$
holds
\sn
\item "{$\circledast_8$}"
$\neg(b_{\zeta,(1,\varepsilon),(1,\varepsilon +1)} \le_{I_*}
b_{\zeta,(1,\varepsilon +2),(1 + \varepsilon + 3)})$.
\ermn
[Why?  Otherwise by indiscernibility, if $\zeta \in (6,\partial)$ then
$b_{\zeta,(1,\zeta),(-1,3)} <_I b_{\zeta,(-1,5),(-1,4)}$.  Hence
$\langle b_{\zeta,(-1,5),(-1,4)}:\zeta \in (6,\partial)\rangle$ is
monotonic in $I_*$, all members realizing the fix cut of $I^2_u$ and
is unbounded in it (by the inequality above) so contradiction to
maximality of $|w_j|$.]
\mr
\item "{$\circledast_9$}"  $\neg(b_{\zeta,(1,\varepsilon
+2),(1,\varepsilon +3)} <_I b_{\zeta,(1,\varepsilon),(1,\varepsilon
+1)})$.
\ermn
[Similarly, as otherwise if $\zeta \in (6,\partial)$ then
$b_{\zeta,(1,\zeta),(-1,\zeta)} <_I b_{\zeta,(1,4),(1,5)})$.  Hence
$\langle b_{\zeta,(1,4),(1,5)}:\zeta \in (6,\partial)\rangle$
contradict the maximality of $(w_1)$.]

So we have proved $(*)_7$
\mr
\item "{$(*)_8$}"  if $u$ is $\mu$-reasonable, $x \in I \backslash
I_u$ then cf$(I_{u,x}) = \aleph_0$.
\ermn
[Otherwise by $(*)_7$ it has a last element say $b =
f_*(a_{t_0},\dotsc,a_{t_{n-1}})$ where $t_0,\dotsc,t_{n-1} \in I_2$ and 
$f_*$ a definable function, \wilog \,  with $n$ minimal hence
$\{a_{t_0},\dotsc,a_{t_{n-1}}\}$ is transcendentally independent and with no
repetitions and $b$ is not algebraic over
$\{a_{t_0},\dotsc,a_{t_{n-1}}\} \backslash \{a_{t_\ell}\}$ for $\ell < n$.
So $\{t_0,\dotsc,t_{n-1}\} \nsubseteq I^2_u$ and let $\ell < n$ be such 
that $t_\ell \notin I^2_u$ hence there are $s_0 <_{I_2} s_1$ such that
$t_\ell \in (s_0,s_1)_{I_2}$ and $(s_0,s_1)_{I_2} \cap I^2_u 
= \emptyset$ (recall $\boxtimes_4$(a),(b) and 
$(*)_6$ about cofinality $\aleph_0$ and $I_2$ being
dense).  Also \wilog \, $\{t_0,\dotsc,t_{n-1}\} \cap (s_0,s_1)_{I_2} =
\{t_\ell\}$, now the function $c \mapsto f^M_*(a_{t_0},\dotsc,
a_{t_{\ell-1}},c,a_{t_{\ell+1}},\dotsc,a_{t_{n-1}})$ for $c \in
(a_{s_0},a_{s_1})_I$ is increasing or decreasing (cannot be constant by
the minimality on $n$ and the elimination of quantifiers for real
closed fields and the transcendental independence of
$\{t_0,\dotsc,t_{n-1}\}$).   So we can find $s'_0,s'_1$ such that $s_0 <_{I_2}
s'_0 <_{I_2} t_\ell <_{I_2} s'_1 <_{I_2} s_1$ such that $X :=
\{f^M_*(a_{t_0},\dotsc,a_{t_{\ell-1}},c,a_{t_{\ell
+1}},\dotsc,a_{t_{n-1}}):c \in (a_{s'_0},a_{s'_1})_I\}$ is
included in $I_{u,x}$.  Again as the function defined above is
monotonic on $(a_{s'_0},a_{s'_1})_I$ so for some 
value $b' \in (a_{s'_0},a_{s'_1})$ we have $b <_I b'$.  But $b$ is
last in $I_{u,x}$ by our assumption toward contradiction hence 
$(b,b')_{I_u} \cap I_u = \emptyset$.  But this is impossible as all members of
$\{f(a_{t_0},\dotsc,a_{t_{\ell -1}},c,a_{t_{\ell +1}},\dotsc,
a_{t_{n-1}}):c \in (a_{s'_1},a_{s'_2})_I\}$ realize the same cut of $I_u$
so $(*)_8$ holds.]
\mr
\item "{$(*)_9$}"  if $u$ is $\mu$-reasonable, $x \in I \backslash
I_u$ then also the inverse of $I_{u,x}$ has cofinality $\aleph_0$.
\ermn
[Why?  Similarly to the proof of $(*)_7 + (*)_8$ or note that the
mapping $y \mapsto -y$ (defined in $M$) maps $I_u$ onto itself and is
an isomorphism from $I$ onto its inverse.]
\mr
\item "{$(*)_{10}$}"  if $u$ is $\mu$-reasonable, then $I_u$ is unbounded in
$I$ from below and from above.
\ermn
[Why?  Easy.]
\mr
\item "{$(*)_{11}$}"  if $h,u_1,u_2$ are as in clauses (a),(b),(c)
below \ub{then} the function $h_4$ defined below is (well defined and)
is, recalling $(*)_2(g)$, an order preserving function from 
$\hat I_{u_1}$ onto $\hat I_{u_2}$ mapping $u_1$ onto $u_2$ and also
the functions $h_0,h_1,h^*_2,h_2,h_3$ are as stated where
{\roster
\itemitem{ $(a)$ }  $u_1,u_2 \subseteq \zeta$ are $\mu$-reasonable
\sn
\itemitem{ $(b)$ }  $h$ is an order preserving function
from $u_1$ onto $u_2$
\sn
\itemitem{ $(c)$ }   $(\alpha) \quad$ for $\alpha \in u_1$, we have 
cf$(\alpha) \ge \theta_1 \Leftrightarrow \text{ cf}(h(\alpha)) \ge \theta_1$
\sn
\itemitem{ ${{}}$ }  $(\beta) \quad$ if $\gamma \in u_1$ then
$(\forall \alpha < \gamma)(\exists \beta \in u_1)(\alpha \le \beta <
\gamma)$ iff $(\forall \alpha < h(\gamma))(\exists \beta \in
u_2)(\alpha \le \beta < h(\gamma)))$
\sn
\itemitem{ $(d)$ }  $(\alpha) \quad h_1$ is the induced 
order preserving function from $I^1_{u_1}$ onto \nl

\hskip25pt $I^1_{u_2}$, i.e., $h_1((\ell,\beta')) = (\ell,\beta'')$
when $h(\beta') = \beta'' < \zeta$ or
\nl

\hskip25pt $\beta' = \beta'' \in [\zeta,\zeta + \omega)$;
\sn
\itemitem{ ${{}}$ }  $(\beta) \quad$ let $h_0$ be the 
partial function from $\zeta + \omega$ into $\zeta + \omega$ 
\nl

\hskip35pt  such that $h_0(\alpha) = \beta \Leftrightarrow (\exists \ell)
[h_1((\ell,\alpha)) = (\ell,\beta)]$
\sn
\itemitem{ $(e)$ }   $h^*_2$ is the order preserving function from
$I^{*,2}_{u_1}$ onto $I^{*,2}_{u_2}$ defined by: for 
$\eta \in I^{*,2}_{\eta_1},h^*_2(\eta) = \langle
h_1(\eta(\ell)):\ell < \ell g(\eta)\rangle = 
\langle(\ell^{\eta(\ell)},h_0(\alpha^{\eta(\ell)})):\ell < \ell
g(\eta)\rangle$, recalling (d)
\sn
\itemitem{ $(f)$ }   $h_2 = h^*_2 \restriction I^2_{u_1}$ is an order
preserving function from $I^2_{u_1}$ onto $I^2_{u_2}$
\sn
\itemitem{ $(g)$ }   $h_3$ is the unique isomorphism from the real
closed field $M_{I^2_{u_1}}$ onto 
\nl

\hskip25pt the real closed field $M_{I^2_{u_2}}$ mapping $a_t$ to
$a_{h_2(t)}$ for $t \in I^2_{u_1}$, 
\nl

\hskip25pt where for $I' \subseteq I_2$ we 
let $M_{I'} \subseteq M$ be the real closure
\nl

\hskip25pt  of $\{a_t:t \in I'\}$ inside $M$
\sn
\itemitem{ $(h)$ }   $h_4$ is the map defined by:
\nl
$h_4(x)=y$ iff
\sn
\itemitem{ ${{}}$ }   $(\alpha) \quad x \in I_{u_1} \wedge y = h_3(x)$
or
\sn
\itemitem{ ${{}}$ }   $(\beta)\quad$ for some $a \in I \backslash
I_{u_1},b \in I \backslash I_{u_2}$ we have $x = I_{u,a},y \in
I_{u,b}$ and \nl

\hskip35pt $(\forall c \in I_u)(c <_I a \equiv h_3(c) <_I b)$
\sn
\itemitem{ $(i)$ }   $\hat I_{u_1} = \text{ Dom}(h_4)$ and $\hat
I_{u_2} = \text{ Rang}(h_4)$ ordered naturally.
\endroster}
\ermn
[Why?  Trivially $h_1$ is an order preserving function from
$I^1_{u_1}$ onto $I^1_{u_2}$.  Recall $I^{2,*}_{u_\ell} = \{\eta \in
I^*_2:\eta(\ell) \in I^1_{u_\ell}$ for $\ell < \ell g(\eta)\}$.
So obviously $h^*_2$ is an order preserving function from $I^{*,2}_{u_1}$ onto
$I^{*,2}_{u_2}$.  Now $h_2 = h^*_2 \restriction I^2_{u_1}$, but does
it map $I^2_{u_1}$ onto $I^2_{u_2}$? we have excluded some members of 
$I^{*,2}_{u_2}$ by $\circledast$ above.  
But by clauses (c) and (d)$(\alpha)$ of the assumption being
excluded/not excluded is preserved by the natural mapping, i.e.,
$h^*_2$ maps $I^2_{u_1}$ onto $I^2_{u_2}$ hence $h_2 = h^*_2
\restriction I^1_{u_1}$ is an isomorphism from $I^1_{u_1}$ onto
$I^1_{u_2}$.  Also by $(*)_1$ being the real closure of the ordered
field $M_0$, and the uniqueness of ``the real closure" $h_3$ is 
the unique isomorphism from the real closed field $M_{I^2_{u_1}}$
onto $M_{I^2_{u_2}}$ mapping $a_t$ to $a_{h_2(t)}$ for $t \in
I^2_{u_1}$.

Let $\langle({\Cal U}^1_\varepsilon,{\Cal
U}^2_\varepsilon):\varepsilon < \varepsilon^*\rangle$ list the pairs
$({\Cal U}_1,{\Cal U}_2)$ such that:
\mr
\item "{$\circledast_{10}$}"  $(a) \quad {\Cal U}_\ell$ has the form
$I_{u_\ell,x}$ for some $x \in I \backslash I_{u_\ell}$ for $\ell =1,2$
\sn
\item "{${{}}$}"  $(b) \quad$ for every $a \in I_{u_1},(\exists y \in
{\Cal U}_1)(a <_I y) \Leftrightarrow (\exists y \in {\Cal U}_2)(h_2(a) <_I y)$.
\ermn
Now
\mr
\item "{$\circledast_{11}$}"  $\langle {\Cal U}^\ell_\varepsilon:
\varepsilon < \varepsilon^*\rangle$ is a partition
of $I \backslash I_{u_\ell}$ for $\ell=1,2$.
\ermn
[Why?  First, note the parallel claim for $I_1$.  For this note that
$h_1((\ell,0)) = (\ell,0)$ as $0 \in u_1 \cap u_2$ as $u_1,u_2$ are
$\mu$-reasonable, see clause (e) of $(*)_4$ and $h_1((\ell,\alpha)) =
(\ell,\beta) \Leftrightarrow h_1((\ell,\alpha +1)) = (\ell,\beta +1)$,
by clause (b) of $(*)_4$ and if $h((\ell,\delta_1)) =
(\ell,\delta_2),\delta_1$ is a limit (equivalently $\delta_2$ is
limit) then

$$
\delta_1 = \sup\{\alpha < \delta:(\ell,\alpha) \in I^1_{u_1}\}
\Leftrightarrow \delta_2 = \sup\{\alpha < \delta:(\ell,\alpha) \in
I^1_{u_2}\}.
$$
\mn
Second, note the parallel claim for $h_2,I^{*,2}_{u_\ell},h^*_2$.

Third, note the parallel claim for $I^2_{u_\ell},h_2$.

Fourth, note the parallel claim for $I_{u_\ell},h_3$ (which is the
required one).]

So it follows that
\mr
\item "{$\circledast_{12}$}"  $h_4$ is as promised.
\ermn
So we are done proving $(*)_{11}$.
\nl
[Why?  By clauses (b),(c) of $(*)_{11}$.]
\mr
\item "{$(*)_{12}$}"  if $u_1,u_2$ are $\mu$-reasonable, $h$ is an
order preserving mapping from $\hat I_{u_1}$ onto $\hat I_{u_2}$ which
maps $I_{u_1}$ onto $I_{u_2}$ then there is an automorphism $h^+$ of
the linear order $I$ extending $h \restriction I_{u_1}$.
\ermn
[Why?  Let $\langle {\Cal U}^1_\varepsilon:\varepsilon <
\varepsilon^*\rangle$ list $\hat I_{u_1} \backslash I_{u_1}$ and
${\Cal U}^2_\varepsilon = h({\Cal U}^1_\varepsilon)$.  Now for every
$\varepsilon$ we choose $\langle a^\ell_{\varepsilon,n}:n \in \Bbb Z
\rangle$ such that
\mr
\item "{$\circledast_{13}$}"  $(a) \quad a^\ell_{\varepsilon,n} \in {\Cal
U}^\ell_\varepsilon$
\sn
\item "{${{}}$}"  $(b) \quad a^\ell_{\varepsilon,n} <_I
a^\ell_{\varepsilon,n+1}$ for $n \in \Bbb Z$
\sn
\item "{${{}}$}"  $(c) \quad \{a^\ell_{\varepsilon,n}:n \in \Bbb Z,n
\ge 0\}$ is unbounded from  above in ${\Cal U}^\ell_\varepsilon$
\sn
\item "{${{}}$}"  $(d) \quad \{a^\ell_{\varepsilon,n}:n \in \Bbb Z,n
< 0\}$ is unbounded from below in ${\Cal U}^\ell_\varepsilon$.
\ermn
This is justified by $u_\ell$ being $\mu$-reasonable by $(*)_6$, $\boxtimes_4$.
Now define $h_5:I \rightarrow I$ by:
\sn
\block
$h_5(x) = h_4(x)$ \ub{if} $x \in I_{u_1}$ and otherwise
\nl
$h_5(x) = a^2_{\varepsilon,n} + (a^2_{\varepsilon,n+1} -
a^2_{\varepsilon,n})(x - a^1_{\varepsilon,n})/(a^1_{\varepsilon,n+1} -
a^1_{\varepsilon,n})$
\endblock
\nl
\hskip55pt \ub{if} $a^1_{\varepsilon,n} 
\le_{I_2} x < a^1_{\varepsilon,n+1}$ and $n \in \Bbb Z$.
\sn
Now check using linear algebra.]
\mr
\item "{$(*)_{13}$}"  $({}^\mu I)/{\Cal E}^{\text{aut}}_{I,\mu}$ has $\le
2^\mu$ members recalling that $f_1 {\Cal E}^{\text{aut}}_{I,h} f_2$ 
\ub{iff} $f_1,f_2$ are functions \nl

\hskip25pt from $\mu$ into $I$ and for some automorphism $h$ of $I$ we have
\nl

\hskip25pt $(\forall \alpha < \mu)(h(f_1(\alpha)) = f_2(\alpha))$
\nl

\hskip25pt [Why?  Should be clear recalling $|I^1_u| \le \mu$,
recalling $(*)_5,(*)_{11},(*)_{12}$.]
\ermn
So we have finished proving part (1) of \scite{734-am3.2.3}.
\nl
2) Really the proof is included in the proof of part (1).  That is,
given $I' \subseteq I$ of cardinality $< \theta_2$ 
by $(*)_5$ there is a $\mu$-reasonable $u \subseteq \zeta$
such that $I' \subseteq I_u$ and $|u| = \mu + |I'|$.  Now clearly
\mr
\item "{$(*)_{14}$}"   for $\mu$-reasonable $u \subseteq \zeta$, the
family $\{I^2_{u,x}:x \in I_2 \backslash I^2_u\}$ has $\le \mu + |u|$
members.
\ermn
[Why?  By $\boxtimes_5$.]
\mr
\item "{$(*)_{15}$}"  for a $\mu$-reasonable $u \subseteq \zeta$, the
family $\{I_{u,x}:x \in I \backslash I_u\}$ has $\le \mu$ members.
\ermn
[Why?  By $(*)_{16}$ below.]
\mr
\item "{$(*)_{16}$}"  if $u$ is $\mu$-reasonable then
$I_{u,b_1} = I_{u,b_2}$ when
{\roster
\itemitem{ $(a)$ }  $b_k = f(a_{t_{k,0}},\dotsc,a_{t_{k,n-1}})$ for
$k=1,2$
\sn
\itemitem{ $(b)$ }  $f$ a definable function in $M$
\sn
\itemitem{ $(c)$ }  $t_{k,0} <_{I_2} \ldots <_{I_2} t_{k,n-1}$ for
$k=1,2$
\sn
\itemitem{ $(d)$ }  $t_{1,\ell} \in I^2_u \vee t_{2,\ell} \in I^2_u
\Rightarrow t_{1,\ell} = t_{2,\ell}$
\sn
\itemitem{ $(e)$ }  if $t_{1,\ell} \notin I^2_u$ then
$I^2_{u,t_{1,\ell}} = I^2_{u,t_2,\ell}$ for $\ell=0,\dotsc,n-1$.
\endroster}
\ermn
[Why?  Use the proof of $(*)_{11}$, for $u_1 = u=u_2,h = \text{\rm
id}_{u_2}$ so ${\Cal U}^1_\varepsilon = {\Cal U}^2_\varepsilon$ for
$\varepsilon < \varepsilon^*$.

By the assumptions for each $\ell$ there is $\varepsilon$ such that
$a_{t_{\varepsilon,1,\ell}},a_{t_{2,\ell}} \in {\Cal U}^1_\varepsilon =
{\Cal U}^2_\varepsilon$.  Now for each $\varepsilon < \varepsilon^*$
there is an automorphism $\pi_\varepsilon$ of ${\Cal U}^1_\varepsilon$
as a linear order mapping $t_{1,\ell}$ to $t_{2,\ell}$ if $t_{1,\ell}
\in {\Cal U}^1_\varepsilon$.  Let $\pi =
\cup\{\pi_\varepsilon:\varepsilon < \varepsilon^*\} \cup \text{ id}_{I_u}$.]
\mr
\item "{$(*)_{17}$}"  if $n < \omega,t^\ell_0 <_I t^\ell_1 <_I \ldots
<_I t^\ell_{n-1}$ for $\ell=1,2,I_{u,t^1_k} = I_{u,t^2_k}$ for
$k=0,1,\dotsc,n-1$ \ub{then} for some automorphism $g$ of $I$ over
$I_u$ we have $k<n \Rightarrow g(t^1_k) = t^2_k$.
\ermn
[Why?  We shall use $g$ such that $g \restriction I_u = \text{
id}_{I_u}$ and $g \restriction I_{u,x}$ is an automorphism of
$I_{u,x}$ for each $x \in I \backslash I_u$.  Clearly it suffices to
deal with the case $\{t^\ell_k:\ell < n$ and $\ell \in \{1,n\}\}
\subseteq I_{u,x}$ for one $x \in I \backslash I_u$.  We choose $s_1 <
s_2$ from $I_{u,x}$ such that $s_1 <_I t^\ell_k < s_2$ for $\ell=1,2$.
We choose $g \restriction I_{u,x}$ such that it is the identity on
$\{s \in I_{u,x}:s \le_I s_1$ or $s_2 \le_I s$, now 
stipulates $t_{-1} = s_1,t_n =
s_2$ and maps $(t^1_k,t^1_{k+1})_I$ onto $(t^2_k,t^2_{k+1})_I$ for
$k=-1,0,\dotsc,n-1$ as in the definition above.]
\nl
So we have completed the proof of part (2) of \scite{734-am3.2.3}.
\nl
3) Obvious from the Definition (\scite{734-0n.2}(9)) and the construction.
\nl
4) First
\mr
\item "{$\odot_1$}"    there is $J^*_1 \subseteq I$ of cardinality
$\mu^+$ such that: for every $J^*_2 \subseteq I$ of cardinality $\le
\mu$ there is an automorphism $\pi$ of $I$ which maps $J^*_2$ into
$J^*_1$.
\ermn
[Why?  Let $u = \mu^+ \times \mu^+ \subseteq \zeta$ and let $J^*_1 
= I_u$.  Clearly $u$ has cardinality $\mu^+$ and so does $J^*_1 =
I_u$.  So suppose $J^*_2 \subseteq I$  has cardinality $\le \mu$.
There is $u_2 \subseteq \zeta$ of cardinality $\mu$ such that $J^*_2
\subseteq I_{u_2}$ and \wilog \, $u_2$ is reasonable.  We define an
increasing function $h$ from $u_2$ into $u_1$, by defining $h(\alpha)$
by induction on $\alpha$:
\mr
\item "{$(*)_{17}$}"  if cf$(\alpha) \le \mu$ then $h(\alpha) =
\cup\{h(\beta)+1:\beta \in u_2 \cap \alpha\}$
\sn
\item "{$(*)_{18}$}"  if cf$(\alpha) > \mu$ then $h(\alpha) =
\cup\{h(\beta)+1:\beta \in u_2 \cap \alpha\} + \mu^+$.
\ermn
Let $u_1 := \{h(\alpha):\alpha \in u_2\}$ so $u_1 \subseteq u$.  Now
$h,u_1,u_2$ satisfies clauses (a),(b),(c) of $(*)_{11}$ hence
$h_1,h^*_2,h_2,h_3,h_4,\hat I_{u_1},\hat I_{u_2}$ are as there.

By $(*)_{12}$ there is an isomorphism
$h^+$ of $I$ which extends $h_4$; now does $h^+$ map $J^*_2$
into $J^*_1$?  Yes, as $J^*_2 \subseteq I_{u_2}$ and
$h^+ \restriction I_{u_2}$ is an isomorphism from $I_{u_2}$ onto
$I_{u_1}$ but $I_{u_1} \subseteq I_u,I_u = J^*_1$, so we are
done proving $\odot_1$.]

Finally
\mr
\item "{$\odot_2$}"   part (4) of \scite{734-am3.2.3} holds, i.e. if
$I^*_0 \subseteq I,|I^*_0| < \theta_2$ \ub{then} for some $I^*_1
\subseteq I$ of cardinality $\le \mu^+ + |I^*_0|$ we have: for every $J
\subseteq I$ of cardinality $\le \mu$ there is an automorphism of $I$
over $I^*_0$ mapping $J$ into $I^*_1$.
\ermn
Why?  Given $I^*_0 \subseteq I$ of cardinality $< \theta_2$ we can find
$u_1 \subseteq \zeta$ of cardinality $\mu + |I^*_0|$ such
that $I^*_0 \subseteq I_{u_1}$.  By $(*)_5$ we can find a
$\mu$-reasonable set $u_2 \subseteq \zeta$ of cardinality $\mu +
|u_1|$ such that $u_1 \subseteq u_2$.

Let $\langle {\Cal U}_\varepsilon:\varepsilon < \varepsilon^*\rangle$
list the sets of the form $I_{u_2,x},x \in I_2 \backslash I_{u_1}$, so
by ($\boxdot_5$) $\varepsilon^* \le \mu + |I^*_0|$.  For each $\varepsilon$
we choose $\langle a_{\varepsilon,n}:n \in \Bbb Z\rangle$ as in
$\circledast_{13}$ from the proof of $(*)_{12}$.  For each $\varepsilon <
\varepsilon^*$ and $n \in \Bbb Z$ let $\pi_{\varepsilon,n}$ be an
isomorphism from $I$ onto 
$(a_{\varepsilon,n},a_{\varepsilon,n+1})_I$, exists by the properties
of ordered fields.  Let $J^*_1 \subseteq I$
be as in $\odot_1$ above and let $I^*_2 = I^*_1
\cup\{a_{\varepsilon,n}:\varepsilon < \varepsilon^*$ and $n < \omega\}
\cup \{\pi_{\varepsilon,n}(J^*_1):\varepsilon < \varepsilon^*$ and $n
\in \Bbb Z\}$.  Easily, $I^*_2$ is as required.
\nl
5) By \scite{734-0n.em.7}.   \hfill$\square_{\scite{734-am3.2.3}}$
\bigskip

\remark{Remark}  Concerning $(*)_{11}$, we could have used more time
\mr
\item "{$(*)'_{11}$}"  $h_2$ is an order preserving function from
$I^2_{u_1}$ onto $I^2_{u_2}$ and $h_3$ is an isomorphism from
$I_{u_1}$ onto $I_{u_2}$ and $h_1$ is an order preserving mapping from
$\hat I_{u_2}$ onto $\hat I_{u_2}$.
\endroster
\endremark
\goodbreak

\head {\S6 Linear orders and equivalence relations} \endhead  \resetall \sectno=6
 \spuriousreset
\bigskip

This section deals with a relative of the stability spectrum.  
We ask: what can be the number of equivalence classes in ${}^\mu I$
for an equivalence realtion on ${}^\mu I$ which is 
so called ``invariant", in fact definable 
(essentially by a quantifier
free infinitary formula, mainly for well ordered $I$).

It is done in a very restricted context, but via EM-models has useful
conclusions, for a.e.c. and also for a.e.c. with amalgamation; i.e. it is
used in \scite{734-am2.18}.
\nl
There are two versions; one for well ordering and one for the class of
linear orders both expanded  by unary relations.
\nl
On $\tau^*_{\alpha(*)},K^{\text{lin}}_{\tau^*_{\alpha(*)}}$ see
\scite{734-0n.2}(4).  We may replace sequences, i.e. inc$_J(I)$ by subsets
of $I$ of cardinality $|J|$, this may help to eliminate $2^{|J|}$
later, but at present it seems not to help in the final 
bounds in \S7.  We do here only enough for \S7.
\bigskip

\demo{\stag{734-lin.0} Context}  We fix $\alpha(*),\bar u^* = (u^-,u^+)$ such
that
\mr
\item "{$(a)$}"   $\alpha(*)$ is an ordinal $\ge 1$
\sn
\item "{$(b)$}"   $u^- \subseteq \alpha(*)$
\sn
\item "{$(c)$}"   $u^+ \subseteq \alpha(*)$.
\endroster
\enddemo
\bigskip

\remark{\stag{734-lin.0.3} Remark}  1) The main cases are
\mr
\item "{$(A)$}"  $\alpha(*) = 1$, so
$K^{\text{lin}}_{\tau^*_{\alpha(*)}}$ is the class of linear orders
\sn
\item "{$(B)$}"  $\alpha(*) = 2,u^+ = \emptyset,u^- = \{0\}$.
\ermn
2) Usually the choice of the parameters does not matter.
\endremark
\bigskip

\definition{\stag{734-lin.1} Definition}  1) For $I,J \in
K^{\text{lin}}_{\tau^*_{\alpha(*)}}$, i.e. both linear orders expanded
by a partition $P_\alpha(\alpha < \alpha(*))$, pedantically the
interpretation of the $P_\alpha$'s, let inc$'_J(I)$
 be the set of embedding of $J$ into $I$; see below, we 
denote members by $h$.
\nl
2) Recalling $\bar u^* = (u^-,u^+)$ where $u^- \cup u^+ \subseteq
\alpha(*)$ let inc$^{\bar u^*}_J(I)$ be the set of $h$ such that
\mr
\item "{$(a)$}"  $h$ is an embedding of $J$ into $I$, i.e. one-to-one,
order preserving function mapping $P^J_\alpha$ into $P^I_\alpha$ for $\alpha <
\alpha(*)$ 
\sn
\item "{$(b)$}"  if $\alpha \in u^-$ and $t \in P^J_\alpha$ and $s <_I
h(t)$ \ub{then} for some $t_1 <_J t$ we have $s \le_I h(t_1)$
\sn
\item "{$(c)$}"  if $\alpha \in u^+$ and $t \in P^J_\alpha$ and $h(t)
<_I s$ \ub{then} for some $t_1$ we have $t <_J t_1$ and $h(t_1) \le_I s$.
\nl
\endroster
\enddefinition
\bn
Concerning $\bar u^*$
\demo{\stag{734-lin.1.5} Observation}  1) For any $h \in 
\text{ inc}^{\bar u^*}_J(I)$
\mr
\item "{$(a)$}"  if $t$ is the successor of $s$ in $J$ (i.e. $s <_J t$
and $(s,t)_J = \emptyset$) and $t \in P^J_\alpha,\alpha \in u^-$ \ub{then}
$h(t)$ is the successor of $h(s)$ in $I$ 
\sn
\item "{$(b)$}"  if $\langle t_i:i < \delta\rangle$ is
$<_J$-increasing with limit $t_\delta \in J$ (i.e. $i < \delta
\Rightarrow t_i <_J t_\delta$ and $\emptyset =
\cap\{(t_i,t_\delta)_J:i < \delta\}$) and $t_\delta \in P^J_\alpha,\alpha \in
u^-$ \ub{then} $\langle h(t_i):i < \delta\rangle$ is $<_I$-increasing
with limit $h(t_\delta)$ in $I$
\sn
\item "{$(c)$}"  if $t$ is the first member of $J$ and $t \in
P^J_\alpha,\alpha \in u^-$ \ub{then} $h(t)$ is the first member of $I$.
\ermn
2) If $h_1,h_2 \in \text{ inc}^{\bar u^*}_J(I)$ \ub{then}
\mr
\item "{$(a)$}"  if $t$ is the successor of $s$ in $J$ and
$t \in P^J_\alpha,\alpha
\in u^-$ \ub{then} $h_1(s) = h_2(s) \Leftrightarrow h_1(t) = h_2(t)$
and $h_1(s) <_I h_2(s) \Leftrightarrow h_1(t) <_I h_2(t)$ and $h_1(s) >_I
h_2(s) \Leftrightarrow h_1(t) >_I h_2(t)$
\sn
\item "{$(b)$}"  if $\langle t_i:i < \delta\rangle$ is
$<_J$-increasing with limit $t_\delta$
and $t_\delta \in P^J_\alpha,\alpha \in u^-$, \ub{then} $(\forall i <
\delta)(h_1(t_i) = h_2(t_i)) \Rightarrow h_1(t_\delta) =
h_2(t_\delta)$ moreover $(\forall i < \delta)(\exists j <\delta)
(h_1(t_i) <_I h_2(t_j) \wedge h_2(t_i) <_I h_1(t_j)) \Rightarrow 
h_1(t_\delta) = h_2(t_\delta)$ and also $(\exists j < \delta)(\forall
i < \delta)(h_1(t_i) <_I h_2(t_j)) \Rightarrow h_1(t_\delta) <_I 
h_2(t_\delta)$.
\ermn
3) Similar to parts (1) + (2) for $\alpha \in u^+$ (inverting the
orders of course).
\nl
4) inc$'_I(J) = \text{ inc}^{(\emptyset,\emptyset)}_I(J)$.
\enddemo
\bigskip

\demo{Proof}  Straight (and see the proof of \scite{734-lin.1.14}).  
\hfill$\square_{\scite{734-lin.1.5}}$
\enddemo
\bigskip

\demo{\stag{734-lin.1.7} Convention}  1) $\alpha(*),\bar u^*$ 
will be constant so usually we shall not
mention them, e.g. write inc$_J(I)$ for inc$^{\bar u^*}_I(I)$ and
pedantically below we should have written $\bold e^{\bar u^*}(J,I),\bold
e^{\bar u^*}_*(J)$ and also in notions like reasonable and wide in
Definition \scite{734-lin.1C} mention $\bar u^*$.
\nl
2) $I,J$ denote members of $K^{\text{lin}}_{\tau^*_{\alpha(*)}}$.
\enddemo
\bn
Below we use mainly ``$e$-pairs" (and weak $e$-pairs and the reasonable case).
\definition{\stag{734-lin.1A} Definition}  1)  let $\bold e(J)$ be the set
of equivalence relations on some subset of $J$ such that each
equivalence class is a convex subset of $J$.
\nl
2) For $h_1,h_2 \in \text{ inc}_J(I)$ we say that $(h_1,h_2)$ is a
strict $e$-pair (for $(I,J))$
 \ub{when} $e \in \bold e(J)$ and $(h_1,h_2)$ satisfies
\mr
\item "{$(a)$}"  $s \in J \backslash \text{ Dom}(e)$ iff $h_1(s) =
 h_2(s)$
\sn
\item "{$(b)$}"   if $s <_J t$ and $s/e \ne t/e$ (so $s,t \in \text{
Dom}(e)$) \ub{then} $h_1(s) <_I h_2(t)$ and $h_2(s) <_I h_1(t)$
\sn
\item "{$(c)$}"  if $s <_J t$ and $s/e = t/e$ (so $s,t \in \text{
Dom}(e)$) \ub{then} $h_1(t) <_I h_2(s)$.
\ermn
2A) We say that $(h_1,h_2)$ is a strict $(e,{\Cal Y})$-pair where $e
\in \bold e(J)$ and ${\Cal Y} \subseteq \text{\rm Dom}(e)/e$ \ub{when} clauses
(a)+(b) from part (2) hold and
\mr
\item "{$(c)'$}"  if $s <_J t$ and $s/e=t/e$ (so $s,t \in \text{\rm
Dom}(e))$ \ub{then} $(h_1(t) <_I h_2(s)) \equiv (s/e \in {\Cal Y}) \equiv
(h_1(s) < h_2(t))$.
\ermn
2B)  We say that $(h_1,h_2)$ is an $e$-pair when $(h_1,h_2)$ is a
strict $(e,{\Cal Y})$-pair for some ${\Cal Y}$ (this relation is
symmetric, see below).
\nl
3) We say that $(h_1,h_2)$ is a weak $e$-pair where 
$h_1,h_2 \in \text{\rm inc}_J(I)$ \ub{when} clauses (a),(b)
hold (this, too, is symmetric!)
\nl
4) For $h_1,h_2 \in \text{ inc}_J(I)$, let $e = \bold e(h_1,h_2)$ 
be the (unique) $e \in \bold e(J)$ such that (see \scite{734-lin.1B}(1) below)
\mr
\item "{$(a)$}"  Dom$(e) = \{s \in J:h_1(s) \ne h_2(s)\}$
\sn
\item "{$(b)$}"   $(h_1,h_2)$ is a weak $e$-pair
\sn
\item "{$(c)$}"  if $e' \in \bold e(J)$ and $(h_1,h_2)$ is a weak $e'$-pair
\ub{then} Dom$(e) \subseteq \text{\rm Dom}(e')$ and $e$ refines $e'
\restriction \text{ Dom}(e)$.
\ermn
5) If $e \in \bold e(J)$ and ${\Cal Y} \subseteq \text{\rm Dom}(e)/e$
\ub{then} we let set$({\Cal Y})  = \{s \in J:s/e \in {\Cal Y}\}$ 
and $e \restriction {\Cal Y} =e \restriction  \text{\rm set}({\Cal Y})$.
\nl
6) Let $\bold e(J,I)$ be the set of $e \in \bold e(J)$ such that there
is an $e$-pair.
\nl
7) Let $\bold e_*(J) = \cup\{\bold e(J,I):I \in
 K^{\text{lin}}_{\tau^*_{\alpha(*)}}\}$. 
\enddefinition
\bn
Concerning $\bar u^*$
\demo{\stag{734-lin.1.14} Observation}  Assume that $e \in \bold e(J,I)$.
\nl
0) 
\mr
\item "{$(a)$}"  If $t$ is the first member of $J$ and $t \in
P^J_\alpha,\alpha \in u^-$ \ub{then} $t \notin \text{ Dom}(e)$.
\sn
\item "{$(b)$}"  If $t \in$ Dom$(e)$ and $t$ is the first member of
$t/e$ and $t \in P^J_\alpha$ \ub{then} $\alpha \notin u^-$.
\ermn
1) If $t$ is the $<_J$-successor of $s$ and $t \in P^J_\alpha,\alpha \in
 u^-$ \ub{then} $s \in \text{ Dom}(e) \Leftrightarrow t \in 
\text{ Dom}(e)$ and $s \in \text{ Dom}(e) \Rightarrow s \in t/e$.
\nl
2) If $\langle t_i:i < \delta\rangle$ is $<_J$-increasing with limit
$t_\delta$ and $t_\delta \in P^J_\alpha$ and $\alpha \in u^-$ \ub{then}:
\mr
\item "{$(a)$}"  if $(\forall i < \delta)(t_i \notin \text{ Dom}(e))$
\ub{then} $t_\delta \notin \text{ Dom}(e)$
\sn
\item "{$(b)$}"  if $(\forall i < \delta)(\neg t_i e t_{i+1})$
or just $(\forall i < \delta)(\exists j < \delta)(i<j \wedge \neg t_i
e t_j)$ \ub{then} $t_\delta \notin \text{ Dom}(e)$
\sn
\item "{$(c)$}"  if $(\forall i < \delta)(t_i \in t_0/e)$ \ub{then}
$t_\delta \in t_0/e$.
\ermn
3) Similar to parts (0),(1),(2) when $\alpha \in u^+$ (inverting the
   order, of course).
\nl
4) $\bold e_*(J)$ is the family of $e \in \bold e(J)$ satisfying the
requirements in parts (0),(1),(2),(3) above so if $\bar u^* =
   (\emptyset,\emptyset)$ then $\bold e_*(J) = \bold e(J)$.
\enddemo
\bigskip

\demo{Proof}  Easy by \scite{734-lin.1.5}, e.g.
\mn
\ub{Part (1)}:  We are assuming $e \in \bold e(J,I)$ hence by
Definition \scite{734-lin.1A} there is an $e$-pair $(h_1,h_2)$ where
$h_1,h_2 \in \text{ inc}_J(I)$.  Now for $\ell=1,2$, 
clearly $h_\ell(s),h_\ell(t) \in I$ and as $s <_J t$ we
have $h_\ell(s) < h_\ell(t)$.  Now if $h_\ell(t)$ is not the
$<_I$-successor of $h_\ell(s)$ then there is $s'_\ell \in
(h_\ell(s),h_\ell(t))_I$ hence by clause (b) of Definition
\scite{734-lin.1}(2) there is $s^*_\ell \in [s,t)_J$ such that $s'_\ell
\le_I h_\ell(s^*_\ell) <_I h_\ell(t)$ so as $h_\ell(s) <_I s'_\ell$ we
have $h_\ell(s) <_I h_\ell(s^*_\ell) <_I h_\ell(t)$ hence $s <_I
s^*_\ell <_J t$, contradiction to the assumption ``$t$ is the
successor of $s$ in $J$".  So indeed $h_\ell(t)$ is the successor of
$h_\ell(s)$ in $I$.

As this holds for $\ell=1,2$, clearly $h_1(s) = h_2(s) \Leftrightarrow
h_1(t) = h_2(t)$ but by Definition \scite{734-lin.1}(2) we know $s \in
\text{ Dom}(e) \Leftrightarrow (h_1(s) \ne h_2(s))$ and similarly for
$t$ hence $s \in \text{ Dom}(e) \Leftrightarrow t \in \text{ Dom}(e)$.
Lastly,  assume $s,t \in \text{ Dom}(e)$, but $s,t$ are nor $e$-equivalent
so by Definition \scite{734-lin.1A}(2) clause (b) we have $h_1(s) <_I h_2(t)
\wedge h_2(s) <_I h_1(t)$ clear contradiction.
\mn
\ub{Part 2}: We leave clauses (a),(b) to the reader.

 For clause (c) of part (2), 
if $t_\delta \notin t_0/e$ then choose $h_1,h_2 \in 
\text{ inc}^{\bar u^*}_J(I)$ such that $(h_1,h_2)$ is an 
$e$-pair, hence an $(e,{\Cal Y})$-pair for some
${\Cal Y} \subseteq \text{ Dom}(e)/e$.  If $(t_0/e) \in {\Cal Y}$ then
$h_2(t_0)$ is above $\{h_1(t_i):i < \delta\}$ by $<_I$ so we have
$h_1(t_\delta) \le_I h_2(t_0)$ but if $t_\delta \notin t_0/e$ this contradicts
clause (b) in Definition \scite{734-lin.1A}(2),(2A).  The proof when
$t_0/e \notin {\Cal Y}$ is similar.  \hfill$\square_{\scite{734-lin.1.14}}$
\enddemo
\bigskip

\demo{\stag{734-lin.1B} Observation}  Let $h_1,h_2 \in \text{ inc}_J(I)$
and $e \in \bold e(J)$.
\nl
1) $\bold e(h_1,h_2)$ is well defined.
\nl
2) $(h_1,h_2)$ is a strict $(e,{\Cal Y}_1)$-pair iff $(h_2,h_1)$ is a
strict $(e,{\Cal Y}_2)$-pair when $({\Cal Y}_1,{\Cal Y}_2)$ is a
partition of Dom$(e)/e$.
\nl
3) $(h_1,h_2)$ is a strict $e$-pair \ub{iff} 
$(h_2,h_1)$ is a strict $(e,\emptyset)$-pair.
\nl
4) $(h_1,h_2)$ is an $e$-pair \ub{iff} $(h_2,h_1)$ is an $e$-pair.
\nl
5) $(h_1,h_2)$ is a weak $e$-pair iff $(h_2,h_1)$ is a weak $e$-pair.
\nl
6) If $(h_1,h_2)$ is a strict $e$-pair \ub{then} $(h_1,h_2)$
is an $e$-pair which implies $(h_1,h_2)$ being a weak $e$-pair.
\nl
7) If $e_\alpha \in  \bold e(J)$ for $\alpha < \alpha^*$, \ub{then} $e
:= \cap\{e_\alpha:\alpha < \alpha^*\} = \{(s,t):s,t$ are
$e_\alpha$-equivalent for every $\alpha < \alpha^*\}$ 
belongs to $\bold e(J)$ with Dom$(e) =
\cap\{\text{Dom}(e_\alpha):\alpha < \alpha^*\}$.  
\nl
8) If $e \in \bold e(J,I)$ \ub{then} for every ${\Cal Y} \subseteq 
\text{ Dom}(e)/e$ also $e \restriction \text{ set}({\Cal Y})$ belongs to 
$\bold e(J,I)$ and there is 
a strict $(e \restriction \text{ set}({\Cal Y}))$-pair $(h'_1,h'_2)$;
moreover, for every ${\Cal Y}_1 \subseteq {\Cal Y}$ there is a strict $(e
\restriction \text{ set}({\Cal Y}),{\Cal Y}_1)$-pair.  
  \hfill$\square_{\scite{734-lin.1B}}$
\enddemo
\bigskip

\demo{Proof}  Easy, e.g.:
\nl
1) Let

$$
\align
e = \{(s_1,s_2):&h_1(s_\ell) \ne h_2(s_\ell) \text{ for } \ell=1,2  
\text{ and if } s_1 \ne s_2 \text{ then} \\
  &\text{ for some } t_1 <_J t_2 \text{ we have }
\{s_1,s_2\} = \{t_1,t_2\} \\
  &\text{ and there is no initial segment } J' \text{ of } J \text{
  such that} \\
  &J' \cap \{t_1,t_2\} = \{t_1\} \text{ and} \\
  &(\forall t' \in J')(\forall t'' \in J \backslash J')
[h_1(t') <_I h_2(t'') \wedge h_2(t') <_I h_1(t'')]\}.
\endalign
$$
\mn
Clearly $e$ is an equivalence relation on $\{t \in J:h_1(t) \ne
h_2(t)\}$ and each equivalence class is convex hence $e_1 \in \bold
e(J)$, so clauses (a),(b) of \scite{734-lin.1A}(1),(4) holds.  
Easily $e$ is as required.
\nl
8) Let $(h_1,h_2)$ be an $e$-pair and ${\Cal Y}_1,{\Cal Y}_2,{\Cal Y}_3$ be a
partition of Dom$(e)/e$.  We define $h'_1,h'_2 \in 
\text{ inc}_J(I)$ as follows, for $\ell \in \{1,2\}$
\mr
\item "{$(a)$}"  if $t \in J \backslash \text{ Dom}(e)$ then
$h'_\ell(t) = h_1(t) \, (= h_2(t))$
\sn
\item "{$(b)$}"  if $t \in \text{ set}({\Cal Y}_1)$ then
$h'_\ell(t) = h_1(t)$
\sn
\item "{$(c)$}"  if $t \in \text{ set}({\Cal Y}_2)$ then $h'_\ell(t)$
is min$\{h_1(t),h_2(t)\}$ if $\ell=1$, and is max$\{h_1(t),h_2(t)\}$ if
$\ell=2$
\sn
\item "{$(d)$}"  if $t \in \text{ set}({\Cal Y}_3)$ then $h'_\ell(t)$
is max$\{h_1(t),h_2(t)\}$ if $\ell=1$ and is min$\{h_1(t),h_2(t)\}$ if
$\ell=2$. 
\ermn
Now $(h'_1,h'_2)$ is a strict $(e \restriction (\text{set}({\Cal Y}_2) \cup
\text{ set}({\Cal Y}_3)),{\Cal Y}_2)$-pair, so we are done.
\hfill$\square_{\scite{734-lin.1B}}$
\enddemo
\bigskip

\definition{\stag{734-lin.6} Definition}  1) For a subset $u$ of $J \in
K^{\text{lin}}_{\tau^*_{\alpha(*)}}$ we define 
$e = e_{J,u} \in \bold e(J)$ on $J \backslash u$ as follows:

$$
s_1 e s_2 \text{ iff } (\forall t \in u)(t <_J s_1 \equiv t <_J s_2).
$$
\mn
2) For $I,J \in K^{\text{lin}}_{\alpha(*)}$, we say that the pair
$(I,J)$ is non-trivial \ub{when}: $\bold e(J,I) \ne \emptyset$.  
\enddefinition
\bigskip

\definition{\stag{734-lin.1C} Definition}  1) For $h_0,\dotsc,h_{n-1} \in
\text{ inc}_J(I)$ let

$$
\qftp^J(\langle h_0,\dotsc,h_{n-1}\rangle,I) 
= \{(\ell,m,s,t):s,t \in J
\text{ and } h_\ell(s) < h_m(t)\}.
$$
\mn
We may write $\qftp^J(h_0,\dotsc,h_{n-1};I)$ and we usually omit $J$
as it is clear from the context.
\nl
2) For $h_1,h_2 \in \text{ inc}_J(I)$ let eq$(h_1,h_2) = \{s \in J:
h_1(s) = h_2(s)\}$.
\nl
3) We say that the pair $(I,J)$ is a reasonable $(\mu,\alpha(*)))$-base
\ub{when}:
\mr
\item "{$(a)$}"   $I,J \in K^{\text{lin}}_{\tau^*_{\alpha(*)}},
|J| \le \mu$ and the pair $(I,J)$
is non-trivial
\sn
\item "{$(b)$}"   if $e \in \bold e(J,I)$  and $h_1,h_2 \in
\text{\rm inc}_J(I)$ and $(h_1,h_2)$ is an $e$-pair then we can
find $h'_1,h'_2,h'_3 \in \text{\rm inc}_J(I)$ and 
${\Cal Y} \subseteq \text{\rm Dom}(e)/e$ such that
{\roster
\itemitem{ $(\alpha)$ }    $\qftp((h'_1,h'_2),I) = \qftp((h_1,h_2),I)$
\sn
\itemitem{ $(\beta)$ }    $(h'_1,h'_3)$ and $(h'_2,h'_3)$ 
are strict $(e,{\Cal Y})$-pairs.
\endroster}
\ermn
4) We say that the pair $(I,J)$ is a wide
$(\lambda,\mu,\alpha(*))$-base \ub{when}:
\mr
\item "{$(a)$}"   $I,J \in K^{\text{lin}}_{\tau^*_{\alpha(*)}},|J| \le
\mu$ and the pair $(I,J)$ is non-trivial
\sn
\item "{$(b)$}"  for every $e \in \bold e(J,I)$ there is a sequence 
$\bar h = \langle h_\alpha:\alpha < \lambda \rangle$ such that 
{\roster
\itemitem{ $(\alpha)$ }  $h_\alpha$ is an embedding of $J$ into $I$
\sn
\itemitem{ $(\beta)$ }  if $\alpha < \beta < \lambda$ then
$(h_\alpha,h_\beta)$ is an $e$-pair.
\endroster}
\ermn
5) We say that the pair $(I,J)$ is a
strongly wide $(\lambda,\mu,\alpha(*))$-base \ub{when}:
\mr
\item "{$(a)$}"  $I,J \in K^{\text{lin}}_{\tau^*_{\alpha(*)}}$, the
pair $(I,J)$ is non-trivial and $J$ has cardinality $\le \mu$
\sn
\item "{$(b)$}"  for every $e \in \bold e(J,I)$ and ${\Cal Y} \subseteq
\text{\rm Dom}(e)/e$ there is $\bar h = \langle h_\alpha:\alpha <
\lambda \rangle$ such that
{\roster
\itemitem{ $(\alpha)$ }  $h_\alpha \in \text{\rm inc}_J(I)$
\sn  
\itemitem{ $(\beta)$ }  if $\alpha < \beta$ then
$(h_\alpha,h_\beta)$ is a strict $(e,{\Cal Y})$-pair.
\endroster}
\ermn
6) Above we may omit $\mu$ meaning $\mu = |J|$ and we may omit $\alpha(*)$,
 as it is determined by $J$ (and by $I$), and then may omit ``base" so
in part (3) we say $(I,J)$ is reasonable and in part (4) we say
$\lambda$-wide and in part (5) say strongly $\lambda$-wide.
\enddefinition
\bigskip

\demo{\stag{734-lin.1F} Observation}  1) If $(I,J)$ is a 
reasonable $(\mu,\alpha(*))$-base \ub{then} $(I,J)$ is a reasonable
$(\mu',\alpha(*))$-base for $\mu' \ge \mu$.
\nl
2) If $(I,J)$ is a wide $(\lambda,\mu,\alpha(*))$-base and $\lambda'
\le \lambda,\mu' \ge \mu$ \ub{then} $(I,J)$ is a wide
$(\lambda',\mu',\alpha(*))$-base.
\nl
3) If $(I,J)$ is a strongly wide $(\lambda,\mu,\alpha(*))$-base,
\ub{then} $(I,J)$ is a wide $(\lambda,\mu,\alpha(*))$-base.
\enddemo
\bigskip

\demo{Proof}  Obvious.  \hfill$\square_{\scite{734-lin.1F}}$
\enddemo
\bigskip

\proclaim{\stag{734-lin.2} Claim}  1) If $\alpha(*) =1$ and 
$\mu \le \zeta(*) < \mu^+ \le \lambda$, \ub{then}
the pair $(\lambda \times \zeta(*),\zeta(*))$ is a reasonable
$(\mu,\alpha(*))$-based which is a wide $(\lambda,\mu,\alpha(*))$-base. 
\nl
2) If $\alpha(*)=2$ and $\bar u^* = (\{0\},\emptyset)$ as in
\scite{734-lin.0.3} and $\mu \le \zeta(*) < \mu^+ < \lambda$ and
$\zeta'(*) = \zeta(*) \times 3$ and $w \subseteq \zeta(*),w \ne \zeta(*)$
\ub{then} the pair $(I^{\text{lin}}_{\mu,\lambda \times
\zeta(*)},I^{\text{lin}}_{\mu,\zeta(*),w})$ is a reasonable
$(\mu,\alpha(*))$-base which is a wide $(\lambda,\mu,\alpha(*))$-base where
\mr
\item "{$(*)$}"  for any ordinal $\beta$ and $w \subseteq \beta$ we
define $I = I^{\text{lin}}_{\mu,\beta,w}$, a $\tau^*_{\alpha(*)}$-model
(if $w = \emptyset$ we may omit it)
{\roster
\itemitem{ $(\alpha)$ }  its universe is $\beta$
\sn
\itemitem{ $(\beta)$ }  the order is the usual one
\sn
\itemitem{ $(\gamma)$ }  $P^I_1 = \{\alpha < \beta:\text{\rm cf}(\alpha) >
\mu$ or $\alpha \in w\}$, (if we write
$I^{\text{lin}}_{\ge\mu,\beta,w}$ we mean here {\rm cf}$(\alpha) \ge \mu$).
\endroster}
\endroster
\endproclaim
\bigskip

\demo{Proof}  1) First: $(I,J)$ = 
\ub{$(\lambda \times \zeta(*),\zeta(*))$ is a  
wide $(\lambda,\mu,\alpha(*))$-base}

Easily $\bold e(J,I) \ne \emptyset,|J| \le \mu$ and $I,J \in
K^{\text{lin}}_{\tau^*_{\alpha(*)}}$ so clause (a) of Definition
\scite{734-lin.1C}(4) holds (recalling Definition \scite{734-lin.6}(2)), so it
suffices to deal with clause (b). 

Let $e \in \bold e(J,I)$ and define

$$
\align
u = \{\zeta < \zeta(*):&\zeta \in \text{ Dom}(e)
\text{ is minimal in } \zeta/e \\
  &\text{ or } \zeta \in \zeta(*) \backslash \text{ Dom}(e)\}.
\endalign
$$
\mn
Now for every $\alpha < \lambda$ we define $h_\alpha \in \text{
inc}_J(I)$ as follows:
\mr
\item "{$(a)$}"  if $\zeta \in \zeta(*) \backslash \text{ Dom}(e)$
then $h_\alpha(\zeta) = \lambda \times \zeta$
\sn
\item "{$(b)$}"  if $\zeta \in \text{ Dom}(e)$ and $\varepsilon =
\text{ min}(\zeta/e)$ then $h_\alpha(\zeta) = \lambda \times
\varepsilon + \zeta(*) \times \alpha + \zeta$.
\ermn
Second: $(I,J)$ = \ub{$(\lambda \times \zeta(*),\zeta(*))$ is a 
reasonable $(\mu,\alpha(*))$-base}

Again clause (a) of Definition \scite{734-lin.1C}(3) holds so we deal with
clause (b).

So assume $e \in \bold e(J,I)$ and $h_1,h_2 \in \text{ inc}_J(I)$ and
$(h_1,h_2)$ is just a weak $e$-pair and ${\Cal Y} \subseteq \text{
Dom}(e)/e$.  Let $u = \text{ Rang}(h_1) \cup \text{ Rang}(h_2)$.
For $\ell=1,2$ let $h^*_\ell \in \text{ inc}_J(I)$ be
$h^*_\ell(\zeta) = \text{ otp}(u \cap h_\ell(\zeta))$, so Rang$(h^*_\ell)
\subseteq \xi(*) := \text{ otp}(u) \le \zeta(*) \times 3$.
\nl
[Why?  If $\zeta(*)$ is finite this is trivial, so assume $\zeta(*)
\ge \omega$.  Let $n < \omega$ and $\alpha$ be such that $\omega^\alpha n
\le \zeta(*) < \omega^\alpha(n+1)$, so $\alpha \ge 1,n \ge 1$.  As
$\omega^\alpha$ is additively indecomposable otp$(u) \le 
\omega^\alpha(2n+1)$, alternatively 
use natural sums \cite{MiRa65} which gives a better bound $\zeta(*)
\oplus \zeta(*)$, [actually $< \mu^+$ 
sufices using $\zeta(*) < \mu^+$ large enough below, still.]

For $\ell=1,2,3$ we define $h'_\ell \in \text{ inc}_J(I)$ as follows:
\mr
\item "{$(a)$}"  if $\zeta \in \zeta(*) \backslash \text{ Dom}(e)$
then $h'_\ell(\zeta) = (\zeta(*) \times 4) \times \zeta$
\sn
\item "{$(b)$}"  if $\zeta \in \text{ Dom}(e)$ and $\varepsilon =
\text{ min}(\zeta/e)$ and $\zeta/e \in {\Cal Y}$ then
{\roster
\itemitem{ $(\alpha)$ }   if $\ell=3$ then $h'_\ell(\zeta) = (\zeta(*)
\times 4) \times \varepsilon + \zeta(*) \times 3 + \zeta$
\sn  
\itemitem{ $(\beta)$ }  if $\ell=1,2$ then $h'_\ell(\zeta) = (\zeta(*)
\times 4) \times \varepsilon + h^*_\ell(\zeta)$
\endroster}  
\item "{$(c)$}"  if $\zeta \in \text{ Dom}(e)$ and $\varepsilon =
\text{ min}(\zeta/e)$ and $\zeta/e \notin {\Cal Y}$ then
{\roster
\itemitem{ $(\alpha)$ }  if $\ell = 3$ then $h'_\ell(\zeta) =
(\zeta(*) \times 4) \times \varepsilon + \zeta$
\sn 
\itemitem{ $(\beta)$ }  if $\ell=1,2$ then $h'_\ell(\zeta) = (\zeta(*)
\times 4) \times \varepsilon + \zeta(*) + h^*_\ell(\zeta)$.
\endroster}  
\ermn
Now check.
\nl
2) \ub{First: $(I,J) = (I^{\text{lin}}_{\mu,\lambda \times
\zeta(*)},I^{\text{lin}}_{\mu,\zeta(*),w})$ is a wide
$(\lambda,\mu,\alpha(*))$-base}.

Note that $P^J_1 = w$ because $\zeta(*) < \mu^+$ and
$P^I_1 = \{\alpha \in I:\text{ cf}(\alpha) >
\mu\}$.  As above clause (a) of the Definition \scite{734-lin.1C} holds so
we deal with clause (b).

Let

$$
u = \{\zeta < \zeta(*):\zeta \in \text{ Dom}(e) 
\text{ is minimal in } \zeta/e \text{ or } \zeta 
\in \zeta(*) \backslash \text{ Dom}(e)\}.
$$
\mn
Clearly $u$ is a closed subset of $\zeta(*)$ and $0 \in u$.

Given $\zeta < \zeta(*)$ let
$\varepsilon_\zeta := \text{ max}(u \cap (\zeta +1))$, clearly well
defined by the choice of $u$ and $\varepsilon_\zeta \le \zeta$.

For every $\alpha < \lambda$ we define $h_\alpha \in 
\text{ inc}_J(I)$ as follows: 

We define $h_\alpha(\zeta)$ by induction on $\zeta < \zeta(*)$ such
that $h_\alpha(\zeta) < \lambda \times (\varepsilon_\zeta +1)$.
\bn
\ub{Case A}:  for $\zeta \in \zeta(*) \backslash \text{ Dom}(e)$
\mn
\hskip20pt \ub{Subcase A1}:  $\zeta \in P^J_1$

\hskip35pt Let $h_\alpha(\zeta)$ be $\lambda \times \varepsilon_\zeta
+ \mu^+$.
\mn
\hskip20pt \ub{Subcase A2}:  $\zeta \in P^J_0$ and $\zeta = 0$

\hskip35pt Let $h_\alpha(\zeta) = 0$.
\mn
\hskip20pt \ub{Subcase A3}:  $\zeta \in P^J_0,\zeta = \xi+1$

\hskip35pt Let $h_\alpha(\zeta) = h_\alpha(\xi) +1$.
\mn
\hskip20pt \ub{Subcase A4}:  $\zeta \in P^J_0,\zeta$ is a limit ordinal,
$\zeta = \sup(u \cap \zeta)$

\hskip35pt Let $h_\alpha(\zeta) = \lambda \times \varepsilon_\zeta$
which is equal to $\cup\{h_\alpha(\zeta'):\zeta' < \zeta\}$.
\mn
\hskip20pt \ub{Subcase A5}:  $\zeta \in P^J_0,\zeta$ is a limit ordinal
and $\xi = \sup(u \cap \zeta) < \zeta$.

\hskip35pt So $(\xi +1)/e$ is an end-segment of $\zeta$, but this is
impossible by \scite{734-lin.1.14}(2)(c).
\bn
\ub{Case B}: $\zeta \in \text{ Dom}(e)$:
\mn
\hskip20pt \ub{Subcase A1}:  $\zeta = \text{ min}(\zeta/e)$ hence
$\zeta \in P^J_1$ (see \scite{734-lin.1.14}(0)(b))

\hskip35pt Let $h_\alpha(\zeta) = \lambda \times \varepsilon_\zeta +
\mu^+ \times \zeta(*) \times \alpha + \mu^+$.
\mn
\hskip20pt \ub{Subcase A2}:  $\zeta \in P^J_0$ hence $\zeta >
\text{ min}(\zeta/e)$

\hskip35pt Let $h_\alpha(\zeta)= \cup\{h_\alpha(\zeta') +1:\zeta' < \zeta\}$.
\mn
\hskip20pt \ub{Subcase A3}:  $\zeta \in P^J_1$ and $\zeta >
\text{ min}(\zeta/e)$

\hskip35pt Let $h_\alpha(\zeta) = \cup\{h_\alpha(\zeta'):\zeta' <
\zeta\} + \mu^+$.

So clearly we can show by induction on $\zeta < \zeta(*)$ that:

$$
h_\alpha(\zeta) < \lambda \times \varepsilon_\zeta + \mu^+ \times \zeta(*)
\times (\alpha 2+2).
$$
\mn
Now check.

Also recalling $\mu^+ < \lambda$ clearly for $\alpha < \lambda,\zeta <
\zeta(*)$ we have $h_\alpha(\zeta) < \lambda \times \varepsilon_\zeta
+ \lambda$.

Now check.  
\mn
\ub{Second $(I^{\text{lin}}_{\mu,\lambda \times
\zeta(*)},I^{\text{lin}}_{\mu,\zeta(*),w})$ 
is a reasonable $(\mu,\alpha(*))$-base}

Combine the proof of ``first" with the parallel proof in part (1).
  \hfill$\square_{\scite{734-lin.2}}$
\enddemo
\bigskip

\definition{\stag{734-lin.4} Definition}  1) Let $I,J \in
K^{\text{lin}}_{\tau^*_{\alpha(*)}}$.  We say that ${\Cal E}$ is an
invariant $(I,J)$-equivalence relation \ub{when}:
\mr
\item "{$(a)$}"  ${\Cal E}$ is an equivalence relation on inc$_J(I)$,
so ${\Cal E}$ determines $I$ and $J$
\sn
\item "{$(b)$}"   if $h_1,h_2,h_3,h_4 \in \text{ inc}_J(I)$ and
$\qftp(h_1,h_2;I) = \qftp(h_3,h_4;I)$ 
then $h_1 {\Cal E} h_2 \Leftrightarrow h_3 {\Cal E} h_4$.
\ermn
2) We add non-trivial \ub{when}:
\mr
\item "{$(c)$}"  if eq$(h_1,h_2) = \{t \in J:h_1(t) = h_2(t)\}$  
is co-finite then $h_1 {\Cal E} h_2$
\sn
\item "{$(d)$}"  there are $h_1,h_2 \in \text{ inc}_J(I)$ such that
$\neg(h_1 {\Cal E} h_2)$.
\ermn
3) Let $J,I_1,I_2 \in K^{\text{lin}}_{\tau^*_{\alpha(*)}}$.  \ub{Then} $I_1
\le^1_J I_2$ means that:
\mr 
\item "{$(a)$}"   $I_1 \subseteq I_2$
\sn
\item "{$(b)$}"   for every $h_1,h_2,h_3 \in \text{\rm inc}_J(I_2)$ we can find
$h'_1,h'_2,h'_3 \in \text{ inc}_J(I_1)$ such that
\qftp$(h'_1,h'_2,h'_3;I_1) = \text{ \qftp}(h_1,h_2,h_3;I_2)$. 
\endroster
\enddefinition
\bigskip

\proclaim{\stag{734-lin.5} Claim}  Assume $J,I_1,I_2 \in
K^{\text{lin}}_{\tau^*_{\alpha(*)}}$.
\nl
1) If $I_1 \subseteq I_2,{\Cal E}$ is an
invariant $(I_2,J)$-equivalence relation \ub{then} ${\Cal E} \restriction 
\text{\rm inc}_J(I_1)$ is an invariant $(I_1,J)$-equivalence relation.
\nl
2) If $I_1 <^1_J I_2$ and ${\Cal E}_1$ is an invariant $(I_1,J)$-equivalence
relation \ub{then} there is one and only one invariant
$(I_2,J)$-equivalence relation ${\Cal E}_2$ such that ${\Cal E}_2 \restriction 
\text{\rm inc}_J(I_1) = {\Cal E}_1$.
\nl
3) Assume $e \in \bold e(J)$ and ${\Cal Y} \subseteq \text{\rm
Dom}(e)/e$.  
If $(h'_1,h'_2)$ is a strict $(e,{\Cal Y})$-pair for 
$(I_1,J)$ and $(h''_1,h''_2)$ is
a strict $(e,{\Cal Y})$-pair for $(I_2,J)$ \ub{then} 
$\qftp(h'_1,h'_2;I_1) = \qftp(h''_1,h''_2;I_2)$.
\nl
4) Assume $\alpha(*) = 1,J = \zeta(*),I_\ell = \beta_\ell$ with the
usual order (for $\ell=1,2$), $\mu \le \zeta(*) < \mu^+$ and $\mu^+ \le
\beta_1 \le \beta_2$.  \ub{Then} $I_1 <^1_J I_2$ (see Definition 
\scite{734-lin.4}(3)).
\nl
5) Assume $\alpha(*)=2,J = I^{\text{lin}}_{\mu,\zeta(*),w},I_\ell =
I^{\text{lin}}_{\mu,\beta_\ell}$ for $\ell=1,2$ and $\mu^{++} \le 
\beta_1 \le \beta_2$.  \ub{Then} $I_1 <^1_J I_2$ (see Definition 
\scite{734-lin.4}(3)).
\endproclaim
\bigskip

\demo{Proof}  1) Obvious.
\nl
2)  We define

$$
\align
{\Cal E}^*_2 = \bigl\{(h_1,h_2):&h_1,h_2 \in \text{ inc}_J(I_2) \text{
and for some} \\
  &h'_1,h'_2 \in \text{ inc}_J(I_1) \text{ we have} \\
  &\qftp(h'_1,h'_2;I_1) = \qftp(h_1,h_2;I_2) \text{ and} \\
  &h'_1 {\Cal E}_1 h'_2\bigr\}.
\endalign
$$
\mn
Now
\mr
\item "{$(*)_1$}"  ${\Cal E}^*_2$ is a set of pairs of members of
inc$_J(I_2)$.
\ermn
[Why?  By its definition]
\mr
\item "{$(*)_2$}"  $h_1 {\Cal E}^*_2 h_1$ if $h_1 \in \text{ inc}_J(I_2)$.
\ermn
[Why?  Let $h' \in \text{ inc}_J(I_1)$ so clearly $h' {\Cal E}_1 h'$ and
$\qftp(h',h';I_1) = \qftp(h,h;I_2)$]
\mr
\item "{$(*)_3$}"  ${\Cal E}^*_2$ is symmetric.
\ermn
[Why?  As ${\Cal E}_1$ is.]
\mr
\item "{$(*)_4$}"  ${\Cal E}^*_2$ is transitive.
\ermn
[Why?  Assume $h_1 {\Cal E}^*_2 h_2$ and $h_2 {\Cal E}^*_2 h_3$ and
let $h'_1,h'_2 \in \text{ inc}_J(I_1)$ witness $h_1 {\Cal E}^*_2 h_2$
and $h''_2,h''_3 \in \text{ inc}_J(I_1)$ witness $h_2 {\Cal E}^*_2 h_3$.

Apply clause (b) of part (3) of Definition \scite{734-lin.4} 
to $(h_1,h_2,h_3)$ so there are $g_1,g_2,g_3 \in
\text{ inc}_J(I_1)$ such that $\qftp(g_1,g_2,g_3;I_1) =
\qftp(h_1,h_2,h_3;I_2)$.  Now $h'_1 {\Cal E}_1 h'_2$ by the
choice of $(h'_1,h'_2)$ and $\qftp(g_1,g_2;I_1) = 
\qftp(h_1,h_2;I_2) = \qftp(h'_1,h'_2;I_1)$ so as
${\Cal E}_1$ is invariant we get $g_1 {\Cal E}_1 g_2$.  Similarly $g_2
{\Cal E}_1 g_3$, so as ${\Cal E}_1$ is transitive we have $g_1 {\Cal
E}_1 g_3$.  But clearly $\qftp(g_1,g_3;I_1) = \qftp(h_1,h_3;I_2)$ 
hence $g_1,g_2$ witness that $h_1 {\Cal E}_2 h_3$ is as required.] 
\mr
\item "{$(*)_5$}"   ${\Cal E}^*_2$ is invariant.
\ermn
[Why?  See its definition.]  
\mr
\item "{$(*)_6$}"  ${\Cal E}^*_2 \restriction \text{ inc}_I(I_1) =
{\Cal E}_1$.
\ermn
[Why?  By the way ${\Cal E}^*_2$ is defined and ${\Cal E}_1$ being
invariant.]

So together ${\Cal E}^*_2$ is as required.  The uniqueness (i.e. if
${\Cal E}_2$ is an invariant equivalent relation on inc$_J(I)$ such
that ${\Cal E}_2 \restriction \text{ inc}_J(I_1) = {\Cal E}_1$ then
${\Cal E}_2 = {\Cal E}^*_2$) is also easy.
\nl
3) Straight. 
\nl
4) See \footnote{Actually instead ``$\mu^+ \le \beta_1$" it suffice to
have $\zeta(*) \times 4 \le \beta_1$ because if $\zeta(*) = \dsize
\sum_{i < \gamma} \zeta_i$ then $\dsize \sum_{i < \gamma} \zeta_i
\times 4 \le \zeta(*) \times 4$ or just the natural sum $\zeta(*)
\oplus \zeta(*) \oplus \zeta(*)$.} the proof of ``Second" in the proof
of \scite{734-lin.2}(1).
\nl
5) Combine \footnote{Here $(\mu^+ + 1) \times (\zeta(*) \times 4)$ will
suffice.} the proof of part (4) and of ``First" in the proof of
\scite{734-lin.2}(2).    \hfill$\square_{\scite{734-lin.5}}$
\enddemo
\bn
Below mostly it suffices to consider ${\Cal D}_{{\Cal E},e}$.
\definition{\stag{734-lin.7} Definition}  1) Let ${\Cal E}$ be an invariant
$(I,J)$-equivalence relation; we define

$$
\align
{\Cal D}_{\Cal E} = \{u \subseteq J:&\text{if } h_1,h_2 \in \text{\rm inc}_J(I)
\text{ satisfies eq}(h_1,h_2) \supseteq u \\
  &\text{then } h_1 {\Cal E} h_2\}
\endalign
$$
\mn
recalling

$$
\text{ eq}(h_1,h_2) := \{t \in J:h_1(t) = h_2(t)\}.
$$
\mn
2) If in addition $e \in \bold e(J,I)$ then we let

$$
\align
{\Cal D}_{{\Cal E},e} = 
\{u \subseteq \text{\rm Dom}(e)/e:&\text{if } h_1,h_2 
\in \text{\rm inc}_J(I) \text{ and } (h_1,h_2) \text{ is an } \\
  &(e \restriction (\text{\rm Dom}(e) \backslash 
\text{set}(u)))\text{-pair then } h_1 {\Cal E} h_2\}.
\endalign
$$
\mn
\enddefinition
\bigskip

\proclaim{\stag{734-lin.9} Claim}  Assume $I,J \in
K^{\text{lin}}_{\tau^*_{\alpha(*)}}$ and $(I,J)$ is reasonable (see
Definition \scite{734-lin.1C}(3),(6))
and ${\Cal E}$ is an invariant $(I,J)$-equivalence relation.
\nl
1) For $u \subseteq J$ such that $e_{J,u} \in \bold e(J,I)$
we have:
$u \in {\Cal D}_{\Cal E}$ \ub{iff} $h_1 {\Cal E} h_2$ for every $e_{J,u}$-pair
$(h_1,h_2)$ \ub{iff} $h_1 {\Cal E} h_2$ for some $e_{J,u}$-pair
$(h_1,h_2)$; see Definition \scite{734-lin.6}(1).
\nl
2) Assume $e \in \bold e(J,I)$, then for any $u \subseteq 
\text{\rm Dom}(e)/e$ we have: $u \in {\Cal D}_{{\Cal E},e}$ \ub{iff} 
$h_1 {\Cal E} h_2$ for any $(e \restriction$ {\rm set}$(u))$-pair 
\ub{iff} $h_1 {\Cal E} h_2$ for some ($e \restriction$ {\rm set}$(u))$-pair.
\nl
3) If $e \in \bold e(J,I)$ and $u_1,u_2 \subseteq \text{\rm Dom}(e)/e$
\ub{then} we can find $h_1,h_2,h_3 \in \text{\rm inc}_J(t)$ such that
$(h_1,h_2)$ is a strict $(e \restriction \text{\rm set}(u_1))$-pair,
$(h_2,h_3)$ is a strict $(e \restriction \text{\rm set}(u_2))$ pair and
$(h_1,h_3)$ is a strict $(e \restriction (\text{\rm set}(u_1 \cup
u_2))$-pair.
\nl
4) Assume $e \in \bold e(J,I)$ and that in clause (b) of Definition
\scite{734-lin.1C}(3) we allow $(h_1,h_2)$ to be a weak $e$-pair,
\ub{then} for any $u \subseteq \text{\rm Dom}(e)/e$ we have: {\rm
Dom}$(e) \backslash u \in {\Cal D}_{{\Cal E},e}$ 
iff $h_1 {\Cal E} h_2$ for every weak $e$-pair $(h_1,h_2)$.
\endproclaim
\bigskip

\demo{Proof} 1) Like part (2).
\nl
2) In short, by transitivity of equivalence and the definitions +
mixing, but we elaborate.
\nl
The ``first implies the second" holds by Definition
\scite{734-lin.7}(2) and ``the second implies the third" holds trivially
as there is such a pair $(h_1,h_2)$ by the assumption 
$e \in \bold e(J,I)$.   So it is 
enough to prove ``the third implies the first"; hence 
suppose that $g_1 {\Cal E} g_2$,
where $(g_1,g_2)$ is an $e_1 := e \restriction \text{\rm
set}(u)$-pair (recalling that $e_1 \in \bold e(J,I)$ by \scite{734-lin.1B}(8)), 
and let $(h_1,h_2)$ be an $e_1$-pair, we need to show that 
$h_1 {\Cal E} e_2$.  By Definition \scite{734-lin.1A}(2B) for
some sets ${\Cal Y}_g,{\Cal Y}_h \subseteq \text{\rm Dom}(e_1)/e_1$ the pair
$(g_1,g_2)$ is a strict $(e_1,{\Cal Y}_g)$-pair and the pair $(h_1,h_2)$ is a
strict $(e_1,{\Cal Y}_h)$-pair.  Recalling clause (b) of \scite{734-lin.1C}(3)
there are $g'_1,g'_2,g'_3$ and ${\Cal Y}$ such that:
\mr
\item "{$(*)_1$}"  $(a) \quad g'_\ell \in \text{\rm inc}_J(I)$ for $\ell=1,2,3$
\sn
\item "{${{}}$}"  $(b) \quad \qftp(g_1,g_2) = \qftp(g'_1,g'_2)$
\sn
\item "{${{}}$}"   $(c) \quad {\Cal Y} \subseteq 
\text{\rm Dom}(e_1)/e_1$
\sn
\item "{${{}}$}"   $(d) \quad (g'_1,g'_3)$ and $(g'_2,g'_3)$ are strict
$(e_1,{\Cal Y})$-pairs.
\ermn
Now for each $s \in \text{ Dom}(e_1)$, we can find a permutation $\bar
\ell_s = (\ell_{s,1},\ell_{s,2},\ell_{s,3})$ of $\{1,2,3\}$ such that
$I \models g'_{\ell_{s,1}}(s) < g'_{\ell_{s,2}}(s) <
g'_{\ell_{s,3}}(s)$.  By $(*)_1(d)$ and $(*)_1(b)$ and $(g_1,g)$ being
an $e_1$-pair, clearly $\bar \ell_s$ depends only on $s/e_1$
and every member of $\{(g'_{\ell_{s,1}}(t):t \in s/e_1\}$ is below
every member of $\{g'_{\ell_{s,2}}(t):t \in s/e_1\}$ and similarly
for the pair $(g'_{\ell_{s,2}},g'_{\ell_{s,3}})$.  Now
we can find $(g''_1,g''_2,g''_3)$ such that:
\mr
\item "{$(*)_2$}"   $(a) \quad g''_\ell \in \text{\rm inc}_J(I)$ for
$\ell=1,2,3$ 
\sn
\item "{${{}}$}"  $(b) \quad (g''_1,g''_2)$ is a strict
$(e_1,{\Cal Y}_h)$
\sn
\item "{${{}}$}"  $(c) \quad (g''_1,g''_3)$ and $(g''_2,g''_3)$ are
strict $(e_1,{\Cal Y}_g)$-pairs
\ermn
[Why?  We do the choice for each $s/e_1$ separately such that:
$\{g''_1 \restriction (s/e_1),g''_2 \restriction (s/e_1),g''_3
\restriction (s/e_1)\} = \{g'_1 \restriction (s/e_1),g'_2 \restriction
(s/e_1),g'_3 \restriction (s/e_1)\}$.]

Clearly $\qftp(g''_1,g''_3;I) = \qftp(g_1,g_2;I)
= \qftp(g''_2,g''_3;I)$ so as ${\Cal E}$ is invariant and 
$g_1 {\Cal E} g_2$ clearly $g''_1 {\Cal E} g''_3 \wedge 
g''_2 {\Cal E} g''_3$ which implies $g''_1 {\Cal E}
g''_2$.  For ${\Cal Y}' = {\Cal Y}_h$ by clause (b) of $(*)_2$ we
conclude that $\qftp(g''_1,g''_2;I) = \qftp(h_1,h_2;I)$ so as 
${\Cal E}$ is invariant we are done.
\nl
3),4)  Similarly.    \hfill$\square_{\scite{734-lin.9}}$
\enddemo
\bigskip

\proclaim{\stag{734-lin.10} Claim}  Assume $I,J \in
K^{\text{lin}}_{\tau^*_{\alpha(*)}}$ and ${\Cal E}$ is an invariant
$(I,J)$-equivalence relation.
\nl
0) If $e \in \bold e(J,I)$ and ${\Cal E}$ is non-trivial \ub{then}
${\Cal D}_{{\Cal E},e}$ contains all co-finite subsets of {\rm Dom}$(e)/e$.
\nl
1) If the pair $(I,J)$ is reasonable and $e \in \bold e(I,J)$
\ub{then} ${\Cal D}_{{\Cal E},e}$ is a filter on {\rm Dom}$(e)/e$ but
possibly $\emptyset \in {\Cal D}_{{\Cal E},e}$.
\nl
2) \quad $(a) \quad {\Cal D}_{\Cal E}$ is a filter on $J$
\sn
\hskip25pt $(b) \quad$ if 
${\Cal E}$ is non-trivial then all cofinite subsets of $J$
belongs to ${\Cal D}_{\Cal E}$ but $\emptyset \notin {\Cal D}_{\Cal E}$.
\endproclaim
\bigskip

\demo{Proof}  0) Easy, see Definition \scite{734-lin.4}(2).
\nl
1) By \scite{734-lin.9}(2) and \scite{734-lin.9}(3).
\nl
2) Trivial by Definition \scite{734-lin.7}(1).  \hfill$\square_{\scite{734-lin.10}}$
\enddemo
\bigskip

\proclaim{\stag{734-lin.11} Main Claim}  Assume
\mr
\item "{$(a)$}"  $I,J \in K^{\text{lin}}_{\tau^*_{\alpha(*)}}$
\sn
\item "{$(b)$}"   ${\Cal E}$ is an invariant $(I,J)$-equivalence relation
\sn
\item "{$(c)$}"  $(I,J)$ is a reasonable
$(\mu,\alpha(*))$-base which is a wide $(\lambda,\mu,\alpha(*))$-base
\sn
\item "{$(d)$}"  $e \in \bold e(J,I)$
\sn
\item "{$(e)$}"  $g$ is a function from {\rm Dom}$(e)/e$ into some
cardinal $\theta$
\sn
\item "{$(f)$}"  ${\Cal D}^* = \{Y \subseteq \theta:g^{-1}(Y) 
\in {\Cal D}_{{\Cal E},e}\}$
is a filter, i.e., $\emptyset \notin {\Cal D}^*$.
\ermn
\ub{Then} ${\Cal E}$ has at least $\chi := \lambda^\theta/{\Cal D}^*$  
equivalence classes.
\endproclaim
\bigskip

\demo{Proof}  Let $\langle f_\alpha:\alpha < \chi \rangle$ be a set
of functions from $\theta$ to $\lambda$ exemplifying $\chi :=
\lambda^\theta/{\Cal D}^*$ so $\alpha \ne \beta \Rightarrow \{i <
\theta:f_\alpha(i) = f_\beta(i)\} \notin {\Cal D}^*$.

Let $\langle h_\zeta:\zeta < \lambda \rangle$ exemplify the pair
$(I,J)$ being a wide $(\lambda,\mu,\alpha(*))$-base, see Definition
\scite{734-lin.1C}(4), so $h_\zeta \in \text{\rm inc}_J(I)$.

Lastly for each $\alpha < \chi$ we define $h^\alpha \in \text{\rm
inc}_J(I)$ as follows:
\nl
$h^\alpha(t)$ is: $h_0(t)$ \ub{if} $t \in J \backslash \text{\rm Dom}(e)$

\hskip25pt $h_{f_\alpha(g(t/e))}(t)$ \ub{if} $t \in \text{\rm
Dom}(e)$.
\mn
Now
\mr
\item "{$(*)_1$}"   $h^\alpha$ is a function from $J$ to $I$.
\ermn
[Why?  Trivially recalling each $h_\zeta$ is.]
\mr
\item "{$(*)_2$}"   $h^\alpha$ is increasing.
\ermn
[Why?  Let $s <_J t$ and we split the proof to cases.
\nl
If $s,t \in J \backslash \text{\rm Dom}(e)$ use ``$h_0 \in \text{\rm
inc}_J(I)$".
\nl
If $s \in J \backslash \text{\rm Dom}(e)$ and $t \in \text{\rm
Dom}(e)$, then $h^\alpha(t) = h_{f_\alpha(g(t/e))}(t),h^\alpha(s) =
h_0(s) = h_{f_\alpha(g(t/e))}(s)$ because $\langle h_\alpha \restriction
(J \backslash \text{\rm Dom}(e)):\alpha < \lambda\rangle$ is constant
(recalling $(h_0,h_\alpha)$ is an $e$-pair (for $\alpha >0$)), so
as $h_{f_\alpha(g(t/e))} \in \text{\rm inc}_J(I)$ we are done.
\nl
If $s \in \text{\rm Dom}(e),t \in J \backslash \text{\rm Dom}(e)$, the
proof is similar.
\nl
If $s,t \in\text{\rm Dom}(e),s/e \ne t/e$, we again use Definition
\scite{734-lin.1A}(2B), clause (b)$(\beta)$ of Definition \scite{734-lin.1C}(4).

Lastly, if $s,t \in \text{\rm Dom}(e),s/e=t/e$ we get $g(s/e) =
g(t/e)$ hence $f_\alpha(g(s/e)) = f_\alpha(g(t/e))$ call it $\gamma$
so $h^\alpha(s) = h_\gamma(s),h^\alpha(t) = h_\gamma(t)$ and of course
$h_\gamma \in \text{\rm inc}_J(I)$ hence $h_\gamma(s) <_I h_\gamma(t)$
so necessarily $h^\alpha(s) <_I h^\alpha(t)$ as required.  
So $(*)_2$ holds.]
\mr
\item "{$(*)_3$}"    $h^\alpha \in \text{ inc}_J(I)$.
\ermn
[Why?  Clearly if $i < \alpha(*)$ and $t \in P^J_i$ then $(\forall
\beta < \lambda)h_\beta(t) \in P^J_i$ hence $\alpha < \chi \Rightarrow
h_{f_\alpha(g(t/e))}(t) \in P^J_i$ which means $\alpha < \chi
\Rightarrow h^\alpha(t) \in P^J_i$; so recalling $(*)_2$, clause (a) of
Definition \scite{734-lin.1}(2) holds.
We should check clauses (b),(c) of Definition \scite{734-lin.1}(2)
which is done as in the proof of \scite{734-lin.1.14} and of $(*)_2$ above.]
\mr
\item "{$(*)_4$}"  if $\alpha < \beta$ and we let $u = u_{\alpha,\beta} :=
\cup\{g^{-1}(\zeta):\zeta < \theta$ and $f_\alpha(\zeta) \ne 
f_\beta(\zeta)\}$ so $u \subseteq \text{\rm Dom}(e)/e$ \ub{then}
$(h^\alpha,h^\beta)$
is a $(e \restriction \text{\rm set}(u))$-pair.
\ermn
[Why?  \ub{Case 1}:  If $s \in J \backslash \text{\rm Dom}(e)$ then
$h^\alpha(s) = h_0(s) = h^\beta(s)$.
\sn
\ub{Case 2}:  If $s \in \text{\rm Dom}(e) 
\backslash \text{\rm set}(u)$ then $h^\alpha(s) =
h_{f_\alpha(g(s/e))}(s) = h_{f_\beta(g(s/e))}(s) = h^\beta(s)$.
\sn
\ub{Case 3}:  If $s,t \in \text{\rm set}(u),s/e \ne t/e,s <_J t$ \ub{then}
$h^\alpha(s) <_I h^\beta(t) \wedge h^\beta(s) <_I h^\alpha(t)$ because
\sn
\ub{Subcase 3A}:  If $f_\alpha(g(s/e)) = 
f_\beta(g(t/e))$ we use $h_{f_\alpha(g(t/e))} \in \text{ inc}_J(I)$ hence

$$
h^\alpha(s) = h_{f_\alpha(g(s/e))}(s) <_I h_{f_\alpha(g(s/e))}(t) 
= h_{f_\beta(g(t/e))}(t) = h^\beta(t)
$$
\mn
and similarly $h^\beta(s) <_I h^\alpha(t)$.
\sn
\ub{Subcase 3B}:  $f_\alpha(g(s/e)) \ne f_\beta(g(t/e))$ we 
use $``(h_{f_\alpha(g(s/e))},h_{f_\beta(g(t/e))})$
is an $e$-pair".
\sn
\ub{Case 4}:  And lastly, if 
$s,t \in$ set$(u),s/e = t/e$ and $s <_J t$ \ub{then}
$h^\alpha(t) <_I h^\beta(s) \equiv (s/e \in u) \equiv h^\alpha(s) <_I
h^\beta(t)$.

Why?   Recalling 
$f_\alpha(g(s/e)) \ne f_\beta(g(t/e))$ as $s,t \in \text{ set}(u)$ by
the definition of $u$, see $(*)_4$ and we just use
$``(h_{f_\alpha(g(s/e))},h_{f_\beta(g(s/e))})"$ is an $e$-pair and
clause (c)$'$ of Definition \scite{734-lin.1A}.]
\mr
\item "{$(*)_5$}"    if $\alpha < \beta$ then $u_{\alpha,\beta} \ne
\emptyset$ mod ${\Cal D}_{{\Cal E},e}$.
\ermn
[Why?  By the choice of $\langle f_\alpha:\alpha < \lambda\rangle$.]
\mr
\item "{$(*)_6$}"  if $\alpha < \beta$ then $h^\alpha,h^\beta$ are not
${\Cal E}$-equivalent.
\ermn
[Why?  By $(*)_4 + (*)_5$  and \scite{734-lin.9}(2).] 

Together we are done.  \hfill$\square_{\scite{734-lin.11}}$
\enddemo
\bigskip

\proclaim{\stag{734-lin.11.8} Claim}   Assume ${\Cal E}$ is an invariant
$(I,J)$-equivalence relation, $I,J$ are well ordered and 
$|\text{\rm inc}_J(I)/{\Cal E}| \ge \lambda =
\text{\rm cf}(\lambda) > \mu = |I| > |2 + \alpha(*)|^{|J|}$.
\ub{Then} for some $e \in \bold e(I,J)$ there is an ultrafilter ${\Cal D}$ on
{\rm Dom}$(e)/e$ extending ${\Cal D}_{{\Cal E},e}$ which is not principal.
\endproclaim
\bigskip

\remark{Remark}    This is close to \cite[\S7]{Sh:620}.
\endremark
\bigskip

\demo{Proof}  Without loss of generality as linear orders, $J$ is
$\zeta(*)$ and $I$ is $\xi(*) \in [\mu,\mu^+)$.

Toward contradiction assume the conclusion fails.  Let $g$ be a 
one-to-one function from $\mu$ onto $[\xi(*)]^{< \aleph_0}$ and $\chi$
be large enough and $\kappa = |J|$ and $\partial = |2 +
\alpha(*)|^{|J|}$ so $\partial^\kappa = \partial$.

We now choose $\langle N_\eta:\eta \in {}^n \mu\rangle$ by induction
on $n < \omega$ such that
\mr
\item "{$\circledast_1$}"  $(a) \quad N_\eta \prec ({\Cal H}(\chi),\in)$
\sn
\item "{${{}}$}"  $(b) \quad \|N_\eta\| = \partial$ and $\partial +1
\subseteq N_\eta$
\sn
\item "{${{}}$}"  $(c) \quad A \subseteq N_\eta \wedge |A| \le \kappa
\Rightarrow A \in N_\eta$
\sn
\item "{${{}}$}"  $(d) \quad I,J$ and $g$ as well as $\eta$ belong to $N_\eta$
\sn
\item "{${{}}$}"  $(e) \quad \nu \triangleleft \eta \Rightarrow N_\nu
\in N_\eta$ (hence $N_\nu \subseteq N_\eta$ so $N_\nu \prec N_\eta$).
\ermn
There is no problem to do this.  Now it suffices to prove that for every
$h \in \text{ inc}_J(I)$, for some $h' \in \cup\{N_\eta:\eta \in
{}^{\omega >}\mu\} \cap \text{ inc}_J(I)$ we have $h {\Cal E} h'$.

Fix $h_* \in \text{ inc}_J(I)$ such that $h_* \notin \cup\{h/{\Cal E}:h
 \in \text{ inc}_J(I) \cap N_\eta$ for some $\eta \in {}^{\omega >} \mu\}$
 and for each $\eta \in 
{}^{\omega >}\mu$ we define $\bar \alpha_\eta,e_\eta$ as follows:
\mr
\item "{$\circledast_2$}"  $(a) \quad \bar\alpha_\eta = \langle
\alpha_{\eta,t}:t \in J\rangle$
\sn
\item "{${{}}$}"  $(b) \quad \alpha_{\eta,t} = 
\text{ min}((\xi(*) +1) \cap N_\eta \backslash h_*(t))$
\sn
\item "{${{}}$}"  $(c) \quad e_\eta := \{(s,t):s,t \in J$ and
$\alpha_{\eta,s} = \alpha_{\eta,t}$ and $\alpha_{\eta,s} > h_*(s)$ and
\nl

\hskip25pt $\alpha_{\eta,t} > h_*(t)\}$ 
\sn
\item "{${{}}$}"  $(d) \quad$ for $\alpha \in N_\eta$ let $X_{\eta,\alpha} :=
\{t \in J:\alpha_{\eta,t} = \alpha > h_*(t)\}$.
\ermn
Note
\mr
\item "{$(*)_1$}"  $\bar \alpha_\eta \in N_\eta$.
\ermn
[Why?  As $[N_\eta]^{\le \kappa} \subseteq N_\eta$ and $|J| = \kappa$
and $\alpha_{\eta,t} \in N_\eta$ for every $t \in J$.]
\mr
\item "{$(*)_2$}"  $(a) \quad e_\eta \in \bold e(J)$, i.e. $e_\eta$ is an
equivalence relation on some subset of $J$ with each equivalence class
a convex subset of $J$, see Definition \scite{734-lin.1A}(1)
\sn
\item "{${{}}$}"  $(b) \quad \langle X_{\eta,\alpha}:\alpha \in
\{\alpha_{\eta,t}:t \in \text{ Dom}(e)\}$ hence $X_{\eta,\alpha} \ne
\emptyset\rangle$ list the 
\nl

\hskip25pt $e_\eta$-equivalence classes.
\ermn
[Why?  Think.]
\mr
\item "{$(*)_3$}"  $h_\eta := h_* \restriction (J \backslash \text{
Dom}(e_\eta)) \in N_\eta$.
\ermn
[Why?  By the definition of $e_\eta$ we have 
$t \in J \wedge t \notin \text{ Dom}(e_\eta) 
\Rightarrow h_*(t) \in N_\eta$ and recall $[N_\eta]^{\le \kappa} 
\subseteq N_\eta$.]
\mr
\item "{$(*)_4$}"  if $t \in \text{ Dom}(e_\eta)$ then
cf$(\alpha_{\eta,t}) > \partial$.
\ermn
[Why?  As $\alpha_{\eta,t} \in N_\eta \prec ({\Cal H}(\chi),\in)$ 
if cf$(\alpha_{\eta,t})
= \theta \le \partial$ then there is a cofinal set $B$ of
$\alpha_{\eta,t}$ of cardinality $\theta$ in $N_\eta$ but 
$\theta \le \partial +1 \subseteq N_\eta$ therefore $B
\subseteq N_\eta$.  In particular as $h_*(t) < \alpha_{\eta,t}$ there
is $\beta \in B$ so that $h_*(t) < \beta$, but this contradicts the
choice of $\alpha_{\eta,t}$.]
\mr
\item "{$(*)_5$}"  $e_\eta \in \bold e(J,I)$.
\ermn
[Why?  Choose $h' \in \text{ inc}_J(I) \cap N_\eta$ similar enough to
$h_*$, specifically: $t \in J \backslash \text{ Dom}(e_\eta) \Rightarrow
h'(t) = h_*(t)$ and $t \in \text{ Dom}(e_\eta) \Rightarrow 
\sup\{\alpha_{\eta,s}:s \in J,s <_J t$ and $s \notin t/e_\eta\} < 
h'(t) < \alpha_{\eta,t}$.  The point being that
sup$\{\alpha_{\eta,s}:s \in J,s <_J t$ and $s \notin t/e_n\} \in
N_\eta$.  Now $(h',h_*)$ is a strict $e$-pair.]
\mr
\item "{$(*)_6$}"  there is $\ell_\eta < \omega$ and a finite sequence $\langle
\beta_{\eta,\ell}:\ell < \ell_\eta\rangle$ of members of Rang$(\bar
\alpha_\eta \restriction \text{ Dom}(e_\eta))$ so $X_{\eta,\beta_{\eta,\ell}}
\in \text{ Dom}(e_\eta)/e_\eta$ for $\ell < \ell_\eta$ such that
$\cup\{X_{\eta,\beta_{\eta,\ell}}:\ell < \ell_\eta\} \in 
{\Cal D}_{{\Cal E},e_\eta}$.
\ermn
[Why?  Otherwise there is an ultrafilter as desired, but toward
contradiction we have assumed this does not occur; in trying to get
generalizations we should act differently.]

Now we choose $(\eta_n,h_n)$ by induction on $n < \omega$ such that
\mr
\item "{$\boxdot$}"  $(a) \quad \eta_n \in {}^n \mu$
\sn
\item "{${{}}$}"  $(b) \quad$ if $n = m+1$ then $\eta_m = \eta_n
\restriction m$
\sn
\item "{${{}}$}"  $(c) \quad h_n \in \text{ inc}_J(I)$
\sn
\item "{${{}}$}"  $(d) \quad h_0 = h_*$
\sn
\item "{${{}}$}"  $(e) \quad$ if $n = m+1$ then:
{\roster
\itemitem{ ${{}}$ }  $(\alpha) \quad h_n {\Cal E} h_m$ hence $h_n
{\Cal E} h_*$ and Dom$(e_{\eta_n}) \subseteq \text{ Dom}(e_{\eta_m})$
\sn
\itemitem{ ${{}}$ }  $(\beta) \quad h_m \restriction (J \backslash \text{
Dom}(e_{\eta_m})) \subseteq h_n$
\sn
\itemitem{ ${{}}$ }  $(\gamma) \quad (h_m \restriction
\cup\{X_{\eta_m,\beta_{\eta_m,\ell}}:\ell < \ell_{\eta_m}\}) \subseteq h_n$
\sn
\itemitem{ ${{}}$ }  $(\delta) \quad h_n 
\restriction (\text{Dom}(e_{\eta_m}) \backslash
\cup\{X_{\eta_m,\beta_{\eta_m,\ell}}:\ell < \ell_{\eta_m}\})$ belongs to
$N_{\eta_m}$
\sn
\itemitem{ ${{}}$ }  $(\varepsilon) \quad$ moreover $t \in 
\text{ Dom}(e_{\eta_m}) \backslash
\cup\{X_{\eta_m,\beta_{\eta_m,\ell}}:\ell < \ell_{\eta_m}\}$ implies
\nl

\hskip40pt $h_n(t) < h_m(t)$
\sn
\itemitem{ ${{}}$ }  $(\zeta) \quad \ell_{\eta_m} > 0$
\endroster}
\item "{${{}}$}"  $(f) \quad Y_{m+1} \subseteq Y_m$ where $Y_m :=
\cup\{X_{\eta_m,\beta_{\eta_m,\ell}}:\ell < \ell_\eta\}$. 
\ermn
Why can we carry out the construction?  For $n=0$ we obviously can
(choose $h_0 = h_*$).  For $n=m+1$ first choose $h'_m \in N_{\eta_m}$
as we choose in the proof of $(*)_5$.  Now recalling $\langle
X_{\eta_m},\beta_{\eta_m,\ell}:\ell < \ell_{\eta_m}\rangle$ was chosen in
$(*)_6$, and define 
$h_n$ by $h_n \restriction (\text{Dom}(e_{\eta_m}) \backslash
\cup\{X_{\eta_m,\beta_{\eta_m,\ell}}:\ell < \ell_{\eta_m}\}) = 
h'_m \restriction (\text{Dom}(e_{\eta_m}) \backslash
\cup\{X_{\eta_m,\beta_{\eta_m,\ell}}:\ell < \ell_\eta\})$ and
$h_n \restriction (J \backslash \text{ Dom}(e_{\eta_m})) = h_m
\restriction (J \backslash  \text{ Dom}(e_{\eta_m})$ and
$h_n \restriction (\cup\{X_{\eta_m,\beta_{\eta_m,\ell}}:\ell <
n_{\eta_m}\}) = h_m \restriction 
(\cup\{X_{\eta_m,\beta_{\eta_m,\ell}}:\ell < \ell_{\eta_m}\})$.  Why
$h_n {\Cal E} h_m$?  Because
\mr
\widestnumber\item{$(iii)$}
\item "{$(i)$}"  as in the proof of $(*)_5,(h_n,h_m)$ form a strict
$\ell_\eta$-pair 
\sn
\item "{$(ii)$}"   they agree on
$\cup\{X_{\eta_m,\beta_{\eta_m,\ell}}:\ell < \ell_\eta\}$
\sn
\item "{$(iii)$}"  $\{X_{\eta_m,\beta_{\eta_m},\ell}:\ell < n\} \in
{\Cal D}_{{\Cal E},\bold e_\eta}$.
\ermn
Lastly, choose $\eta_n = \eta_m \char 94 \langle \gamma_m\rangle$ where
$\gamma_m$ is chosen such that $g(\gamma_m) =
\{\sup(\beta_{\eta_m,\ell} \backslash \sup\{h_m(t):t \in
X_{\beta_{\eta_m,\ell}}\}):\ell < \ell_{\eta_m}\}$ recalling that $g$
is a function from $\mu$ onto $[\xi(*)]^{< \aleph_0} = [I]^{< \aleph_0}$.

Now check that $\eta_n,h_n$ are as required.

Note that this induction never stops in the sense that $h_n \notin
N_{\eta_n}$ recalling the choice of $h_*$ and $h_n {\Cal E} h_*$.  
Now ${\Cal U}_n := \{\beta_{\eta_m,\ell}:\ell < n_\eta\}$
is a finite non-empty set of ordinals, and if $n = m+1$, then easily
$(\forall \ell < \ell_{\eta_n})(\exists k < 
\ell_{\eta_m})(\beta_{\eta_n,\ell} < \beta_{\eta_m,k})$ because for
$\ell < \ell_{\eta_n}$ letting $t \in X_{\eta_n,\ell}$ we know that
for some $k \le \ell_{\eta_m}$ we have $t \in X_{\eta_m,k}$ and
$\eta_n(m)$ was chosen above such that as $\gamma_m$, now $h_*(t) \le
\gamma_n \in N_{\eta_n},\gamma_m \le \alpha_{\eta_m,t}$ and the
inequality is strict as cf$(\alpha_{\eta_m,t}) > 0$.  So
$\langle\text{max}({\Cal U}_n):n < \omega\rangle$ is a decreasing sequence of
ordinals, contradiction, so we are done.   \hfill$\square_{\scite{734-lin.11.8}}$
\enddemo
\bn
\margintag{734-lin.eX}\ub{\stag{734-lin.eX} Example}:  For $e \in \bold e(J,I),J \in
K^{\text{lin}}_{\tau^*_{\alpha(*)}}$ and $I \in
K^{\text{lin}}_{\tau^*_{\alpha(*)}}$ we define ${\Cal E}^*_e 
= {\Cal E}^*_{e,I}$; it is
an invariant equivalent relation on inc$_J(I)$, by: $h_1 {\Cal E}^*_{e,I}
h_2$ iff:
\mr
\item "{$(a)$}"  if $t \in J \backslash \text{\rm Dom}(e)$ then
$h_1(t) = h_2(t)$
\sn
\item "{$(b)$}"   if $t \in \text{ Dom}(e)$ then cnv$_{I,h_1}(t) 
= \text{ cnv}_{I,h_2}(t)$ where cnv$_{I,h}(s) :=$ the convex hull
(in $I$) of the set $\{h_1(s)\} \cup \bigcup\{[h_1(s),h_1(t)]_I:s <_J
t$ and $t \in s/e\} \cup \bigcup \{[h(t),h(s)]_I:t <_J s$ and $t \in s/e\}$.
\ermn
1) If $J,I \in K^{\text{lin}}_{\tau^*_{\alpha(*)}}$ are well ordered
 and $e = J \times J$ then ${\Cal E}^*_{e,I}$ from part (1) 
has $\le |I| + \aleph_0$ equivalence classes.
\nl
2) If $J \in K^{\text{lin}}_{\tau^*_{\alpha(*)}}$ and $e$ as in part (2),
$\theta = \text{\rm cf}(J)$ and $|J| < \lambda = \lambda^{< \theta} <
\lambda^\theta$ \ub{then} there is 
$I \in K^{\text{lin}}_{\tau^*_{\alpha(*)}}$ of cardinality $\lambda$ such
that ${\Cal E}^*_{e,I}$ has $\lambda^\theta$ equivalence classes.
\bigskip

\remark{Remark}  We can define the stability spectrum for some
classes, essentially this is done in \S7, generally we intend to look
at it in \cite{Sh:F782}.
\endremark
\goodbreak

\head {\S7 Categoricity for a.e.c. with bounded amalgamation} \endhead  \resetall \sectno=7
 \spuriousreset
\bigskip

Recall that \scite{734-f.20} is the main result of this chapter; 
we think that it will lead to understanding
the categoricity spectrum of an a.e.c.  In particualr we hope
eventually to prove that this spectrum contains or is disjoint to some
 end segments of the class of cardinals.  Still here we like to show that what
we have is enough at least for restricted enough families of a.e.c. 
${\frak K}$'s, those definable by $\Bbb L_{\kappa,\omega},\kappa$ a
measurable cardinal or with enough amalgamation (concerning them and
earlier results see \chaptercite{E53}).  We could have relied on
\footnote{In the references to \cite{Sh:394}, e.g. 1.6tex is to 1.6 in
the published version and 1.8 is in the e-version.}
\cite{Sh:394}, but though we mention connections, we do not rely on
it, preferring self-containment.

We can say much even if we replace categoricity by strong solvability, 
but do this only when it is cheap; we can work even with weak and even
pseudo-solvability but not here.
\bigskip

\demo{\stag{734-loc.0} Hypothesis}  1) ${\frak K}$ is an a.e.c., so ${\Cal
S}(M) = {\Cal S}_{{\frak K}_\lambda}(M)$ for $M \in K_\lambda$, see
\marginbf{!!}{\cprefix{600}.\scite{600-0.12}}.
\nl
2) Let $K^x_\mu$ be the class $K_\mu$ if $K$ is categorical in $\mu$
and the class of superlimit models in ${\frak K}_\mu$ if there is one,
(the two definitions are compatible).
\enddemo
\bn
The following is a crucial claim because lack of locality is the problem
in \cite{Sh:394}.
\proclaim{\stag{734-loc.1} Claim}  Assume
\mr
\item "{$(a)$}"   {\rm cf}$(\mu) > \kappa \ge \text{\rm LS}({\frak K})$
\sn
\item "{$(b)$}"   ${\frak K}_{< \mu}$ has amalgamation 
\sn
\item "{$(c)$}"  $\Phi \in \Upsilon^{\text{or}}_\kappa[{\frak K}]$
satisfies:  if $I$ is $\theta$-wide and $\theta \in (\kappa,\mu)$ 
\ub{then} {\rm EM}$_{\tau({\frak K})}(I,M)$ is 
$\theta$-saturated (see \scite{734-0n.2}(1), \marginbf{!!}{\cprefix{600}.\scite{600-0.15}}(2) and
\marginbf{!!}{\cprefix{600}.\scite{600-0.19}}).
\ermn
\ub{Then}
\mr
\item "{$(\alpha)$}"  for some $\mu_* < \mu$, the class 
$\{M \in K_{< \mu}:M$ is saturated$\}$ is $[\mu_*,\mu)$-local, see Definition
\scite{734-loc.1.3}(3) below
\sn
\item "{$(\alpha)^+$}"  this applies not only to 
${\Cal S}(M) = {\Cal S}^1(M)$ but also for ${\Cal S}^\partial(M)$ 
if {\rm cf}$(\mu) > \kappa^\partial$.
\endroster
\endproclaim
\bn
Recall
\definition{\stag{734-loc.1.2} Definition}  ${\frak K}$ is 
$\mu$-stable if $\mu \ge \text{ LS}({\frak K})$ and
$M \in K_{\le \mu} \Rightarrow |{\Cal S}(M)| \le \mu$.
\enddefinition
\bn
Recall (\cite[Def.1.8=1.6tex]{Sh:394}(1),(2).
\definition{\stag{734-loc.1.3} Definition}  1) For $M \in {\frak K},\mu \ge
\text{ LS}({\frak K})$, satisfying $\mu \le \|M\|$ and 
$\alpha$, let $\Bbb E_{M,\mu,\alpha}$ be the following
equivalence relation on ${\Cal S}^\alpha(M):p_1 
\Bbb E_{M,\mu,\alpha} p_2$ iff for every $N \le_{\frak K} M$ of
cardinality $\mu$ we have $p_1 \restriction N = p_2 \restriction N$.
We may suppress $\alpha$ if it is 1, similarly below; let $\Bbb
E_{\mu,\alpha}$ be $\bigcup\{\Bbb E_{M,\mu,\alpha}:M \in K\}$ and so 
$\Bbb E_\mu = \Bbb E_{\mu,1}$. 
\nl
2) We say that $M \in {\frak K}$ is $\mu-\alpha$-local when
$\Bbb E_{M,\mu,\alpha}$ is the equality; we say that 
$p \in {\Cal S}^\alpha(M)$ is
$\mu$-local if $p/\Bbb E_{M,\mu,\alpha}$ is a singleton and we say
e.g. $K' \subseteq {\frak K}$ is $\mu-\alpha$-local (in ${\frak K}$,
if not clear from the context) when every $M \in K'$ is.
\nl
3) We say $K' \subseteq {\frak K}$ is $[\mu_*,\mu)-\alpha$-local if
every $M \in K' \cap {\frak K}_{[\mu^*,\mu)}$ is $\mu_*-\alpha$-local.  
\nl
4) We say that $\bar a \in N$ realizes $\bold p \in 
{\Cal S}^\alpha_{\frak K}(M)/\Bbb E_{\mu,\alpha}$ if 
$M \le_{\frak K} N$ and for every $M'
\le_{\frak K} N$ of cardinality $\mu$ the sequence $\bar a$
realizes $\bold p \restriction M'$ in $N$ or pedantically realizes $q
\restriction M'$ for some, equivalently every $q \in \bold p$. 
\enddefinition
\bigskip

\remark{Remark}  If $M \in {\frak K}_\mu$, then $M$ is $\mu-\alpha$-local.
\endremark
\bigskip

\demo{Proof of \scite{734-loc.1}}  Recall $\Phi \in
\Upsilon^{\text{or}}_\kappa[{\frak K}]$, see Definition 
\scite{734-11.1.3A}(2) and Claim \scite{734-0n.8}.
Easily there is $\langle I_\theta:\theta \in [\kappa,\mu)\rangle$, an 
increasing sequence of wide linear orders which are strongly
$\aleph_0$-homogeneous (that is dense with neither first nor last
element such that if $n < \omega$ and $\bar s,\bar t \in
{}^n(I_\theta)$ are $<_I$-increasing then some automorphism of
$I_\theta$ maps $\bar s$ to $\bar t$, e.g. the order of any real
closed field/or just ordered field) satisfying $|I_\theta|=\theta$.

Recalling $\Bbb Q$ here is the rational order, we let 
$J_\theta = \Bbb Q + I_\theta,M_\theta = \text{ EM}_{\tau({\frak
K})}(I_\theta,\Phi)$ and $N_\theta = \text{ EM}_{\tau({\frak K})}
(J_\theta,\Phi)$.  So 
\mr
\item "{$\circledast$}"  $(a) \quad M_\theta \le_{{\frak K}_\theta} N_\theta$
\sn
\item "{${{}}$}"   $(b) \quad M_{\theta_1} \le_{\frak K}
M_{\theta_2}$ and $N_{\theta_1} \le_{\frak K} N_{\theta_2}$ when $\kappa
 \le \theta_1 < \theta_2 < \mu$
\sn
\item "{${{}}$}"  $(c) \quad M_\theta$ is saturated (for ${\frak K}$,
of course) when $\theta > \kappa$
\sn
\item "{${{}}$}"  $(d) \quad$ every type from ${\Cal S}(M_\theta)$ is
realized in $N_\theta$
\sn 
\item "{${{}}$}"  $(e) \quad$ if $n < \omega,\bar a \in 
{}^n(N_\theta)$ then for some $\bar a' \in {}^n(N_\kappa)$ and 
automorphism  \nl

\hskip25pt $\pi$ of $N_\theta,\pi(\bar a) = \bar a'$ and 
$\pi$ maps $M_\theta$ onto itself.
\ermn
[Why?  Clauses (a),(b) holds by clause (c) of Claim \scite{734-0n.8}(1)
recalling Definition \scite{734-11.1.3A}(2).
\nl
Clause (c) holds by Clause (c) of the assumption of \scite{734-loc.1};
you may note \cite[6.7=6.4tex]{Sh:394}(2).
\nl
Clause (d) holds as EM$_{\tau({\frak K})}(\theta^+ + J_\theta,\Phi)
\in {\frak K}_{\theta^+}$ is saturated, and use the definition of a
type (or like the proof of claue (e) below using appropriate $I' +
I_\theta$ instead $\theta^+ + J_\theta$); you 
may note \cite[6.8=6.5tex]{Sh:394}.
\nl
Clause (e) holds as for every finite sequence $\bar t$ from $J_\theta$
there is an automorphism $\pi$ of $J_\theta$ such that: $\pi$ is the
identity on $\Bbb Q$, it maps $I_\theta$ onto itself and it maps $\bar
t$ to a sequence from $J_\kappa = \Bbb Q + I_\kappa$, such $\pi$
exists as $I_\theta$ is strongly $\aleph_0$-homogeneous and $I_\kappa
\subseteq I_\theta$ is infinite.]
\nl
For any $a \ne b$ from $N_\kappa$ let

$$
\align
\mu(a,b) = \text{ Min}\{\theta:&\theta \ge \kappa \text{ and if } 
\theta < \mu \\
  &\text{ then } \ortp_{\frak  K}(a,M_\theta,N_\theta) 
\ne \text{ \ortp}_{\frak K}(b,M_\theta,N_\theta)\}.
\endalign
$$
\mn
So $\mu(a,b) \le \mu$.  Let

$$
\mu_* = \sup\{\mu(a,b):a,b \in N_\kappa \text{ and }
\mu(a,b) < \mu\}.
$$
\mn
So $\mu_*$ is defined as the supremum on a set of $\le \kappa
\times \kappa$ cardinals $< \mu$, which is a cardinal of cofinality 
$\text{cf}(\mu) > \kappa$, hence clearly $\mu_* < \mu$.
Also $\mu_* \ge \kappa$ as there are $a \ne b$ from $M_\kappa$ hence
$\mu(a,b) = \kappa$.  Now suppose that $\theta \in [\mu_*,\mu),M \in {\frak
K}_\theta$ is saturated and $p_1 \ne p_2 \in {\Cal S}(M)$
and we shall find $M' \le_{\frak K}
M,M' \in {\frak K}_{\mu_*}$ such that $p_1 \restriction M' \ne p_2
\restriction M'$, this suffice.

Clearly $M_\theta \in K_\theta$ is saturated (by clause (c) of
 $\circledast$) hence the models
$M,M_\theta$ are isomorphic so without loss of generality 
$M=M_\theta$.  But by clause (d) of $\circledast$ every type 
from ${\Cal S}(M_\theta)$ is realized in $N_\theta$, so let $b_\ell$
 be such that $p_\ell = \text{ \ortp}_{\frak K}
(b_\ell,M_\theta,N_\theta)$ for $\ell=1,2$.  
Now there is an automorphism $\pi$ of
$N_\theta$ which maps $M_\theta$ onto itself and maps $b_1,b_2$ into
$N_\kappa$ (by clause (e) of $\circledast$) and let $a_\ell =
\pi(b_\ell)$ for $\ell=1,2$, so $a_1,a_2 \in N_\kappa$.

Now $$\align
\ortp(a_1,M_\theta,N_\theta) & =
\ortp(\pi(b_1),\pi(M_\theta),\pi(N_\theta)) =
\pi(\ortp(b_1,M_\theta,N_\theta)) \ne \cr
&\ne \pi(\ortp(b_2,M_\theta,N_\theta))
= \ortp(\pi(b_2),\pi(M_\theta),\pi(N_\theta)) =
\ortp(a_2,M_\theta,N_\theta).
\endalign$$
Hence by the 
definition of $\mu(a_1,a_2)$ we
have $\mu(a_1,a_2) \le \theta < \mu$.  Hence by the definition of $\mu_*$
we have $\mu(a_1,a_2) \le \mu_*$ which implies that \ortp$_{\frak
K}(a_1,M_{\mu_*},N_{\mu_*}) \ne \text{ \ortp}_{\frak K}
(a_2,M_{\mu_*},N_{\mu_*})$.
As $\pi$ is an automorphism of $N_\theta$ and 
$M_{\mu^*} \le_{\frak K} M_\theta$
it follows that $\ortp_{\frak K}(\pi^{-1}(a_1),
\pi^{-1}(M_{\mu^*}),\pi^{-1}(N_\theta)) \ne \ortp_{\frak
K}(\pi^{-1}(a_2),\pi^{-1}(M_{\mu^*}),\pi^{-1}(N_\theta))$ which means
\nl
$\ortp_{\frak K}(b_1,\pi^{-1}(M_{\mu^*}),N_\theta) \ne \ortp_{\frak
K}(b_2,\pi^{-1}(M_{\mu^*}),N_\theta)$,  but
$\pi^{-1}(M_{\mu^*}) \le_{\frak K} M_\theta$ as $\pi$ maps $M_\theta$
onto itself and recall that $p_\ell = \ortp_{\frak
 K}(b_\ell,M_\theta,N_\theta)$ so $p_\ell \restriction
 \pi^{-1}(M_{\mu^*})$ is well defined for $\ell=1,2$.  Hence
$p_1 \restriction \pi^{-1}(M_{\mu^*}) \ne p_2 \restriction 
\pi^{-1}(M_{\mu^*})$ and clearly $\pi^{-1}(M_{\mu^*})$
 has cardinality $\mu^*$ and is $\le_{\frak K} M_\theta$, so 
we are done proving
 clause $(\alpha)$.  The proof of clause $(\alpha)^+$ is the same
except that
\mr
\item "{$(*)_1$}"  if $\theta \in [\kappa,\mu),\bar t \in
{}^\partial(I_\theta)$ then some automorphism $\pi$ of $I_\theta$ maps
$\bar t$ to some $\bar t' \in {}^\partial(I_\kappa)$, justified by
\scite{734-am3.2.3} 
\sn
\item "{$(*)_2$}"  we replace $\Bbb Q$ by $\partial^+$
\sn
\item "{$(*)_3$}"   ${}^\partial(N_\kappa)$ has cardinality $\le
(\partial^+ + \kappa)^\partial \le \kappa^\partial < \text{ cf}(\mu)$.
\endroster
\nl
${{}}$  \hfill$\square_{\scite{734-loc.1}}$ 
\enddemo
\bn
Implicit in non-$\mu$-splitting is
\definition{\stag{734-am2.10} Definition}  Assume $\alpha < \mu^+,
N \in K_{\le \mu},N \le_{\frak K} M$ and 
$p \in {\Cal S}^\alpha(M)$ does not $\mu$-split over
$N$, see Definition \marginbf{!!}{\cprefix{705}.\scite{705-gr.1}}(1).  
The scheme of the non-$\mu$-splitting, ${\frak p} = 
\text{ sch}_\mu(p,N)$ is $\{(N'',c,\bar b)_{c \in N}/\cong :$ we
have 
$N \le_{\frak K} N' \le_{\frak K} M$ and 
$N' \le_{\frak K} N'',\{N',N''\} \subseteq K_\mu$ and the sequence 
$\bar b$ realizes $p \restriction N'$ in the model $N''\}$.
\enddefinition
\bigskip

\definition{\stag{734-am3.4} Definition}  For a cardinal $\mu$ and model $M$ let
\nl
1) 
$$
\align
\text{ps}-{\Cal S}_\mu(M) = {\Cal S}_{{\frak K},\mu}(M) = \{\bold
p:&\bold p \text{ is a function with
domain } \{N \in K_\mu:N \le_{\frak K} M\} \\
  &\text{ such that }\bold p(N) \in {\Cal S}(N) 
\text{ and } N_1 \le_{\frak K} N_2 \in
\text{ Dom}(\bold p) \\
  &\Rightarrow \bold p(N_1) = \bold p(N_2) \restriction N_1\}.
\endalign
$$
\mn
2) For $p \in {\Cal S}(M)$ let $p \restriction (\le \mu)$ be the
function ${\bold p}$ with domain $\{N \in K_\mu:N \le_{\frak K}
M\}$ such that ${\bold p}(N) = p \restriction N$.
\enddefinition
\bn

\demo{\stag{734-am3.5} Observation}  1) The function $p \mapsto p \restriction (\le
\mu)$ is a function from ${\Cal S}(M)$ into ps-${\Cal S}_\mu(M)$ such
that for $p_1,p_2 \in {\Cal S}(M)$ we have $p_1 \restriction (\le \mu)
= p_2 \restriction (\le \mu) \Leftrightarrow p_1 \Bbb E_\mu p_2$.
\nl
2) The subset $\{p \restriction (\le \mu):p \in{\Cal S}(M)\}$ of
 ps-${\Cal S}_\mu(M)$ has cardinality $|{\Cal S}(M)/\Bbb E_\mu|$.
\enddemo
\bigskip

\demo{Proof}  Should be clear.  \hfill$\square_{\scite{734-am3.5}}$
\enddemo
\bigskip

\proclaim{\stag{734-am2.18} Claim}  Every (equivalently some) $M \in K^x_\mu$ is
$\lambda^+$-saturated \ub{when}:
\mr
\item "{$(a)$}"  $(\alpha) \quad {\frak K}$ is categorical in $\mu$
\nl
or just
\item "{${{}}$}"  $(\beta) \quad {\frak K}$ is strongly solvable in $\mu$
\sn
\item "{$(b)$}"  {\rm LS}$({\frak K}) \le \lambda < \chi \le \mu$ and
$2^{2^\lambda} \le \mu$ (actually $2^\lambda \le \mu$ suffice)
\sn
\item "{$(c)$}"  $(\alpha) \quad \aleph_{\lambda^{+4}} =
\lambda^{+ \lambda^{+4}} \le \chi$
\nl

\hskip25pt \ub{or} at least
\sn
\item "{${{}}$}"  $(\beta) \quad$ if $\theta = 
{ \text{\rm cf\/}}(\theta) \le \lambda$ is $\aleph_0$ or a
measurable cardinal \ub{then} for some 
\nl

\hskip25pt $\partial \in (\lambda,\chi)$ we have:
$\partial = \partial^{< \theta} < \partial^\theta$ \ub{or} at least
$\partial^{<\theta>_{\text{tr}}} > \partial$ (i.e. there 
\nl

\hskip25pt is a tree ${\Cal T}$ with 
$\theta$ levels, $\partial$ nodes and the number of
$\theta$-branches 
\nl

\hskip25pt  of ${\Cal T}$ is $> \chi$, see \cite{Sh:589})
\sn
\item "{$(d)$}"  ${\frak K}_{\ge \partial} \ne \emptyset$ for every $\partial$,
equivalently $K_{\ge \theta} \ne \emptyset$ for arbitrarily large
$\theta < \beth_{1,1}(\text{\rm LS}({\frak K}))$
\sn
\item "{$(e)$}"  $(\alpha) \quad {\frak K}_{< \mu}$ has 
amalgamation and {\rm JEP}
\nl

\hskip25pt or just 
\sn
\item "{${{}}$}"  $(\beta) \quad$ if $\text{\rm LS}({\frak K}) \le \partial <
\chi$ then
{\roster
\itemitem{ ${{}}$ }  $\quad (i) \quad {\frak K}_\partial$ has 
amalgamation and {\rm JEP} and
\sn
\itemitem{ ${{}}$ }  $\quad (ii) \quad {\frak K}$ has $(\partial,\le
\partial^+,\mu)$-amalgamation \footnote{It suffices to have: if $M_0
\le_{\frak K} M_1 \in K_{\partial^+},M_1 \le_{\frak K} M_2 \in
K^x_\mu$ and $M_0 \in K_\partial$ then $M_1$ can be $\le_{\frak
K}$-embedded into some $M_3 \in K^x_\mu$.  Similarly in \scite{734-am3.6}.}
(see \marginbf{!!}{\cprefix{88r}.\scite{88r-2.5}}(2)) hence \footnote{Why?  Assume $M \in
K_{\partial^+}$ let $M_2 \in K^x_\mu$, let $M_0 \le_{\frak K} M_2$ be of
cardinality $\partial$, let $M_1 \in K_{\partial^+}$ be a $\le_{\frak
K}$-extension of $M_0$ which there is an $\le_{\frak K}$-embedding $f$
of $M$ into $M_1$ (exists as ${\frak K}_\partial$ has amalgamation and
JEP).  Lastly, use ``${\frak K}$ has 
$(\partial,\le\partial^+,\mu)$-amalgamation}
\sn
\itemitem{ ${{}}$ }  $\quad (iii) \quad$ every 
$M \in K_{\partial^+}$ has a $\le_{\frak K}$-extension in $K^x_\mu$
\nl

\hskip45pt (actually (i) + (iii) suffices).
\endroster}
\endroster 
\endproclaim
\bigskip

\remark{Remark}  1) $M$ is $\lambda^+$-saturated is well defined as
${\frak K}_{\le \lambda}$ has amalagamation.
\nl
2) We assume $2^{2^\lambda} \le \mu$ because
the proof is simpler with not much loss
(at least as long as other parts of the analysis are not much tighter).
\nl
3) We can weaken the assumptions.  In particular using solvability
instead categoricity, but for non-essential reasons this is delayed;
similarly in \scite{734-am3.6}.  
\nl
4) If $\mu = \mu^\lambda$ the claim is easy (as in \S1).
\endremark
\bigskip

\demo{Proof}  Note that by \cite[IX,\S2]{Sh:g}, \cite[II,3.1]{Sh:g} if
clause $(c)(\alpha)$ holds then clause $(c)(\beta)$ holds, hence we
can assume $(c)(\beta)$.

Let $\Phi \in \Upsilon^{\text{or}}_{\frak K}$ see Definition
\scite{734-11.1.3A}(2), exist by \scite{734-0n.8} and clause $(d)$ of the
assumption and $I \in K^{\text{lin}}_\mu \Rightarrow \text{
EM}_{\tau({\frak K})}(I,\Phi) \in K^x_\mu$ (trivially if $K$ is
categorical in $\mu$, otherwise by the definition of solvable). 

Clearly 
\mr
\item "{$(*)_0$}"  if $\partial \in [\text{LS}({\frak K}),\chi)$ \ub{then}
${\frak K}$ is stable in $\partial$.
\ermn
[Why?  We prove assuming clause $(e)(\beta)$, as the case of clause
$(e)(\alpha)$ is easier.  Otherwise as 
${\frak K}_\partial$ has amalgamation there are $M_0
\le_{\frak K} M_1$ such that $M_0 \in K_\partial,M_1 \in K_{\partial^+}$ and
$\{\ortp_{\frak K}(a,M_0,M_1):a \in M_1\}$ has cardinality $\partial^+$.
By assumption $(e)(\beta)(iii)$ there is $N_1$ such that 
$M_1 \le_{\frak K} N_1 \in {\frak K}_\mu$ and \wilog \, $N_1 \in
K^x_\mu$.  Let $I$ be as in
\scite{734-am3.2.3} with $(\lambda,\theta_2,\theta_1,\mu)$ there standing for
$(\mu,\partial^{++},\partial^+,\partial)$ here and 
$N_2 := \text{ EM}_{\tau({\frak K})}(I,\Phi)$.  Now by 
\scite{734-am3.2.3}(2), $N_1 \ncong N_2$,
contradiction to ``${\frak K}$ categorical in $\mu$".  Or you may
see \cite[1.7=1.5tex]{Sh:394}.]
\enddemo
\bn
The proof now splits to two cases.
\nl
\ub{Case 1}:  For every $M \in K^x_\mu$ we have $\mu 
\ge |{\Cal S}(M)/\Bbb E_\lambda|$.  

For every $M \in K^x_\mu$ there is $M'$ such that: $M \le_{\frak K} M'
\in K_\mu$ and for every $\bold p \in {\Cal S}(M)/\Bbb E_\lambda$
either $\bold p$ is realized in $M'$ or there are no $M'',a$ such that
$M' \le_{\frak K} M'' \in K_\mu$ and $a \in M''$ realizes $p$ in $M''$.
\nl
[Why?  Let $\langle p_i/\Bbb E_\lambda:i < \mu\rangle$ list ${\Cal
S}(M)/\Bbb E_\lambda$, exists by the assumptions and choose $M_i$ for
$i \le \mu,\le_{{\frak K}_\mu}$-increasing continuous such that
$M_{i+1}$ satisfies the demand for $\bold p = p_i/\Bbb E_\lambda$, 
possibly no $p \in p_i/\Bbb E_\lambda$ has an
extension in ${\Cal S}(M_{i+1})$ (hence is not realized in it), so
then the desired demand holds trivially; note that it is not
unreasonable to assume ${\frak K}_\mu$ has
amalgamation and it clarifies but it is not necessary.]

Also \wilog \, $M' \in K^x_\mu$ as any model $M$ from $K_\mu$ has a
$\le_{\frak K}$-extension in $K^x_\mu$ (at least if $M$ does
$\le_{\frak K}$-extend some $M' \in K^x_\mu$).

Now we can choose by induction on $i \le \lambda^+$ a
model $M_i \in K^x_\mu,\le_{\frak K}$-increasing continuous with
$i$, such that for every $p \in {\Cal S}(M_i)$ \ub{either} 
there is $q \in {\Cal S}(M_i)$
realized in $M_{i+1}$ which is $\Bbb E_\lambda$-equivalent to $p$ \ub{or}
there is no $\le_{\frak K}$-extension of $M_{i+1}$ satisfying this.
Now we shall prove that $M_{\lambda^+}$ is $\lambda^+$-saturated recalling
Definition \marginbf{!!}{\cprefix{600}.\scite{600-0.15}}.  Now if $N \le_{\frak K} M_{\lambda^+},
\|N\| \le \lambda$ and $p \in {\Cal S}(N)$ then there
is $i < \lambda^+$ such that $N \le_{\frak K} M_i$ and we can find $p' \in
{\Cal S}(M_{\lambda^+})$ extending $p$.  (Why?  If clause $(e)(\alpha)$
holds then this follows by ${\frak K}_{< \mu}$ having amalgamation, 
see \marginbf{!!}{\cprefix{88r}.\scite{88r-2.8}}.  If clause $(e)(\beta)$ holds, use ``${\frak K}$ has
the $(\lambda,\le \lambda^+,\mu)$-amalgamation property" recalling
LS$({\frak K}) \le \lambda < \chi$.)  Hence
there is $a \in M_{i+1}$ such that \ortp$(a,M_i,M_{i+1}) 
\Bbb E_\lambda(p' \restriction M_i)$, hence $a$ realizes 
$p$ in $M_{i+1}$ hence in $M_{\lambda^+}$.
\bn
\ub{Case 2}:  Not Case 1.

Let $I$ be as in \scite{734-am3.2.3} with $(\lambda,\theta_2,\theta_1,\mu)$ there
standing for $(\mu,\lambda^{++},\lambda^+,\lambda)$ here, so $|I|=\mu$.  
Let $M = \text{ EM}_{\tau({\frak K})}(I,\Phi)$,
so by not Case 1 we can find $p_i \in {\Cal S}(M)$ for $i < \mu^+$ pairwise
non-$\Bbb E_\lambda$-equivalent.  As ${\frak K}_\lambda$ is a
$\lambda$-a.e.c. with amalgamation and is stable in $\lambda$ (by $(*)_0$) we
can deduce, see \marginbf{!!}{\cprefix{705}.\scite{705-gr.6}}(2), that: if $p \in {\Cal S}(M)$ then
for some $N \le_{\frak K} M$ of cardinality $\lambda$ 
the type $p$ does not $\lambda$-split over
$N$ (or see \cite[3.2 = 3.2tex]{Sh:394}(1)).
For each $i$ choose $N_i \le_{\frak K} M$ of
cardinality $\lambda$ such that $p_i$ does not $\mu$-split over $N_i$.
As there is no loss in increasing $N_i$ (as long as it
is $\le_{\frak K} M$ and has cardinality $\lambda$) \wilog 
\mr
\item "{$(*)_1$}"  $N_i = \text{ EM}_{\tau({\frak K})}(I_i,\Phi)$ 
where $I_i \subseteq I$ and $|I_i| = \lambda$ 
and let $\bar t_i = \langle t^i_\varepsilon:\varepsilon < \lambda \rangle$
list $I_i$ with no repetitions.
\ermn
As $2^\lambda \le \mu$ \wilog \, the $I_i$'s are
pairwise isomorphic, so without loss of generality 
for $i,j < \mu^+$, the mapping 
$t^i_\varepsilon \mapsto t^j_\varepsilon$ is such an isomorphism. 
Moreover, without loss of generality 
\mr
\item "{$(*)_2$}"  for every $i,j < \mu^+$ there is an
automorphism $\pi_{i,j}$ of $I$ mapping $t^i_{\varepsilon}$ to
$t^j_{\varepsilon}$ for $\varepsilon < \lambda$.
\ermn
[Why?  By \scite{734-am3.2.3}(1) as we can replace 
$\langle p_i:i < \mu^+ \rangle$ by
$\langle p_i:i \in {\Cal U} \rangle$ for every unbounded ${\Cal U}
\subseteq \mu^+$.]

Let ${\frak p}_i$ be the non-$\lambda$-splitting scheme of $p$ over $N_i$
(see Definition \scite{734-am2.10}).  Without loss of generality:
\mr
\item "{$(*)_3$}"  for $i,j < \mu^+$, the isomorphism $h_{i,j}$
from $N_j = \text{ EM}_{\tau({\frak K})}(I_j,\Phi)$ onto $N_i = \text{
EM}_{\tau({\frak K})}(I_i,\Phi)$ induced by the mapping
$t^j_\zeta \mapsto t^i_\zeta$ (for $\zeta < \lambda$)  satisfies
{\roster
\itemitem{ $(i)$ }  it is an isomorphism from $N_j$ onto $N_i$
\sn
\itemitem{ $(ii)$ }  it maps ${\frak p}_j$ to ${\frak p}_i$.
\endroster}
\ermn
[Why?  For (i) this holds by the definition of EM$(I_i,\Phi)$.  For
(ii) let $h_{i,0}$ map ${\frak p}_i$ to ${\frak p}'_i$.
The number of schemes is $\le 2^{2^\lambda}$; so if $\mu \ge
2^{2^\lambda}$ then \wilog \, $i < \mu^+ \Rightarrow {\frak p}'_i =
{\frak p}'_1$ hence we are done (with no real loss).  If we weaken the
assumption $\mu \ge 2^{2^\lambda}$ to $\mu \ge 2^\lambda$ (or even
$\mu > \lambda$ so waive $(*)_2$) using \scite{734-am3.2.3}(4) we can find
$I^+_i$ such that $I_i \subseteq I^+_i \subseteq I,|I^+_i| \le \lambda^+$
and for every $J \subseteq I$ of cardinality $\le \lambda$ there is an
automorphism of $I$ over $I_i$ mapping $J$ into $I^+_i$.  So only
$\langle {\frak p}'_i((\text{EM}_{\tau({\frak K})}(I^+_0,\Phi),
c,\bar b))_{c \in \text{EM}_{\tau({\frak K})}(I_0,\Phi)}/\cong):\bar b \in
{}^\lambda(\text{EM}_{\tau({\frak K})}(I^+_0,\Phi))\rangle$ matters
(an overkill) but this
is determiend by $p_i \restriction \text{ EM}_{\tau({\frak
K})}(I^+_0,\Phi))$ which $\in {\Cal S}(\text{\rm EM}_{\tau({\frak K})}
(I^+_0,\Phi))$ by $(*)_0$ and as
${\frak K}$ is stable in $\lambda^+$ \wilog \, 
${\frak p}'_{1+i} = {\frak p}'_1$ and we are done.]

Now we translate our problem to one on expanded (by unary predicates) 
linear orders which was treated
in \S6.  Recall that by \scite{734-am3.2.3}(3), we can use $I = 
\text{\rm EM}_{\{<\}}(I^*,\Psi)$ where 
$\Psi \in \Upsilon^{\text{lin}}_{\aleph_0}[2]$, see Definition
\scite{734-X1.2}(5), and 
$I^* = I^{\text{lin}}_{\lambda,\mu \times \lambda^+}$ from \scite{734-lin.2}(2)
with $\alpha(*) = 2$.  Recall that $I^* = I^{\text{lin}}_{\lambda,\mu
\times \lambda^{++}}$ is $\mu \times \lambda^{++}$ expanded by
$P_1 = \{\alpha \in I^*:\text{\rm cf}(\alpha) \ge \lambda^+\},P_0 = I_*
\backslash P_0$ so $I^*$ is a well ordered $\tau^*_2$-model,
i.e. $\in K^{\text{lin}}_{\tau^*_2}$, see Definition \scite{734-X1.2}(5).
Without loss of generality $I_i = 
\text{\rm EM}_{\{<\}}(I^*_i,\Psi)$ where $I^*_i \subseteq I^*$ has
cardinality $\lambda$ and the pair $(I^*,I^*_i)$ is a reasonable 
$(\lambda,\alpha(*))$-base which is a wide
$(\mu,\lambda,\alpha(*))$-base, 
see Definition \scite{734-lin.1C}(3)(4), Claim \scite{734-lin.2}(2).  Without
loss of generality for every $i < \mu^+$ there is $h_i$, 
an isomorphism from $I^*_0$ onto $I^*_i$ such that (see below) the
induced function $h^{[1]}_1$ maps
$\bar t_0$ to $\bar t_i$.  Let $J^* = I^*_0$ and $J=I_0$.  We like
to apply \S6 for $J^*,I^*$ fixing $\alpha(*) = 2,\bar u^* = (u^-,u^+)
= (\{0\},\emptyset)$.  So recalling Definition \scite{734-lin.1}(2) 
for every $h \in \text{ inc}^{\bar u^*}_{J^*}(I^*)$ we 
can naturally define the function $h^{[1]}$ by 
$h^{[1]}(\sigma^{\text{EM}(J^*,\Psi)}(t_0,\dotsc,t_{n-1})) =
\sigma^{\text{EM}(J^*,\Psi)}(a_{h(t_0)},\dotsc,a_{h(t_{n-1})})$ whenever
$\sigma(x_0,\dotsc,x_{n-1})$ is a $\tau(\Psi)$-term and $J^* \models
``t_0 < \ldots < t_{n-1}"$ so it is an isomorphism from
EM$_{\{<\}}(J^*,\Psi)$ onto EM$_{\{<\}}(I^* \restriction \text{
Rang}(h),\Psi)$ so as $J^* \subseteq I^*$ by \scite{734-am3.2.3}(5) there
is an automorphism $h^{[2]}$ of $I$ extending $h^{[1]}$ and so there is
an automorphism $h^{[3]}$ of EM$(I,\Phi)$ such that $h^{[3]}(a_t) =
a_{h^{[2]}(t)}$ for $t \in I$ and
$h^{[3]}(\sigma^{\text{EM}(I,\Phi)}(a_{t_0},\dotsc,a_{t_{n-1}})) =
\sigma^{\text{EM}(I,\Phi)}(a_{h^{[2]}(t_0)},
\dotsc,a_{h^{[2]}(t_n)})$ where $t_0 <_I
\ldots <_I t_{n-1}$ and $\sigma(x_0,\dotsc,x_{n-1})$ is a
$\tau(\Phi)$-term.

Note that
\mr
\item "{$(*)_4$}"  if $h',h''$ are automorphisms of EM$_{\tau[{\frak
K}]}(I,\Phi)$ extending $h^{[3]} \restriction 
\text{ EM}_{\tau{[\frak K}]}(I_0)$ then 
$h'(p_0/\Bbb E_\lambda) = h''(p_0 /\Bbb E_\lambda)$.
\ermn
[Why?  Because $p_0$ does not $\lambda$-split over EM$_{\tau[{\frak
K}]}(I_0,\Phi)$.] 

We define a two-place relation ${\Cal E}$ on inc$_{J^*}(I^*)$ by:
$h_1 {\Cal E} h_2$ if $h^{[3]}_1(p_0/\Bbb E_\lambda) = 
h^{[3]}_2(p_0/\Bbb E_\lambda)$.
(Note that $h \mapsto h^{[3]}$ is a function so this is well defined
and $h^{[3]}$ is an automorphism of EM$_{\tau({\frak K})}(I,\Phi))$.
By $(*)_4$ clearly ${\Cal E}$ is an invariant equivalence 
relation on inc$^{\bar u^*}_{J^*}(I^*)$ with $> \mu$ equivalence
classes as exemplified by $\langle h_i:i < \mu^+\rangle$.  

By \scite{734-lin.11.8} there is $e \in \bold e(J^*,I^*)$ 
such that (recalling Definition \scite{734-lin.9})
 the filter ${\Cal D}_{{\Cal E},e}$ 
has an extension to a non-principal ultrafilter ${\Cal D}$ so for some
regular $\theta \le \lambda$ there is a function $g$ from Dom$(\bold e)/e$
onto $\theta$ which maps ${\Cal D}$ to a uniform ultrafilter 
$g({\Cal D})$ on $\theta$, so $\partial^{<\theta>_{\text{tr}}} \le
\partial^{\text{Dom}(e)/e}/{\Cal D}_{{\Cal E},e}$ for every cardinal 
$\partial$.
Choose such a pair $(g,\theta)$ with minimal 
$\theta$ so ${\Cal D}$ is $\theta$-complete hence
$\theta = \aleph_0$ or $\theta$ is a measurable cardinal $\le \lambda$.
By clause $(c)(\beta)$ of our assumption justified in the beginning of
the proof there is $\partial \in (\lambda^+,\chi)$ such that 
$\partial < \partial^{<\theta>_{\text{tr}}}$ 
hence $\partial^+ \le \partial^{<\theta>_{\text{tr}}} \le
\partial^{\text{Dom}(e)/e}/{\Cal D}_{{\Cal E},e}$.
So letting $I^0_\partial = I^{\text{lin}}_{\lambda,\partial \times 
\lambda^{++}} \subseteq I^*$ the set $\{\bar t/{\Cal E}:
\bar t \in \text{ incr}_{J^*}(I^*)$ and Rang$(\bar t) \subseteq 
I^0_\partial\}$ has cardinality $> \partial$.  Now for each
$\bar t \in \text{ inc}_{J^*}{}^{\bar u^*}(I^*)$ 
let $\pi_{\bar t} \in \text{ Aut}(I)$ be such that $\pi_{\bar t}(\bar t_0) 
= \bar t$ and let $\hat \pi_{\bar t}$ be the automorphism of
EM$_{\tau({\frak K})}(I,\Phi)$ which $\pi_{\bar t}$ induce,
 and let $p_t = \hat \pi_{\bar t}(p_0) \in 
{\Cal S}(M)$.   Hence $\{\hat \pi_{\bar t}(p_0) \restriction 
\text{ EM}_{\tau({\frak K})}(I^{\text{lin}}_{\lambda,\partial \times
\lambda^+},\Phi):\bar t \in \text{ inc}_{J^*}{}^{\bar u^*}(I^*)$ and
$\text{Rang}(\bar t) \subseteq I^{\text{lin}}_{\lambda,\partial \times 
\lambda^{++}}\}$ is of cardinality $> \partial$,
contradicting ``${\frak K}$ stable in $\partial$" from $(*)_0$.
\hfill$\square_{\scite{734-am2.18}}$ 
\bn
Note but we shall not use
\demo{\stag{734-am2.20} Conclusion}  1) Under the assumptions of
\scite{734-am2.18} we have $\kappa({\frak K}_\mu) = \aleph_0$, see below.
\nl
2) Moreover, $\bold \kappa_{\text{st}}({\frak K}_\mu) = \emptyset$.  
\enddemo
\bn
Recall
\definition{\stag{734-am2.20.7} Definition}  If ${\frak K}_\mu$ is an 
$\mu$-a.e.c. with amalgamation which is stable, \ub{then}:
\mr
\item "{$(a)$}"   $\kappa({\frak K}_\mu) = \aleph_0 + 
\sup\{\kappa^+:\kappa$ regular
$\le \mu$ and there is an $\le_{{\frak K}_\mu}$-increasing continuous sequence
$\langle M_i:i \le \kappa \rangle$ and $p \in {\Cal S}(M_\kappa)$ such
that $M_{2i+2}$ is universal over $M_{2i+1}$ and $p
\restriction M_{2i+2}$ does $\mu$-split over $M_{2i+1}\}$
\sn
\item "{$(b)$}"  $\kappa_{\text{sp}}({\frak K}_\mu) := 
\{\kappa:\kappa$ regular $\le \mu$ and there is
an $\le_{{\frak K}_\mu}$-increasing continuous sequence $\langle M_i:i \le
\kappa \rangle$ and $p \in {\Cal S}(M_\kappa)$ which $\mu$-splits over
$M_i$ for each $i < \kappa$ and $M_{2i+2}$ is universal over
$M_{2i+1}\}$.
\endroster
\enddefinition
\bigskip

\demo{Proof of \scite{734-am2.20}}  By playing with 
EM$(I,\Phi)$, (or see Claim \cite[5.7=5.7tex]{Sh:394} and Definition
\cite[4.9=4.4tex]{Sh:394}).  \hfill$\square_{\scite{734-am2.20}}$
\enddemo
\bn
\margintag{734-am3.3}\ub{\stag{734-am3.3} Question}:  Can we omit assumption \scite{734-am2.18}(c) 
(see below so $\chi = \text{ LS}({\frak K})$)?
\bigskip

\proclaim{\stag{734-am3.6} Theorem}  For some cardinal $\lambda_* < \chi$ and
 a cardinal $\lambda_{**} < \beth_{1,1}(\lambda^{+ \omega}_*)$ above 
$\lambda_*,{\frak K}$ is categorical in every cardinal $\lambda \ge
\lambda_{**}$ but in no $\lambda \in (\lambda_*,\lambda_{**})$ provided
that:
\mr
\item "{$\circledast^{\mu,\chi}_{\frak K}$}"  $(a) \quad K$ is an
a.e.c. cateogorical in $\mu$
\sn
\item "{${{}}$}"  $(b) \quad {\frak K}$ has amalgamation and {\rm JEP} in every
$\lambda < \aleph_\chi,\lambda \ge \text{\rm LS}({\frak K})$
\sn
\item "{${{}}$}"  $(c) \quad \chi$ is a limit cardinal, 
{\rm cf}$(\chi) > \text{\rm LS}({\frak K})$,
and for arbitrarily large $\lambda < \chi$ 
\nl

\hskip25pt the sequence $\langle 2^{\lambda^{+n}}:
n < \omega \rangle$ is increasing
\sn
\item "{${{}}$}"  $(d) \quad \mu > \beth_{1,1}(\lambda)$ for every
$\lambda < \chi$ hence $\mu \ge \aleph_\chi$
\sn
\item "{${{}}$}"  $(e) \quad$ every $M \in K_{< \aleph_\chi}$ has a
$\le_{\frak K}$-extension in $K_\mu$. 
\endroster
\endproclaim
\bigskip

\remark{Remark}  1) Concerning \cite{Sh:394} note
\mr
\item "{$(a)$}"  there the central case was ${\frak K}$ with full
amalgamation (not just below $\chi \ll \mu$!), trying to concentrate on
the difficulty of lack of localness,
\sn
\item "{$(b)$}"  when we use clause (e) this is just to
get the $``M \in K_\mu$ is $\lambda$-saturated", this is where we use
\scite{734-am2.18}
\sn
\item "{$(c)$}"  we demand ``cf$(\chi) > \text{ LS}({\frak K})$" to
prove locality.
\ermn
2) We rely on \chaptercite{600} and \chaptercite{705} in the end.
\nl
3) The assumption (e) of \scite{734-am3.6} follows if ${\frak K}$ has 
amalgamation in every $\lambda' \le \beth_{1,1}(\lambda)$ 
for $\lambda < \chi$ which is a reasonable assumption.
\nl
4) Most of the proof works even if we weaken the assumption (a) to
``${\frak K}$ is strongly solvable in $\mu$" and even weakly solvable,
i.e. up to $\boxdot_7$, we continue in and see more \cite{Sh:F782}.
\nl
5) Theorem \scite{734-am3.6} also continue 
Kolman-Shelah \cite{KlSh:362}, \cite{Sh:472}, as its assumptions are
proved there.
\endremark
\bigskip

\demo{Proof}  Let $\kappa = \text{\rm LS}({\frak K})$ and let $\Phi
\in \Upsilon^{\text{or}}_\kappa[{\frak K}]$ be as guaranteed by
\scite{734-0n.8}(1) hence
\mr
\item "{$(*)_1$}"  if $I \in K^{\text{lin}}_\lambda$ then
EM$_{\tau({\frak K})}(I,\Phi)$ belongs to $K_\lambda$ 
for $\lambda \ge \text{\rm LS}({\frak K})$ (and in the strongly
solvable case, $I \in K^{\text{lin}}_\mu \Rightarrow 
\text{ EM}_{\tau({\frak K})} (I,\Phi) \in K^x_\mu$)
\ermn
and
\mr
\item "{$(*)_2$}"  if $I \subseteq J$ are from $K^{\text{lin}}$ then
 EM$_{\tau({\frak K})}(I,\Phi) \le_{\frak K} 
\text{\rm EM}_{\tau({\frak K})}(J,\Phi)$.
\ermn
Also
\mr
\item "{$(*)_3$}"  $\langle{\Cal S}_{\frak K}(M):
M \in {\frak K}_{< \aleph_\chi}\rangle$ has the reasonable basic properties.
\ermn
[Why?  See \marginbf{!!}{\cprefix{600}.\scite{600-0.12}} and \marginbf{!!}{\cprefix{600}.\scite{600-0.12A}} 
because ${\frak K}_{< \aleph_\chi}$ has the
amalgamation property by clause (b) of the assumption 
$\circledast^{\mu,\chi}_{\frak K}$).]
\mr
\item "{$(*)_4$}"  if $M \in K_\mu$ then $M$ is $\chi$-saturated
(hence $\chi$-model homogeneous).
\ermn
[Why?  We shall prove that: if LS$({\frak K}) \le \lambda <
\chi$ and $M \in K^x_\mu$ \ub{then} $M$ is $\lambda^+$-saturated.  We
shall show that all the assumptions of \scite{734-am2.18} with 
$(\mu,\chi,\lambda)$ there standing for
$(\mu,\aleph_\chi,\lambda)$ here hold.  Let us check; clause (a) of
\scite{734-am2.18} means ``${\frak K}$ is categorical in $\mu$" (or is
strongly solvable) which
holds by clause (a) of $\circledast^{\mu,\chi}_{\frak K}$.  Clause (b)
of \scite{734-am2.18} says
that LS$({\frak K}) \le \lambda < \aleph_\chi \le \mu$ and $2^{2^\lambda} \le
\mu$; the first holds because of the way $\lambda$ was chosen above and
the second holds as clause (d) of $\circledast^{\mu,\chi}_{\frak K}$
says that $\mu > \beth_{1,1}(\lambda)$ and $\mu  \ge \aleph_\chi$.  
Clause $(c)(\alpha)$ of \scite{734-am2.18} holds as $\lambda^{+ \lambda^{+4}} <
\aleph_{\lambda^{+5}}$ which is $< \aleph_\chi$ as $\chi$ is a limit
cardinal and $\aleph_\chi$ here plays the role of $\chi$ there.  
Clause (d) of \scite{734-am2.18} says ${\frak K}_{\ge \partial} \ne
\emptyset$ for every cardinal $\partial$, holds by $(*)_1$ above.
Lastly, clause (e) of \scite{734-am2.18} holds more exactly clauses
$(e)(\beta)(i)+(iii)$ hold by clauses (b) + (e) of
$\circledast^{\mu,\chi}_{\frak K}$ and they suffice.  

We have shown that all the assumptions of \scite{734-am2.18} holds,
hence its conclusion, which says, as $M \in K_\mu$, that
$M$ is $\lambda^+$-saturated.  The ``$\chi$-model homogeneous" holds by
\marginbf{!!}{\cprefix{600}.\scite{600-0.19}}.]
\mr
\item "{$(*)_5$}"  if $M \le_{\frak K} N$ are from $K^x_\mu$ then $M
\prec_{\Bbb L_{\infty,\chi}[{\frak K}]} N$.
\ermn
[Why?  Obvious by $(*)_4$.]
\mr
\item "{$(*)_6$}"  if $\lambda \in (\kappa,\chi)$ and $I \in
K^{\text{lin}}_{\ge \lambda}$ is $\lambda$-wide 
then EM$_{\tau({\frak K})}(I,\Phi)$ is
$\lambda$-saturated; moreover, if 
$I^+ \in K^{\text{lin}}_\lambda$ is wide over $I$ then
every $p \in {\Cal S}(\text{EM}_{\tau(K)}(I,\Phi))$ is realized in
EM$_{\tau({\frak K})}(I^+,\Phi)$.
\ermn
[Why?  By \scite{734-11.3}, its assumption ``$\Phi$ satisfies the
conclusion of \scite{734-11.2}" holds by $(*)_5$,
 (or as in \cite[6.8=6.5tex]{Sh:394}).  The ``moreover" is immediate by
$(*)_4$ as in the proof of $\circledast(d)$ inside the proof of
\scite{734-loc.1} above or see the proof of $(*)_{10}$ below.]
\mr
\item "{$(*)_7$}"   ${\frak K}$ is stable in $\lambda$ when $\kappa
\le \lambda < \chi$.
\ermn
[Why?  Recalling clause (e) of the assumption of \scite{734-am3.6},
by Claim \scite{734-am2.18} or more accurately $(*)_0$ in its proof as we
have proved (in the proof of $(*)_4$) that the assumptions of
\scite{734-am2.18} holds with $(\mu,\chi,\lambda)$ there standing for
$(\mu,\aleph_\chi,\lambda)$ here.]
\mr
\item "{$(*)_8$}"  if $\lambda \in [\kappa,\chi)$ and $M \in
K^x_\lambda$ \ub{then} there is $N \in {\frak K}_\lambda$ which is
$(\lambda,\aleph_0)$-brimmed over $M$
\ermn
[Why?  By $(*)_7$ and \marginbf{!!}{\cprefix{600}.\scite{600-0.22}}(1)(b) remembering the
amalgamation, clause (b) of the assumption of the theorem.]
\mr
\item "{$(*)_9$}"    if $\langle M_\alpha:\alpha \le
\lambda\rangle$ is $\le_{\frak K}$-increasing continuous, 
$\kappa \le \|M_\lambda\| \le \lambda < \chi$, \ub{then} no 
$p \in {\Cal S}_{\frak K}(M_\lambda)$ satisfies 
$p \restriction M_{i+1}$ does $\lambda$-split
over $M_i$ for every $i < \lambda$.
\ermn
[Why?  Otherwise we get contradiction to stability in $\lambda$,
i.e. $(*)_7$, see in \marginbf{!!}{\cprefix{705}.\scite{705-gr.6}}(1B), using amalgamation (using 
the tree ${}^{\theta >} 2$ when
$\theta = \text{ min}\{\partial:2^\partial > \lambda\}$; also we can
prove it as in the proof of case 2 inside the proof of \scite{734-am2.18}.]

We could use more
\mr
\item "{$(*)_{10}$}"  if $I_1,I_2$ are wide linear orders of
cardinality $\lambda \in (\kappa,\chi)$ and $I_2$ is wide over
$I_1$ so $I_1 \subseteq I_2$ and 
$M_\ell = \text{ EM}_{\tau({\frak K})}(I_\ell,\Phi)$,
\ub{then} $M_2$ is universal over $M_1$ and even brimmed
over $I_1$, even $(\lambda,\partial)$-brimmed for any regular $\partial <
\lambda$.
\ermn
[Why?  As $I_2$ is wide over $I_1$, we can find a sequence $\langle
J_\gamma:\gamma < \lambda\rangle$ of pairwise disjoint subsets of $I_2
\backslash I_1$ such that each $J_\gamma$ is a convex subset of $I_2$
and in $J_\gamma$ there is a monotonic sequence $\langle
t_{\gamma,n}:n < \omega\rangle$ of members.  Let $\langle
\gamma_\varepsilon:\varepsilon < \lambda \times \partial\rangle$ list
$\lambda$, and let $I_{2,0} = I_1$ and $I_{2,1 + \varepsilon} = I_2
\backslash \cup\{J_{\gamma_\zeta}:\zeta \in [1 + \varepsilon,\lambda
\times \partial)\}$ and $M'_\zeta = \text{ EM}_{\tau({\frak
K})}(I_{2,\varepsilon},\Phi)$.  So $\langle M'_\zeta:\zeta \le \lambda
\times \partial \rangle$ is $\le_{\frak K}$-increasing continuous
sequence of members of $K_\lambda$; first member $M_1$, last member
$M_2$.

By \marginbf{!!}{\cprefix{600}.\scite{600-0.22}}(4)(b) it is enough to prove that if 
$\varepsilon < \lambda \times \partial$ and $p
\in {\Cal S}(M_\varepsilon)$ then $p$ is realized in $M_{\varepsilon
+1}$.  As $I_1$ is wide of cardinality $\lambda$ so is
$I_{2,\varepsilon}$ hence $M'_\varepsilon$ is saturated.  Also for
each $\varepsilon$ we can find a linear order $I^+_{2,\varepsilon}$ of
cardinality $\lambda$ such that $I_{2,\varepsilon +1} \subseteq
I^+_{2,\varepsilon}$ and $J^+_\varepsilon = I^+_{2,\varepsilon +1}
\backslash I_{2,\varepsilon}$ is a convex subset of
$I^+_{2,\varepsilon +1}$ and is a wide linear order of cardinality
$\lambda$ which is strongly $\aleph_0$-homogeneous, (recall
$J_{\gamma_\varepsilon} \subseteq J^+_{\gamma_\varepsilon}$ is 
infinite).  So in $M^+_{\varepsilon +1} = \text{ EM}_{\tau({\frak
K})}(I^+_{2,\varepsilon +2},\Phi)$ every $p \in {\Cal
S}(M^1_\varepsilon)$ is realized (as $I^+_{2,\varepsilon +1}$ is wide
over $I_{2,\varepsilon}$ as $J^+_\varepsilon$ is wide of cardinality
$\lambda$), moreover realized in $M'_{\varepsilon +1}$
(why? by the strong $\aleph_0$-homogeneous every element and even finite
sequence from $M^+_{\varepsilon +1}$ can be mapped by
some automorphism of $M^+_{\varepsilon +1}$ over $M_\varepsilon$ into
$M_{\varepsilon +1}$).  As said above, this suffices.]
\mr
\item "{$\circledast_1$}"  $\chi_*$ is well defined $\in (\kappa,\chi)$ where
$$
\align
\chi_* = \text{\rm Min}\{\theta:&\kappa < \theta < \chi \text{ and for
every saturated} \\
  &M \in {\frak K}, \text{ if } \theta \le \|M\| < \chi, \text{ every} \\
  &p \in {\Cal S}(M) \text{ is } \theta\text{-local, see Definition
\scite{734-loc.1.3}(2)}\}.
\endalign
$$
\ermn
[Why?  By \scite{734-loc.1} which we apply with $(\mu,\kappa)$ there
standing for $(\chi,\kappa)$ here recalling $\kappa = \text{
LS}({\frak K})$; this is O.K. as: clause (a) in
\scite{734-loc.1} holds by clause (c) of the assumption here, clause (b)
in \scite{734-loc.1} holds by clause (b) of the assumption here as $\chi \le
\aleph_\chi$.  Lastly, 
clause (c) in \scite{734-loc.1} easily follows by $(*)_6$ above.]
\mr
\item "{$\circledast_2$}"   if $\lambda \in (\kappa,\chi)$ and 
$\langle M_i:i \le \delta\rangle$ is 
$\le_{{\frak K}_\lambda}$-increasing continuous,
$M_{i+1}$ is $\le_{\frak K}$-universal over $M_i$ for $i < \delta$
\ub{then}
$M_\delta$ is saturated and moreover
every $p \in {\Cal S}(M_\delta)$ does not $\lambda$-split over
$M_\alpha$ for some $\alpha < \delta$.
\ermn
[Why?  For $i \le \delta$ let $I_i$ be the linear order $\lambda
\times \lambda \times (1+i)$ and $M'_i = \text{ EM}_{\tau({\frak
K})}(I_i,\Phi)$.  So $\langle M'_i:i \le \delta\rangle$ is
$\le_{{\frak K}_\lambda}$-increasing continuous.  Also for $i \le
\delta,\zeta \le \lambda$ let $I_{i,\zeta} = \lambda \times \lambda
\times (1+i) + \lambda  \times \zeta$ and $M'_{i,\zeta} = \text{
EM}_{\tau({\frak K})}(I_{i,\zeta},\Phi)$, so for each $i < \delta$ the
sequence $\langle M'_{i,\zeta}:\zeta \le \lambda\rangle$ is 
$\le_{{\frak K}_\lambda}$-increasing 
continuous, $M'_{i,0} = M'_i,M'_{i,\lambda} =
M'_{i+1}$.  Now for $i< \delta,\zeta < \lambda$ every $p \in {\Cal
S}(M_{i,\zeta})$ is realized in $M'_{i,\zeta +1}$ by $(*)_6$ and the
definition of type, varying the linear order.  By
\marginbf{!!}{\cprefix{600}.\scite{600-0.22}}(4)(b) the model
$M'_{i+1}$ is $\le_{{\frak K}_\lambda}$-universal
over $M'_i$ and by Definition \marginbf{!!}{\cprefix{600}.\scite{600-0.21}} the models $M'_\delta$ and
$M_\delta$ are $(\lambda,\text{cf}(\delta))$-brimmed hence by
\marginbf{!!}{\cprefix{600}.\scite{600-0.22}}(3) are isomorphic.  But $M'_\delta$ is saturated by
$(*)_6$, hence $M_\delta$ is saturated.

What about the ``moreover"?  (Note that if $\lambda =
\lambda^{\text{cf}(\delta)}$ then $(*)_9$ does not cover it.)
We can find easily $\langle I''_\alpha:\alpha \le \lambda \times \delta
+1\rangle$ such that:
\mr
\item "{$(a)$}"    $I''_\alpha$ is a linear order of cardinality
$\lambda$ into which $\lambda$ can be embedded
\sn
\item "{$(b)$}"    $I''_\alpha$ is increasing continuous with $\alpha$
\sn
\item "{$(c)$}"    $I''_\alpha$ is an initial segment of $I''_\beta$
for $\alpha < \beta \le \delta +1$
\sn
\item "{$(d)$}"    $I''_{\alpha +1}$ has a subset of order types $\lambda
\times \lambda$ whose convex hull is disjoint to $I''_\alpha$
\sn
\item "{$(e)$}"    if $\alpha \le \beta < \lambda \times \delta$ and
$s \in I''_{\lambda \times \delta +1} \backslash I''_{\lambda \times \delta}$
then there is an automorphism $\pi_{\alpha,\beta,s}$ of $I''_{\lambda
\times \delta +1}$ mapping $I''_{\beta +1}$ onto 
$I''_{\lambda \times \delta}$ and is over $I''_\alpha \cup \{t \in 
I''_{\lambda \times \delta +1}:s \le_{I''_{\lambda \times \delta +1}} t\}$.
\ermn
Let $M''_\alpha = \text{ EM}_{\tau({\frak K})}(I''_\alpha,\Phi)$, so
$\langle M''_{\lambda \times \alpha}:\alpha \le \delta\rangle$ has 
the properties of $\langle M'_\alpha:\alpha \le \delta\rangle$, i.e.
every $p \in {\Cal S}(M''_\alpha)$ is realized in $M''_{\alpha +1}$
hence $M''_{\alpha + \lambda}$ is $\le_{{\frak K}_\lambda}$-universal
over $M''_\alpha$. So (easily or see \marginbf{!!}{\cprefix{600}.\scite{600-0.22}},\marginbf{!!}{\cprefix{600}.\scite{600-0.21}})
there is an isomorphism $f$ from $M_\delta$ onto $M''_{\lambda \times
\delta}$ such that $M''_{\lambda \alpha} \le_{\frak K} f(M_{\alpha + 1}) \le 
M''_{\lambda \alpha + 2}$.  So it suffices to prove 
the ``moreover" for $\langle M''_{\lambda \times \alpha}:\alpha \le
\delta\rangle$, equivalently for $\langle M''_\alpha:
\alpha \le \lambda \times \delta\rangle$.  Let $p \in {\Cal S} 
(M''_{\lambda \times \delta})$ so some $a \in M''_{\lambda \times
\delta +1}$ realizes it, hence for some $t_0 < \ldots < t_{n-1}$ from
$I''_{\lambda \times \delta +1}$ and $\tau_\Phi$-term
$\sigma(x_0,\dotsc,x_{n-1})$ we have $a = 
\sigma^{\text{EM}(I''_{\lambda \times \delta +1},\Phi)}(a_{t_0},
\dotsc,a_{t_{n-1}})$, it follows that for some $m \le n$ we
have $t_\ell \in I''_{\lambda \times \delta} \Leftrightarrow \ell < m$ and let
$\alpha < \lambda \times \delta$ be such that $\{t_\ell:\ell < m\} \subseteq
I''_\alpha$; if $m=n$ choose any $t_n \in I''_{\lambda \times \delta +1}
\backslash I''_{\lambda \times \delta}$.  
If $\beta \in (\alpha,\lambda \times \delta)$ and
$\ortp_{\frak K}(a,M''_\delta,M''_{\delta +1})$ does $\lambda$-split
over $M''_\beta$ then $\pi' := \pi_{\beta,\beta,t_m}$ is an
automorphism of $I''_{\lambda \times \delta +1}$ mapping 
$I''_{\beta +1}$ onto $I''_{\lambda \times \delta}$ and is over $I''_\beta
\cup \{s \in I''_{\lambda \times \delta +1}:t_m \le_{I''_{\lambda \times
\delta +1}} s\}$ hence it is the identity on $\{t_\ell:\ell < n\}$;
now $\pi'$ induces an automorphism $\hat \pi'$ of
EM$_{\tau({\frak K})}(I''_{\lambda \times \delta +1},\Phi)$, so 
clearly it maps $a$ to itself and maps
$\ortp_{\frak K}(a,M''_{\beta +1},M''_{\lambda \times \delta +1})$ 
to $\ortp_{\frak K}(a,M''_{\lambda \times \delta},
M''_{\lambda \times \delta +1})$ and it maps 
$M''_\beta$ onto itself, hence also 
$\ortp_{\frak K}(a,M''_{\beta +1},M''_{\delta +1})$ does
$\lambda$-split over $M''_\beta$.  So if for some $\beta \in
(\alpha,\lambda \times \delta)$, the type
$\ortp_{\frak K}(a,M''_\delta,M''_{\delta + 1})$ does not $\lambda$-split over
$M''_\beta$ we get the desired conclusion, but 
otherwise this contradicts $(*)_9$.]
\mr
\item "{$\circledast_3$}"    If $\lambda \in [\chi_*,\chi)$ and
$M \in K_\lambda$ is saturated and $p \in {\Cal S}(M)$
\ub{then} for some $N$ we have:
{\roster
\itemitem{ $(a)$ }  $N \le_{\frak K} M$
\sn
\itemitem{ $(b)$ }  $N \in K_{\chi_*}$ is saturated
\sn
\itemitem{ $(c)$ }  $p$ does not $\chi_*$-split over $N$
\sn
\itemitem{ $(d)$ }   $p$ does not $\lambda$-split over $N$ (follows by
(a),(b),(c)).
\endroster}
\ermn
[Why $\circledast_3$ holds?  For clauses (a),(b),(c) 
use $\circledast_2$ or just $(*)_9$; for clause (d) use localness, i.e. 
recall $\circledast_1$ and Definition \scite{734-loc.1.3}.]
\mr
\item "{$\circledast_4$}"   Assume $\lambda \in [\kappa,\chi)$ and 
$M_1 \le_{\frak K} M_2 \le_{\frak K} M_3$ are members of
$K,M_2$ is $\lambda^+$-saturated and $p \in {\Cal S}(M_3)$.
If $N_\ell \le_{\frak K} M_\ell$ is from $K_{\le \lambda}$ and $p
\restriction M_{\ell+1}$ does not $\lambda$-split over $N_\ell$ for $\ell
= 1,2$ \ub{then} $p$ does not $\lambda$-split over $N_1$.
\ermn
[Why?  Easy manipulations.  Without loss of generality $N_1 \le_{\frak
K} N_2$ as we can increase $N_2$.  So for some pair $(M_4,a)$ we have
$M_3 \le_{\frak K} M_4,a \in M_4$ and $p = \ortp_{\frak K}(a,M_3,M_4)$.  
Assume $\alpha < \lambda^+$ and let
$\bar b,\bar c \in {}^\alpha(M_3)$ be such that $\ortp_{\frak K}
(\bar b,N_1,M_3) = \ortp_{\frak K}(\bar c,N_1,M_3)$.   
As $M_2$ is $\lambda^+$-saturated and $N_2 \le_{\frak K}
M_2 \le_{\frak K} M_3$ we can find $\bar b',
\bar c' \in {}^\alpha(M_2)$ such that $\ortp_{\frak K}(\bar b' \char 94 \bar
c',N_2,M_3) = \ortp_{\frak K}(\bar b \char 94 \bar c,N_2,M_3)$ using
\marginbf{!!}{\cprefix{600}.\scite{600-0.19}}.  Hence $\ortp_{\frak K}(\bar b',N_1,M_3) = 
\ortp_{\frak K}(\bar b,N_1,M_3) = \ortp_{\frak K}(\bar c,N_1,M_3) 
= \ortp_{\frak K}(\bar c',N_1,M_3)$.

By the choice of $(M_4,a)$ and the assumption on $N_1$ that $p
\restriction M_2$ does not $\lambda$-split over $N_1$ we get

$$
\ortp_{\frak K}(\langle a\rangle \char 94 \bar b',N_1,M_4) = 
\ortp_{\frak K}(\langle a \rangle \char 94 \bar c',N_1,M_4).
$$
\mn
Clearly $\ortp_{\frak K}(\bar b',N_2,M_3) =
\ortp_{\frak K}(\bar b,N_2,M_3)$ hence by the choice of $(M_4,a)$ and
the assumption on $N_2$ that $p$ does not $\lambda$-split over $N_2$
we have $\ortp_{\frak K}(\langle a\rangle \char 94 \bar b',N_2,M_4) =
\ortp(\langle a\rangle \char 94 \bar b,N_2,M_4)$
hence by monotonicity 

$$
\ortp_{\frak K}(\langle a\rangle \char 94 \bar
b',N_1,M_4) = \ortp_{\frak K}(\langle a\rangle \char 94 \bar b,N_1,M_4).
$$
\mn
Similarly 

$$
\ortp_{\frak K}(\langle a\rangle \char 94 \bar c',N_1,M_4) =
\ortp_{\frak K}(\langle a\rangle \char 94 \bar c,N_1,M_4).
$$
\mn
As equality of types is transitive 
$\ortp_{\frak K}(\langle a\rangle \char 94 \bar c,N_1,M_4)
= \ortp_{\frak K}(\langle a\rangle \char 94 \bar c',N_1,M_4)
= \ortp_{\frak K}(\langle a\rangle \char 94 \bar b',N_1,M_4) =
\ortp_{\frak K}(\langle a\rangle \char 94 \bar b,N_1,M_4)$ as required.]
\mr
\item "{$\circledast_5$}"   Assume $I_3 = I_0 + I'_1 + I'_2$ are wide  
linear orders of cardinality $\lambda$ where $\chi > \lambda > \kappa$
and let $I_\ell = I_0 + I'_\ell$ for $\ell=1,2$ and
$M_\ell = \text{ EM}_{\tau({\frak K})} (I_\ell,\Phi)$ for
$\ell=0,1,2,3$.  If $\ell \in \{1,2\}$ and 
$\bar a \in {}^{\lambda >}(M_\ell)$ \ub{then}
$\ortp_{{\frak K}_\lambda}
(\bar a,M_{3-\ell},M_3)$ does not $\lambda$-split over $M_0$, (moreover
if $\ortp_{{\frak K}_\lambda}(\bar a,M_0,M_3)$ does not
$\lambda$-split over $N \in K_{\le \lambda}$ then also $\ortp_{{\frak
K}_\lambda}(\bar a,M_{3 -\ell},M_3)$ does not $\lambda$-split over $N$).
\ermn
[Why?  For $\ell=2$, if the desired conclusion fails 
we get a contradiction as in the proof of
$\circledast_2$ so for $\ell=2$ we get the conclusion.  
For $\ell=1$ if the desired conclusion fails (but it holds for
$\ell=2$) we get a contradiction to categoricity in $\mu$ by the
order property (by \scite{734-11.1.2}).]
\mr
\item "{$\circledast_6$}"  If $\lambda \in (\chi_*,\chi),\delta <
\lambda^+,\langle M_i:i \le \delta\rangle$ is 
$\le_{{\frak K}_\lambda}$-increasing continuous and 
$i < \delta \Rightarrow M_i$ saturated \ub{then} $M_\delta$ is saturated.
\ermn
[Why?  Let $N \le_{\frak K} M_\delta,\|N\| < \lambda$ and $p \in {\Cal
S}(N)$.  If cf$(\delta) > \|N\|$ this is easy so assume cf$(\delta)
\le \|N\|$ hence cf$(\delta) < \lambda$ 
and \wilog \, $\delta = \text{ cf}(\delta)$ and choose a
cardinal $\theta$ such that  LS$({\frak K}) < \chi_* + |\text{cf}(\delta)| +
\|N\| \le \theta < \lambda$ and $\|N\|^+ < \lambda \Rightarrow \|N\| <
\theta$ and let $q \in {\Cal S}(M_\delta)$ extend $p$, exist as
${\frak K}_{\le \lambda}$ has amalgamation.

Now for every $X \subseteq M_\delta$ of cardinality $\le \theta$ we
can choose $N_i \le_{\frak K} M_i$ by induction 
on $i \le \delta$ such that $N_i
\in K_\theta$ is saturated, is $\le_{\frak K}$-increasing continuous
with $i$ and $N_i$ is $\le_{\frak K}$-universal over $N_j$ and
includes $(X \cup N) \cap M_i$ when $i=j+1$.  So by $\circledast_2$
(we justify the choice of $N_i$ for limit $i$ and)
the model $N_\delta$ is saturated, so if 
$\|N\|^+ < \lambda$ then $N \le_{\frak
K} N_\delta,N_\delta$ is saturated of cardinality $\theta > \|N\|$ so  
we are done as $N_\delta \le_{\frak K} M_\delta$, 
so \wilog \, $\lambda = \|N\|^+$ hence $\lambda = \theta^+$.

Also for some $\alpha_* < \delta$ and $N_* \le_{\frak K} M_{\alpha_*}$
of cardinality $\theta$, the type $q$ does not $\theta$-split over
$N_*$.  [Why?  Otherwise we choose $(N_i,N^+_i)$ by induction on $i \le
\delta$ such that $N_i \le_{\frak K} N^+_i$ are from $K_\theta,N_i
\le_{\frak K} M_i,N^+_i \le_{\frak K} M_\delta,N_i$ is $\le_{\frak
K}$-increasing continuous, $N_i$ is $\le_{\frak K}$-universal over
$N_j$ if $i=j+1$ and $q \restriction N^+_i$ does $\theta$-split over
$N_i$ and $\cup\{N^+_j \cap M_i:j<i\} \subseteq N_i$.  In the end we
get a contradiction to $\circledast_2$.]

We can find $N' \le_{\frak K} M_{\alpha_*}$ from $K_{\chi_*}$ such
that $q \restriction M_{\alpha_*}$ does not $\theta$-split over $N'$,
(why? by $\circledast_3$) and \wilog \, $N' \le_{\frak K} N_*$ and $N'
\le_{\frak K} N$.  Also
$q$ does not $\theta$-split over $N'$ (why? by applying
$\circledast_4$, with $\theta,N_*,M_{\alpha_*},M_\delta$ here standing
for $\lambda,M_1,M_2,M_3,N_1,N_2$ there; or use $N' = N_*$).

By $(*)_6$ as $M_{\alpha_*}$ is saturated \wilog \, $M_{\alpha_*} = \text{
EM}_{\tau({\frak K})}(\lambda,\Phi)$ and for $\varepsilon < \lambda$
let $M_{\alpha_*,\varepsilon}
= \text{ EM}_{\tau({\frak K})}(\theta \times \theta \times (1 +
\varepsilon),\Phi)$, so $M_{\alpha_*,\varepsilon} \in K_\theta$ is
saturated and is brimmed over $M_{\alpha^*,\zeta}$ when $\varepsilon =
\zeta +1$ by $(*)_{10}$.  So for each $\varepsilon < \lambda$ there is
$a_\varepsilon \in M_{\alpha^*,\varepsilon +1}$ realizing $q
\restriction M_{\alpha_*,\varepsilon}$.  Also \wilog \, $M_\delta
\le_{\frak K} \text{ EM}_{\tau({\frak K})}(\lambda + \lambda,\Phi)$ as
in the proof of $\circledast_2$ or by $(*)_{10}$,
now for some $\varepsilon(*) < \lambda$ we have $N \le_{\frak K}
\text{ EM}_{\tau({\frak K})}(I_2,\Phi)$ and $N_* \le_{\frak K} 
\text{ EM}_{\tau({\frak K})}(I_0,\Phi)$ where $I_0 = \theta \times \theta
\times (1 + \varepsilon(*))$ and $I_2 = [\lambda,\lambda +
\varepsilon(*)) \cup I_0$.  Let $I_1 = \theta \times \theta \times
\zeta(*)$ where $\zeta(*) \in (\varepsilon(*),\lambda)$ is large
enough such that $a_{\varepsilon(*)} \in 
\text{ EM}_{\tau({\frak K})}(I_1,\Phi)$, e.g. $\zeta(*) = 1 +
\varepsilon(*)+1$ and let $I_3 = I_1 \cup I_2 \subseteq \lambda +
\lambda$.  Let $M'_\ell = \text{ EM}_{\tau({\frak K})}(I_\ell,\Phi)$
for $\ell=0,1,2,3$.

Now we apply $\circledast_5$, the ``moreover" with
$\theta,I_0,I_1,I_2,I_1 \backslash I_0,I_2 \backslash
I_0,a_{\varepsilon(*)},N'$ here standing for
$\lambda,I_0,I_1,I_2,I'_1,I'_2,\bar a,N$ there and we conclude that
$\ortp_{{\frak K}_\lambda}(a_{\varepsilon(*)},M'_2,M'_3)$
does not $\theta$-split over $N'$.

As $N' \le_{\frak K} M'_0 \le_{\frak K} M'_2$ also the type 
$q' := \ortp_{{\frak K}_\lambda}
(a_{\varepsilon(*)},M'_2,M'_3)$ does not
$\theta$-split over $N'$.  Let us sum up: $q \restriction M'_2,q'$ 
belong to ${\Cal S}_{{\frak K}_\lambda}(M'_2)$, does not $\theta$-split over
$N',N' \in K_{\chi_*}$ and $\chi_* \le \theta$.  Also $N' \le_{{\frak K}_*}
M'_0 \le_{{\frak K}_*} M'_2$, the model $M'_0$ is
$\theta$-saturated and $q \restriction M_{\alpha_*} = q'
\restriction M_{\alpha_*}$.  By the last two sentences obviously $q =
q'$ (it may be more transparent to consider $q \restriction (\le \chi_*) = q'
\restriction (\le \chi_*))$, so we are done proving $\circledast_6$.]
\mr
\item "{$\circledast_7$}"  If $\lambda \in (\chi_*,\chi)$ \ub{then} 
the saturated $M \in {\frak K}_\lambda$ is superlimit.
\ermn
[Why?  By $\circledast_6$, (existence by $(*)_6$,
the non-maximality by $(*)_6 + $ uniqueness; you may
look at \cite[6.7=6.4tex]{Sh:394}(1).]

Now we have arrived to the main point
\mr
\item "{$\odot_1$}"  If $\lambda \in (\chi_*,\chi)$ \ub{then} 
${\frak s}_\lambda$ is a full good 
$\lambda$-frame, $K_{{\frak s}_\lambda}$ categorical  
where ${\frak s}_\lambda$ is defined by
{\roster
\itemitem{ $(a)$ }  ${\frak K}_{{\frak s}_\lambda} = {\frak K}_\lambda
\restriction \{M \in K_\lambda:M$ saturated$\}$
\sn
\itemitem{ $(b)$ }   ${\Cal S}^{\text{bs}}_{{\frak s}_\lambda}(M) 
= {\Cal S}^{\text{na}}_{{\frak s}_\lambda}(M) := \{\ortp_{\frak s}(a,M,N):M
\le_{{\frak K}_\lambda} N$ and $a \in N \backslash M\}$ 
for $M \in K_{{\frak s}_\lambda}$
\sn
\itemitem{ $(c)$ }  $p \in {\Cal S}^{\text{bs}}_{{\frak s}_\lambda}
(M_2)$ does not fork over $M_1$ \ub{when} $M_1 \le_{{\frak s}_\lambda} M_2$
and for some $M \le_{\frak K} M_1$ of cardinality $\chi_*$, the
type $p$ does not $\chi_*$-split over $N$.
\endroster}
\ermn
[Why?  We check the clauses of Definition \marginbf{!!}{\cprefix{600}.\scite{600-1.1}}.
\mn
\ub{$K_{{\frak s}_\lambda}$ is categorical}:

By \marginbf{!!}{\cprefix{600}.\scite{600-0.34}}(1) and $\circledast_7$.
\mn
\ub{Clause (A),Clause (B)}:  By $\circledast_7$ recalling that 
there is a saturated $M \in K_{{\frak s}_\lambda}$ (and it is not
$<_{{\frak s}_\lambda}$-maximal) by $(*)_6$ and trivially 
recalling \marginbf{!!}{\cprefix{600}.\scite{600-0.34}}, of course.
\mn
\ub{Clause (C)}:  By categoricity and $(*)_6$ clearly no $M \in
K_{{\frak s}_\lambda}$ is maximal; amalgamation and JEP holds by
clause (b) of the assumption of the claim.
\mn
\ub{Clause (D)(a),(b)}:  By the definition.
\mn

\ub{Clause (D)(c)}:  Density is obvious; in fact ${\frak s}_\lambda$ is full.
\mn
\ub{Clause (D)(d)}:  (bs - stability).

Easily ${\Cal S}_{{\frak s}_\lambda}(M) = {\Cal S}_{{\frak
K}_\lambda}(M)$ which has cardinality $\le \lambda$ by the moreover in
$(*)_6$.
\mn
\ub{Clause (E)(a)}:  By the definition.
\mn
\ub{Clause (E)(b)}:  Monotonicity (of non-forking).

By the definition of ``does not $\chi_*$-split".
\mn
\ub{Clause (E)(c)}:  Local character.

Why?  Let $\langle M_\alpha:\alpha \le \delta\rangle$ be $\le_{{\frak
s}_\lambda}$-increasing continuous, $\delta < \lambda^+$ and $q \in
{\Cal S}^{\text{bs}}_{{\frak s}_\lambda}(M_\delta)$.  
Using the third paragraph of the
proof of $\circledast_6$ for $\theta = \chi_*$, for some $\alpha_* <
\delta$ and $N_* \le_{{\frak s}_\lambda} M_{\alpha_*}$ of cardinality
$\theta$ the type $q$ does not $\theta$-split over $N_*$.  So clearly
$q$ does not fork over $M_{\alpha_*}$ (for ${\frak s}_\lambda$), as required.
\mn
\ub{Clause (E)(d)}:  Transitivity of non-forking.

By $\circledast_4$.
\mn
\ub{Clause (E)(e)}:  Uniqueness.

Holds by the choice of $\chi_*$, i.e. by $\circledast_1$.
\mn
\ub{Clause (E)(f)}:  Symmetry.

Why?  Let $M_\ell$ for $\ell\le 3$ and $a_0,a_1,a_2$ be as in (E)(f)$'$ in
\marginbf{!!}{\cprefix{600}.\scite{600-1.16E}}.  We can find a $\le_{\frak K}$-increasing
continuous sequence $\langle M_{0,\alpha}:\alpha \le \lambda^+\rangle$
such that $M_{0,0} = M_0,M_{0,\alpha +1}$ is $\le_{{\frak
s}_\lambda}$-universal over $M_{0,\alpha}$ and \wilog \, $M_{0,\alpha}
= \text{ EM}_{\tau({\frak K})}(\gamma_\alpha,\Phi)$ so is 
$\le_{\frak K}$-increasing continuous, and $\lambda$ divides $\gamma_\alpha$.

By (E)(g) proved below we 
can find $a^\ell_\alpha \in M_{0,\alpha +1}$ realizing
$\ortp_{{\frak s}_\lambda}(a_\ell,M_0,M_{\ell +1})$ such that 
$\ortp_{{\frak s}_\lambda}(a^\ell_\alpha,M_{0,\alpha},
M_{0,\alpha +1})$ does not fork over $M_0 =
M_{0,0}$, for $\ell=1,2$.  We can find
$N_* \le_{\frak K} M_0$ of cardinality $\chi_*$ such that
$\ortp_{{\frak s}_\lambda}(\langle a_1,a_2 \rangle,M_0,M_3)$ does 
not $\chi_*$-split over $N_*$ so $N_* \le_{\frak K} M_{0,0}$.

Then as in \scite{734-11.1.2} we get a contradiction (recalling
\marginbf{!!}{\cprefix{600}.\scite{600-1.16E}}).
\mn
\ub{Clause (E)(g)}:  Extension existence.

If $M \le_{{\frak s}_\lambda} N$ and $p \in {\Cal S}^{\text{bs}}_{{\frak
s}_\lambda}(M) = {\Cal S}^{\text{na}}_{\frak K}(M)$, then $p$ does not
$\chi_*$-split over $M_*$ for some $M_* \le_{\frak K} M$ of
cardinality $\chi_*$ by $\circledast_3$.  Let $M^* \in K_{\chi^*}$ be  
such that $M_* \le_{\frak K} M^* \le_{\frak K} M$ and $M^*$ is $\le_{\frak
K}$-universal over $M_*$.
As $M,N \in K_{{\frak s}_\lambda} \subseteq K_\lambda$ are 
saturated there is an isomorphism
$\pi$ from $M$ onto $N$ over $M^*$ and let $q = \pi(p)^+$.

Now $q \restriction M = p$ by $\circledast_1$ as both are from ${\Cal
S}^{\text{na}}_{\frak K}(M)$, does not $\chi_*$-split over $M_*$ and
has the same restriction to $M^*$.
\mn
\ub{Clause (E)(h)}:  Follows by \marginbf{!!}{\cprefix{600}.\scite{600-1.16A}}(3),(4) recalling
${\frak s}_\lambda$ is full.
\mn
\ub{Clause (E)(i)}:  Follows by \marginbf{!!}{\cprefix{600}.\scite{600-1.15}}.

So we have finished proving ``${\frak s}_\lambda$ is a good $\lambda$-frame.]
\mr
\item "{$\odot_2$}"    If $\lambda \in (\chi_*,\chi)$ 
then ${\frak K}^{{\frak s}_\lambda}$ is 
${\frak K} \restriction \{M:M$ is $\lambda$-saturated$\}$.
\ermn
[Why?  Should be clear.]
\mr
\item "{$\odot_3$}"    $\lambda_*$ is well defined where
\nl
$\lambda_* = \text{\rm Min}\{\lambda:\chi_* < \lambda < \chi$ and
$2^{\lambda^{+n}} < 2^{\lambda^{+n+1}}$ for every $n < \omega\}$.
\ermn
[Why?  By clause (c) of the assumption.]

Let $\Theta = \{\lambda^{+n}_*:n < \omega\}$.
\mr
\item "{$\odot_4$}"   ${\frak s}_\lambda$ is weakly succsesful for
$\lambda \in \Theta$.
\ermn
[Why?  Recalling that ``${\frak s}_\lambda$ categorical", by Definition
\marginbf{!!}{\cprefix{705}.\scite{705-stg.0A}}, Definition \marginbf{!!}{\cprefix{600}.\scite{600-nu.1}} and Observation
\marginbf{!!}{\cprefix{600}.\scite{600-nu.13.1}}(b) this means that if
$(M,N,a) \in K^{3,\text{bs}}_{{\frak s}_\lambda}$ then for some
$(M_1,N_1,a) \in K^{3,\text{uq}}_{{\frak s}_\lambda}$ we have
$(M,N,a)  \le^{\text{bs}}_{{\frak s}_\lambda} (M_1,N_1,a)$ (see
Definition \marginbf{!!}{\cprefix{600}.\scite{600-nu.1A}}).  
Toward contradiction, assume that this fails.
Let $\langle M_\alpha:\alpha < \lambda^+\rangle$ be $\le_{{\frak
s}_\lambda}$-increasing continuous, $M_{\alpha +1}$ is brimmed over
$M_\alpha$ for $\alpha < \lambda^+$ such that $M_0 = M$.  
Now directly by the definitions 
(as in \sectioncite[\S5]{600}, see more in \xCITE{838}) 
we can find $\langle M_\eta,f_\eta:\eta \in
{}^{\lambda^+ >}2\rangle$ such that:
\mr
\item "{$(a)$}"   if $\eta \triangleleft \nu \in {}^{\lambda^+ >}2$ then
$M_\eta \le_{{\frak s}_\lambda} M_\nu$
\sn
\item "{$(b)$}"   if $\eta \in {}^{\lambda^+ >}2$ then 
$f_\eta$ is a one-to-one function from $M_{\ell g(\eta)}$ to $M_\eta$
over $M_0 = M$ such that $\rho \triangleleft \eta \Rightarrow f_\rho \subseteq
f_\eta$ and $f_\eta(M_{\ell g(\eta)}) \le_{{\frak s}_\lambda} M_\eta$ 
in fact $f_0 = \text{ id}_M$ and $(M,N,a)
\le^{\text{bs}}_{{\frak s}_\lambda}(f_\eta(M_{\ell g(\eta)}),M_\eta,a) \in
K^{\text{bs}}_{\frak s}$
\sn
\item "{$(c)$}"   if $\nu = \eta \char 94 \langle \ell\rangle
\in {}^{\lambda >}2$ then $M_\nu$ is brimmed over $M_\eta$
\sn
\item "{$(d)$}"   if $\eta \in {}^{\lambda^+>}2$ then $f_{\eta\char 94
<0>} (M_{\ell g(\eta)+1}) = f_{\eta \char 94 <1>}(M_{\ell g(\eta)+1})$
\sn
\item "{$(e)$}"   if $\eta \in {}^{\lambda >}2$ then there is no triple
$(N,f_0,f_1)$ such that $f_{\eta \char 94 \langle 1 \rangle}
(M_{\ell g(\eta)+1}) \le_{\frak s} N$, and
$f_\ell$ is a $\le_{{\frak s}_\lambda}$-embedding of 
$M_{\eta \char 94 <\ell>}$ into $N$ over 
$f_{\eta \char 94 <\ell>}(M_{\ell g(\eta)+1})$ for $\ell=0,1$ and $f_0
\restriction M_\eta = f_1 \restriction M_\eta$.
\ermn
Having carried the induction  by renaming \wilog \, $\eta \in
{}^{\lambda^+ >} 2 \Rightarrow f_\eta = \text{ id}_{M_{\ell g(\eta)}}$.
Now $M_* := \cup\{M_\alpha:\alpha < \lambda^+\}$; it belongs to ${\frak
s}_{\lambda^+}$ and is saturated and for $\eta \in {}^{\lambda^+}2$ let
$M_\eta := \cup\{M_{\eta \restriction \alpha}:\alpha <\lambda^+\}$ so
$M_* \le_{{\frak s}_{\lambda^+}} M_\eta \in K_{{\frak s}_{\lambda^+}}$.
But $\chi$ is a limit cardinal so
also $\lambda^+ \in (\kappa,\chi)$ so let 
$N_* \in K_{{\frak s}_{\lambda^+}}$ be 
$\le_{{\frak s}_{\lambda^+}}$-universal over $M_*$, so for
every $\eta \in {}^{\lambda^+}2$ there is an 
$\le_{{\frak s}^+}$-embedding $h_\eta$ of $M_\eta$ into 
$N_*$ over $M_*$.  But $2^\lambda <
2^{\lambda^+}$ by the choice of $\lambda_*$ so by
\marginbf{!!}{\cprefix{88r}.\scite{88r-0.wD}} we get a contradiction to clause (e).]
\mr
\item "{$\odot_5$}"  for $\lambda \in \Theta$, if $M \in 
K^{{\frak s}_\lambda}_{\lambda^+}$ is saturated above $\lambda$ for
$K^{{\frak s}_\lambda}$, \ub{then} $M$ is saturated for ${\frak K}$.
\ermn
[Why?  Should be clear and implicitly was proved above.]
\mr
\item "{$\boxdot_1$}"   NF$_{{\frak s}_\lambda}$ is well defined
and is a non-forking relation on ${\frak K}_{{\frak s}_\lambda}$
respecting ${\frak s}_\lambda$ (for $\lambda \in \Theta$).
\ermn
[Why?  By \sectioncite[\S6]{600} as ${\frak s}_\lambda$ is a weakly
successful good $\lambda$ frame.]
\mr
\item "{$\boxdot_2$}"   ${\frak s}_\lambda$ is a good$^+ \,
\lambda$-frame (for $\lambda \in \Theta$).
\ermn
[Recalling Definition \marginbf{!!}{\cprefix{705}.\scite{705-stg.1}}, assume that this fails
so there are $\langle
M_i,N_i:i < \lambda^+\rangle$ and $\langle a_{i+1}:i <
\lambda^+\rangle$, as there, i.e. $a_{i+1} \in M_{i+2} \backslash
M_{i+1},\ortp_{{\frak s}_\lambda}(a_{i+1},
M_{i+1},M_{i+2})$ does not fork over $M_0$ for
${\frak s}_\lambda$, but $\ortp_{{\frak
s}_\lambda}(a_{i+1},N_0,M_{i+1})$ forks over $M_0$.  Also, recalling
Definition \marginbf{!!}{\cprefix{705}.\scite{705-stg.1}} the model 
$M = \cup\{M_i:i < \lambda^+\}$ is saturated for ${\frak K}^{{\frak
s}_\lambda}_{\lambda^+}$ hence by $\odot_5$ for ${\frak K}$, so  it belongs to
$K_{{\frak s}_{\lambda^+}}$.
\nl
We can find an isomorphism $f_0$ from $M$ onto EM$_{\tau({\frak
K})}(\lambda^+,\Phi)$, by $(*)_6$.  
By the ``moreover" from $(*)_6$, more exactly by $(*)_{10}$
we can find a $\le_{\frak K}$-embedding $f_1$ of
$N = : \cup\{N_i:i <\lambda^+\}$ into EM$_{\tau({\frak K})}
(\lambda \times \lambda,\Phi)$ extending $f_0$.  As we can increase the $N_i$'s
without loss of generality $f_1$ is onto EM$_{\tau({\frak K})}(\lambda
\times \lambda,\Phi)$.  We can find $\delta < \lambda^+$ such that $N_\delta =
\text{EM}_{\tau({\frak K})}(u,\Phi)$ where $u = \{\lambda \alpha +
\beta:\alpha,\beta < \delta\}$.  
By $a_{\delta +1}$ we get a contradiction to $\circledast_5$.]
\mr
\item "{$\boxdot_3$}"  Let $\lambda \in \Theta$
\sn
\item "{${{}}$}"  $(\alpha) \quad \le^*_{{\frak s}_\lambda}$ is a
partial order on $K^{\text{nice}}_{\lambda^+}[{\frak s}_\lambda] = 
K_{{\frak s}_{\lambda^+}}$ and $(K_{{\frak s}_{\lambda^+}},
\le^*_{{\frak s}_\lambda})$ satisfies
\nl

\hskip25pt  the demands on a.e.c. except possibly smoothness, see
\sectioncite[\S7]{600} 
\sn
\item "{${{}}$}"  $(\beta) \quad$ if $M \in K_{\lambda^+}$ is
saturated and $p \in {\Cal S}_{\frak K}(M)$ \ub{then} for some pair
$(N,a)$ 
\nl

\hskip25pt we have $M \le^*_{{\frak s}_\lambda} N$ and $a \in N$
realizes $p$
\sn
\item "{${{}}$}"  $(\gamma) \quad$ if $M \in K_{\lambda^+}$ is
saturated \ub{then} some $N$ satisfies:
{\roster
\itemitem{ ${{}}$ }  $(a) \quad N \in K_{\lambda^+}$ is saturated
\sn
\itemitem{ ${{}}$ }  $(b) \quad N$ is $\le_{\frak K}$-universal over $M$
\sn
\itemitem{ ${{}}$ }  $(c) \quad M \le^*_{{\frak s}_\lambda} N$
\endroster}
\item "{${{}}$}"  $(\delta) \quad {\frak s}_\lambda$ is successful.
\ermn
[Why?  \ub{Clause $(\alpha)$}:  

We know that both $K^{\text{nice}}_{\lambda^+}[{\frak s}_\lambda]$
and $K_{{\frak s}_{\lambda^+}}$ are the class of saturated $M \in
K_\lambda$.  The rest holds by \sectioncite[\S7,\S8]{600}.
\bn
\ub{Clause $(\beta)$}:

By $\circledast_3$ we 
can find $M_* \le_{\frak K} M$ of cardinality $\chi_*$ such that $p$
does not $\chi_*$-split over it (equivalently does not
$\lambda^+$-split over it).

Let $\langle M_\alpha:\alpha < \lambda^+\rangle$ be $\le_{{\frak
s}_\lambda}$-increasing continuous such that $M_{\alpha +1}$
is brimmed over $M_\alpha$ for ${\frak s}_\lambda$ for every $\alpha < 
\lambda^+$ and $M_* \le_{\frak K} M_0$ (so $\|M_*\| < \|M_0\|$
otherwise we would require $M_0$ is brimmed over $M_*$).  Hence 
$\cup\{M_\alpha:\alpha < \lambda^+\} \in
K_{\lambda^+}$ is saturated (by $\odot_5$) so \wilog \, is equal to
$M$.  We can choose
$a_*,N_\alpha(\alpha < \lambda)$ such that $\langle N_\alpha:\alpha <
\lambda^+\rangle$ is $\le_{{\frak s}_\lambda}$-increasing continuous,
$M_\alpha \le_{{\frak s}_\lambda} M_\alpha$, NF$_{{\frak
s}_\lambda}(M_\alpha,N_\alpha,M_\beta,M_\beta)$ for $\alpha < \beta <
\lambda^+,N_{\alpha +1}$ is brimmed over $M_{\alpha +1} \cup N_\alpha$
and $\ortp_{{\frak s}_\lambda}(a,N_0,M_0) = p \restriction M_0$ so $a \in
N_0$.  Let $N = \cup\{N_\alpha:\alpha < \lambda^+\}$ so again $N \in
K_{\lambda^+}$ is saturated (equivalently $N \in
K^{\text{nice}}_{\lambda^+}[{\frak s}_\lambda]$) and $M \le_{\frak K}
N$ and even $M \le^*_{{\frak s}_\lambda} N$ (by the definition of
$\le^*_{{\frak s}_\lambda}$).  For each $\alpha < \lambda^+$ we have
NF$_{{\frak s}_\lambda}(M_0,N_0,M_\alpha,N_\alpha)$ but NF$_{{\frak
s}_\lambda}$ respect ${\frak s}_\lambda$ hence $\ortp_{{\frak
s}_\lambda}(a,M_\alpha,N_\alpha)$ does not fork over $M_0$ hence by
the definition of ${\frak s}_\lambda$ the type 
$\ortp_{{\frak s}_\lambda}(a,M_\alpha,N_\alpha)$ does 
not $\lambda$-split over $M_*$ hence
$\ortp_{{\frak s}_\lambda}(a,M_\alpha,N_\alpha)= p \restriction M_\alpha$.  As
this holds for every $\alpha < \lambda^+$, by the choice of $\chi_*$,
i.e. by $\circledast_1$ clearly $a$ realizes $p$.
\bn
\ub{Clause $(\gamma)$}:

By clause $(\beta)$ as in the proofs in \sectioncite[\S4]{600}; that is, we
choose $N\in K_{\lambda^+}$ which is $\le_{{\frak
K}_\lambda}$-universal over $M$.  We now try to choose
$(M_\alpha,f_\alpha,N_\alpha)$ by induction on $\alpha < \lambda^+$
such that: $M_0 = M,N_0 = N,f_0 = \text{ id}_M,M_\alpha$ is
$\le^*_{{\frak s}_\lambda}$-increasing continuous, $N_\alpha$ is
$\le_{\frak K}$-increasing continuous, $f_\alpha$ is a $\le_{\frak
K}$-embedding of $M_\alpha$ into $N_\alpha,f_\alpha$ is
$\subseteq$-increasing continuous with $\alpha$ and $\alpha = \beta +1
\Rightarrow f_\alpha(M_\alpha) \cap N_\beta \ne f_\beta(M_\beta)$.

For $\alpha=0,\alpha$ limit no problems.  If $\alpha = \beta +1$ and
$f_\alpha(M_\alpha) = N_\alpha$ we are done and otherwise use clause
$(\beta)$.  But by Fodor lemma we cannot carry the induction for every
$\alpha < \lambda^+$, so we are done proving $(\gamma)$.
\bn
\ub{Clause $(\delta)$}:

We should verify the conditions in Definition \marginbf{!!}{\cprefix{705}.\scite{705-stg.0A}}.
Now clause (a) there, being weakly successful, holds
by $\odot_4$.  As for clause (b) there, it suffices to prove 
that if $M_1,M_2 \in K^{\text{nice}}_{\lambda^+}[{\frak s}_\lambda] = 
K_{{\frak s}^+_\lambda}$ and
$M_1 \le_{\frak K} M_2$ then $M_1 \le^*_{{\frak s}_\lambda} M_2$ which
means: \ub{if} $\langle M^\ell_\alpha:\alpha < \lambda^+\rangle$ is
$\le_{{\frak s}_\lambda}$-increasing continuous, $M^\ell_{\alpha +1}$
is brimmed over $M^\ell_\alpha$ with $M_\ell =
\cup\{M^\ell_\alpha:\alpha < \lambda^+\}$, \ub{then} for some club
$E$ of $\lambda^+$ for every $\alpha < \beta$ from $E$, 
NF$_{{\frak s}_\lambda}(M^1_\alpha,M^2_\alpha,M^1_\beta,M^2_\delta)$.

By clause $(\gamma)$ there is $N \in K_{{\frak s}^+_\lambda}$ such
that $M_1 \le^*_{{\frak s}_{\lambda^+}} N$ 
(hence $M_1 \le_{\frak K} N$) and $N$ is
$\le_{{\frak K}^{{\frak s}_\lambda}}$-universal over $M_1$.  So \wilog $M_2
\le_{\frak K} N$ but by \marginbf{!!}{\cprefix{600}.\scite{600-ne.3}}(3) all this implies $M_1
\le^*_{\lambda^+} M_2$.  So we are done proving $\boxdot_3$. 
\mr
\item "{$\boxdot_4$}"   ${\frak s}_{\lambda^+}$ is the successor of
${\frak s}_\lambda$ for $\lambda \in \Theta$.
\ermn
[Why?  Now by $\boxdot_3$ the good frame ${\frak s}_\lambda$ is  
successful; by \marginbf{!!}{\cprefix{705}.\scite{705-stg.3}} we know that ${\frak s}^+_\lambda$
is a well defined good $\lambda^+$-frame.  
Clearly $K_{{\frak s}_\lambda(+)}$ is the class of saturated $M
\in {\frak K}_{\lambda^+}$, by $\odot_5$, see the definitions in
\marginbf{!!}{\cprefix{600}.\scite{600-ne.1}}, \marginbf{!!}{\cprefix{600}.\scite{600-rg.7}}(5).  But
${\frak s}_\lambda$ is good$^+$ by
$\boxdot_2$ so by \marginbf{!!}{\cprefix{705}.\scite{705-stg.3B}} we know that
$\le_{{\frak s}_\lambda(+)} = <^*_{\lambda^+}[{\frak s}_\lambda]$ is equal to
$\le_{\frak K} \restriction K_{{\frak s}_\lambda(+)}$, so ${\frak
K}_{{\frak s}_\lambda(+)} = {\frak K}_{{\frak s}_{\lambda^+}}$.
As both ${\frak s}_\lambda(+)$ and ${\frak s}_{\lambda^+}$ are full,
clearly ${\Cal S}^{\text{bs}}_{{\frak s}_\lambda(+)} = {\Cal
S}^{\text{bs}}_{{\frak s}_{\lambda^+}}$.  For $M_1 \le_{{\frak
s}_\lambda(+)} M_2 \le_{{\frak s}_\lambda(+)} M_3$ and $a \in M_3
\backslash M_2$, comparing the two definitions of ``$\ortp_{{\frak
K}_{{\frak s}_\lambda(+)}}(a,M_2,M_1)$ does not fork over $M_1$" they
are the same.  So we are done.]
\mr
\item "{$\boxdot_5$}"  ${\frak s}_{\lambda^{+\omega}_*}$ is the limit of
$\langle {\frak s}_{\lambda_*}^{+n}:n < \omega\rangle$.
\ermn
[Why?  Should be clear.]
\mr
\item "{$\boxdot_6$}"   ${\frak s}_\lambda$ satisfies the hypothesis
\marginbf{!!}{\cprefix{705}.\scite{705-12.1}} of \sectioncite[\S12]{705} if $\lambda \in \Theta
\backslash \lambda^{+3}_*$ holds.
\ermn
[Why?  By $\boxdot_2,\boxdot_3,\boxdot_4$ and \marginbf{!!}{\cprefix{705}.\scite{705-12D.1}}.]

Hence
\mr
\item "{$\boxdot_7$}"  ${\frak s}_{\lambda_*}$ is beautiful $\lambda^{+
\omega}_*$-frame.
\ermn
[Why?  By \marginbf{!!}{\cprefix{705}.\scite{705-12b.14}} and \marginbf{!!}{\cprefix{705}.\scite{705-12f.16A}}.]
\mr
\item "{$\boxdot_8$}"   $K[{\frak s}_{\lambda^{+ \omega}_*}]$ is
categorical in one $\chi > \lambda^{+ \omega}_*$ iff it is categorical
in every $\chi > \lambda^{+ \omega}$.
\ermn
[Why?  By \marginbf{!!}{\cprefix{705}.\scite{705-12f.16A}}(d),(e).]
\mr
\item "{$\boxdot_9$}"  if $\lambda \ge \beth_{1,1}(\lambda^{+
\omega}_*)$ then ${\frak K}_\lambda = 
{\frak K}_\lambda[{\frak s}_{\lambda^{+ \omega}_*}]$.
\ermn
[Why?  The conclusion $\supseteq$ is obvious.  For the other inclusion let $M
\in K_\lambda$, now by the definition of class in the left, it is 
enough to prove that $M$ is $(\lambda^{+ \omega}_*)^+$-saturated.
But otherwise by the omitting type theorem for
a.e.c., i.e. by \scite{734-0n.8}(1),(d), (or see
\cite[8.6=X1.3A]{Sh:394}) there is such a 
model $M' \in K_\mu$, contradiction to $(*)_4$.]
\sn
By $\boxdot_8 + \boxdot_9$ we are done.  \hfill$\square_{\scite{734-am3.6}}$ 
\enddemo


\nocite{ignore-this-bibtex-warning} 
%
  REFERENCES.  
  \bibliographystyle{lit-plain}
  \bibliography{lista,listb,listx,listf,liste}

\enddocument